\documentclass[a4paper,12pt]{amsart}
\usepackage[utf8x]{inputenc}
\usepackage{amsmath}
\usepackage{amsthm}
\usepackage{amscd}
\usepackage{amssymb}
\usepackage{mathrsfs}
\usepackage{amsfonts}
\usepackage{graphicx}
\usepackage{longtable}
\usepackage[all]{xy}
\usepackage{verbatim}
\usepackage{hyperref} 

\theoremstyle{plain}
\newtheorem{theorem}{Theorem}
\newtheorem{lemma}{Lemma}[section]
\newtheorem{assumption}[lemma]{Assumption}

\newtheorem{corollary}[lemma]{Corollary}
\newtheorem{definition}[lemma]{Definition}
\newtheorem{example}{Example}
\newtheorem{fact}[lemma]{Fact}
\newtheorem{proposition}[lemma]{Proposition}

\newtheorem{remark}[lemma]{Remark}
\newtheorem{atheorem}{Theorem}[section]
\newtheorem{aprop}{Proposition}[section]
\newtheorem{alemma}[aprop]{Lemma}
\newtheorem{adefinition}[aprop]{Definition}
\newtheorem{aexample}[aprop]{Example}
\newtheorem{aremark}[aprop]{Remark}

\def\A{\operatorname{A}}
\def\B{\operatorname{B}}
\def\BC{\operatorname{BC}}
\def\C{\operatorname{C}}
\def\D{\operatorname{D}}
\def\E{\operatorname{E}}
\def\F{\operatorname{F}}
\def\G{\operatorname{G}}
\def\GL{\operatorname{GL}}

\def\Pin{\operatorname{Pin}}

\def\PSO{\operatorname{PSO}}

\def\O{\operatorname{O}}

\def\SO{\operatorname{SO}}
\def\Sp{\operatorname{Sp}}
\def\SU{\operatorname{SU}}
\def\U{\operatorname{U}}

\def\ad{\operatorname{ad}}
\def\Ad{\operatorname{Ad}}

\def\Aut{\operatorname{Aut}}
\def\Ave{\operatorname{Ave}}

\def\opd{\operatorname{d}}

\def\depth{\operatorname{depth}}
\def\der{\operatorname{der}}

\def\diag{\operatorname{diag}}

\def\Fold{\operatorname{Fold}}

\def\Gal{\operatorname{Gal}}

\def\Hom{\operatorname{Hom}}
\def\id{\operatorname{id}}

\def\Ima{\operatorname{Im}}
\def\Ind{\operatorname{Ind}}
\def\Int{\operatorname{Int}}

\def\Lie{\operatorname{Lie}}

\def\nd{\operatorname{nd}}

\def\Out{\operatorname{Out}}

\def\Pin{\operatorname{Pin}}

\def\rank{\operatorname{rank}}

\def\span{\operatorname{span}}

\def\Spin{\operatorname{Spin}}
\def\st{\operatorname{st}}
\def\Stab{\operatorname{Stab}}
\def\supp{\operatorname{supp}}
\def\Sym{\operatorname{Sym}}

\def\tr{\operatorname{tr}}

\newcommand{\scrR}{\mathscr{R}}
\newcommand{\rmd}{\mathrm{d}}
\newcommand{\rmm}{\mathrm{m}}
\newcommand{\rms}{\mathrm{s}}

\newcommand{\bbR}{\mathbb{R}}
\newcommand{\bbZ}{\mathbb{Z}}

\newcommand{\fre}{\mathfrak{e}}

\newcommand{\frg}{\mathfrak{g}}

\newcommand{\frs}{\mathfrak{s}}

\begin{document}

\title{The functorial source problem via dimension data}
\date{February 24, 2022}

\author{Jun Yu}
\address{Beijing International Center for Mathematical Research and School of Mathematical Sciences, Peking
University, No. 5 Yiheyuan Road, Beijing 100871, China.}
\email{junyu@bicmr.pku.edu.cn}

\abstract{}
For an automorphic representation $\pi$ of Ramanujan type, there is a conjectural subgroup $\mathcal{H}_{\pi}$
of the Langlands L-group $^{L}G$ associated to $\pi$ (\cite{Langlands}), called the {\it functional source} of
$\pi$. The functorial source problem proposed by Langlands and refined by Arthur (\cite{Langlands},
\cite{Arthur-note}, \cite{Arthur-Howe}, \cite{Arthur-functoriality}) intends to determine $\mathcal{H}_{\pi}$
via analytic and arithmetic data of $\pi$. In this paper we consider the functorial source problem of automorphic
representations of a split group, a unitary group or an orthogonal group which do not come from endoscopy and have
minimal possible ramification. In this setting, $\mathcal{H}_{\pi}$ must be an S-subgroup of $^{L}G$ (Definition
\ref{D:Sgroup}). We approach the functorial source problem by proving distinction and linear independence among
dimension data of S-subgroups. Nice results along this direction are shown in this paper (e.g, Theorems
\ref{T1}-\ref{T6}). We define a notion of quasi root system and use it as the key tool for studying S-subgroups
and their dimension data.
\endabstract

\noindent

\noindent

\maketitle

\tableofcontents

\section*{Introduction}

This paper is motivated by some questions of Langlands and Arthur concerning the functorial source of automorphic
forms (\cite{Langlands}, \cite{Arthur-note}, \cite{Arthur-Howe}, \cite{Arthur-functoriality}). Let $\mathbb{G}/F$
be a connected reductive linear algebraic group over a number field $F$ and $\mathbb{G}(A_{F})$ be the corresponding
group over the adelic ring $A_{F}$ of $F$. Let $\pi$ be an automorphic representation of $\mathbb{G}(A_{F})$ of
Ramanujan type. Then there is a conjectural subgroup $\mathcal{H}_{\pi}$ of the Langlands L-group $^{L}G$ associated
to $\pi$ (\cite{Langlands}), called the {\it functional source} of $\pi$. The mysterious group $\mathcal{H}_{\pi}$
should satisfy the following simple condition: for each finite-dimensional irreducible complex linear representation
$\rho$ of $^{L}G$, the pole order at $s=1$ of the Langlands L-function $L(s,\pi,\rho)$ is equal to the dimension of
the space of $\mathcal{H}_{\pi}$ invariant vectors in $\rho$. The right hand side of the above condition is called
the dimension datum of $\mathcal{H}_{\pi}$ as a subgroup of $^{L}G$. Langlands thought that the group
$\mathcal{H}_{\pi}$ is the key to show functoriality beyond that the endoscopy theory has achieved. In
\cite{Langlands} Langlands proposed a strategy of showing this general form of functoriality through a refined form
of stable trace formulas by inserting information of pole orders at $s=1$ of Langlands L-functions $L(s,\pi,\rho)$.
This is called Langlands' beyond endoscopy program nowadays. If succeeded, this general functoriality will give a
tremendous advancement in automorphic form theory. In \cite{Arthur-note}, \cite{Arthur-Howe},
\cite{Arthur-functoriality}, \cite{Arthur-stratification} and some other papers, Arthur further refined Langlands'
ideas motivated by his great successful strategy of classifying automorphic representations of classical groups using
the stabilization of trace formulas (\cite{Arthur-endoscopy}, \cite{Mok}). Particularly, Arthur wrote out more concrete
forms of stable trace formulas after inserting information of L-functions and he
proposed a primitization process of the stable trace formulas (\cite{Arthur-Howe}, \cite{Arthur-functoriality}).
Arthur also proposed ideas/questions to determine $\mathcal{H}_{\pi}$ using dimension datum of $\mathcal{H}_{\pi}$
and Satake parameters of $\pi$ (\cite[Section 2]{Arthur-note}, \cite[Section 4]{Arthur-functoriality}). This paper
concerns the determination of $\mathcal{H}_{\pi}$ via its dimension datum alone. We show an affirmative answer to the
``refined question" in \cite[p. 17]{Arthur-functoriality} when $\mathbb{G}$ is a split group, a unitary group or an
orthogonal group and $\pi$ does not come from endoscopy and has minimal possible ramification.

The group $\mathcal{H}_{\pi}$ should be a reductive subgroup of $^{L}G$. Thus, it is equivalent to consider
its maximal compact subgroup as a subgroup of a maximal compact subgroup of $^{L}G$. Working with compact
groups provides simplicity while taking topological argument. Let $G$ be a compact Lie group with cyclic
component group. We call a closed subgroup $H$ of $G$ an {\it $S$-subgroup} if $H$ is generated by $H^{0}$
and an element $h_{0}\in g_{0}G^{0}$ such that $h_{0}^{d}\in Z_{G}(G^{0})H^{0}$ and $Z_{G}(H)=Z(G)$, where
$d$ is the order of $\Ad(g_{0})|_{\mathfrak{g}}\in\Out(\mathfrak{g})=\Aut(\mathfrak{g})/\Int(\mathfrak{g})$
and $g_{0}\in G$ is an element generating $G/G^{0}$. The definition of $S$-subgroups is motivated by properties
of the group $\mathcal{H}_{\pi}$ should satisfied while the automorphic representation $\pi$ has minimal
possible ramification and is not the image of any endoscopic functoriality. It is named after Dynkin's
S-subalgebra (\cite{Dynkin},\cite{Minchenko}). Our aim in this paper is to consider the functorial source
problem while $\mathbb{G}/F$ is either split or is a unitary group or an orthogonal group. Turning to the
L-group side, we are led to study dimension data of S-subgroups of the following groups: (1)a connected
compact Lie group; (2)$\U(N)\rtimes\langle\tau\rangle$, where $\tau^{2}=I$ and $\tau X\tau^{-1}=\bar{X}$
($\forall X\in\U(N)$); (3)$\O(2N)$. Besides these groups, we also study dimension data of S-subgroups of
$\Spin(8)\rtimes\langle\tau\rangle$ (triality form) and $\E_{6}\rtimes\langle\tau\rangle$ (disconnected
$\mathbf{E}_{6}$). Moreover, we reduce the study of S-subgroups of isogeny forms of these groups to
S-subgroups of them.

We define a notion of quasi root system and use it as the key tool for studying S-subgroups and their dimension
data. We call a compact abelian Lie group $\tilde{S}$ a {\it quasi torus} if its component group
$\tilde{S}/\tilde{S}^{0}$ is a cyclic group. A {\it quasi root system} $\Phi$ on a quasi torus $\tilde{S}$
(or in $X^{\ast}(\tilde{S})$) is a finite subset of $X^{\ast}(\tilde{S})$ with certain properties (Definition
\ref{D:QRS1}). In particular, if $\Phi$ is a quasi root system in $X^{\ast}(\tilde{S})$, then the set $p(\Phi)$
(called the restricted root system of $\Phi$) defined by \[p(\Phi)=\{\alpha|_{\tilde{S}^{0}}:\alpha\in\Phi\}\]
is a root system in $X^{\ast}(\tilde{S}^{0})$. When $\tilde{S}$ is a generalized Cartan subgroup of a compact Lie
group $H$, there is a quasi root system $R(H,\tilde{S})$ constructed from the action of $\tilde{S}$ on the
complexified Lie algebra $\mathfrak{h}$ of $H$ which is moreover reduced in the sense that $\alpha\neq 2\beta$
for any two roots $\alpha,\beta\in R(H,\tilde{S})$. Quasi root systems are complicated in general and the quasi
root system $R(H,\tilde{S})$ does characterize the group $H$ up to isomorphism in general. It is miraculous that
the quasi root system of an S-subgroup $H$ of $\U(N)\rtimes\langle\tau\rangle$ or $\O(2N)$ characterizes the
subgroup well (cf. Lemmas \ref{L:ARS-conjugacy-TA1} and \ref{L:ARS-conjugacy-TD}) and is subject to strong
constraints (cf. \S \ref{SS:TA1-ARS}, \S \ref{SS:TD1} and \S \ref{SS:TD2}).

Let's state main results of this paper. Two closed subgroup $H_{1},H_{2}$ of a group $G$ are said to be element
conjugate if there exists an isomorphism $\phi: H_{1}\rightarrow H_{2}$ such that \[\phi(x)\sim x,\quad
\forall x\in H_{1}\] (cf. \cite{Larsen2}). Element conjugate subgroups must have the same dimension datum. For
S-subgroups of a connected compact Lie group, we show the following Theorems \ref{T1} and \ref{T2} in
\S \ref{S:connected}. When $G=\U(N)$ or $\SU(N)$, Theorem \ref{T1} recovers \cite[Theorem 2]{Larsen-Pink}
(concerning distinction of dimension data) and \cite[Theorem 1.3]{Yu-dimension} (concerning linear independence
of dimension data).

\begin{theorem}[]\label{T1}
Let $G$ be a connected compact Lie group. Then the dimension data of any tuple of pairwise non-element conjugate
connected S-subgroups of $G$ are linearly independent.
\end{theorem}

\begin{theorem}\label{T2}
Let $H_{1}$ and $H_{2}$ be two non-element conjugate S-subgroups of a connected compact Lie group $G$. Then
\[\mathscr{D}_{H_{1}}\neq\mathscr{D}_{H_{2}}.\]
\end{theorem}

For S-subgroups of $\U(N)\rtimes\langle\tau\rangle$, we show the following Theorems \ref{T3} and \ref{T4} in
\S \ref{S:TA1}.

\begin{theorem}\label{T3}
Let $H_{1},\dots,H_{s}$ be a set of pairwise non-conjugate S-subgroups of $\U(N)\rtimes\langle\tau\rangle$. Suppose
that: (1)for any quasi Cartan subgroup $\tilde{S}$ of any of $H_{1},\dots,H_{s}$,
$(\Psi_{\tilde{S}},\chi_{\rho}|_{\tilde{S}^{0}})$ does not contain an irreducible factor $(\Psi_{i},\chi_{i})$
isomorphic to $(\B_{n}^{(2)},([\frac{1}{2}]+[-\frac{1}{2}])^{\otimes n})$ ($n\geq 2$); (2)$H_{1},\dots,H_{s}$
are all semisimple or all non-semisimple. Then, the dimension data \[\mathscr{D}_{H_{1}},\dots,\mathscr{D}_{H_{s}}\]
are linearly independent.
\end{theorem}

\begin{theorem}\label{T4}
Let $H_{1}$ and $H_{2}$ be two non-conjugate S-subgroups of $\U(N)\rtimes\langle\tau\rangle$. Then
\[\mathscr{D}_{H_{1}}\neq\mathscr{D}_{H_{2}}.\]
\end{theorem}

In Theorem \ref{T3}, $\Psi_{\tilde{S}}$ is a canonical quasi root system associated to the quasi torus $\tilde{S}$
in $\U(N)\rtimes\langle\tau\rangle$, $\rho=\mathbb{C}^{N}$ is the natural representation of $\U(N)$ and
$\chi_{\rho}$ means its character.

For S-subgroups of $\O(2N)$, we show the following Theorems \ref{T5} and \ref{T6} in \S \ref{S:TD}.

\begin{theorem}\label{T5}
Let $H_{1},\dots,H_{s}$ be a set of pairwise non-conjugate S-subgroups of $\O(2N)$. Suppose that
they are all semisimple or all non-semisimple. Then, the dimension data \[\mathscr{D}_{H_{1}},\dots,
\mathscr{D}_{H_{s}}\] are linearly independent.
\end{theorem}

\begin{theorem}\label{T6}
Let $H_{1}$ and $H_{2}$ be two non-conjugate S-subgroups of $\O(2N)$. Then \[\mathscr{D}_{H_{1}}\neq
\mathscr{D}_{H_{2}}.\]
\end{theorem}

From Theorems \ref{T2}, \ref{T4} and \ref{T6} we see that: non-element conjugate S-subgroups have different
dimension datum in the respective setting. Such a statement is remarkable while applying dimension datum in
the beyond endoscopy program. In \cite{Larsen-Pink}, \cite{An-Yu-Yu} and \cite{Yu-dimension}, we have
seen examples of non-conjugate (or non-isomorphic) connected closed subgroups with the same dimension datum.
From Theorems \ref{T1}, \ref{T3} and \ref{T5} we see that: we have linear independence statement for dimension
data of S-subgroups under mild assumptions. The linear independence statement is important while applying stable
trace formula in Langlands' beyond endoscopy program. Note that the assumptions in Theorems \ref{T3} and \ref{T5}
are also necessary as Theorems \ref{T:counter1}-\ref{T:counter3} indicated. In \cite{Yu-dimension}, we classified
connected closed subgroups with the same dimension datum and found generators of linear relations among different
dimension data of connected closed subgroups. In summary, we could say that: the S-subgroup condition is a natural
condition which guarantees distinction (or linear independence) among dimension data under no (or mild) assumption.

Considering more general algebraic groups $\mathbb{G}/F$ and automorphic representations $\pi$ with wider
ramification, perhaps we must use Satake parameters (beside dimension data) to treat the functional source
problem. It looks to the author that the idea of considering automorphic representations not coming from endoscopy
and the tool of quasi root system might still be useful.

Let's explain the proofs of Theorems \ref{T1}-\ref{T6}. The most technical part is to show constraints on
the quasi root system $\Phi=R(H,\tilde{S})$ of an S-subgroup $H$. Precisely to say, by Lemmas \ref{L:root3},
\ref{L:TA1-Psi0} and \ref{L:TD-Psi0}, we know that the quasi root system $\Phi$ is contained in a canonical
reduced quasi root system $\Psi$ depending on $\tilde{S}\subset G$ and $p(\Phi)\supset p(\Psi)^{\circ}$ (=the set
of short roots in $p(\Psi)$). When $G$ is split, this suffices to show Theorems \ref{T1}-\ref{T2} and a
generalization of them: Theorem \ref{T:S3}. When $G=\U(N)\rtimes\langle\tau\rangle$ or $\O(2N)$, the hard part is
to show constraints on folding indices (cf. remark after Lemma \ref{L:ARS3}) and the fractional factor of $\Phi$
(cf. remark ahead of Lemma \ref{L:ARS5}) while $p(\Phi)$ is given. To bound the folding indices of $\Phi$ we employ
two tools of Jantzen's twining character formula (Theorem \ref{T:twining1}) and irreducible restriction character
(Definition \ref{D:resChar1}). With the condition $p(\Phi)\supset p(\Psi)^{\circ}$, the intersection of $p(\Phi)$
with an irreducible factor of $p(\Psi)$ has very few possibilities except when the irreducible factor is isomorphic
to $\B_{n}$ or $\BC_{n}$ ($n\geq 2$). This exceptional case is treated in Subsection \ref{SS:fully}. The fractional
factor of $\Phi$ is determined in Lemmas \ref{L:TA1-fracFactor3}, \ref{L:TDi-2} and \ref{L:TD-fracFactor} using the
the S-subgroup condition while the restricted root system $p(\Phi)$ and folding indices of $\Phi$ are given.

In the literature, dimension data of non-finite closed subgroups was first studied by Larsen and Pink in their
pioneering work \cite{Larsen-Pink} where they studied dimension data of connected closed semisimple subgroups of
$\U(N)$. In \cite{An-Yu-Yu} and \cite{Yu-dimension}, the author of this paper and his collaborators took a
systematic study of dimension data of connected closed subgroups (without semisimplicity constraint) of a compact
Lie group. There are applications of dimension data to monodromy groups (\cite{Larsen-Pink3}) and isospectral
manifolds (\cite{Pesce}, \cite{Sunada}, \cite{Sutton}, \cite{An-Yu-Yu}) via Pesce-Sunada method and its
generalization. However, the most important and most promising application of dimension data is that in the beyond
endoscopy program (\cite{Langlands}, \cite{Arthur-Howe}). This paper is the first one which studies dimension data
of non-connected closed subgroups. The idea of S-subgroup and the tool of quasi root system are the main innovation
of this paper.

This paper is organized as follows. In Subsection \ref{SS:ARS}, we define and classify quasi root systems. In
Subsection \ref{SS:qrs-char}, we associate a character to a quasi root system on a quasi torus. In Subsections
\ref{SS:MCC}-\ref{SS:dim-ARS}, we construct a quasi root system $R(H,\tilde{S})$ for a generalized Cartan subgroup
$\tilde{S}$ in a compact Lie group $H$, discuss the character associated to $R(H,\tilde{S})$ and its application
to dimension data of closed subgroups in a given compact Lie group. Subsection \ref{SS:actual} concerns the so-called
actual leading weight of the character associated to a quasi root system. Subsection \ref{SS:BCn} makes a thorough
study of quasi sub-root systems of a quasi root system $\Psi$ with $p(\Psi)=\BC_{n}$ and associated characters with
respect to a large group action. In Section \ref{S:twining} we recall Jantzen's twining character formula and gives
a proof of it. In Section \ref{S:Res-character} we discuss irreducible restriction characters of a quasi root system.
In Section \ref{S:S-group} we discuss S-subgroups and applications to the functorial source problem. In Section
\ref{S:connected}, we treat dimension data of S-subgroups of a connected compact Lie group. In Section \ref{S:TA1},
we treat dimension data of S-subgroups of $\U(N)\rtimes\langle\tau\rangle$. In Section \ref{S:TD}, we treat dimension
data of S-subgroups of $\O(2N)$. We remark that quasi root system could be regarded as an analogue of affine root
system (\cite{Bourbaki}, \cite{Bruhat-Tits}, \cite{Macdonald}, \cite{Reeder-Yu}). The appendix is a note written by
Professor Jiu-kang Yu which discusses concrete relations between quasi root system and affine root system.

\smallskip

\noindent{\it Notation and conventions.} Let $G$ be a compact Lie group. Write $G^{0}$ for the neutral subgroup of
$G$, $\mathfrak{g}_{0}$ for the Lie algebra of $G$, and $\mathfrak{g}=\mathfrak{g}_{0}\otimes_{\mathbb{R}}\mathbb{C}$
for the complexified Lie algebra of $G$. Write $\mathfrak{g}_{\der}=[\mathfrak{g},\mathfrak{g}]$ for the derived
subalgebra of $\mathfrak{g}$ and write $z(\mathfrak{g})$ for the center of $\mathfrak{g}$. When $G$ is connected,
write $G_{\der}=[G,G]$ for the derived subgroup of $G$ and write $Z(G)$ for the center of $G$. Then we have
$G_{\der}Z(G)^{0}=G$. Write $[x]=xN$ ($x\in G$) for an element in a quotient group $G/N$.

Let $H$ be an abelian compact Lie group. Write $X^{\ast}(H)=\Hom(H,\U(1))$ for the character group of $H$, and write
$X_{\ast}(H)=\Hom(\U(1),H)$ for the cocharacter group of $H$. We call elements in $X^{\ast}(H)$ linear characters
of $H$. For a linear character $\lambda\in X^{\ast}(H)$, write $[\lambda]: H\rightarrow\U(1)$ for the corresponding
evaluation function. For two linear characters $\lambda,\mu\in X^{\ast}(H)$ and $k_1,k_2\in\mathbb{Z}$, define a
linear character $k_{1}\lambda+k_{2}\mu\in X^{\ast}(H)$ by \[[k_{1}\lambda+k_{2}\mu](x)=
\lambda(x)^{k_{1}}\mu(x)^{k_{2}}\ (\forall x\in H).\] Write $\mathbf{1}$ for the identity character of
$\U(1)$. For any character $\lambda\in X^{\ast}(H)$ and cocharacter $\check{\mu}\in X_{\ast}(H)$, there is a unique
integer $k$ such that $\lambda\circ\check{\mu}=k\mathbf{1}$. Define \[\langle\lambda,\check{\mu}\rangle=k.\] Then,
\[\langle\cdot,\cdot\rangle: X^{\ast}(H)\times X_{\ast}(H)\rightarrow\mathbb{Z}\] is a bi-additive pairing, which
is a perfect pairing when $H$ is a torus. For a linear character $\lambda\in X^{\ast}(H)$, let $\opd\lambda\in
\mathfrak{h}_{0}^{\ast}$ denote its differential. When $H$ is a torus, a linear character $\lambda$ and its
differential $\opd\lambda$ determine each other, in this case we don't distinguish them.

Write $\E_6$ or $\E_6^{sc}$ for a connected and simply-connected compact simple Lie group of type $\E_6$; write
$\E_6^{ad}$ for a connected adjoint type compact simple Lie group of type $\E_6$; write $\fre_6$ for a compact
simple Lie algebra of type $\E_6$, and $\fre_6(\mathbb{C})$ for a complex simple Lie algebra of type $\E_6$. Similar
notations apply to simple Lie groups and simple Lie algebras of other types.

For any $n\geq 1$, write $J_{n}=\left(\begin{array}{cc}0_{n}&I_{n}\\-I_{n}&0_{n}\\\end{array}\right)$.

\smallskip

\noindent{\it Acknowledgements.} I would like to thank Jiu-Kang Yu and Xuhua He for helpful communications. Particularly,
I thank Jiu-kang Yu for suggesting the names of quasi root system, quasi-connected compact Lie group, etc and for very
helpful discussions with him concerning quasi root system and affine root system. The author learned a lot from reading
\cite{Arthur-note}, \cite{Arthur-Howe} and \cite{Arthur-functoriality}. I am grateful for that. Sections \ref{S:ARS} and
\ref{S:twining} were mostly done when the author visited MPI Bonn in the summer of 2016. The author would like to thank
MPI Bonn for the support and hospitality. This research is partially supported by the NSFC Grant 11971036.

\section{Quasi root systems}\label{S:ARS}

We call a compact abelian Lie group with cyclic component group a {\it quasi torus} and we define a {\it quasi root system}
on a quasi torus to be a finite subset of the character group of it with some properties. We give a classification of quasi
root systems. Given a quasi torus, we associate a character to each reduced quasi root system on it. For a generalized
Cartan subgroup $\tilde{S}$ in a compact Lie group $H$, we show that {\it infinite roots} from the action of $\tilde{S}$
on the complexified Lie algebra of $H$ is naturally a reduced quasi root system. We express the Weyl density function
in terms of the character associated to this quasi root system, which is the starting point of using quasi root systems
and associated characters to study dimension data of non-connected closed subgroups. Given a quasi torus $\tilde{S}$
in a compact Lie group $G$, we show that the union $\Psi'_{\tilde{S}}$ of quasi root systems of closed subgroups of
$G$ with $\tilde{S}$ a generalized Cartan subgroup is also a quasi root system. In Subsection \ref{SS:actual} we
discuss the so-called actual leading weight in the character associated to a quasi root system. In Subsection
\ref{SS:BCn} we present a precise description of reduced quasi sub-root systems in a quasi root system $\Psi$ with
$\bar{\Psi}\cong\BC_{n}$ and their associated characters with respect to a large group averaging.

\subsection{Quasi root systems: definition and classification}\label{SS:ARS}

Following usual conventions in the literature (cf. \cite{Bruhat-Tits}, \cite{Prasad-Raghunathan}), we divide roots
in a root system into several types.

\begin{definition}\label{D:root1}
Let $\Phi$ be a root system in a Euclidean space $V$. \begin{enumerate}
\item[(i)]We call a root $\alpha$ a short root if $\frac{2(\alpha,\beta)}{(\beta,\beta)}\in\{-1,0,1\}$ ($\forall\beta\in
\Phi-\{\pm{\alpha}\}$).
\item[(ii)]Call $\alpha$ a long root if $\frac{2(\alpha,\beta)}{(\alpha,\alpha)}\in\{-1,0,1\}$ ($\forall\beta\in
\Phi-\{\pm{\alpha}\}$).
\item[(iii)]Call $\alpha$ a multipliable root if $2\alpha\in\Phi$.
\item[(iv)]Call $\alpha\in\Phi$ a divisible root if $\frac{1}{2}\alpha\in\Phi$.
\item[(v)]Call $\alpha$ a non-divisible and non-multipliable root (ndm root for short) if there exists two
other roots $\beta_{1},\beta_{2}\in\Phi-\{\pm{\alpha}\}$ such that \[\frac{2(\alpha,\beta_{1})}{(\beta_{1},\beta_{1})}
=2\textrm{ and }\frac{2(\beta_{2},\alpha)}{(\alpha,\alpha)}=2.\]
\end{enumerate}
Moreover, \begin{enumerate}
\item[(i)]We write $\Phi^{\circ}$ for the set of short roots in $\Phi$.
\item[(ii)]Write $\Phi^{\nd}$ for the set of roots $\alpha\in\Phi$ which are not divisible.
\end{enumerate}
\end{definition}

Note that divisible root, multipliable root, ndm root exist only in non-reduced irreducible factors of $\Phi$; a
divisible root is also regarded as a long root and a multipliable root is also regarded as a short long; a root in
a simply-laced irreducible factor of $\Phi$ is regarded both as a long root and as a short root.
The following lemma is well-known.

\begin{lemma}\label{L:root2}
Let $\Phi\subset\Psi$ be two root systems in a Euclidean space $V$. If the subgroup of $V$ generated by roots
in $\Phi$ and that generated by roots in $\Psi$ are equal, then $\Phi\supset\Psi^{\circ}$.
\end{lemma}

Let $L$ be a finitely generated abelian group. We call $L$ a {\it quasi-lattice} if its torsion subgroup \[L_{tor}
:=\{x\in L:\exists n\in\mathbb{Z}_{>0}\textrm{ s.t. } nx=0\}\] is a cyclic group and it is given a positive definite
inner product $(\cdot,\cdot)$ on $L/L_{tor}$. Write \[p:L\rightarrow L/L_{tor}\] for the natural projection. Put
\[(\lambda,\mu)=(p(\lambda),p(\mu)),\ \forall\lambda,\mu\in L.\] Then, it gives a bi-additive form
$(\cdot,\cdot): L\times L\rightarrow\mathbb{R}$. Put \[L^{\check{}}=\Hom(L,\mathbb{Z})\] and write
\[\langle\cdot,\cdot\rangle: L^{\check{}}\times L\rightarrow\mathbb{Z}\] for the natural pairing between $L^{\check{}}$
and $L$. Let $\alpha\in L-L_{tor}$ be an element satisfying that \[\frac{2(\lambda,\alpha)}{(\alpha,\alpha)}\in
\mathbb{Z},\ \forall\lambda\in L.\] Then we define an element $\check{\alpha}\in L^{\check{}}$ by
\[\langle\check{\alpha},\lambda\rangle=\frac{2(\lambda,\alpha)}{(\alpha,\alpha)},\ \forall\lambda\in L\] and define
$s_{\alpha}:L\rightarrow L$ by  \[s_{\alpha}(\lambda)=\lambda-\langle\check{\alpha},\lambda\rangle\alpha=\lambda-
\frac{2(\lambda,\alpha)}{(\alpha,\alpha)}\alpha,\ \forall\lambda\in L.\] We call $s_{\alpha}$ the reflection
associated to $\alpha$.

\begin{definition}\label{D:QRS1}
\begin{itemize}
\item[(i)]Let $L$ be a quasi-lattice with a positive definite inner product $(\cdot,\cdot)$ on $L/L_{tor}$.
We call a finite subset $\Phi$ of $L-L_{tor}$ a quasi root system in $L$ if it satisfies the following
conditions.
\begin{enumerate}
\item[(QRS1)]$\frac{2(\lambda,\alpha)}{(\alpha,\alpha)}\in\mathbb{Z},\ \forall\lambda\in L, \forall\alpha\in\Phi$;
\item[(QRS2)] $s_{\alpha}(\Phi)=\Phi,\ \forall\alpha\in\Phi$.
\end{enumerate}
\item[(ii)]Moreover, we call $\Phi$ a reduced quasi root system if it satisfies one more condition:
\begin{enumerate}
\item[(QRS3)] for any $\alpha\in\Phi$, $2\alpha\not\in\Phi$.
\end{enumerate}
\end{itemize}
\end{definition}

The following lemma is clear.

\begin{lemma}\label{L:root1}
If $\Phi$ is a quasi root system in a quasi-lattice $L$, then $p(\Phi)$ is a root system in the lattice $L/L_{tor}$
in the sense of \cite[Definition 2.2]{Yu-dimension}.
\end{lemma}

Let $\Phi$ be a quasi root system in a quasi-lattice $L$. We call an element $\alpha\in\Phi$ a root and call
$\check{\alpha}$ the coroot associated to $\alpha$; call $p(\Phi)$ the {\it restricted root system} of $\Phi$
and call $p(\alpha)$ the {\it restricted root} associated to a root $\alpha\in\Phi$. Let \[W_{\Phi}:=
\langle s_{\alpha}:\alpha\in\Phi\rangle\subset\Aut(L)\] and call it the {\it Weyl group} of $\Phi$ (acting on $L$).
We call a set of roots $\{\alpha_{i}:1\leq i\leq l\}$ in $\Phi$ a {\it simple system} of $\Phi$ if
$\{p(\alpha_{i}):1\leq i\leq l\}$ is a simple system of $p(\Phi)^{\nd}$; call $l$ the {\it rank} of $\Phi$ and
denote it by $\rank\Phi$; call $\Phi$ {\it irreducible} if $p(\Phi)^{\nd}$ is irreducible.

\begin{definition}\label{D:QRS2}
Let $\Phi_{1}$ (resp. $\Phi_{2}$) be a quasi root system in a quasi-lattice $L_{1}$ (resp. $L_{2}$).
\begin{itemize}
\item[(i)]We say that $(\Phi_{1},L_{1})$ is isomorphic to $(\Phi_{2},L_{2})$ (and write $(\Phi_{1},L_{1})\cong
(\Phi_{2},L_{2})$) if there exists an isomorphism $\phi: L_{1}\rightarrow L_{2}$ of abelian groups such that
$\phi(\Phi_{1})=\Phi_{2}$.
\item[(ii)]We say that $\Phi_{1}$ is isomorphic to $\Phi_{2}$ (and write $\Phi_{1}\cong\Phi_{2}$) if there exist
another quasi root system $\Phi_{3}$ in a quasi-lattice $L_{3}$ and two injective homomorphisms
$\phi_{j}: L_{3}\rightarrow L_{j}$ which map $\Phi_{3}$ onto $\Phi_{j}$ ($j=1,2$).
\end{itemize}
\end{definition}

Clearly, the isomorphism between pairs of quasi root system and quasi lattice is an equivalent relation. By the
following Lemma \ref{L:QRS1} the isomorphism between quasi root systems is also an equivalent relation. Let
$\Phi$ be a quasi root system in a quasi-lattice $L$. Write $L_{\Phi}$ for the subgroup of $L$ generated by
elements in $\Phi$. Then, $L_{\Phi}$ is a quasi-lattice with an inner product on $L_{\Phi}/(L_{\Phi})_{tor}$
inherited from that on $L/L_{tor}$ by restriction.

\begin{lemma}\label{L:QRS1}
Let $\Phi_{1}$ (resp. $\Phi_{2}$) be a quasi root system in a quasi-lattice $L_{1}$ (resp. $L_{2}$). Then $\Phi_{1}
\cong\Phi_{2}$ if and only if there is an isomorphism of abelian groups $\phi: L_{\Phi_{1}}\rightarrow L_{\Phi_{2}}$
such that $\phi(\Phi_{1})=\Phi_{2}$.
\end{lemma}

\begin{proof}
Necessarity. Suppose that $\Phi_{1}\cong\Phi_{2}$. Then, by Definition \ref{D:QRS1} there exist another quasi root
system $\Phi_{3}$ in a quasi-lattice $L_{3}$ and two injective homomorphisms $\phi_{j}: L_{3}\rightarrow L_{j}$ which
map $\Phi_{3}$ onto $\Phi_{j}$ ($j=1,2$). Write $\phi'_{j}=\phi_{j}|_{L_{\Phi_{3}}}$. Then, we have isomorphisms
\[\phi'_{j}: L_{\Phi_{3}}\rightarrow L_{\Phi_{j}}.\] Put $\phi=\phi'_{2}\circ\phi_{1}^{'-1}$. Then, $\phi:
L_{\Phi_{1}}\rightarrow L_{\Phi_{2}}$ is an isomorphism and $\phi(\Phi_{1})=\Phi_{2}$.

Sufficiency. Put $L_{3}=L_{\Phi_{1}}$, $\Phi_{3}=\Phi_{1}$ and let inner product on $L_{3}/(L_{3})_{tor}$ be inherited
from that on $L_{1}/(L_{1})_{tor}$. Let $\phi_{1}: L_{3}\rightarrow L_{1}$ be the inclusion map and let $\phi_{2}:
L_{3}\rightarrow L_{2}$ be the composition of $\phi$ with $\phi_{1}$.
\end{proof}

\begin{lemma}\label{L:QRS2}
Let $\Phi$ be a quasi root system in a quasi-lattice $L$. For two roots $\alpha,\beta\in\Phi$, put
$p=\max\{k\geq 0: \beta-k\alpha\in\Phi\}$ and $q=\max\{k\geq 0: \beta+k\alpha\in\Phi\}$. Then,
\[\frac{2(\beta,\alpha)}{(\alpha,\alpha)}=p-q.\]
\end{lemma}

\begin{proof}
Taking the reflection $s_{\alpha}$ on $\Phi$, we must have $s_{\alpha}(\beta-p\alpha)=\beta+q\alpha$. Then,
$\beta-s_{\alpha}(\beta)=(p-q)\alpha$. Thus, \[\frac{2(\beta,\alpha)}{(\alpha,\alpha)}=p-q.\qedhere\]
\end{proof}

\begin{lemma}\label{L:QRS3}
Let $\Phi_{1}$ be an irreducible quasi root system in a quasi-lattice $L_{1}$ such that $\Phi_{1}$ generates
a finite index subgroup of $L_{1}$. Assume that $\phi: L_{1}\rightarrow L_{2}$ is an isomorphism of abelian
groups such that $\phi(\Phi_{1})=\Phi_{2}$. Then, there exists a positive integer $c$ such that
\[(\phi(x),\phi(y))=c(x,y),\ \forall x,y\in L_{1}.\]
\end{lemma}

\begin{proof}
By Lemma \ref{L:QRS2} and the assumption of this lemma, we have
\[\frac{2(\phi(\beta),\phi(\alpha))}{(\phi(\alpha),\phi(\alpha))}=\frac{2(\beta,\alpha)}{(\alpha,\alpha)},
\ \forall\alpha,\beta\in\Phi_{1}.\] By the irreducibility of $\Phi_{1}$, it follows that there exists a positive
integer $c$ such that \[(\phi(\alpha),\phi(\beta))=c(\alpha,\beta),\ \forall \alpha,\beta\in\Phi_{1}.\] Since
$\Phi_{1}$ generates a finite index subgroup of $L_{1}$, then \[(\phi(x),\phi(y))=c(x,y),
\ \forall x,y\in L_{1}.\qedhere\]
\end{proof}

Using Lemma \ref{L:QRS3} one can show that: if $\phi: (\Phi_{1},L_{1})\rightarrow(\Phi_{2},L_{2})$ is an isomorphism,
then $\phi$ is an isometry up to scalar while restricted to each irreducible factor of $\Phi_{1}$.

Put \[d=|L_{tor}|\] and choose a generator $\beta_{0}$ of $L_{tor}$, which is an element of order $d$.

\begin{lemma}\label{L:ARS1}
Let $\Phi$ be a quasi root system in a quasi-lattice $L$. Then for any root $\alpha\in\Phi$, there
exists $k\geq 1$ with $k|d$ such that \[\{\beta\in\Phi:p(\beta)=p(\alpha)\}=\{\alpha+jk\beta_{0}:
0\leq j\leq\frac{d}{k}-1\}.\]
\end{lemma}

\begin{proof}
Let $\beta=\alpha+j\beta_0\in\Phi$ and $\beta'=\alpha+j'\beta_0\in\Phi$ ($0\leq j,j'\leq d-1$). By calculation
one shows that $s_{\beta}s_{\alpha}(\beta')=\alpha+(j'+2j)\beta_0$. This implies the conclusion of the lemma.
\end{proof}

As in Lemma \ref{L:ARS1}, the number $\frac{d}{k}$ depends only on the restricted root $p(\alpha)$. We call it the
{\it folding index} of $p(\alpha)$ (and of $\alpha$), and denote it by $m'_{p(\alpha)}$ (and $m'_{\alpha}$).
Apparently, the folding index $m'_{p(\alpha)}$ takes the same value on each $W_{p(\Phi)}$ orbit in $p(\Phi)$.

\begin{lemma}\label{L:ARS2}
Let $\Phi$ be a reduced quasi root system in a quasi-lattice $L$ and $\alpha\in\Phi$ be a root such that $p(\alpha)$
is a multipliable root in $p(\Phi)$. Then there exists $k\geq 1$ with $2k|d$ such that \begin{eqnarray*}&&
\{\beta\in\Phi:p(\beta)\in\mathbb{Z}_{>0}p(\alpha)\}\\&=&\{\alpha+jk\beta_{0}:0\leq j\leq \frac{d}{k}-1\}\sqcup
\{2\alpha+(2j+1)k\beta_{0}:0\leq j\leq\frac{d}{2k}-1\}.\end{eqnarray*}
\end{lemma}

\begin{proof}
By Lemma \ref{L:ARS1}, we may take $k\geq 1$ with $k|d$ such that \[\{\beta\in\Phi:p(\beta)=p(\alpha)\}=\{\alpha+
jk\beta_{0}: 0\leq j\leq\frac{d}{k}-1\}.\] Let $\beta\in\Phi$ be a root such that $p(\beta)=tp(\alpha)$
($t\in\mathbb{Z}_{>1}$). By the condition (QRS1), we have \[\frac{2}{t}=\frac{2(\alpha,\beta)}{(\beta,\beta)}\in
\mathbb{Z}.\] Then, $t=2$ due to $t>1$. Write $\beta=2\alpha+i\beta_{0}$ ($1\leq i\leq d-1$). By calculation one
shows $s_{\alpha}s_{\beta}(\alpha)=\alpha-i\beta_0$. Then, $k|i$. Moreover, by the condition (QRS3) we get that:
$\frac{i}{k}$ and $\frac{d+i}{k}$ are both odd. Then, $2k|d$ and $\beta\in
\{2\alpha+(2j+1)k\beta_{0}:0\leq j\leq\frac{d}{2k}-1\}$. Assume that $2\alpha+(2j+1)k\beta_{0}\in\Phi$. Then, \[s_{2\alpha+(2j+1)k\beta_{0}}s_{\alpha+j\beta_{0}}(2\alpha+(2j+1)k\beta_{0})=2\alpha+(2j+3)k\beta_{0}\in\Phi.\]
By iteration one shows the conclusion of the lemma.
\end{proof}

One shows the following lemma by a similar argument as for Lemma \ref{L:ARS2}.

\begin{lemma}\label{L:ARS3}
Let $\Phi$ be a quasi root system in a quasi-lattice $L$ and $\alpha\in\Phi$ be a root such that $2\alpha\in\Phi$.
Then, the set $\{\beta\in\Phi:p(\beta)\in\mathbb{Z}_{>0}p(\alpha)\}$ is equal to one of the following: \begin{enumerate}
\item[(1)] $\{\pm{}\alpha+jk\beta_{0}:0\leq j\leq \frac{d}{k}-1\}\sqcup\{\pm{}2\alpha+jk\beta_{0}:0\leq j\leq\frac{d}{k}-1\}$
where $k\geq 1$ and $k|d$.
\item[(2)] $\{\pm{}\alpha+jk\beta_{0}:0\leq j\leq \frac{d}{k}-1\}\sqcup\{\pm{}2\alpha+2jk\beta_{0}:0\leq j\leq\frac{d}{2k}-1\}$
where $k\geq 1$ and $2k|d$.
\item[(3)] $\{\pm{}\alpha+jk\beta_{0}:0\leq j\leq \frac{d}{k}-1\}\sqcup\{\pm{}2\alpha+4jk\beta_{0}:0\leq j\leq\frac{d}{4k}-1\}$
where $k\geq 1$ and $4k|d$.
\end{enumerate}
\end{lemma}

\begin{lemma}\label{L:ARS0}
Let $\Phi\subset\Psi$ be quasi root systems in a quasi-lattice $L$. Assume that the folding index of each
multipliable root of $p(\Psi)$ is not a multiple of 4, $\Phi$ is reduced, and $p(\Phi)\supset p(\Psi)^{\circ}$.
Then $\Phi$ does not contain any root $\beta\in\Psi$ such that $\beta=2\alpha$ for some root $\alpha\in\Psi$.
\end{lemma}

\begin{proof}
Suppose that $\Phi$ contains a root $\beta\in\Psi$ such that $\beta=2\alpha$ for some root $\alpha\in\Psi$. Since
$p(\Phi)\supset p(\Psi)^{\circ}$, there exists a root $\alpha'\in\Phi$ such that $p(\alpha')=p(\alpha)$. Put $r=
m'_{\alpha}$. By assumption $r$ is not a multiple of 4. By Lemma \ref{L:ARS1}, there exists an order $r$ element
$\gamma_{0}\in L_{tor}$ such that \[\{\gamma\in\Psi:p(\gamma)=p(\alpha)\}=\{\alpha'+j\gamma_{0}:0\leq j\leq r-1\}.\]
By Lemma \ref{L:ARS2}, there exists a positive integer $k$ with $2k|r$ such that \begin{eqnarray*}&&\{\gamma\in\Phi:
p(\gamma)\in\mathbb{Z}_{>0}p(\alpha)\}\\&=&\{\alpha'+jk\gamma_{0}:0\leq j\leq\frac{r}{k}-1\}\sqcup\{2\alpha'+
(2j+1)k\gamma_{0}:0\leq j\leq\frac{r}{2k}-1\}.\end{eqnarray*} Since $r$ is not a multiple of 4, both $k$ and
$\frac{r}{2k}$ are odd. Write $\alpha=\alpha'+i\gamma_{0}$ where $0\leq i\leq r-1$. Then, there exists an integer
$j$ with $0\leq j\leq\frac{r}{2k}-1$ such that \[2\alpha'+2i\gamma_{0}=2\alpha=\beta=2\alpha'+(2j+1)k\gamma_{0}.\]
Then, $2i=(2j+1)k$ or $(2j+1)k+r$. However, $2i$ is even and $(2j+1)k$ and $(2j+1)k+r$ are both odd, we get a
contradiction.
\end{proof}

\begin{remark}\label{R:ARS1}
The assumption of ``the folding index of each multipliable root of $p(\Psi)$ is not a multiple of 4" in Lemma
\ref{L:ARS0} is necessary. For example, when $d=4$ and $\rank L=1$, we could have
\[\Psi=\{\pm{\alpha},\pm{}(\alpha+\beta_{0}),\pm{}(\alpha+2\beta_{0}),\pm{}(\alpha+3\beta_{0}),
\pm{2\alpha},\pm{}(2\alpha+2\beta_{0})\}\] and \[\Phi=\{\pm{\alpha},\pm{}(\alpha+2\beta_{0}),
\pm{}(2\alpha+2\beta_{0})\}.\] In this case, $\Phi\subset\Psi$ are quasi root systems, $\Phi$ is reduced and
$p(\Phi)\supset p(\Psi)^{\circ}$. However, $\alpha+\beta_{0}\in\Psi$ and $2(\alpha+\beta_{0})\in\Phi$.
\end{remark}

Let $\Phi$ be an irreducible quasi root system in a quasi-lattice $L$. We call the folding index of a long root
(resp. a short root) in $p(\Phi)$ the {\it ordinary folding index} (resp. the {\it total folding index}) of
$\Phi$; call the folding index of an ndm root in $p(\Phi)$ (exists only when $p(\Phi)$ is non-reduced) the
{\it ndm folding index} of $\Phi$. We call the ratio between the total folding index (resp. the ndm folding
index) and the ordinary folding index of $\Phi$ the {\it twisted number} (resp. the {\it ndm twisted number})
of $\Phi$. By Lemma \ref{L:fold} below, the twisted number of an irreducible quasi root system takes values
in $\{1,2,3,4\}$.

\begin{lemma}\label{L:fold}
Let $\Phi$ be an irreducible quasi root system in a quasi-lattice $L$ and let $m_0$ be the total folding index
of $\Phi$.
\begin{enumerate}
\item[(i)]when $p(\Phi)$ is of type $\B_{n}$ ($n\geq 2$), $\C_{n}$ ($n\geq 2$) or $\F_{4}$, the ordinary folding
index of $\Phi$ is equal to $m_0$ or $\frac{m_{0}}{2}$.
\item[(ii)]when $p(\Phi)$ is of type $\G_{2}$, the ordinary folding index of $\Phi$ is equal to $m_0$ or
$\frac{m_{0}}{3}$.
\item[(iii)]when $p(\Phi)$ is of type $\BC_{1}$ and $\Phi$ is reduced, the ordinary folding index of $\Phi$ is equal
to $\frac{m_{0}}{2}$.
\item[(iv)]when $p(\Phi)$ is of type $\BC_{n}$ ($n\geq 2$) and $\Phi$ is reduced, the ordinary folding
index of $\Phi$ is equal to $\frac{m_{0}}{2}$ and the ndm folding index of $\Phi$ is equal to $m_0$.
\item[(v)]when $p(\Phi)$ is of type $\BC_{1}$ and $\Phi$ is non-reduced, the ordinary folding index of $\Phi$ is equal
to $m_0$, $\frac{m_{0}}{2}$ or $\frac{m_{0}}{4}$.
\item[(vi)]when $p(\Phi)$ is of type $\BC_{n}$ ($n\geq 2$) and $\Phi$ is non-reduced, let $m_{1}$ (resp. $m_{2}$) be
the ndm folding index (resp. ordinary folding index) of $\Phi$. Then the pair $(m_{1},m_{2})$ is equal to one of
$(m_{0},m_{0})$, $(m_{0},\frac{m_{0}}{2})$, $(\frac{m_{0}}{2},\frac{m_{0}}{2})$, $(\frac{m_{0}}{2},\frac{m_{0}}{4})$.
\end{enumerate}
\end{lemma}

\begin{proof}
We show when $p(\Phi)\cong\B_{n}$ ($n\geq 2$). The proof for other types is similar. Assume that $p(\Phi)\cong\B_{n}$.
Take liftings of simple roots in a simple system of $p(\Phi)$ and take the action of associated reflections, we may
assume that $\Phi$ contains $\{\pm{e_{i}}\pm{e_{j}},\pm{e_{k}}:1\leq i<j\leq n, 1\leq k\leq n\}$ for some elements
$e_{1},\dots,e_{n}\in L$ which are orthogonal to each other and have the same length. By Lemma \ref{L:ARS1} there is
an element $\gamma_{0}\in L_{tor}$ such that \[\{\beta\in\Phi: p(\beta)=p(e_{1})\}=\{e_{1}+j\gamma_{0}:0\leq j\leq m_{0}
-1\}.\] Suppose that $e_{1}+e_{2}+\beta'\in\Phi$ where $\beta'\in L_{tor}$. Then, \[s_{e_{1}}s_{e_{1}+e_{2}+\beta'}(e_{2})
=e_{1}-\beta'\] is a root. Thus, $\beta'=j\gamma_{0}$ for some $j\in\{0,\dots,m_{0}-1\}$. On the other hand, suppose
that $e_{1}+e_{2}+j\beta\in\Phi$. Then, \[e_{1}+e_{2}+(j+2)\beta=s_{e_{1}+\beta}s_{e_{1}}(e_{1}+e_{2}+j\beta)\in\Phi.\]
Thus, the folding index of $p(e_{1}+e_{2})$ is equal to $m_0$ or $\frac{m_{0}}{2}$.
\end{proof}

Choose a section $q: L/L_{tor}\rightarrow L$ of the projection $p: L\rightarrow L/L_{tor}$. Let $\{\bar{\alpha}_{i}:
1\leq i\leq l\}$ be a simple system of $p(\Phi)^{\nd}$. For each $i$, write $\alpha'_{i}=q(\bar{\alpha}_{i})$
and $m_{i}=m'_{\bar{\alpha}_{i}}$. Write a root $\alpha_{i}\in\Phi$ with $p(\alpha_{i})=\bar{\alpha}_{i}$ as
the form $\alpha_{i}=\alpha'_{i}+k_{i}\beta_0$ where $k_{i}\in\mathbb{Z}$. By Lemma \ref{L:ARS1}, the fractional
numbers $\frac{k_{i}m_{i}}{m}\pmod{1}$ do not depend on the choice of $\alpha_{i}$. We call the tuple
\[\{\frac{k_{i}m_{i}}{m}\pmod{1}:1\leq i\leq l\}\] the {\it fractional factor} of $\Phi$ relative to the section
$q$. The following lemma follows from Lemmas \ref{L:ARS1}-\ref{L:ARS3}.

\begin{lemma}\label{L:ARS5}
Let $L$ be a given quasi-lattice. Fix a root system $\bar{\Phi}$ in the lattice $L/L_{tor}$ and a section $q: L/L_{tor}
\rightarrow L$. A quasi root system $\Phi$ in $L$ with $p(\Phi)=\bar{\Phi}$ is determined by the following data: the
fractional factor of $\Phi$; the ordinary folding indices, twisted numbers, ndm twisted numbers and the
reducedness/non-reducedness of irreducible factors of $\Phi$.
\end{lemma}

We introduce some notations for quasi root systems. Let $\Phi$ be an irreducible quasi root system with restricted
root system $\bar{\Phi}$. \begin{enumerate}
\item[(i)](1)We write $\Phi\cong\bar{\Phi}$ if the total folding index of $\Phi$ is equal to 1. (2) Write $\Phi\cong
\bar{\Phi}^{(2)}$ if the ordinary folding index and the total folding index of $\Phi$ are both equal to 2.
\item[(ii)]When $\bar{\Phi}\cong\B_{n}$ ($n\geq 2$), we write $\Phi\cong\B_{n}^{(2,1)}$ if the ordinary folding
index of $\Phi$ is equal to 1 and the total folding index of $\Phi$ is equal to 2.
\item[(iii)]When $\bar{\Phi}\cong\C_{n}$ ($n\geq 2$), we write $\Phi\cong\C_{n}^{(2,1)}$ if the ordinary folding
index of $\Phi$ is equal to 1 and the total folding index of $\Phi$ is equal to 2.
\item[(iv)]When $\bar{\Phi}\cong\F_{4}$, we write $\Phi\cong\F_{4}^{(2,1)}$ if the ordinary folding index of
$\Phi$ is equal to 1 and the total folding index of $\Phi$ is equal to 2.
\item[(v)]When $\bar{\Phi}\cong\BC_{1}$. (1)We write $\Phi\cong\BC_{1}^{'(2,1)}$ if $\Phi$ is non-reduced,
the ordinary folding index of $\Phi$ is equal to 1 and the total folding index of $\Phi$ is equal to 2.
(2)Write $\Phi\cong\BC_{1}^{(2,1)}$ if $\Phi$ is reduced, the ordinary folding index of $\Phi$ is equal to 1
and the total folding index of $\Phi$ is equal to 2.
\item[(vi)]When $\bar{\Phi}\cong\BC_{n}$ ($n\geq 2$). (1)We write $\Phi\cong\BC_{n}^{(2,1,1)}$ if the ordinary
and ndm folding indices of $\Phi$ are equal to 1 and the total folding index of $\Phi$ is equal to 2. (2)Write
$\Phi\cong\BC_{n}^{'(2,2,1)}$ if $\Phi$ is non-reduced, the ordinary folding index of $\Phi$ is equal to 1 and
the ndm and total folding indices of $\Phi$ are equal to 2. (3)Write $\Phi\cong\BC_{n}^{(2,2,1)}$ if $\Phi$ is
reduced, the ordinary folding index of $\Phi$ is equal to 1 and the ndm and total folding indices of $\Phi$ are
equal to 2.
\end{enumerate}
By Lemma \ref{L:fold} and Lemma \ref{L:ARS5}, these represent all isomorphism classes of irreducible quasi root
systems with total folding index equal to 1 or 2. We remark that: if $\bar{\Phi}\cong\BC_{n}$ ($n\geq 2$) and
it is reduced, then we must have $\Phi\cong \BC_{n}^{(2,2,1)}$; if $\bar{\Phi}\cong\BC_{1}$ and it is reduced,
then we must have $\Phi\cong\BC_{1}^{(2,1)}$.

There is strong analogy between quasi root system and affine root system. More concrete connections between them
are explained in the appendix written by Professor Jiu-kang Yu. Affine root systems have a vast of applications
in a lot of directions in mathematics. We are curious if quasi root systems have other applications besides that have
been exploited in this paper.

\subsection{Associating a character to a quasi root system}\label{SS:qrs-char}

We call a compact abelian Lie group $\tilde{S}$ a {\it quasi torus} if its component group $\tilde{S}/\tilde{S}^{0}$
is a cyclic group. Choose and fix a positive definite inner product on the Lie algebra $\mathfrak{s}_{0}=\Lie\tilde{S}$,
denoted by $(\cdot,\cdot)'$. Then, it induces an isomorphism $\phi:\mathfrak{s}_{0}\rightarrow\mathfrak{s}_{0}^{\ast}$
determined by \[\phi(X)(Y)=(X,Y)'\ (\forall X,Y\in\mathfrak{s}_{0}).\] Define a positive definite inner product
$(\cdot,\cdot)$ on $\mathbf{i}\mathfrak{s}_{0}^{\ast}$ by \[(x,y)=(\phi^{-1}(\mathbf{i}x),\phi^{-1}(\mathbf{i}y))'
\ (\forall x,y\in\mathbf{i}\mathfrak{s}_{0}^{\ast}).\] Note that $\opd X^{\ast}(\tilde{S}^{0})$ is a lattice of
$\mathbf{i}\mathfrak{s}_{0}^{\ast}$. For two characters $\lambda,\mu\in X^{\ast}(\tilde{S}^{0})$ (or $\lambda,\mu\in
X^{\ast}(\tilde{S})$), define \[(\lambda,\mu)=(\opd\lambda,\opd\mu).\] Then, $(\cdot,\cdot)$ gives a positive definite
inner product on the lattice $X^{\ast}(\tilde{S}^{0})=X^{\ast}(\tilde{S})/X^{\ast}(\tilde{S})_{tor}$ and makes
$X^{\ast}(\tilde{S})$ a quasi-lattice. We call a quasi root system in $X^{\ast}(\tilde{S})$ a {\it quasi root system
on $\tilde{S}$}. Let $\Phi$ be a quasi root system on $\tilde{S}$. For a root $\alpha\in\Phi$, the induced action of
the refection $s_{\alpha}$ on $\tilde{S}$ is given by \[s_{\alpha}(x)=x\check{\alpha}(\alpha(x)^{-1}),\ \forall x\in
\tilde{S},\] where the coroot $\check{\alpha}$ is regarded as an element in \[X_{\ast}(\tilde{S})\cong
\Hom(X^{\ast}(\tilde{S}),\mathbb{Z}).\] Write $W(\Phi)\subset\Aut(\tilde{S})$ for the group generated by reflections
$s_{\alpha}$ ($\alpha\in\Phi$) and call it the Weyl group of $\Phi$ (acting on $\tilde{S}$). Put \[d=
|\tilde{S}/\tilde{S}^{0}|.\] Then, $|X^{\ast}(\tilde{S})_{tor}|=d$. Write $p: X^{\ast}(\tilde{S})
\rightarrow X^{\ast}(\tilde{S}^{0})$ for the projection \[p(\lambda)=\lambda|_{\tilde{S}^{0}},\ \forall\lambda
\in X^{\ast}(\tilde{S}).\] In \cite[Definition 3.1]{Yu-dimension}, we have defined \[\Psi_{\tilde{S}^{0}}=
\{0\neq\alpha\in X^{\ast}(\tilde{S}^{0}):\frac{2(\lambda,\alpha)}{(\alpha,\alpha)}\in\mathbb{Z},\forall\lambda
\in X^{\ast}(\tilde{S}^{0})\},\] which is a root system in the lattice $X^{\ast}(\tilde{S}^{0})$
(\cite[Definition 2.2]{Yu-dimension}). Define \[\Psi''_{\tilde{S}}=\{\alpha\in X^{\ast}(\tilde{S}):
\alpha|_{\tilde{S}^{0}}\in\Psi_{\tilde{S}^{0}}\}.\] Then, $\Psi''_{\tilde{S}}$ is a quasi root system on
$\tilde{S}$ and it contains all quasi root systems on $\tilde{S}$.

\begin{remark}\label{R:Psi}
Different with that in \cite{Yu-dimension}, to save notations we use the notation $\Psi''_{\tilde{S}}$ instead
of $\Psi_{\tilde{S}}$. The notation $\Psi_{\tilde{S}}$ is reserved for another quasi root system which is defined
while $\tilde{S}$ is a quasi torus in a fixed compact Lie group.
\end{remark}

Fix a connected component $S$ of $\tilde{S}$ generating $\tilde{S}/\tilde{S}^{0}$.

\begin{definition}\label{D:ARS4}
We call an element $s\in S$ a pinned element of $\Phi$ if \[\{p(\alpha):\alpha\in\Phi, \alpha(s)=1\}\] is equal to
$p(\Phi)^{\nd}$.
\end{definition}

If $\{\alpha_{i}:1\leq i\leq l\}$ is a simple system of $\Phi$, then any point $s\in S$ such that \[\alpha_{1}(s)=
\cdots=\alpha_{l}(s)=1\] is a pinned element of $\Phi$. Choose a pinned element $s\in S$ of $\Phi$. Put
\[\Phi_{s}=\{\alpha\in\Phi:\alpha(s)=1\}.\] For a restricted root $\bar{\alpha}\in p(\Phi)^{\nd}$, write \[R_{\bar{\alpha}}=\{\alpha\in\Phi:\alpha|_{\tilde{S}^{0}}\in\mathbb{Z}_{>0}\cdot\bar{\alpha}\}.\] Put
\[m_{\bar{\alpha}}=m'_{\bar{\alpha}}+2m'_{2\bar{\alpha}}.\] For a root $\alpha\in\Phi$ with
$\alpha|_{\tilde{S}^{0}}\in p(\Phi)^{\nd}$, put \[R_{\alpha}=R_{\alpha|_{\tilde{S}^{0}}}\textrm{ and }
m_{\alpha}=m_{\alpha|_{\tilde{S}^{0}}}.\]

The following Lemma follows from Lemmas \ref{L:ARS1}-\ref{L:ARS2} directly.

\begin{lemma}\label{L:ARS6}
Let $\Phi$ be a reduced quasi root system on $\tilde{S}$. For any restricted root $\bar{\alpha}\in\Phi^{\nd}$,
we have $m'_{\bar{\alpha}}|d$ and $m_{\bar{\alpha}}=m'_{\bar{\alpha}}$ or $2m'_{\bar{\alpha}}$; for any root
$\alpha\in\Phi$ with $\alpha|_{\tilde{S}^{0}}\in p(\Phi)^{\nd}$, the character $m_{\alpha}\alpha$ depends
only on $\alpha|_{\tilde{S}^{0}}$.
\end{lemma}

Choose and fix a positive system $\Psi_{\tilde{S}^{0}}^{+}$ of $\Psi_{\tilde{S}^{0}}$. Set \[ p(\Phi)^{\nd,+}=
 p(\Phi)^{\nd}\cap\Psi_{\tilde{S}^{0}}^{+},\] \[\delta_{s,\Phi}=\frac{1}{2}\sum_{\alpha\in\Phi_{s},
\alpha|_{\tilde{S}^{0}}\in p(\Phi)^{\nd,+}}m_{\alpha}\alpha,\] $$A_{\Phi}=\sum_{w\in W_{\Phi_{s}}}\epsilon(w)
[\delta_{s,\Phi}-w\delta_{s,\Phi}].$$ For a subgroup $W$ of $\Aut(\tilde{S})$ containing $W_{\Phi}$, set
\[F_{\Phi,W}=\frac{1}{|W|}\sum_{\gamma\in W}\gamma(A_{\Phi}).\] More generally for any $\mu\in X^{\ast}(\tilde{S})$,
set $$A_{\Phi,\mu}=\sum_{w\in W_{\Phi_{s}}}\epsilon(w)[\mu+\delta_{s,\Phi}-w\delta_{s,\Phi}]$$ and \[F_{\Phi,\mu,W}=
\frac{1}{|W|}\sum_{\gamma\in W}\gamma(A_{\Phi,\mu}).\] For a weight $\lambda\in X^{\ast}(\tilde{S})$ and a finite
subgroup $W$ of $\Aut(\tilde{S})$, set \[\chi^{\ast}_{\lambda,W}=\frac{1}{|W|}\sum_{\gamma\in W}[\gamma\mu]\] and call
it the W-averaging of $\mu$. The following lemma follows Lemma \ref{L:ARS6}.

\begin{lemma}\label{L:AR1}
Let $\Phi$ be a reduced quasi root system on $\tilde{S}$. Then,
\begin{enumerate}
\item[(i)] the characters $2\delta_{s,\Phi}$ and $\delta_{s,\Phi}-w\delta_{s,\Phi}$ ($w\in W_{\Phi_{s}}$) are in
$X^{\ast}(\tilde{S})$. They are independent of the choice of the pinned element $s$, and so is $A_{\Phi}\in
\mathbb{Q}[X^{\ast}(\tilde{S})]$.
\item[(ii)] Each character $F_{\Phi,W}$ (or $F_{\Phi,\mu,W}$) is a finite linear combination of
$\{\chi^{\ast}_{\lambda,W}:\lambda\in X^{\ast}(\tilde{S})\}$ with integral coefficients.
\end{enumerate}
\end{lemma}

The following lemma follows from Lemmas \ref{L:ARS1}-\ref{L:ARS2}.

\begin{lemma}\label{L:AR2}
Let $\Phi$ be a reduced quasi root system on $\tilde{S}$ and $\alpha\in\Phi$ be a root such that
$\alpha|_{\tilde{S}^{0}}\in\Phi^{\nd}$. Then, \begin{equation}\label{Eq:density3}
\prod_{\beta\in R_{\alpha}}(1-\beta(x))=1-\alpha(x)^{m_{\alpha}},\ \forall x\in S.
\end{equation}
\end{lemma}

\begin{lemma}\label{L:AR3}
Let $\Phi$ be a reduced quasi root system on $\tilde{S}$. Then, \[\frac{1}{|W_{\Phi_{s}}|}\prod_{\alpha\in\Phi}(1-\alpha(x))
=F_{\Phi,W_{\Phi_{s}}}(x),\ \forall x\in S.\]
\end{lemma}

\begin{proof}
By Lemma \ref{L:AR2}, we have \[\prod_{\alpha\in\Phi}(1-\alpha(x))=\prod_{\alpha\in\Phi_{s}}(1-\alpha(x)^{m_{\alpha}}).\]
Since $\{m_{\alpha}\alpha:\alpha\in R_{s}\}$ forms an abstract root system by Lemma \ref{L:fold}, the calculation
in the proof of \cite[Proposition 3.7]{Yu-dimension} shows that \[\frac{1}{|W_{\Phi_{s}}|}\prod_{\alpha\in\Phi_{s}}
(1-\alpha(x)^{m_{\alpha}})=F_{\Phi,W_{\Phi_{s}}}(x).\qedhere\]
\end{proof}

\subsection{Generalized Cartan subgroup and quasi root system}\label{SS:MCC}

As defined in \cite[Definition IV.4.1]{BtD}, a closed abelian subgroup $\tilde{S}$ of a compact Lie group $G$ is
called a {\it generalized Cartan subgroup} of $G$ if $\tilde{S}$ has a dense cyclic subgroup and $\tilde{S}^{0}=
Z_{G}(\tilde{S})^{0}$. We call a closed commutative connected subset $S$ of a compact Lie group $G$ a {\it maximal
commutative connected subset} if $s^{-1}S=Z_{G}(S)^{0}$ for any $s\in S$. For an element $g\in G$, choose a maximal
torus $T$ of $Z_{G}(g)^{0}$. Put $\tilde{S}=\langle g,T\rangle$ and $S=gT$. Then, $\tilde{S}$ is a generalized
Cartan subgroup of $G$ and $S$ is a maximal commutative connected subset of $G$. Moreover, all generalized Cartan
subgroups (resp. maximal commutative connected subsets) of $G$ arise in this way. Let $S$ be a maximal commutative
connected subset of $G$. Set \[W(G,S)=N_{G^{0}}(S)/s^{-1}S,\ s\in S\] and we call $W(G,S)$ the {\it Weyl group of
$G$ with respect to $S$}.

Let $\tilde{S}$ be a generalized Cartan subgroup of $G$. Then $\tilde{S}^{0}\subset\tilde{S}\cap G^{0}\subset
Z_{G^{0}}(\tilde{S}).$ As the following example indicates, we have $\tilde{S}^{0}\neq\tilde{S}\cap G^{0}$ and
$\tilde{S}\cap G^{0}\neq Z_{G^{0}}(\tilde{S})$ in general.

\begin{example}\label{E:A3}
Write $L=\left(\begin{array}{cccc}0&0&1&0\\0&0&0&1\\1&0&0&0\\0&1&0&0\end{array}\\\right)$ and
\[T=\{\diag\{\lambda_{1},\lambda_{2},\lambda_{1}^{-1},\lambda_{2}^{-1}\}:\lambda_{1},\lambda_{2}\in
\mathbb{C}^{\times},|\lambda_{1}|=|\lambda_{2}|=1\}.\] Set $G_{1}=(\SU(4)/\langle-I\rangle)\rtimes\langle
\sigma_{1}\rangle,$ where $\sigma_{1}^2=[\mathbf{i}I]$ and \[\sigma_{1}[X]\sigma_{1}^{-1}=[L\overline{X}L^{-1}],
\ \forall X\in\SU(4).\] Put $\tilde{S}_{1}=\langle T,\sigma_{1}\rangle.$ Then, $\tilde{S}_{1}$ is a generalized
Cartan subgroup of $G_{1}$. In this case, $[\mathbf{i}I]\in\tilde{S}_{1}\cap G_{1}^{0}-\tilde{S}_{1}^{0}.$
Thus, $\tilde{S}_{1}^{0}\neq\tilde{S}_{1}\cap G_{1}^{0}$.

Set $G_{2}=(\SU(4)/\langle-I\rangle)\rtimes\langle\sigma_{2}\rangle,$ where $\sigma_{2}^2=1$ and
\[\sigma_{2}[X]\sigma_{2}^{-1}=[L\overline{X}L^{-1}],\ \forall X\in\SU(4).\] Put $\tilde{S}_{2}=
\langle T,\sigma_{2}\rangle.$ Then, $\tilde{S}_{2}$ is a quasi Cartan subgroup of $G_{2}$. In this case,
$[\mathbf{i}I]\in Z_{G_{2}^{0}}(\tilde{S}_{2})-\tilde{S}_{2}\cap G_{2}^{0}.$ Thus, $\tilde{S}_{2}\cap G_{2}^{0}
\neq Z_{G_{2}^{0}}(\tilde{S}_{2})$.
\end{example}

Lemmas \ref{L:compact-conjugacy}-\ref{L:integration} below are from \cite[Chapter IV]{BtD} and \cite{Wendt}.

\begin{lemma}[\cite{BtD}]\label{L:compact-conjugacy}
We have the following properties of maximal commutative connected subsets.
\begin{enumerate}
\item[(i)]Any two maximal commutative connected subsets in a connected component $g'G^{0}$ of $G$ are $G^{0}$ conjugate.
\item[(ii)]If $S$ is a maximal commutative connected subset in $g'G^{0}$, then every $G^{0}$ conjugacy class in $g'G^{0}$
intersects with $S$.
\item[(iii)]Two elements in $S$ are in a $G^{0}$ conjugacy class if and only if they are in the same $W(G,S)$ orbit.
\end{enumerate}
\end{lemma}

\begin{lemma}[\cite{BtD},\cite{Wendt}]\label{L:S-q}
The map $$q: G^{0}/\tilde{S}^{0}\times S\rightarrow g'G^{0},\quad (g\tilde{S}^{0},s)\mapsto gsg^{-1}$$ is a surjective
finite degree map with degree equal to $|W(G,S)|$.
\end{lemma}

Write $\opd g$ for a right $G^{0}$ invariant measure of volume $1$ on a connected component $g'G^{0}$, $\opd\bar{g}$ for a
$G^{0}$ invariant measure on $G^{0}/\tilde{S}^{0}$ of volume $1$, and $\opd s$ for a left $\tilde{S}^{0}$ invariant
measure on $S$ of volume $1$.

\begin{lemma}(\cite[Lemma 2.1]{Wendt})\label{L:det}
We have $q^{\ast}\opd g=\det(q)\opd\bar{g}\wedge\opd s$ and \[\det(q)(g\tilde{S}^{0},s)=
|\det(1-\Ad(s))|_{\frg/\frs}|,\ \forall (g,s)\in G\times S.\]
\end{lemma}

Write \begin{equation}\label{Eq:density1}D_{G}(gsg^{-1})=\frac{1}{|W(G,S)|}|\det(1-\Ad(s))|_{\frg/\frs}|,
\ \forall (g,s)\in G^{0}\times S\end{equation} can call the {\it Weyl density function} on the connected component
$g'G^{0}$. By Lemma \ref{L:compact-conjugacy}, $D_{G}(\cdot)$ is a well-defined function on $g'G^{0}$. It is actually
a smooth function.

Choose a maximal commutative connected subset in each connected component of $G$ and write $S_1,S_2,\dots,S_{m}$ for
the chosen subsets. For each $i$, write $D_{i}(s)$ for the Weyl density function on $S_{i}$. Let $\tilde{\opd} g$ be
a normalized Haar measure on $G$.

\begin{lemma}(\cite[Prop. 2.3]{Wendt})\label{L:integration}
For any $G^{0}$ conjugation invariant continuous function $f$ on $g'G^{0}$, we have $$\int_{g'G^{0}}f(g)\opd g=
\int_{S}f(s)D_{G}(s)\opd s.$$ For any $G^{0}$ conjugation invariant continuous function $f$ on $G$, one has
$$\int_{G}f(g)\tilde{\opd} g=\frac{1}{|G/G^{0}|}\sum_{1\leq i\leq m}\int_{S_{i}}f(s)D_{G}(s)\opd s.$$
\end{lemma}

In the following lemma, we compare dimension and minimal order of elements of maximal commutative connected subsets
in different connected components of $G$ while the component group $G/G^0$ is cyclic.

\begin{lemma}\label{L:MCC1}
Assume that $G/G^0$ is a cyclic group. Let $g\in G$ be an element generating $G/G^{0}$ and $S$ be a maximal commutative
connected subset in $gG^{0}$. Let $S'$ be a maximal commutative connected subset in another connected component $g'G^{0}$.
Then we have the following assertions:
\begin{enumerate}
\item[(i)](1)$\dim S'\geq\dim S$; (2)the minimal order of elements of $S'$ is less than or equal to that of $S$.
\item[(ii)]$\dim S'=\dim S$ holds only when $g'$ generates $G/G^{0}Z_{G}(G^{0})$.
\item[(iii)]The minimal order of elements of $S'$ is equal to that of $S$ holds only when $g'$ generates $G/G^{0}$.
\end{enumerate}
\end{lemma}

\begin{proof}
(i)(1)Take a maximal torus $T$ of $Z_{G}(g)^{0}$. We may take $S=gT$. Since $G/G^0$ is cyclic, we may assume that
$g'=g^{k}$ for some integer $k$. Take a maximal torus $T'$ of $Z_{G}(g')^{0}$ containing $T$ and we may assume
that $S'=g'T'$. Then, \[\dim S'=\dim T'\geq\dim T=\dim S.\]

(2)Let $gx$ be a minimal order element in $S$, where $x\in\tilde{S}^{0}$. Write $d$ for the order of $gx$. Then,
$(g'x^{k})^{d}=(gx)^{dk}=1$. Thus, the minimal order of elements of $S'$ $\leq d$.

(ii)Follow the proof of (i)(1). Consider the adjoint homomorphism \[\pi: G\rightarrow\Out(\mathfrak{g}_{0}):=
\Aut(\mathfrak{g}_{0})/\Int(\mathfrak{g}_{0}).\] Apparently, $\dim S$ (resp. $\dim S'$) is determined by $\pi(g)$
(resp. $\pi(g')$). It is clear that the homomorphism $\pi$ factors through $G/G^{0}Z_{G}(G^{0})$ and there is an
injective homomorphism \[G/G^{0}Z_{G}(G^{0})\rightarrow\Out(\mathfrak{g}_{0}).\] Write $(\mathfrak{g}_{0})_{\der}
=[\mathfrak{g}_{0},\mathfrak{g}_{0}]$ for the derived subalgebra of $\mathfrak{g}_{0}$ and write $z(\mathfrak{g}_{0})$
for the center of $\mathfrak{g}_{0}$. Then, \[\Out(\mathfrak{g}_{0})=\Out((\mathfrak{g}_{0})_{\der})\times
\GL(z(\mathfrak{g}_{0})).\] For an automorphism $\sigma$ of $(\mathfrak{g}_{0})_{\der}$, let $s(\sigma)$ be the
rank of $(\mathfrak{g}_{0})_{\der}^{\sigma}$. By Lemma \ref{L:compact-conjugacy}, the number $s(\sigma)$ depends
only on the class $[\sigma]$ of $\sigma$ in $\Out((\mathfrak{g}_{0})_{\der})$. For $\sigma\in\GL(z(\mathfrak{g}_{0}))$,
let $s(\sigma)=\dim z(\mathfrak{g}_{0})^{\sigma}$. It suffices to show the following assertions:
\begin{enumerate}
\item[(1)]For any $\sigma\in\Aut((\mathfrak{g}_{0})_{\der})$ and an integer $k$, $s(\sigma^{k})=s(\sigma)$ if and
only if $k$ is prime to the order of $[\sigma]\in\Out((\mathfrak{g}_{0})_{\der})$.
\item[(2)]For any $\sigma\in\GL(z(\mathfrak{g}_{0}))$ and an integer $k$, $s(\sigma^{k})=s(\sigma)$ if and
only if $k$ is prime to the order of $\sigma$.
\end{enumerate}
(2) follows by considering eigenvalues of $\sigma:z(\mathfrak{g}_{0})_{\mathbb{C}}\rightarrow
z(\mathfrak{g}_{0})_{\mathbb{C}}$. For (1), we may assume that $\sigma$ is a pinned automorphism which acts
transitively on simple factors of $(\mathfrak{g}_{0})_{\der}$. Then, one shows the assertion by a case by
case verification.

(iii)Follow the proof of (i)(2). The order of $(gx)^{k}=g'x^{k}\in S'$ is equal to $\frac{d}{\gcd(k,d)}\leq d$.
If the minimal order of elements of $S'$ is equal to that of $S$, then we have $\gcd(k,d)=1$. Write $e=|G/G^{0}|$.
Due to $(gx)^{d}=1$ and $gx$ generates $G/G^{0}$, then $e$ divides $d$. Thus, $\gcd(k,e)=1$. Hence, $g'=g^{k}$
generates $G/G^{0}$ as $g$ does.
\end{proof}

Endow $G$ with a biinvariant Riemannian metric. Then, it induces a $W(G,S)$ invariant positive definite inner
product on $\mathfrak{s}_0$. As in the beginning of the last subsection, this induces a $W(G,S)$ invariant positive
definite inner product on $\mathbf{i}\mathfrak{s}_{0}^{\ast}$. For two characters $\lambda,\mu\in X^{\ast}(\tilde{S})$
(or $\lambda,\mu\in X^{\ast}(\tilde{S}^0)$), define \[(\lambda,\mu)=(\opd\lambda,\opd\mu).\] The conjugation
action of $\tilde{S}$ on $\mathfrak{g}$ gives a decomposition $$\mathfrak{g}=\mathfrak{s}\oplus
\bigoplus_{\lambda\in X^{\ast}(\tilde{S})-\{0\}}\mathfrak{g}_{\lambda},$$ where \[\frg_{\lambda}=\{X\in\mathfrak{g}:
\Ad(s)X=\lambda(s)X,\ \forall s\in\tilde{S}\}.\] Put $$\tilde{R}(G,\tilde{S})=\{\alpha\in X^{\ast}(\tilde{S})-\{0\}:
\mathfrak{g}_{\alpha}\neq 0\}.$$ An element $\alpha\in\tilde{R}(G,S)$ is called a {\it root} and
$\mathfrak{g}_{\alpha}$ is called the {\it root space}. A root $\alpha\in \tilde{R}(G,\tilde{S})$ is called an
{\it infinite root} if $\alpha|_{\tilde{S}^{0}}\neq 0$; otherwise it is called a {\it finite root}. Write
$R(G,\tilde{S})$ for the set of infinite roots in $\tilde{R}(G,S)$, and write $R_{0}(G,\tilde{S})$ for the set
of finite roots in $\tilde{R}(G,S)$.

\begin{lemma}\label{L:infinite1}
Let $\alpha$ be an infinite root. Then $\dim\frg_{\alpha}=1$, $-\alpha$ is also an infinite root, and
$2\alpha$ is not a root.
\end{lemma}

\begin{proof}
Put $G^{[\alpha]}=Z_{G}(\ker\alpha)^{0}$ and $G_{[\alpha]}=(G^{[\alpha]})_{\der}.$ Write $T_{[\alpha]}=G_{[\alpha]}
\cap\tilde{S}$. Then, the complexified Lie algebra of $G^{[\alpha]}$ is spanned by $\frs$ and $\frg_{k\alpha}$
($k\in\mathbb{Z}-\{0\}$). From this we see that $\tilde{S}=T_{[\alpha]}\cdot\ker\alpha$, $\dim T_{[\alpha]}=1$,
and $T_{[\alpha]}$ is a maximal torus of $G_{[\alpha]}$. Thus, the complexified Lie algebra of $G_{[\alpha]}$ is
of rank one. Hence, it is isomorphic to $\mathfrak{sl}_{2}(\mathbb{C})$. By this, $\dim\frg_{\alpha}=1$, $-\alpha$
is an infinite root, and $2\alpha$ is not a root.
\end{proof}

There is a unique cocharacter $\check{\alpha}\in X_{\ast}(T_{[\alpha]})\subset X_{\ast}(\tilde{S}^{0})$ such
that $\langle\alpha,\check{\alpha}\rangle=2$. Define \begin{equation}\label{Eq:reflection}s_{\alpha}(x)=
x\check{\alpha}(\alpha(x)^{-1}),\ \forall x\in\tilde{S}.\end{equation} Then, $s_{\alpha}|_{\ker\alpha}=\id$,
$s_{\alpha}|_{\Ima\check{\alpha}}=-1$ and $s_{\alpha}^2=1$. Thus, $s_{\alpha}$ is an automorphism of $\tilde{S}$
which stabilizes $S$. Since $G_{[\alpha]}$ is a connected semisimple Lie group of rank one, the following statement
is clear.

\begin{lemma}\label{L:reflection}
There exists $n\in N_{G_{[\alpha]}}(T_{[\alpha]})$ such that $\Ad(n)|_{S}=s_{\alpha}.$
\end{lemma}

\begin{lemma}\label{L:Cartan-ARS}
$R(G,\tilde{S})$ is a reduced quasi root system on $\tilde{S}$.
\end{lemma}

\begin{proof}
Considering $\mathfrak{g}$ as a representation of the rank one connected semisimple Lie group $G_{[\alpha]}$, then the
condition (QRS1) in the Definition \ref{D:QRS1} follows. The condition (QRS2) follows from Lemma \ref{L:reflection}
and the condition (QRS3) follows from Lemma \ref{L:infinite1}.
\end{proof}

Put $d=|\tilde{S}/\tilde{S}^{0}|$. Write $p: X^{\ast}(\tilde{S})\rightarrow X^{\ast}(\tilde{S}^{0})$ for the projection
map and write \[R(G,\tilde{S}^{0}):=p(R(G,\tilde{S}))\] for the restricted root system of $R(G,\tilde{S})$.

\begin{lemma}\label{L:T}
If $\tilde{S}$ is a generalized Cartan subgroup of a compact Lie group $G$, then $Z_{G^{0}}(\tilde{S}^{0})$ is
a maximal torus of $G^{0}$.
\end{lemma}

\begin{proof}
Put $T=Z_{G^{0}}(\tilde{S}^{0})$. Choose an element $y\in\tilde{S}$ generating $\tilde{S}/\tilde{S}^{0}$. Then,
$\mathfrak{t}^{y}=\mathfrak{s}\subset z(\mathfrak{t})$. From the Levi decomposition $\mathfrak{t}=z(\mathfrak{t})
\oplus\mathfrak{t}_{\der}$, it follows that $(\mathfrak{t}_{\der})^{y}=0$. By a theorem of Borel, this indicates
that $\mathfrak{t}_{\der}=0$. Thus, $\mathfrak{t}$ is abelian. We know that $T$ is connected
(\cite[Corollary 4.51]{Knapp}) and it contains all maximal tori of $G^{0}$ containing $\tilde{S}^{0}$. Thus, $T$
is a maximal torus of $G^{0}$.
\end{proof}

\begin{lemma}\label{L:split}
There is an exact sequence $$1\rightarrow N_{T}(S)/\tilde{S}^{0}\rightarrow W(G,S)\rightarrow W_{R(G,\tilde{S}^{0})}
\rightarrow 1,$$ where $T=Z_{G^{0}}(\tilde{S}^{0})$.
\end{lemma}

\begin{proof}
Restricting the action of $W(G,S)$ on $\tilde{S}$ to $\tilde{S}^{0}$, we get a homomorphism
\[\phi: W(G,S)\rightarrow\Aut(\tilde{S}^0).\] First, we show that \[\ker\phi=N_{T}(S)/\tilde{S}^{0}.\] For any
$g\in N_{G^0}(S)$, in order that $g\tilde{S}^{0}\in\ker\phi$ it is necessary and sufficient that
$g\in Z_{G^{0}}(\tilde{S}^{0})=T$. Hence, $\ker\phi=N_{T}(S)/\tilde{S}^{0}$.

Second, we show that $W_{R(G,\tilde{S}^{0})}\subset\Ima\phi$. For any infinite root $\alpha\in R(G,\tilde{S})$, by
Lemma \ref{L:reflection} there exists $n_{\alpha}\in N_{G^0}(S)$ such that $\Ad(n_{\alpha})|_{\tilde{S}}=s_{\alpha}$.
This implies that $s_{\alpha}|_{\tilde{S}^{0}}\in\Ima\phi$. Thus, $W_{R(G,\tilde{S}^{0})}\subset\Ima\phi$.

Third, we show that $\Ima\phi\subset W_{R(G,\tilde{S}^{0})}$. Choose a positive system $R^{+}(G,\tilde{S}^{0})$ of
$R(G,\tilde{S}^{0})$. Write \[X=\frac{\mathbf{i}}{2}\sum_{\bar{\alpha}\in R^{+}(G,\tilde{S}^{0})}\check{\bar{\alpha}}
\in\mathfrak{s}_0.\] Then, $-\mathbf{i}X$ is regular dominant with respect to $R^{+}(G,\tilde{S}^{0})$. Thus,
$Z_{G^{0}}(X)=Z_{G^{0}}(\tilde{S}^{0}).$ Suppose $\phi(g\tilde{S}^{0})\not\in W_{R(G,\tilde{S}^{0})}$ for some
$g\in N_{G^0}(S)$. Multiplying $g$ by a product of $n_{\alpha}$ ($\alpha\in R(G,\tilde{S})$) if necessary, we may
assume that the induced action of $\phi(g\tilde{S}^0)$ on $R(G,\tilde{S}^{0})$ maps $R^{+}(G,\tilde{S}^{0})$ to
itself. Then, $g\in Z_{G^{0}}(X)=Z_{G^{0}}(\tilde{S}^{0})$ and $\phi(g\tilde{S}^{0})=1$, which is a contradiction.
Thus, $\Ima\phi\subset W_{R(G,\tilde{S}^{0})}$.
\end{proof}

The following lemma is elementary and well-known.

\begin{lemma}\label{L:prime-identity}
If $m=p^{k}$ is a prime power, then $$\prod_{1\leq j\leq m-1,\gcd(m,j)=1}(1-e^{\frac{2j\pi i}{m}})=p.$$ If $m$ is not a
prime power, then $$\prod_{1\leq j\leq m-1,\gcd(m,j)=1}(1-e^{\frac{2j\pi i}{m}})=1.$$
\end{lemma}

\begin{lemma}\label{L:FixedPoint}
Let $H$ be a compact torus and $\phi$ be an automorphism of it with only finitely many fixed points. Write
$\phi_{\ast}\in\GL(\mathfrak{h}_{0})$ for the differential of $\phi$ and write $f(t)=\det(tI-\phi_{\ast})$ for
the characteristic polynomial of $\phi_{\ast}$. Then 1 is not a root of $f(t)$ and the number of fixed points of
$\phi$ is equal to $f(1)$.
\end{lemma}

\begin{proof}
Write $L=X_{\ast}(H)$ and $U=L\otimes_{\mathbb{Z}}\mathbb{R}.$ Then, one can identify $\mathfrak{h}_{0}$ with $U$,
and identify $H$ with $U/L$. Since $\phi$ has only finitely many fixed points, $1$ is not an eigenvalue of
$\phi_{\ast}$. Equivalently, 1 is not a root of $f(t)$. Then, $H^{\phi}$ can be identified with
$(I-\phi_{\ast})^{-1}(L)/L$. Comparing volumes of $U/L$ and $U/(I-\phi_{\ast})^{-1}(L)$, we get \[|H^{\phi}|=
|(I-\phi_{\ast})^{-1}(L)/L|=|\det(I-y_{\ast})|=|f(1)|.\] Moreover, $\phi$ must be of finite order and $f(t)$ is
a product of cyclotomic polynomials. Then, $f(1)>0$ by Lemma \ref{L:prime-identity}. Thus, $|H^{\phi}|=f(1)$.
\end{proof}

\begin{lemma}\label{L:finiteRoot}
We have \begin{equation}\label{Eq:density2}\prod_{\alpha\in R_{0}(G,\tilde{S})}(1-\alpha(x))=|N_{T}(S)/\tilde{S}^{0}|,
\ \forall x\in S.\end{equation}
\end{lemma}

\begin{proof}
The group $N_{T}(S)$ consists of elements $t\in T$ such that $txt^{-1}x^{-1}\in\tilde{S}^{0}$. Considering the quotient
group $T/\tilde{S}^{0}$ and the conjugation action by $x$, we have $$N_{T}(S)/\tilde{S}^{0}\cong\{t\tilde{S}^{0}\in
T/\tilde{S}^{0}:x(t\tilde{S}^{0})x^{-1}(t\tilde{S}^{0})^{-1}=\tilde{S}^{0}\}.$$ The latter is just the fixed point group
of $x$ in $T/\tilde{S}^{0}$. The action of $x$ on $T/\tilde{S}^{0}$ is given by a matrix $X\in\GL(k,\mathbb{Z})$, where
$k=\dim(T/\tilde{S}^{0})$. Write $f(u)=\det(uI-X)$ for the characteristic polynomial of $X$. Then, $f(u)=
\prod_{\alpha\in R_{0}(G,\tilde{S})}(u-\alpha(s)).$ By Lemma \ref{L:FixedPoint} we get the conclusion of the lemma.
\end{proof}

Combining Lemmas \ref{L:split}, \ref{L:finiteRoot} and \ref{L:AR3}, the following lemma follows.

\begin{lemma}\label{L:AR4}
We have \[D_{G}(x)=F_{R(G,\tilde{S}),W_{\Phi_{s}}}(x),\ \forall x\in S.\]
\end{lemma}

Recall that an automorphism $\theta$ of a complex semisimple Lie algebra $\mathfrak{h}$ is called a {\it pinned automorphism}
of $\mathfrak{h}$ if $\theta$ stabilizes a Cartan subalgebra $\mathfrak{t}$ and permutes a set of root vectors
$\{X_{\beta_{1}},\dots,X_{\beta_{l}}\}$ for a simple system $\{\beta_{1},\dots,\beta_{l}\}$ of the root system
$R(\mathfrak{h},\mathfrak{t})$.

\begin{lemma}\label{L:basePoint}
In order that an element $s\in S$ is a pinned element for the quasi root system $R(G,\tilde{S})$ it is necessary and
sufficient that $\Ad(s)|_{\mathfrak{g}_{\der}}$ is a pinned automorphism.
\end{lemma}

\begin{proof}
The sufficiency is clear. For the necessarity, assume that $s\in S$ is a pinned element. Put $\theta=
\Ad(s)|_{\mathfrak{g}_{\der}}$. It is shown in Lemma \ref{L:T} that $T=Z_{G^{0}}(\tilde{S}^{0})$ is a maximal torus
of $G^{0}$. Choose a positive system $R^{+}(G,\tilde{S}^{0})$ of $R(G,\tilde{S}^{0})$. Define a root $\alpha\in
R(\mathfrak{g}_{\der},T)$ positive if $\alpha|_{\tilde{S}^{0}}\in R^{+}(G,\tilde{S}^{0})$. In such a way, one
obtains a positive system of $R(\mathfrak{g}_{\der},T)$. Since $s$ commutes with $\tilde{S}^{0}$. Then, $\theta$
stabilizes $\mathfrak{t}$ and the corresponding simple system of $R(\mathfrak{g}_{\der},T)$. Moreover, $s$ is a
pinned element implies that one can normalize root vectors $X_{\beta_{i}}$ ($1\leq i\leq l$) appropriately such that
$\theta$ permutes $\{X_{\beta_{1}},\dots,X_{\beta_{l}}\}$. Hence, $\Ad(s)|_{\mathfrak{g}_{\der}}=\theta$ is a
pinned automorphism.
\end{proof}

Note that some properties of the set $R(G,\tilde{S})$ in special cases have been studied in the literature
(e.g, \cite{Gantmacher}, \cite{deSie} and \cite[Chapter 3]{GOV}).

\subsection{Pinned automorphism and quasi root system}\label{SS:pinned-QRS}

\begin{definition}\label{D:adjoint}
\begin{enumerate}
\item[(i)]We call a quasi root system $\Phi$ on a quasi torus $\tilde{S}$ of simply-connected type if
coroots $\{\check{\alpha}:\alpha\in\Phi\}$ generate the cocharacter group $X_{\ast}(\tilde{S}^{0})$.
\item[(ii)]Call a quasi root system $\Phi$ on a quasi torus $\tilde{S}$ of adjoint type if
\[\{x\in\tilde{S}: \alpha(x)=1,\ \forall\alpha\in\Phi\}=1.\]
\item[(iii)]Call a root system $\Phi'$ on a torus $T$ of adjoint type if \[\{x\in T: \alpha(x)=1,
\ \forall\alpha\in\Phi'\}=1.\]
\end{enumerate}
\end{definition}

Apparently, a quasi root system $\Phi$ on a quasi torus $\tilde{S}$ being of simply-connected type depends only on
the restricted root system $p(\Phi)$ on the torus $\tilde{S}^{0}$. However, it being of adjoint type depends on
the actual quasi root system $\Phi$ and the quasi torus $\tilde{S}$.

Let $\Phi'$ be a root system on a torus $T$ and $\theta$ be an order $m$ automorphism of $T$ which stabilizes a
positive system $\Phi'^{+}$ of $\Phi'$. Put \[\tilde{S}=(T^{\theta})^{0}\times\langle\theta\rangle.\] We construct
a reduced quasi root system $\Phi$ on $\tilde{S}$. There is a connected compact Lie group $G$ with $T$ a maximal
torus and with root system $R(G,T)=\Phi'$. Then, the automorphism $\theta$ extends to a pinned automorphism of
$G$ with respect to $T$ and a positive system $R^{+}(G,T)=\Phi'^{+}$. Form \[\tilde{G}=G\rtimes\langle\theta\rangle.\]
Then, $\tilde{S}$ is a generalized Cartan subgroup of $\tilde{G}$. Let $\Phi=R(\tilde{G},\tilde{S})$. It is a reduced
quasi root system on $\tilde{S}$. We call $\Phi$ {\it a quasi root system resulting from the action of $\theta$
on $\Phi'$}. In the following lemma, we show a converse of this construction by assuming that $\Phi$ is of adjoint
type.

\begin{lemma}\label{L:ARS-isogeny1}
Given a reduced quasi root system $\Phi$ on a quasi torus $\tilde{S}$ of adjoint type, then there exists a unique
triple $(\Phi',T,\theta)$ up to isomorphism with $\Phi'$ a root system on a torus $T$ of adjoint type and $\theta$
an automorphism of $T$ which stabilizes a positive system $\Phi'^{+}$ of $\Phi'$ such that the resulting quasi
root system is isomorphic to $\Phi$.
\end{lemma}

\begin{proof}
For simplicity we assume that $\Phi$ is irreducible and has ordinary folding index 1. We always take an adjoint type
compact Lie group $G$ with a maximal torus $T$ such that $R(G,T)=\Phi$. Let's describe the pair $(\Phi',\theta)$ case
by case.
\begin{enumerate}
\item[(i)]When the twisted number of $\Phi$ is equal to 1, we have $\Phi'\cong p(\Phi)$ and $\theta=1$.
\item[(ii)]When $\Phi=\B_{n}$ ($n\geq 3$) and the twisted number equals 2, we have $\Phi'\cong\D_{n+1}$ and $\theta$
a diagram automorphism of order 2.
\item[(iii)]When $\Phi=\C_{n}$ ($n\geq 2$) and the twisted number equals 2, we have $\Phi'\cong\A_{2n-1}$ and $\theta$
a diagram automorphism of order 2.
\item[(iv)]When $\Phi=\BC_{n}$ ($n\geq 2$) and the twisted number equals 2, we have $\Phi'\cong\A_{2n}$ and $\theta$
a diagram automorphism of order 2.
\item[(v)]When $\Phi=\F_{4}$ and the twisted number equals 2, we have $\Phi'\cong\E_{6}$ and $\theta$
a diagram automorphism of order 2.
\item[(vi)]When $\Phi=\G_{2}$ and the twisted number equals 3, we have $\Phi'\cong\D_{4}$ and $\theta$
a diagram automorphism of order 3.
\end{enumerate}
\end{proof}

\begin{lemma}\label{L:ARS-isogeny2}
Let $G_{1}$ (resp. $G_{2}$) be a compact semisimple Lie group with cyclic component group and $\tilde{S}_{1}$ (resp.
$\tilde{S}_{2}$) be a generalized Cartan subgroup of it which meets all connected components. Write $\Phi_{j}=
R(G_{j},\tilde{S}_{j})$ ($j=1,2$) for the quasi root system obtained from the action of $\tilde{S}_{j}$ on the
complexified Lie algebra of $G_{j}$. Suppose that $\Phi_{1}\cong\Phi_{2}$. Then, there exists a compact Lie group
$K$ and finite surjective homomorphisms $\phi_{j}: G_{j}\rightarrow K$ ($j=1,2$) such that $\phi_{1}(\tilde{S}_{1})
=\phi_{2}(\tilde{S}_{2})$.
\end{lemma}

\begin{proof}
Taking the quotient of $G_{j}$ by $Z_{G_{j}}(G_{j}^{0})$ ($j=1,2$) replaces the quasi root system $\Phi_{j}$ by an
adjoint type quasi root system, which are also isomorphic. For this reason we assume that $Z_{G_{j}}(G_{j}^{0})=1$.
That is to say, $G_{j}$ is of adjoint type. Then, the isomorphism class of $G_{j}$ is determined by its root system
$\Phi'_{j}$ and a pinned automorphism $\theta_{j}\in\tilde{S}_{j}$ generating $G_{j}/G_{j}^{0}$. By Lemma
\ref{L:ARS-isogeny1}, the pair $(\Phi'_{j},\theta_{j})$ is determined up to isomorphism by $\Phi_{j}$. Hence, the
isomorphism class of $G_{j}$ is determined.
\end{proof}

\begin{example}\label{E:ARS-isogeny3}
Take an odd integer $N=2m+1>1$. Let $G_{1}=\SU(N)\rtimes\langle\tau\rangle$ where $\tau^{2}=I$ and \[\tau X\tau^{-1}
=\overline{X},\ \forall X\in\SU(N).\] Put $G_{2}=G_{1}/Z(\SU(N))$. Write $p:G_{1}\rightarrow G_{2}$ for the
quotient map. Let $\tilde{S}_{1}\subset G_{1}$ be a quasi torus generated by $\tau$ and \[\big\{\left(\begin{array}
{cccccc}\cos\theta_{1}&\sin\theta_{1}&&&&\\-\sin\theta_{1}&\cos\theta_{1}&&&&\\&&\ddots&&&\\&&&\cos\theta_{m}&
\sin\theta_{m}&\\&&&-\sin\theta_{m}&\cos\theta_{m}&\\&&&&&1\end{array}\right):\theta_{j}\in\mathbb{R}\big\}\] and
let $\tilde{S}_{2}=p(\tilde{S}_{1})$. Then, $\tilde{S}_{1}$ (resp. $\tilde{S}_{2}$) is a quasi Cartan subgroup of
$G_{1}$ (resp. $G_{2}$). Since $\tilde{S}_{1}\cap\ker p=1$, the map $p:\tilde{S}_{1}\rightarrow\tilde{S}_{2}$ is
an isomorphism. It is clear that the map $p$ gives an isomorphism \[(\tilde{S}_{1},R(G_{1},\tilde{S}_{1}))\cong
(\tilde{S}_{2},R(G_{2},\tilde{S}_{2})).\] However, $G_{1}\not\cong G_{2}$.
\end{example}

By Lemma \ref{L:ARS-isogeny2} and Example \ref{E:ARS-isogeny3}, we get the following proposition.

\begin{proposition}\label{P:ARS-isogeny4}
The quasi root system of a compact semisimple Lie group $G$ with cyclic component group determines the group up to
isogeny, but even the quasi root system together with the quasi Cartan subgroup do not determine the group $G$ up
to isomorphism in general.
\end{proposition}

\begin{lemma}\label{L:QRS-sc}
Let $G$ be a compact Lie group with simply-connected neutral subgroup $G^{0}$. If $\tilde{S}$ is a generalized Cartan
subgroup of $G$, then the root system $R(G,\tilde{S})$ on the quasi torus $\tilde{S}$ is of simply-connected type.
\end{lemma}

\begin{proof}
Choose an element $s\in\tilde{S}$ generating $\tilde{S}/\tilde{S}^{0}$ such that $\Ad(s)|_{\mathfrak{g}}$ is a pinned
automorphism of $\mathfrak{g}$. Write \[\theta=\Ad(s)|_{\mathfrak{g}}.\] Then, $\tilde{S}^{0}$ is a maximal torus of
$(G^{0})^{\theta}$ (which is connected by a theorem of Steinberg). Note that $R(G,\tilde{S})$ on $\tilde{S}$ being of
simply-connected type depends only on the restricted root system $R(G^{0},\tilde{S}^{0})$. Then, the issue reduces to
the case that $\theta$ acts transitively on irreducible factors of $\mathfrak{g}$. Furthermore, it reduces to the case
that $\mathfrak{g}$ is simple. Write $\Phi'$ for the root system of $\mathfrak{g}$. When $\Phi'$ is simply-laced, the
issue is trivial. When $o(\theta)=2$ and $\Phi'$ is of type $\A_{2n-1}$ ($n\geq 2$), $\D_{n}$ ($n\geq 4$), $\E_{6}$,
or $o(\theta)=3$ and $\Phi'$ is of type $\D_{4}$, it follows from the fact that $R(G^{0},\tilde{S}^{0})$ is equal to
the root system of a connected and simply-connected compact Lie group $G'\subset G$ with $\tilde{S}^{0}$ a maxima torus
where \[G'\cong\left\{\begin{array}{cccc}\Sp(n)\textrm{ when }\Phi\cong\A_{2n-1}\textrm{ and }o(\theta)=2;
\\\Spin(2n-1)\textrm{ when }\Phi\cong\D_{n}\textrm{ and }o(\theta)=2;\\ \F_{4}\textrm{ when }\Phi\cong\E_{6}
\textrm{ and }o(\theta)=2;\\ \G_{2}\textrm{ when }\Phi\cong\D_{4}\textrm{ and }o(\theta)=3.\end{array}\right.\]
When $o(\theta)=2$ and $\Phi'$ is of type $\A_{2n}$ ($n\geq 1$), it is as in Example \ref{E:ARS-isogeny3} and one
can make a direct verification. Alternatively, in this case $R(G^{0},\tilde{S}^{0})\cong\BC_{n}$. Then, the root
lattice coincides with integral weight weight. Hence, coroots generate the cocharacter group $X_{\ast}(\tilde{S}^{0})$.
\end{proof}

In the following theorem, we show that any reduced quasi root system on a quasi torus can be realized as the quasi
root system associated to a generalized Cartan subgroup in a compact Lie group.

\begin{theorem}\label{T:QRS-rs}
Let $\Phi$ be a reduced quasi root system on a quasi torus $\tilde{S}$. Then there exists a compact Lie group $H$ and
a generalized Cartan subgroup $\tilde{S}'$ of $H$ meeting with all connected components of $H$ such that  \[(R(H,\tilde{S}'),\tilde{S}')\cong(\Phi,\tilde{S}).\]
\end{theorem}

\begin{proof}
Let $L$ be the set of cocharacters $\check{\lambda}\in X_{\ast}(\tilde{S}^{0})$ such that $n\check{\lambda}$ is an
integral linear combination of coroots $\{\check{\alpha}:\alpha\in\Phi\}$ for some $n\in\mathbb{Z}_{>0}$. Then,
there is a sub-torus $T_{1}$ of $\tilde{S}^{0}$ such that $X_{\ast}(T_{1})=L$. Let \[\ker\Phi=\{x\in\tilde{S}:
\alpha(x)=1,\ \forall\alpha\in\Phi\}\] and put \[\tilde{S}^{ad}=\tilde{S}/\ker\Phi.\] Then $\Phi$ could be regarded
as an adjoint type quasi root system on the quasi torus $\tilde{S}^{ad}$. By Lemma \ref{L:ARS-isogeny1}, we construct
a root system $\Phi'$ and a diagram automorphism $\theta$ such that the resulting quasi root system from the action
of $\theta$ on $\Phi'$ is isomorphic to $\Phi$. Let $H_{2}^{sc}$ be a connected and simply-connected compact Lie
group with root system isomorphic to $\Phi'$. Then $\theta$ lifts to a pinned automorphism of $H_{2}^{sc}$. Take a
maximal torus $T_{2}$ of $(H_{2}^{sc})^{\theta}$. Then, the restricted root system $R(H_{2}^{sc},T_{2})\cong p(\Phi)$.
By Lemma \ref{L:QRS-sc}, coroots of $R(H_{2}^{sc},T_{2})$ generate the cocharacter group $X_{\ast}(T_{2})$. Then,
there is a finite surjective homomorphism \[\pi: T_{2}\rightarrow T_{1}\] such that the pull-back map $\pi^{\ast}:
X^{\ast}(T_{1})\rightarrow X^{\ast}(T_{2})$ maps $p(\Phi)$ onto $R(H_{2}^{sc},T_{2})$. Then, $\ker\pi
\subset Z(H_{2}^{sc})$. Using Lemma \ref{L:QRS-sc} and the projection $\tilde{S}\rightarrow\tilde{S}^{ad}$, one
constructs a homomorphism $\phi:\tilde{S}\rightarrow\Aut(H_{2}^{sc})$ with three properties: (1)$\ker\phi=\ker\Phi$; (2)$\phi|_{T_{1}}\circ\pi=\Ad|_{T_{2}}$; (3)there is an element $s$ generating $\tilde{S}/\tilde{S}^{0}$ such that
$\phi(s)=\theta$. Then, $\phi(\tilde{S})$ acts trivially on $\ker\pi$. Put \[H_{1}=(H_{2}^{sc}/\ker\pi)
\rtimes_{\phi}\tilde{S}.\] One shows that \[\{([x],\pi(x)^{-1}):x\in T_{2}\}\subset Z(H_{1}).\] Put
\[H=H_{1}/\{([x],\pi(x)^{-1}): x\in T_{2}\}\] and let $\tilde{S}'$ be the image of $1\times\tilde{S}$ in $H$. Then,
one can show that $\tilde{S}'$ is a generalized Cartan subgroup of $H$ which meets all connected components of $H$
and \[(R(H,\tilde{S}'),\tilde{S}')\cong(\Phi,\tilde{S}).\qedhere\]
\end{proof}

\subsection{Connection with dimension data}\label{SS:dim-ARS}

In this subsection let $G$ be a given compact Lie group with a biinvariant Riemannian metric. Let $\tilde{S}\subset G$
be a quasi torus, and $S$ be a connected component of $\tilde{S}$ generating it. Put \[\Gamma^{\circ}=
N_{G}(\tilde{S})/Z_{G}(\tilde{S})\] and \[\Gamma^{'\circ}=N_{G}(S)/Z_{G}(S).\] Note that, the induced action of
$\Gamma^{\circ}$ on $\tilde{S}/\tilde{S}^{0}$ gives a natural homomorphism \[\phi: \Gamma^{\circ}\rightarrow
\Aut(\tilde{S}/\tilde{S}^{0})\] and we have $\Gamma^{'\circ}=\ker\phi$. Hence, $\Gamma^{'\circ}$ is a normal subgroup
of $\Gamma^{\circ}$. As defined in \S \ref{SS:qrs-char}, there is a quasi root systems $\Psi''_{\tilde{S}}$ on
$\tilde{S}$ defined with the inner product on $\mathfrak{s}_{0}$ induced from the Riemannian metric on $G$. Let
$\Psi'_{\tilde{S}}$ be the minimal quasi sub-root system of $\Psi''_{\tilde{S}}$ containing all quasi root systems
$R(H,\tilde{S})$ with $H$ runs over closed subgroups of $G$ with $\tilde{S}$ a generalized Cartan subgroup. We call
$\Psi'_{\tilde{S}}$ the {\it enveloping quasi root system} associated to the quasi torus $\tilde{S}$ in $G$. The
following lemma is analogous to \cite[Proposition 3.3 and Corollary 3.4]{Yu-dimension}.

\begin{lemma}\label{L:Psi'}
We have $W_{\Psi'_{\tilde{S}}}\subset\Gamma^{\circ}$, and $\Psi'_{\tilde{S}}$ equals the union of quasi root systems
$R(H,\tilde{S})$ with $H$ runs over closed subgroups of $G$ with $\tilde{S}$ a generalized Cartan subgroup.
\end{lemma}

\begin{proof}
Write $X$ for the union of quasi root systems $R(H,\tilde{S})$ with $H$ runs over closed subgroups of $G$ with
$\tilde{S}$ a generalized Cartan subgroup. Then, $X\subset\Psi'_{\tilde{S}}$. Clearly, $X$ is $\Gamma^{\circ}$ stable
and $s_{\alpha}\in\Gamma^{\circ}$ for any $\alpha\in X$. Thus, $s_{\alpha}(\beta)\in X$ for any $\alpha,\beta\in X$.
Hence, $X$ satisfies the condition (QRS2) in Definition \ref{D:QRS1}. The condition (QRS1) is trivially satisfied.
Thus, $X$ is a quasi root system on $\tilde{S}$. Then, $\Psi'_{\tilde{S}}=X$ by the minimality of $\Psi'_{\tilde{S}}$.
Thus, \[W_{\Psi'_{\tilde{S}}}=W_{X}\subset\Gamma^{\circ}.\]
\end{proof}

\begin{proposition}\label{P:dim-char}
Let $H_{1},H_{2},\dots,H_{s}$ be a set of closed subgroups of $G$ and $c_{1},\dots,c_{s}\in\mathbb{C}$ be constants.
In order that \[\sum_{1\leq i\leq s}c_{i}\mathscr{D}_{H_{i}}=0\] it is necessary and sufficient that: for any closed
commutative connected subset $S$ in $G$, we have \begin{equation}\label{Eq:dim-char}\sum_{1\leq j\leq t}
c_{i_{j}}m_{i_{j}}^{-1}\sum_{1\leq k\leq t_{j}}F_{\Phi_{i_{j},k},\Gamma^{'\circ}}(x)=0,\ \forall x\in S\end{equation}
where $H_{i_{j}}$ ($i_{1}<\cdots<i_{t}$) runs over all members among $H_{1},H_{2},\dots,H_{s}$ with a maximal
commutative connected subset conjugate to $S$, $t_{j}$ is the number of connected components of $H_{i_{j}}$ with
a maximal commutative connected subset conjugate to $S$, $m_{i}=|H_{i}/H_{i}^{0}|$, $\Gamma^{'\circ}=
N_{G}(S)/Z_{G}(S)$, and $\Phi_{i_{j},k}$ is the quasi root system of $H_{i_{j}}$ with respect to a generalized
Cartan subgroup of it which is generated by a maximal commutative connected subset of $H_{i_{j}}$ conjugate to $S$.
\end{proposition}

\begin{proof}
Based on Lemmas \ref{L:integration} and \ref{L:AR4}, this proposition can be shown in the same way as
\cite[Proposition 3.8]{Yu-dimension}.
\end{proof}

Note that in \eqref{Eq:dim-char}, we identify $S$ with a maximal commutative connected subset of $H_{i_{j}}$ conjugate
to it and regard $\Phi_{i_{j},k}$ as a quasi root system on $\tilde{S}$.

\begin{corollary}\label{C:dim-char2}
In Proposition \ref{P:dim-char}, assume that the component groups \[H_{i}/(H_{i})^{0},\ 1\leq i\leq s\] are all cyclic.
Let $S$ be a closed commutative connected subset in $G$ which is conjugate to a maximal commutative connected subset of
some $H_{i_{0}}$. Suppose that whenever a maximal commutative connected subset of a connected component of some $H_{i}$
is conjugate to $S$, this connected component generates $H_{i}/(H_{i})^{0}$. Then, \eqref{Eq:dim-char} is equivalent to
\begin{equation}\label{Eq:dim-char2}\sum_{1\leq j\leq t}c_{i_{j}}t_{j}m_{i_{j}}^{-1}F_{\Phi_{i_{j}},\Gamma^{\circ}}(x)
=0,\ \forall x\in S\end{equation} where $H_{i_{j}}$ ($i_{1}<\cdots<i_{t}$) runs over all members among $H_{1},H_{2},
\dots,H_{s}$ with a maximal commutative connected subset conjugate to $S$, $t_{j}$ is the number of connected components
of $H_{i_{j}}$ with a maximal commutative connected subset conjugate to $S$, $m_{i}=|H_{i}/H_{i}^{0}|$, $\Gamma^{\circ}=
N_{G}(\tilde{S})/Z_{G}(\tilde{S})$, and $\Phi_{i_{j}}$ is the quasi root system of $H_{i_{j}}$ with respect to a
generalized Cartan subgroup of it which is generated by a maximal commutative connected subset of $H_{i_{j}}$ conjugate
to $S$.
\end{corollary}

\begin{proof}
Assume that $H_{i}$ has a connected component with a maximal commutative connected subset conjugate to $S$. Without
loss of generality we identity this maximal commutative connected subset with $S$. Then, $\tilde{S}$ meets with all
connected components of $H_{i}$ and its intersection with each connected component of $H_{i}$ is $|\tilde{S}
\cap G^{0}/\tilde{S}^{0}|$ left $\tilde{S}^{0}$ cosets in $\tilde{S}$, all are closed commutative connected sets.
Suppose that $S'$ is a maximal commutative connected subset of another connected component of $H_{i}$ which is $G$
conjugate to $S$. By the above, taking a $H_{i}^{0}$ conjugate one if necessary we assume that $S'$ contains a left
$\tilde{S}^{0}$ coset in $\tilde{S}$. Since it is assumed that $S'$ is $G$ conjugate to $S$, then $\dim S'=\dim S$
and the minimal order of elements in $S'$ is equal to that in $S$. By $\dim S'=\dim S$, we get that: $S'$ is equal
to a left $\tilde{S}^{0}$ coset in $\tilde{S}$. By the fact of the minimal order of elements in $S'$ is equal to
that in $S$, we get that: $S'$ generates $\tilde{S}/\tilde{S}^{0}$ as well. Then, a conjugation taking $S'$ to $S$
induces a conjugation on $\tilde{S}$ and hence it comes from an element in $\Gamma^{\circ}$. This shows that:
maximal commutative connected subsets in $H_{i}$ which are $G$ conjugate to $S$ are $H_{i}^{0}$ conjugates of
\[\gamma\cdot S,\ \gamma\in\Gamma^{\circ}.\] Then, \[\sum_{1\leq k\leq t_{j}}F_{\Phi_{i_{j},k},\Gamma^{'\circ}}(x)
=t_{j}F_{\Phi_{i_{j}},\Gamma^{\circ}}(x),\ \forall x\in S.\] This shows the conclusion of this corollary.
\end{proof}

\subsection{Actual leading weight}\label{SS:actual}

Let $\tilde{S}$ be quasi torus and $S$ be a connected component of $\tilde{S}$ generating it.

\begin{lemma}\label{L:muPhi}
Let $\Phi$ be a reduced quasi root system on $\tilde{S}$ with the folding index of each root being equal to 1 or 2.
Fix a positive system $p(\Phi)^{+}$ of $p(\Phi)$. Then, there exists a unique $p(\Phi)^{+}$-dominant weight
$\bar{\mu}_{0}\in 2\mathbb{Z}p(\Phi)$ with the following property: there exists a linear character
$\mu_{0}\in X^{\ast}(\tilde{S})$ appearing in $A_{\Phi}$ such that $\mu_{0}|_{\tilde{S}^{0}}=\bar{\mu}_{0}$; for
any other linear character $\mu$ appearing in $A_{\Phi}$ such that $\mu|_{\tilde{S}^{0}}\in 2\mathbb{Z}p(\Phi)$,
we have $|\mu|_{\tilde{S}^{0}}|\leq|\bar{\mu}_{0}|$.
\end{lemma}

\begin{proof}
Choose a pinned element $s\in S$ of $\Phi$ and put \[\Phi_{s}=\{\alpha\in\Phi:\alpha(s)=1\}.\] Then, $\Phi_{s}$ and $\{m_{\alpha}\alpha:\alpha\in\Phi_{s}\}$ are both abstract root systems and the projection $p: X^{\ast}(\tilde{S})
\rightarrow X^{\ast}(\tilde{S}^{0})$ gives an isomorphism $\Phi_{s}\cong p(\Phi)$. Put \[\Phi_{s}^{+}=\{\alpha\in
\Phi_{s}:\alpha|_{\tilde{S}^{0}}\in p(\Phi)^{+}\}\] and \[\delta_{\Phi}=\frac{1}{2}\sum_{\alpha\in\Phi_{s}^{+}}
m_{\alpha}\alpha.\] Then \begin{eqnarray*}&&A_{\Phi}\\&=&\sum_{w\in W_{\Phi_{s}}}\epsilon(w)[\delta_{\Phi}-
w\delta_{\Phi}]\\&=&\sum_{w\in W_{\Phi_{s}}}\epsilon(ww_{l})[\delta_{\Phi}-ww_{l}\delta_{\Phi}]\\&=&
(-1)^{|\Phi_{s}^{+}|}\sum_{w\in W_{\Phi_{s}}}\epsilon(w)[\delta_{\Phi}+w\delta_{\Phi}],\end{eqnarray*} where
$w_{l}\in W_{\Phi_{s}}$ is the unique element such that $w_{l}\Phi_{s}^{+}=-\Phi_{s}^{+}$. Note that \[|\delta_{\Phi}+w\delta_{\Phi}|^{2}=2|\delta_{\Phi}|^{2}-|\delta_{\Phi}-w\delta_{\Phi}|^{2}.\] It suffices
to find an element $w_{0}\in W_{\Phi_{s}}$ such that: (1)$\delta_{\Phi}+w_{0}\delta_{\Phi}\in 2\mathbb{Z}\Phi_{s}$;
(2)for any $w_{0}\neq w\in W_{\Phi_{s}}$ with $\delta_{\Phi}+w\delta_{\Phi}\in 2\mathbb{Z}\Phi_{s}$, we have
\[|\delta_{\Phi}-w\delta_{\Phi}|>|\delta_{\Phi}-w_{0}\delta_{\Phi}|.\]

Without loss of generality we assume that $\Phi_{s}\cong p(\Phi)$ is irreducible. Let's describe $w_{0}$ and
$\mu_0=\delta_{\Phi}+w_{0}\delta_{\Phi}$ case by case. \begin{enumerate}
\item[(i)]When all $m_{\alpha}=2$, we have $\delta_{\Phi}\in\mathbb{Z}\Phi_{s}$. Then, $w_{0}=1$.
\item[(ii)]When $\Phi\cong\B_{n}^{(2,1)}$, $\BC_{n}^{(2,2,1)}$, $\E_{6}$, $\E_{8}$, $\F_{4}$, $\F_{4}^{(2,1)}$ or
$\G_{2}$, we have $\delta_{\Phi}\in\mathbb{Z}\Phi_{s}$. Then, $w_{0}=1$.
\item[(iii)]When $\Phi=\A_{n-1}$, we have $\delta_{\Phi}=(\frac{n-1}{2},\frac{n-3}{2},\dots,\frac{1-n}{2})$.
If $n$ is odd, then we have $\delta_{\Phi}\in\mathbb{Z}\Phi_{s}$ and $w_{0}=1$. If $n$ is even, then we have
\[w_{0}=(12)(34)\cdots((n-1)n)\] and \[\delta_{\Phi}+w_{0}\delta_{\Phi}=(n-2,n-2,n-6,n-6,\dots,2-n,2-n).\]
\item[(iv)]When $\Phi=\B_{n}$, we have $\delta_{\Phi}=(n-\frac{1}{2},n-\frac{3}{2},\dots,\frac{1}{2})$.
Then, we have \[w_{0}=\left\{\begin{array}{cc}(12)(34)\cdots((n-1)n)\textrm{ if $n$ is even}\\(12)(34)
\cdots((n-2)(n-1))s_{e_{n}}\textrm{ if $n$ is odd}\end{array}\right.\] and \[\delta_{\Phi}+w_{0}\delta_{\Phi}
=\left\{\begin{array}{cc}(2n-2,2n-2,2n-6,2n-6,\dots,2,2)\textrm{ if $n$ is even}\\(2n-2,2n-2,2n-6,2n-6,
\dots,4,4,0)\textrm{ if $n$ is odd}.\end{array}\right.\]
\item[(v)]When $\Phi=\C_{n}$, we have $\delta_{\Phi}=(n,n-1,\dots,1)$. If $n(n+1)$ is a multiple of 4, then
we have $\delta_{\Phi}\in\mathbb{Z}\Phi_{s}$ and $w_{0}=1$. If $n(n+1)$ is not a multiple of 4, then we have
$w_{0}=s_{e_{n}}$ and \[\delta_{\Phi}+w_{0}\delta_{\Phi}=(2n,2n-2,\dots,6,4,0).\]
\item[(vi)]When $\Phi=\C_{n}^{(2,1)}$, we have $\delta_{\Phi}=(2n-1,2n-3,\dots,1)$. If $n$ is even, then we
have $\delta_{\Phi}\in\mathbb{Z}\Phi_{s}$ and $w_{0}=1$. If $n$ is odd, then we have $w_{0}=s_{e_{n}}$ and
\[\delta_{\Phi}+w_{0}\delta_{\Phi}=(4n-2,4n-6,\dots,10,6,0).\]
\item[(vii)]When $\Phi=\D_{n}$, we have $\delta_{\Phi}=(n-1,n-2,\dots,0)$. If $n(n-1)$ is a multiple of 4,
then we have $\delta_{\Phi}\in\mathbb{Z}\Phi_{s}$ and $w_{0}=1$. If $n(n-1)$ is not a multiple of 4, then
we have $w_{0}=s_{e_{n-1}}s_{e_{n}}$ and \[\delta_{\Phi}+w_{0}\delta_{\Phi}=(2n-2,2n-4,\dots,6,4,0,0).\]
\item[(viii)]When $\Phi=\E_{7}$, we have \[2\delta_{\Phi}\equiv \alpha_{2}+\alpha_{5}+\alpha_{7}
\pmod{2\mathbb{Z}\Phi_{s}}.\] Then, we have $w_{0}=s_{\alpha_{2}}s_{\alpha_{5}}s_{\alpha_{7}}$ and
$\delta_{\Phi}+w_{0}\delta_{\Phi}=2\delta_{\Phi}-(\alpha_{2}+\alpha_{5}+\alpha_{7})$.
\end{enumerate}
\end{proof}

In Lemma \ref{L:muPhi}, we call $2\bar{\mu}_{\Phi}:=2\bar{\mu}_{0}$ {\it the actual leading weight of $A_{\Phi}$}.
In the proof of Lemma \ref{L:muPhi}, an explicit formula for the actual leading weight $2\bar{\mu}_{\Phi}$ of
each irreducible quasi root system $\Phi$ with total folding index being equal to 1 or 2 is given.

\subsection{Illustration in the $\BC_{n}$ case}\label{SS:BCn}

Like the study for characters associated to reduced sub-root systems in a given root system taken in
\cite{Yu-dimension}, one could make a thorough study for equalities and linear relations among characters
associated to reduced quasi sub-root systems in a given quasi root system $\Psi$. Besides the obvious
complication due to folding indices and fractional factor, another complication is the averaging operation
by a group $W$. Below we consider the case that \[\dim\tilde{S}=n,\ \Psi_{\tilde{S}^{0}}=\BC_{n}
\textrm{ and }W\supset\Hom(\tilde{S}/\tilde{S}^{0},\tilde{S}^{0}).\]

Fix $d\geq 1$. For each $n\geq 1$, let $\tilde{S}_{n}=\U(1)^{n}\times\langle s_{0}\rangle$, where $s_{0}$ has
order $d$. Then $\tilde{S}_{n}^{0}=\U(1)^{n}$. Write $(x_{1},\dots,x_{n},s_{0}^{j})$ ($x_{i}\in\U(1)$,
$j\in\mathbb{Z}$) for a general element in $\tilde{S}_{n}$. Define $e_{1},\dots,e_{n},\beta_{0}
\in X^{\ast}(\tilde{S}_{n})$ by \[e_{i}(x_{1},\dots,x_{n},s_{0}^{j})=x_{i},\ 1\leq i\leq n\]
and \[\beta_{0}(x_{1},\dots,x_{n},s_{0}^{j})=e^{\frac{2j\pi\mathbf{i}}{d}}.\] We define an inner product on
$X^{\ast}(\tilde{S}_{n}^{0})$ by taking $\{\bar{e}_{i}=e_{i}|_{\tilde{S}_{n}^{0}}: 1\leq i\leq n\}$ an
orthonormal basis of $\mathbf{i}\mathfrak{s}_{0}^{\ast}$. Thus,  \[\Psi_{\tilde{S}_{n}^{0}}=
\{\pm{\bar{e}_{i}}\pm{\bar{e}_{j}},\pm{\bar{e}_{k}},\pm{2\bar{e}_{k}}:1\leq i<j\leq n,1\leq k\leq n\}\] and
it is a root system of type $\BC_{n}$. Write $\mathbb{Z}^{n}$ for the character group of $\tilde{S}_{n}^{0}$,
$W_{n}$ for the Weyl group of $\Psi_{\tilde{S}_{n}^{0}}$, $\mathbb{Z}_{n}=\mathbb{Q}[\mathbb{Z}^{n}]$, and
$Y_{n}=\mathbb{Z}_{n}^{W_{n}}$.

Put $S_{n}=s_{0}\tilde{S}_{n}$, which is a connected component of $\tilde{S}_{n}$. We may view $W_{n}$ as a
subgroup of $\Stab_{\Aut(\tilde{S}_{n})}$ by $W_{n}=\Stab_{\Aut(\tilde{S}_{n})}(s_{0})$. Write $Y_{n,d}$ for
the algebra of the restriction to $S_{n}$ of $W_{n}$ invariant rational coefficient characters on
$\tilde{S}_{n}$. It is clear that $Y_{n,d}$ possesses a set of generators
$\chi^{\ast}_{a_{1}e_{1}+\cdots+a_{n}e_{n}+k\beta_{0},W_{n}}$, where $a_{1},\dots,a_{n}$ are integers with
$a_{1}\geq\cdots a_{n}\geq 0$. Then, the following lemma is clear.

\begin{lemma}\label{L:Sn-Wn}
By mapping $\chi^{\ast}_{a_{1}e_{1}+\cdots+a_{n}e_{n}+k\beta_{0},W_{n}}$ to $e^{\frac{2k\pi\mathbf{i}}{d}}
x_{a_{1}}\cdots x_{a_{n}}$, we get an algebra injection \[E:\lim_{\longrightarrow_{n}}Y_{n,d}\hookrightarrow
\mathbb{C}[x_{0},x_{1},\dots].\]
\end{lemma}

For each positive integer $k$, we define homomorphisms \[\Fold_{[k]},\Ave_{k},\Ave'_{k}:\mathbb{C}[x_{0},x_{1},\dots]
\rightarrow\mathbb{C}[x_{0},x_{1},\dots]\] by \[\Fold_{k}(x_{j})=x_{jk},\ \Ave_{k}(x_{j})=\epsilon_{k,j}x_{j},\
\Ave'_{k}(x_{j})=\epsilon'_{k,j}x_{j},\ \forall j\geq 0,\] where \[\epsilon_{k,j}=\left\{\begin{array}{cc}
1\ \textrm{ if }k|j\\0\textrm{ if }k\not|j\end{array}\right.\textrm{ and }\epsilon'_{k,j}=\left\{\begin{array}{cc}
(-1)^{\frac{j}{k}}\textrm{ if }k|j\\0\quad\textrm{ if }k\not|j\end{array}\right..\] We call
\[\Fold_{k},\ \Ave_{k},\ \Ave'_{k}\] the {\it $k$-folding operation}, the {\it $k$-averaging operation},
the twisted {\it $k$-twisted averaging operation} respectively. The following lemma is clear.

\begin{lemma}\label{L:fold-ave}
For positive integers $k$ and $k'$ with $k|k'$, we have \[\Ave_{k'}\circ\Fold_{k}=\Fold_{k}\circ\Ave_{\frac{k'}{k}}\]
and \[\Ave'_{k'}\circ\Fold_{k}=\Fold_{k}\circ\Ave'_{\frac{k'}{k}}.\]
\end{lemma}

Put \[\Gamma_{n}=\Hom(\tilde{S}_{n}/\tilde{S}_{n}^{0},\tilde{S}_{n}^{0})\rtimes W_{n}.\] Then, the algebra of
$\Gamma_{n}$ invariant characters on $S_{n}$ is a subalgebra of $Y_{n,d}$. For any reduced quasi sub-root system
$\Phi$ of $\Psi''_{\tilde{S}_{n}}$, the character $F_{\Phi,\Gamma_{n}}$ is a $\Gamma_{n}$ invariant character on
$S_{n}$.

We show below how to get the polynomial $E(F_{\Phi,\Gamma_{n}})$ from the ordinary folding index, the twisted number
and the fractional factor of a quasi root system $\Phi$ on $\tilde{S}_{n}$. First, if there is a partition
$\{1,\dots,n\}=I_{1}\bigsqcup I_{2}$ such that $\Phi=\Phi_{1}\sqcup\Phi_{2}$, then \[E(F_{\Phi,\Gamma_{n}})=
E(F_{\Phi_{1},\Gamma_{n}})E(F_{\Phi_{2},\Gamma_{n}}),\] where $\Phi_{j}$ consists of roots $\alpha\in\Phi$ such
that $\alpha|_{\tilde{S}^{0}}$ is a linear combination of $\bar{e}_{i}$ ($i\in I_{j}$). We call a quasi root
system $\Phi$ on $\tilde{S}_{n}$ {\it in-decomposable} if there exists no proper partition $\{1,\dots,n\}=
I_{1}\bigsqcup I_{2}$ such that $\Phi=\Phi_{1}\sqcup\Phi_{2}$. Thus, it suffices to calculate
$E(F_{\Phi,\Gamma_{n}})$ for in-decomposable reduced quasi sub-root systems $\Phi$ on $\tilde{S}_{n}$.

\begin{definition}\label{D:fractional2}
We say that an in-decomposable quasi root system $\Phi$ on $\tilde{S}_{n}$ has an even fractional factor if
$\Phi$ is $\Gamma_{n}$ conjugate to a quasi root system which has a simple system contained in
\[\{\pm{e_{i}}\pm{e_{j}},\pm{e_{k}},\pm{2e_{k}}:1\leq i<j\leq n,1\leq k\leq n\};\] otherwise we say that
$\Phi$ has an odd fractional factor.
\end{definition}

Now assume that $\Phi$ is in-decomposable. Let $e$ be the ordinary folding index of $\Phi$ and $t$ be the
twisted number of $\Phi$. Choose a pinned element $s\in S_{n}$ of $\Phi$, and let $\{\alpha_{i}=
\alpha'_{i}+k_{i}\alpha_{0}:1\leq i\leq l\}$ be a simple system of $\Phi_{s}$, where each $\alpha'_{i}\in
\span_{\mathbb{Z}}\{e_{j}:1\leq j\leq n\}$ and each $k_{i}\in\mathbb{Z}$. In the following lemma, we
describe $E(F_{\Phi,\Gamma_{n}})$ using the averaging and twisted averaging operations. Recall that the
polynomials $a_{k},b_{k},b'_{k},c_{k},d_{k}$ are defined in \cite[\S 7]{Yu-dimension}.

\begin{lemma}\label{L:Ave1}
\begin{enumerate}
\item[(i)] If $p(\Phi)=\A_{k-1}$, then \[E(F_{\Phi,\Gamma_{n}})=\Ave_{d}(\Fold_{e}(a_{k})).\]
\item[(ii)](1)If $p(\Phi)=\B_{k}$ and $t=1$, then \[E(F_{\Phi,\Gamma_{n}})=\Ave_{d}(\Fold_{d}(b_{k})).\]
(2)If $p(\Phi)=\B_{k}$ and $t=2$, then \[E(F_{\Phi,\Gamma_{n}})=\Ave_{d}(\Fold_{e}(c_{k})).\]
\item[(iii)](1)If $p(\Phi)=\C_{k}$, $t=1$ and $\Phi$ has even fractional factor, then \[E(F_{\Phi,\Gamma_{n}})
=\Ave_{d}(\Fold_{e}(c_{k})).\] (2)If $p(\Phi)=\C_{k}$, $t=1$ and $\Phi$ has odd fractional factor, then \[E(F_{\Phi,\Gamma_{n}})=\Ave'_{d}(\Fold_{e}(c_{k})).\] (3)If $p(\Phi)=\C_{k}$, $t=2$ and $\Phi$ has even
fractional factor, then \[E(F_{\Phi,\Gamma_{n}})=\Ave_{d}(\Fold_{2e}(b_{k})).\] (4)If $p(\Phi)=\C_{k}$,
$t=2$ and $\Phi$ has odd fractional factor, then \[E(F_{\Phi,\Gamma_{n}})=\Ave'_{d}(\Fold_{2e}(b_{k})).\]
\item[(iv)](1)If $p(\Phi)=\D_{k}$ and $\Phi$ has even fractional factor, then \[E(F_{\Phi,\Gamma_{n}})=
\Ave_{d}(\Fold_{e}(d_{k})).\] (2)If $p(\Phi)=\D_{k}$ and $\Phi$ has odd fractional factor, then
\[E(F_{\Phi,\Gamma_{n}})=\Ave'_{d}(\Fold_{e}(d_{k})).\]
\item[(v)]If $p(\Phi)=\BC_{k}$, then \[E(F_{\Phi,\Gamma_{n}})=\Ave_{d}(\Fold_{2e}(c_{k})).\]
\end{enumerate}
\end{lemma}

\begin{proof}
We prove (iv). The proofs for (i), (ii), (iii), (v) are similar.

(iv)Assume that $p(\Phi)=\D_{k}$ ($k\geq 2$). Replacing $\Phi$ by a $\Gamma_{n}$ conjugate one we may assume
that $\alpha_{i}=e_{i}-e_{i+1}$ ($1\leq i\leq k-1$) and $\alpha_{n}=e_{n-1}+e_{n}+k\beta_{0}$ where
\[k=\left\{\begin{array}{cc}\quad 0\quad\textrm{ if }d\textrm{ is odd;}\\ 0\textrm{ or }1\textrm{ if }d
\textrm{ is even.}\end{array}\right.\] Then, the twisted number $t=1$ due to $p(\Phi)$ is simply-laced. We
have \[\delta_{\Phi}=\frac{e}{2}\sum_{\alpha\in\Phi_{s}^{+}}\alpha\] and \[A_{\Phi}=\sum_{w\in W_{\Phi_{s}}}
\epsilon(w)[\delta_{\Phi}-w\delta_{\Phi}].\] For each $w\in W_{\Phi_{s}}$, the weight $\delta_{\Phi}-
w\delta_{\Phi}$ is of the form \[a_{1}e_{1}\!+\!\sum_{2\leq j\leq n-2}(a_{j}-a_{j-1})e_{j}\!+\!(a_{n}+a_{n-1}
-a_{n-2})e_{n-1}\!+\!(a_{n}-a_{n-1})e_{n}\!+\!a_{n}k\alpha_0.\] For the term $[\delta_{\Phi}-w\delta_{\Phi}]$
to survive under the averaging by $\Hom(\tilde{S}_{n}/\tilde{S}_{n}^{0},\tilde{S}_{n}^{0})$ it is necessary
and sufficient that \[a_{1},a_{2}-a_{1},\dots,a_{n-2}-a_{n-3},a_{n}+a_{n-1}-a_{n-2},a_{n}-a_{n-1}\] are all
multiples of $d$. Then, \[|a_{1}|+\sum_{2\leq j\leq k-2}|a_{j}-a_{j-1}|+|a_{n}+a_{n-1}-a_{n-2}|+|a_{n}-
a_{n-1}|\equiv 2a_{n}\pmod{2d}.\] By these it follows that: if $k=0$ (i.e., $\Phi$ has even fractional factor),
then \[E(F_{\Phi,\Gamma_{n}})=\Ave_{d}(\Fold_{e}(d_{k}));\] if $k=1$ and $d$ is even (i.e., $\Phi$ has odd
fractional factor), then \[E(F_{\Phi,\Gamma_{n}})=\Ave'_{d}(\Fold_{e}(d_{k})).\]
\end{proof}

The following Lemma \ref{L:Ave2} gives precise description of the $d$-averaging operation on polynomials $a_{k}$,
$b_{k}$, $b'_{k}$, $c_{k}$, $d_{k}$. The proof is by direct application of matrix expressions of these polynomials
(\cite[Proposition 7.2]{Yu-dimension}). We omit the proof. The description for twisted $d$-averaging operation on
$a_{k}$, $b_{k}$, $b'_{k}$, $c_{k}$, $d_{k}$ is similar.

\begin{lemma}\label{L:Ave2}
Let $p,r$ be non-negative integers and $0\leq r\leq d-1$. Then, \begin{enumerate}
\item[(i)] $\Ave_{d}(a_{dp+r})=\Fold_{d}(a_{p+1})^{r}\Fold_{d}(a_{p})^{d-r}$.
\item[](1)When $d$ is even and $0\leq r<\frac{d}{2}$, \[\Ave_{d}(b_{dp+r})=\Fold_{d}(a_{2p+1})^{r}
\Fold_{d}(a_{2p})^{\frac{d}{2}-r}.\] (2)When $d$ is even and $\frac{d}{2}\leq r<d$, \[\Ave_{d}(b_{dp+r})=
\Fold_{d}(a_{2p+2})^{r-\frac{d}{2}}\Fold_{d}(a_{2p+1})^{d-r}.\] (3)when $d$ is odd and $0\leq r<\frac{d+1}{2}$, \[\Ave_{d}(b_{dp+r})=\Fold_{d}(a_{2p+1})^{r}
\Fold_{d}(a_{2p})^{\frac{d-1}{2}-r}\Fold_{d}(b_{p}).\]
\item when $d$ is odd and $\frac{d+1}{2}\leq r<d$, \[\Ave_{d}(b_{dp+r})=\Fold_{d}(a_{2p+2})^{r-\frac{d+1}{2}}
\Fold_{d}(a_{2p+1})^{d-r}\Fold_{d}(b_{p+1}).\]
\item when $d$ is even and $0\leq r<\frac{d}{2}$, \[\Ave_{d}(c_{dp+r})=\Fold_{d}(a_{2p+1})^{r}
\Fold_{d}(a_{2p})^{\frac{d-2}{2}-r}\Fold_{d}(b_{p})\Fold_{d}(c_{p}).\]
\item when $d$ is even and $\frac{d}{2}\leq r<d$, \[\Ave_{d}(c_{dp+r})=\Fold_{d}(a_{2p+2})^{r-\frac{d}{2}}
\Fold_{d}(a_{2p+1})^{d-1-r}\Fold_{d}(b_{p+1})\Fold_{d}(c_{p}).\]
\item when $d$ is odd and $0\leq r\leq\frac{d-1}{2}$, \[\Ave_{d}(c_{dp+r})=\Fold_{d}(a_{2p+1})^{r}
\Fold_{d}(a_{2p})^{\frac{d-1}{2}-r}\Fold_{d}(c_{p}).\]
\item when $d$ is odd and $\frac{d+1}{2}\leq r<d$, \[\Ave_{d}(c_{dp+r})=\Fold_{d}(a_{2p+2})^{r-\frac{d-1}{2}}
\Fold_{d}(a_{2p+1})^{d-1-r}\Fold_{d}(c_{p}).\]
\item when $d$ is even and $0\leq r<\frac{d}{2}$, \[\Ave_{d}(d_{dp+r+1})=\Fold_{d}(a_{2p+1})^{r}
\Fold_{d}(a_{2p})^{\frac{d-2}{2}-r}\Fold_{d}(d_{p+1})\Fold_{d}(\tilde{b}_{p}).\]
\item when $d$ is even and $\frac{d}{2}\leq r<d$, $\Ave_{d}(d_{dp+r+1})$ is equal to
\[\Fold_{d}(a_{2p+2})^{r-\frac{d}{2}}\Fold_{d}(a_{2p+1})^{d-1-r}\Fold_{d}(d_{p+1})\Fold_{d}(\tilde{b}_{p+1}).\]
\item when $d$ is odd and $0\leq r\leq\frac{d-1}{2}$, \[\Ave_{d}(d_{dp+r+1})=\Fold_{d}(a_{2p+1})^{r}
\Fold_{d}(a_{2p})^{\frac{d-1}{2}-r}\Fold_{d}(d_{p+1}).\]
\item when $d$ is odd and $\frac{d+1}{2}\leq r<d$, \[\Ave_{d}(d_{dp+r+1})=\Fold_{d}(a_{2p+2})^{r-\frac{d-1}{2}}
\Fold_{d}(a_{2p+1})^{d-1-r}\Fold_{d}(d_{p+1}).\]
\end{enumerate}
\end{lemma}

\section{Pinned automorphism and twining character formula}\label{S:twining}

We call a compact Lie group $G$ {\it quasi-connected} if its component group $G/G^{0}$ is a cyclic group.
Assume that $G$ is a quasi-connected compact Lie group. We call a generalized Cartan subgroup $\tilde{S}$
of $G$ meeting with all connected components of $G$ a {\it quasi Cartan subgroup}.

\begin{definition}\label{D:resChar1}
Let $\tilde{S}$ be a quasi Cartan subgroup in a quasi-connected compact Lie group $G$. Choose a connected component
$S$ generating $\tilde{S}$. For an (irreducible) finite-dimensional complex linear representation $\rho$ of $G$,
\begin{enumerate}
\item[(i)]we call $\chi_{\rho}|_{S}$ the (irreducible) {\it twisted character} of $\rho$ on $S$;
\item[(ii)]and we call $\chi_{\rho}|_{\tilde{S}^{0}}$ the (irreducible) {\it restriction character} of $\rho$ on
$\tilde{S}^{0}$.
\end{enumerate}
\end{definition}

The twining character formula expresses an irreducible twisted character in terms of the character of an irreducible
representation of a connected compact Lie group of smaller dimension called the {\it twining group}. The twining
character formula was first found by Jantzen (\cite{Jantzen}) and there are many different proofs in the literature,
e.g., \cite{Wendt} and \cite{Kumar-Lusztig-Prasad}. In this section we give a proof of the twining character
formula based on the Weyl orthogonality relation. The accounting in this section follows \cite{Kumar-Lusztig-Prasad}
and the Weyl orthogonality relation is also used in the proof in \cite{Wendt}. Moreover, we generalize the twining
character formula to any quasi-connected compact Lie group.

Let $G$ be a connected and simply-connected compact Lie group, and $T$ be a maximal torus of $G$. Choose a positive
system $R^{+}(\mathfrak{g},T)$ of the root system  $R(\mathfrak{g},T)$, and let $\mathfrak{b}$ be the
corresponding Borel subalgebra of $\mathfrak{g}$. Let $\{\alpha_1,\dots,\alpha_{l}\}$ be the corresponding simple
system. Choose a root vector $X_{\alpha_{i}}$ for each simple root $\alpha_{i}$ ($1\leq i\leq l$). We call the triple  $\mathcal{E}=(T,\mathfrak{b},\{X_{\alpha_{i}}\}_{1\leq i\leq l})$ a {\it pinning} of $G$. An automorphism $\theta$
of $G$ which normalizes $T$ and $\mathfrak{b}$ and permutes $\{X_{\alpha_{i}}\}_{1\leq i\leq l}$ is called a
{\it pinned automorphism} of $G$ with respect to the pinning $\mathcal{E}$. Write $\Aut(G,\mathcal{E})$ for the
group of pinned automorphisms of $G$ with respect to the pinning $\mathcal{E}$. Then $\Aut(G)=\Int(G)\rtimes
\Aut(G,\mathcal{E})$. Moreover, in each connected component of $\Aut(G)$ the set of pinned automorphisms consists
in a single $\Int(G)$ conjugacy class.

Let $\theta\in\Aut(G,\mathcal{E})$. Then, the action of $\theta$ on $\{X_{\alpha_{i}}\}_{1\leq i\leq l}$ induces
a permutation on the index set $I=\{1,\dots,l\}$, still denoted by $\theta$. There is a unique direct product
decomposition $G=G_{1}\times\cdots\times G_{s}$ such that $\theta$ normalizes each $G_{j}$ and $\theta$ acts
transitively on simple factors of $\mathfrak{g}_{j}$ ($1\leq j\leq s$). Write \[Y=\Hom(\U(1),T)\] for the cocharacter
group of $T$, and write \[X=\Hom(T,\U(1))\] for the character group of $T$. There is a perfect pairing \[\langle\cdot,
\cdot\rangle: Y\times X\rightarrow\mathbb{Z}.\] For each simple root $\alpha_{i}$, let $\check{\alpha}_{i}\in Y$ be
the corresponding coroot. The following fact is well-known.

\begin{fact}\label{F:rootDatum}
The tuple $(Y,X,\langle\cdot,\cdot\rangle,\{\check{\alpha}_{i}\},\{\alpha_{i}\})$ is a root datum.
\end{fact}

Put \[Y_{\theta}=Y/(\theta-1)Y\] and \[^{\theta}X=\{\lambda\in X: \theta(\lambda)=\lambda\}.\] Still write
$\langle\cdot,\cdot\rangle$ for the induced perfect pairing \[Y_{\theta}\times^{\theta}X\rightarrow\mathbb{Z}.\]
Write $I_{\theta}$ for the set of $\theta$-orbits in $I$. For each orbit $\mathcal{O}$, let
$\check{\alpha}_{\mathcal{O}}\in Y_{\theta}$ be the image of $\check{\alpha}_{i}$ ($i\in\mathcal{O}$) under
$Y\rightarrow Y_{\theta}$. Let $\alpha_{\mathcal{O}}=2^{h}\sum_{i\in\mathcal{O}}\alpha_{i}\in^{\theta}X$, where
$h=1$ if there exists $i,i'\in\mathcal{O}$ such that $\alpha_{i}+\alpha_{i'}$ is a root, and $h=0$ otherwise. Note
that $h=1$ except when $\mathcal{O}$ is contained in a union of irreducible root systems of type $\A_{2n}$ and
$\alpha_{i}$ ($i\in\mathcal{O}$) is placed at the $n$-th or $(n+1)$-th position of a Dynkin diagram of type
$\A_{2n}$.

\begin{fact}\label{F:rootDatum2}
The tuple $(Y_{\theta},^{\theta}X,\{\alpha_{\mathcal{O}}\},\{\check{\alpha}_{\mathcal{O}}\})$ is a root datum.
\end{fact}

Let $G_{\theta}$ be a connected compact semisimple Lie group associated to the root datum $(Y_{\theta},^{\theta}X,
\{\alpha_{\mathcal{O}}\},\{\check{\alpha}_{\mathcal{O}}\})$. Note that $G_{\theta}$ is simply-connected, and $T_{\theta}:=
T/\{\theta(t)t^{-1}:t\in T\}$ is a maximal torus of $G_{\theta}$. Let $X^{+}$ be the set of dominant characters in $X$, and
$^{\theta}X^{+}$ be the set of dominant characters in $^{\theta}X$. It is clear that $^{\theta}X^{+}=^{\theta}X\cap X^{+}$.
For any $\lambda\in X^{+}$ (resp. $\lambda\in ^{\theta}X^{+}$), write $V_{\lambda}$ (resp. $V'_{\lambda}$) for an irreducible
representation of $G$ (resp. of $G_{\theta}$) with highest weight $\lambda$, and write $0\neq v_{\lambda}$ for a vector in
$V_{\lambda}$ with weight $\lambda$. Let $\varpi:T\rightarrow T_{\theta}$ be the canonical projection.

\begin{theorem}[Twining character formula]\label{T:twining1}
For any $t\in T$, we have \[\tr(t\theta:V_{\lambda}\rightarrow V_{\lambda})=\tr(\varpi(t):V'_{\lambda}\rightarrow
V'_{\lambda}).\]
\end{theorem}

\begin{lemma}\label{L:twining2}
The map $t\theta\mapsto\varpi(t)$ gives a bijection between the set of $G$-conjugacy in $G\theta$ and the set of conjugacy
classes in $G_{\theta}$.
\end{lemma}

\begin{proof}
Write $\tilde{G}=G\rtimes\langle\theta\rangle$ and $S=T^{\theta}\theta$. By Lemma \ref{L:compact-conjugacy}, two elements
in $S$ are in a $G$ conjugacy class if and only if they are in the same $W(G,S)$ orbit. By Lemma \ref{L:split}, there is
an exact sequence \[1\rightarrow N_{T}(S)/\tilde{S}^{0}\rightarrow W(G,S)\rightarrow W_{R(G,T^{\theta})}\rightarrow 1.\]
For two elements $t,t'\in T^{\theta}$, $\varpi(t)$ and $\varpi(t')$ are in the same $G_{\theta}$ conjugacy class if and
only $t'(wt)^{-1}\in(1-\theta)T\cap T^{\theta}$ for some $w\in W(G_{\theta},T_{\theta})=W_{R(G,T^{\theta})}$. Therefore,
it suffices to show that \[N_{T}(S)/\tilde{S}^{0}\cong(1-\theta)T\cap T^{\theta}.\] The map \[N_{T}(S)/\tilde{S}^{0}
\longrightarrow (1-\theta)T\cap T^{\theta},\quad t\mapsto t\theta t^{-1}\theta^{-1}=t\theta(t)^{-1}\] gives this isomorphism.
\end{proof}

\begin{proof}[Proof of Theorem \ref{T:twining1}]
Write $d$ for the order of $\theta$. For any $\lambda\in ^{\theta}X^{+}$, let $\sigma_{\lambda}$ denote the representation
of $G\rtimes\langle\theta\rangle$ which is an extension of $V_{\lambda}$ by defining $\theta v_{\lambda}=v_{\lambda}$. Let
$\chi_0$ be a unitary character of $G\rtimes\langle\theta\rangle$ defined by $\chi_{0}|_{G}=1$ and $\chi_{0}(\theta)=
e^{\frac{2\pi\mathbf{i}}{d}}$. It is easy to show the following: let $\sigma$ be an irreducible complex linear
representation of $G\rtimes\langle\theta\rangle$ such that $V_{\mu}\subset\sigma|_{G}$. If $\theta\mu\neq\mu$, then
$\chi_{\sigma}|_{G\theta}=0$; if $\theta\mu=\mu$, then $\sigma\cong\sigma_{\mu}\otimes\chi_{0}^{k}$ for some $k$
($0\leq k\leq d-1$). Applying the orthogonality relation to $\sigma_{\lambda}\otimes\chi_{0}^{k}$ ($0\leq k\leq m-1$),
we see that \[\int_{G\theta}|\chi_{\sigma_{\lambda}}(x)|^{2}\opd x=1,\] where $\opd x$ means a $G$-biinvariant metric on
$G\theta$ with volume 1. Consider $\chi_{\lambda}(t\theta):=\chi_{\sigma_{\lambda}}(t\theta)$ and $\chi'_{\lambda}(t\theta)
:=\tr(\varpi(t):V'_{\lambda}\rightarrow V'_{\lambda})$ ($t\in T^{\theta}$) for $\lambda\in ^{\theta}X^{+}$. By Lemma
\ref{L:twining2}, $\chi'_{\lambda}(t\theta)$ is also $W(G,S)$ invariant. Note that the leading terms of $\chi'_{\lambda}$
and $\chi_{\lambda}$ are both equal to $[\lambda]$. Hence, $\chi'_{\lambda}(t\theta)$ admits an expansion \[\chi'_{\lambda}=\chi_{\lambda}+\sum_{|\mu|<|\lambda|}a_{\mu}\chi_{\mu}.\] Note that, both of $\chi_{\lambda}$
and $\chi'_{\lambda}$ satisfy the orthogonality relation by Lemma \ref{L:AR4}. Then, we must have $\chi'_{\lambda}
=\chi_{\lambda}$. This is just the twining character formula.
\end{proof}

Now assume that $G$ is a quasi-connected compact Lie group. Choose an element $g_{0}$ generating $G/G^0$ and such
that $\Ad(g_{0})|_{\mathfrak{g}_{\der}}$ is a pinned automorphism. Write $\theta=\Ad(g_0)|_{\mathfrak{g}_{0}}.$
Choose a quasi Cartan subgroup $\tilde{S}$ containing $g_{0}$ and write $S=g_{0}\tilde{S}^{0}$. Since
$\theta$ is a pinned automorphism, it stabilizes a positive system $R^{+}(\mathfrak{g},T)$ of
$R(\mathfrak{g},T)$ where $T=Z_{G^{0}}(\tilde{S}^{0})$ is a maximal torus of $G^{0}$ by Lemma \ref{L:T}.
Then, $\theta$ acts on the set of dominant integral weights with respect to $R^{+}(\mathfrak{g},T)$, denoted
by $\Lambda_{T}^{+}$. A restricted root $p(\alpha)\in R(G,\tilde{S}^{0})$ is defined to be positive if
$(\lambda|_{\tilde{S}^{0}},p(\alpha))\geq 0$ for all $\lambda\in(\Lambda_{T}^{+})^{\theta}$. Let
$R^{+}(G,\tilde{S}^{0})$ be the set of positive restricted roots. Write $d=|G/G^0|$.

\begin{lemma}\label{L:Weyl1}
Let $\sigma$ be an irreducible complex linear representation of $G$, and $\lambda$ be the highest weight of some
irreducible subrepresentation of $\sigma|_{G^0}$. If $\theta\cdot\lambda\neq\lambda$, then
$\chi_{\sigma}|_{g_{0}G^{0}}=0$.
\end{lemma}

\begin{proof}
Let $k=\min\{t\geq 1:\theta^{t}\lambda=\lambda\}$. Then $1<k\leq d$ and $k|d$. Put \[G'=\langle G^{0},g_{0}^{k}
\rangle.\] By assumption, $\sigma$ is isomorphic to an irreducible sub-representation of
$\Ind_{G^{0}}^{G}(V_{\lambda})$. By induction in stages, we have $\Ind_{G^{0}}^{G}(V_{\lambda})=
\Ind_{G'}^{G}(\Ind_{G^{0}}^{G'}(V_{\lambda}))$. Since $\theta^{k}\lambda=\lambda$, we have
$(g_{0}^{k})^{\ast} V_{\lambda}\cong V_{\lambda}$. Thus, $V_{\lambda}$ extends to an irreducible representation
of $G'$. Let $\sigma'$ be such an extension. Then, \[\Ind_{G^{0}}^{G'}(V_{\lambda})\cong
\bigoplus_{0\leq j\leq\frac{d}{k}-1}\sigma'\otimes\chi_{0}^{j},\] where $\chi_0$ is a unitary character of $G$
defined by \[\chi_{0}|_{G^0}=1\textrm{ and }\chi_{0}(g_{0})=e^{\frac{2\pi\mathbf{i}}{d}}.\] Without loss of
generality we assume that $\sigma\subset\Ind_{G'}^{G}(\sigma')$. We have \[\Ind_{G'}^{G}(\sigma')|_{G^{0}}
\cong\bigoplus_{0\leq j\leq k-1}(\theta^{\ast})^{j}V_{\lambda}\] and the $k$ irreducible representations
$(\theta^{\ast})^{j}V_{\lambda}$ ($0\leq j\leq k-1$) are pairwise isomorphic and are permuted transitively
by $\langle g_{0}\rangle$. Thus, $\Ind_{G'}^{G}(\sigma')$ is irreducible. Then, $\sigma=\Ind_{G'}^{G}(\sigma')$.
Thus, $\chi_{\sigma}|_{g_{0}G^{0}}=0$.
\end{proof}

By the assumption on $g_0$, we have \[z_{0}:=g_{0}^{d}\in Z(G).\] Let $\lambda\in(\Lambda_{T}^{+})^{\theta}$.
Choose $\eta_{\lambda}\in\mathbb{C}^{\times}$ such that $|\eta_{\lambda}|=1$ and $z_0$ acts by $\eta_{\lambda}^{d}$
on $V_{\lambda}$. Define $\sigma_{\lambda}$ to be the unique extension of $V_{\lambda}$ to $G$ determined by
$g_{0}\cdot v_{\lambda}=\eta_{\lambda}v_{\lambda}$. Then, all other extensions of $V_{\lambda}$ to $G$ are of the
form $\sigma_{\lambda}\otimes\chi_{0}^{j}$ ($0\leq j\leq d-1$). Define $\mu\in X^{\ast}(\tilde{S})$ by $\mu(xg_{0}^{j})
=\lambda(x)$ ($\forall x\in\tilde{S}^0$). Define \[R^{+}(G,\tilde{S})=\{\alpha\in R(G,\tilde{S}):\alpha|_{\tilde{S}^0}\in
R^{+}(G,\tilde{S}^{0})\}.\] Put \[\delta_{g_0}=\frac{1}{2}\sum_{\alpha\in R^{+}(G,\tilde{S})}\alpha.\]

Write $R=R(G,\tilde{S})$ for simplicity. The following lemma follows from Theorem \ref{T:twining1} and the Weyl
character formula.

\begin{lemma}\label{L:twining3}
Let $\lambda\in(\Lambda_{T}^{+})^{\theta}$. Then $\chi'_{\lambda}:=\chi_{\sigma_{\lambda}}|_{S}$ satisfies \begin{equation}\label{Eq:chi'2}\chi'_{\lambda}\prod_{\alpha\in R^{+}(G,\tilde{S})}(1-[-\alpha])=\eta_{\lambda}
\sum_{w\in W_{R(G,\tilde{S})_{g_0}}}\epsilon(w)[w\mu+w\delta_{g_{0}}-\delta_{g_0}].
\end{equation}
\end{lemma}

Let $\sigma$ be an irreducible complex linear representation of $G$ such that $\sigma|_{G^0}$ is irreducible. By
Lemma \ref{L:twining3}, the twisted character $\chi_{\sigma}|_{S}$ differs from the pull-back of the character
of an irreducible representation of a smaller group (twining group) with root system $\{m_{\alpha}\alpha:\alpha
\in R(G,\tilde{S})_{g_0}\}$ by a scalar. Therefore, there is a constant $\eta\neq 0$ such that \[\chi_{\sigma}|_{S}
=\eta\sum_{\mu}b_{\mu}[\mu],\] where all $b_{\mu}>0$, $\mu\in X^{\ast}(\tilde{S})$ with $[\mu](g_{0})=1$ and such
that $\mu_{\tilde{S}^{0}}$ are distinct. Put \[A_{\lambda,R}=\eta_{\lambda}\sum_{w\in W_{R_{g_{0}}}}
\epsilon(w)[\mu+\delta_{g_{0}}-w\delta_{g_{0}}].\]

\begin{lemma}\label{L:chi-D}
We have \begin{equation}\label{Eq:chi'3}\chi_{\sigma_{\lambda}}(x)D_{G}(x)=\frac{1}{|W_{R_{g_{0}}}|}
\sum_{w\in W_{R_{g_{0}}}}w\cdot A_{\lambda,R},\ \forall x\in S.\end{equation}
\end{lemma}

\begin{proof}
By Lemma \ref{L:twining3}, it reduces to show that \begin{eqnarray*}&&\sum_{\tau\in W_{R_{g_{0}}}}
\epsilon(\tau)[\tau\mu+\tau\delta_{g_{0}}]\sum_{w\in W_{R_{g_{0}}}}\epsilon(w)[-w\delta_{g_{0}}]\\&=&
\sum_{\tau\in W_{R_{g_{0}}}}\tau\cdot\sum_{w\in W_{R_{g_{0}}}}\epsilon(w)[\mu+\delta_{g_{0}}-w\delta_{g_{0}}].
\end{eqnarray*} This follows from the same calculation as that in \cite[Prop. 3.7]{Yu-dimension}.
\end{proof}

\section{Irreducible restriction characters}\label{S:Res-character}

Let $\Phi'$ be a root system on a torus $T$. Write \[\chi_{\Phi'}(\lambda)=\chi_{V_{\lambda}}|_{T}\] for the
character of an irreducible representation with highest weight $\lambda$ of a connected compact Lie group with
$T$ a maximal torus and with root system $\Phi'$. Let $\Phi$ be a quasi root system on a quasi torus $\tilde{S}$
which could be obtained from the action of a diagram automorphism $\theta$ of a root system $\Phi'$ on a torus
$T$. Write \[\chi_{\Phi}(\lambda)=\chi_{\Phi'}(\lambda)|_{\tilde{S}^{0}}\] and call it an {\it irreducible
restriction character} associated to the quasi root system $\Phi$ and the highest weight $\lambda$.
For a weight $\mu'\in X^{\ast}(T)$, write $m_{\Phi',\lambda}(\mu')$ for the coefficient of the linear
character $[\mu']$ appearing in the irreducible character $\chi_{\Phi'}(\lambda)$. Put
\[m_{\Phi'}(\lambda,\mu')=m_{\Phi',\lambda}(\lambda-\mu').\] For a weight $\mu\in X^{\ast}(\tilde{S}^{0})$,
write $m_{\Phi,\lambda}(\mu)$ for the coefficient of the linear character $[\mu]$ appearing in the irreducible
restriction character $\chi_{\Phi}(\lambda)$. Put \[m_{\Phi}(\lambda,\mu)=
m_{\Phi,\lambda}(\lambda|_{\tilde{S}^{0}}-\mu).\]

\subsection{The Levi reduction method}\label{SS:Levi}

Let $\{\bar{\alpha}_{i}:1\leq i\leq l\}$ be a simple system of $p(\Phi)^{\nd}$ and let $\{\bar{\omega}_{i}:
1\leq i\leq l\}$ be the corresponding set of fundamental weights.

\begin{definition}\label{D:LeviWeight2}
Let $\mu=\sum_{1\leq i\leq l}b_{i}\bar{\alpha}_{i}\in X^{\ast}(\tilde{S}^{0})$ be a weight with all
$b_{i}\in\mathbb{Q}_{\geq 0}$.
\begin{enumerate}
\item[(i)]Let $\supp\mu$ be the set of simple roots $\bar{\alpha}_{i}$ such that $b_{i}>0$.
\item[(ii)]Let $\Phi_{\mu}$ be the quasi sub-root system of $\Phi$ consisting of roots $\alpha\in\Phi$ such
that $\alpha|_{\tilde{S}^{0}}$ is a linear combination of simple roots in $\supp\mu$.
\item[(ii)]Let $\Phi'_{\mu}$ be the sub-root system of $\Phi'$ consisting of roots $\alpha\in\Phi'$ such that
$\alpha|_{\tilde{S}^{0}}$ is a linear combination of simple roots in $\supp\mu$.
\end{enumerate}
\end{definition}

\begin{lemma}\label{L:Levi1}
For any weight $\mu\in X^{\ast}(\tilde{S}^{0})$, we have $m_{\Phi}(\lambda,\mu)=m_{\Phi_{\mu}}(\lambda,\mu)$.
\end{lemma}

\begin{proof}
We have \[m_{\Phi}(\lambda,\mu)=\sum_{\mu'\in X^{\ast}(T),\mu'|_{\tilde{S}^{0}}=\mu}m_{\Phi',\lambda}(\lambda-\mu')\]
and \[m_{\Phi_{\mu}}(\lambda,\mu)=\sum_{\mu'\in X^{\ast}(T),\mu'|_{\tilde{S}^{0}}=\mu}
m_{\Phi'_{\mu},\lambda}(\lambda-\mu').\] It suffices to show that: \[m_{\Phi',\lambda}(\lambda-\mu')=
m_{\Phi'_{\mu},\lambda}(\lambda-\mu')\] for any weight $\mu'\in X^{\ast}(T)$ such that $\mu'|_{\tilde{S}^{0}}=\mu$.
Recall Kostant's multiplicity formula \begin{eqnarray}\label{Eq:Kostant}m_{\Phi',\lambda}(\lambda-\mu')=
\sum_{w\in W_{\Phi'}}\mathcal{P}_{\Phi'}(w(\lambda+\delta_{\Phi'})-(\lambda+\delta_{\Phi'}-\mu'))\end{eqnarray}
where $\mathcal{P}_{\Phi'}(\cdot)$ is Kostant 's partition function associated to the positive system $\Phi'^{+}$.
When $w\not\in W_{\Phi'_{\mu}}$, the expansion of $(\lambda+\delta_{\Phi'})-w(\lambda+\delta_{\Phi'})$ into a sum
of simple roots in $\Phi_{s}^{+}$ must contain a term $b\alpha$ with $b>0$ for some simple root $\alpha\not\in
\Phi'_{\mu}$. Then, $\mathcal{P}_{\Phi'}(w(\lambda+\delta_{\Phi'})-(\lambda+\delta_{\Phi'}-\mu'))=0$. Thus, only
those terms corresponding to elements $w\in W_{\Phi'_{\mu}}$ contribute to \eqref{Eq:Kostant}. By Kostant's
multiplicity formula for $m_{\Phi'_{\mu},\lambda}(\lambda-\mu')$, these terms have sum equal to
$m_{\Phi'_{\mu},\lambda}(\lambda-\mu')$. Therefore, \[m_{\Phi',\lambda}(\lambda-\mu')=
m_{\Phi'_{\mu},\lambda}(\lambda-\mu').\qedhere\]
\end{proof}


For a subset $I$ of $\{1,\dots,l\}$, let $\Phi_{I}$ be the quasi sub-root system of $\Phi$ consisting of roots
$\alpha\in\Phi$ such that $\supp\alpha|_{\tilde{S}^{0}}\subset\{\bar{\alpha}_{i}: i\in I\}$. Let
$\chi_{\Phi,I}(\lambda)$ be the sum of terms $m_{\Phi}(\lambda,\mu)[\lambda|_{\tilde{S}^{0}}-\mu]$ in
$\chi_{\Phi}(\lambda)$ where $\supp\mu\subset\{\bar{\alpha}_{i}: i\in I\}$. The following lemma follows from
the Lemma \ref{L:Levi1} directly.

\begin{lemma}\label{L:Levi2}
We have \[\chi_{\Phi,I}(\lambda)=\chi_{\Phi_{I}}(\lambda).\]
\end{lemma}

Write $\{\beta_{i}:1\leq i\leq l'\}$ for a set of simple roots of $\Phi'$ and write $\{\omega_{i}:1\leq i\leq l'\}$
for the corresponding set of fundamental weights. Put \[\delta_{\Phi'}=\sum_{1\leq i\leq l'}\omega_{i}.\] For each
$i$ ($1\leq i\leq l'$), let $\sqrt{k_{i}}$ be the ratio of the length of $\beta_{i}$ and the length of a short root
in the irreducible factor of $\Phi'$ containing $\beta_{i}$. For each root $\alpha=\sum_{1\leq i\leq l'}b_{i}(\alpha)
\beta_{i}\in\Phi'$, define a linear function \[f_{\alpha}(t_{1},\dots,t_{l'})=\sum_{1\leq i\leq l'}b_{i}(\alpha)
k_{i}t_{i}.\] Define \[f_{\Phi'}=\prod_{\beta\in\Phi'^{+}}f_{\beta}.\] The following lemma follows from the
Weyl dimension formula.

\begin{lemma}\label{L:Levi3}
For a highest weight $\lambda=\sum_{1\leq i\leq l'}a_{i}\omega_{i}$ ($a_{i}\in\mathbb{Z}_{\geq 0}$), we have
\[\dim\chi_{\Phi}(\lambda)=\frac{f_{\Phi'}(a_{1}+1,\dots,a_{l'}+1)}{f_{\Phi'}(\underbrace{1,\dots,1}_{l'})}.\]
\end{lemma}

We call the equality in Lemma \ref{L:Levi2} the Levi reduction equality. This method of Levi reduction will be used
effectively together with the dimension formula for $\chi_{\Phi}(\lambda)$ in Lemma \ref{L:Levi3}. Note that from
explicit expression of positive roots in each irreducible factor of $\Phi'$ (cf. \cite[Appendix C.2]{Knapp}), one
obtains a precise expression of $f_{\Phi'}$.

\begin{lemma}\label{L:Levi4}
Fix a root system $\Phi'$ on a torus $T$ and an automorphism $\theta$ of $\Phi$ which stabilizes a positive system
$\Phi^{'+}$. Assume that the order of $\theta$ is equal to 1 or 2 and the resulting quasi root system $\Phi$ is
irreducible. Then the restriction character $\chi_{\Phi}(\lambda)$ determines the highest weight $\lambda$ up to
the action of $\theta$.
\end{lemma}

\begin{proof}
If $o(\theta)=1$, then $\Phi=\Phi'$ and the conclusion is trivial. Assume that $o(\theta)=2$ below. Then, the pair
$(\Phi',\Phi)$ is isomorphic to one in the following list: \begin{enumerate}
\item[(i)]$(\Phi_{0}\sqcup\Phi_{0},\Phi_{0}^{(2)})$ with $\Phi_{0}$ an irreducible root system;
\item[(ii)]$(\A_{2l-1},\C_{l}^{(2,1)})$, $l\geq 2$;
\item[(iii)]$(\A_{2l},\BC_{l}^{(2,2,1)})$, $l\geq 1$;
\item[(iv)]$(\D_{l+1},\B_{l}^{(2,1)})$, $l\geq 3$;
\item[(v)]$(\E_{6},\F_{4}^{(2,1)})$.
\end{enumerate}

When $(\Phi',\Phi)=(\A_{1}\sqcup\A_{1},\A_{1}^{(2)})$, the conclusion follows from the tensor product formula  \[\chi_{\A_{1}}(a\bar{\omega}_{1})\chi_{\A_{1}}(b\bar{\omega}_{1})=\sum_{0\leq j\leq\min\{a,b\}}
\chi_{\A_{1}}((a+b-2j)\bar{\omega}_{1}).\]

When $(\Phi',\Phi)=(\Phi_{0}\sqcup\Phi_{0},\Phi_{0}^{(2)})$ with $\Phi_{0}$ an irreducible root system, define a
distance function $d(\cdot,\cdot)$ on the index set $\{1,\dots,l\}$ by the length of shortest path connecting
corresponding vertices in the Dynkin diagram of $\Phi_{0}$. Write $\lambda=(\lambda_{1},\lambda_{2})$ and $\lambda
=(\lambda'_{1},\lambda'_{2})$, where \[\lambda_{1}=\sum_{1\leq i\leq l}a_{i}\bar{\omega}_{i},\ \lambda_{2}=
\sum_{1\leq i\leq l}b_{i}\bar{\omega}_{i}\] and \[\lambda'_{1}=\sum_{1\leq i\leq l}a'_{i}\bar{\omega}_{i},
\ \lambda'_{2}=\sum_{1\leq i\leq l}b'_{i}\bar{\omega}_{i}.\] Suppose that \begin{equation}\label{Eq:Phi0}
\chi_{\Phi_{0}}(\lambda_{1})\chi_{\Phi_{0}}(\lambda_{2})=\chi_{\Phi_{0}}(\lambda'_{1})\chi_{\Phi_{0}}(\lambda'_{2}),
\end{equation} $(\lambda_{1},\lambda_{2})\neq(\lambda'_{1},\lambda'_{2})$ and $(\lambda_{1},\lambda_{2})\neq
(\lambda'_{2},\lambda'_{1})$. Reducing to a Levi factor of rank one using Lemma \ref{L:Levi2}, we get
$\{a_{i},b_{i}\}=\{a'_{i},b'_{i}\}$ for each $i$ from \eqref{Eq:Phi0}. Let $I_{1}$ be the set of indices $i\in I$
such that $(a_{i}-b_{i})(a'_{i}-b'_{i})>0$ and let $I_{2}$ be the set of indices $i\in I$ such that
$(a_{i}-b_{i})(a'_{i}-b'_{i})<0$. Due to $(\lambda_{1},\lambda_{2})\neq(\lambda'_{1},\lambda'_{2})$ and
$(\lambda_{1},\lambda_{2})\neq(\lambda'_{2},\lambda'_{1})$, we have $I_{1}\neq\emptyset$ and $I_{2}\neq
\emptyset$. Choose $i_{1}\in I_{1}$ and $i_{2}\in I_{2}$ such that \[d(i_{1},i_{2})=
\min_{(i'_{1},i'_{2})\in I_{1}\times I_{2}}d(i'_{1},i'_{2}).\] Then, any other index $i$ such that
$\bar{\alpha}_{i}$ lies in the unique shortest path connecting $\bar{\alpha}_{i_{1}}$ and $\bar{\alpha}_{i_{2}}$
is contained in neither $I_{1}$ nor $I_{2}$. Hence, $a_{i}=b_{i}=a'_{i}=b_{i}$. Let $I$ be the set of indices $i$
such that $\bar{\alpha}_{i}$ lies in the unique shortest path connecting $\bar{\alpha}_{i_{1}}$ and
$\bar{\alpha}_{i_{2}}$ and let $\Phi_{I}$ be the quasi sub-root system consisting of roots $\alpha\in\Phi$ such
that $\supp p(\alpha)\subset\{\bar{\alpha}_{i},i\in I\}$. Without loss of generality we assume that
$a_{i_{1}}>b_{i_{1}}$ and $a_{i_{2}}>b_{i_{2}}$. By the dimension formula in Lemma \ref{L:Levi3} it follows that $\dim\chi_{\Phi_{I}}(\lambda)<\dim\chi_{\Phi_{I}}(\lambda')$. Then, by Lemma \ref{L:Levi2} we get
$\chi_{\Phi}(\lambda)\neq\chi_{\Phi}(\lambda')$, which is a contradiction.

When $(\Phi',\Phi)=(\A_{2l-1},\C_{l}^{(2,1)})$, define a distance function $d(\cdot,\cdot)$ on the index set
$\{1,\dots,l\}$ by the length of shortest path connecting corresponding vertices in the Dynkin diagram of $\C_{l}$.
Write \[\lambda=\sum_{1\leq i\leq 2l-1}a_{i}\omega_{i}\quad\textrm{and}\quad\lambda'=\sum_{1\leq i\leq 2l-1}
b_{i}\omega_{i}.\] Suppose that \begin{equation}\label{Eq:Cm}\chi_{\Phi}(\lambda)=\chi_{\Phi}(\lambda'),
\end{equation} $\lambda\neq\lambda'$ and $\lambda\neq\theta\lambda'$. Reducing to a Levi factor of rank one using
Lemma \ref{L:Levi2}, we get $\{a_{i},a_{2l-i}\}=\{b_{i},b_{2l-i}\}$ for each $i$ ($1\leq i\leq l-1$) from
\eqref{Eq:Cm}. Let $I_{1}$ be the set of indices $i$ ($1\leq i\leq l-1$) such that
$(a_{i}-a_{2l-i})(b_{i}-b_{2l-i})>0$ and let $I_{2}$ be the set of indices $i$ ($1\leq i\leq l-1$) such that $
(a_{i}-a_{2l-i})(b_{i}-b_{2l-i})<0$. Due to $\lambda\neq\lambda'$, we have $I_{2}\neq\emptyset$; due to
$\lambda\neq\theta\lambda'$, we have $I_{1}\neq\emptyset$. The remaining argument is the same as that for
$(\Phi_{0}\sqcup\Phi_{0},\Phi_{0}^{(2)})$.

The proof for each case of $(\A_{2l},\BC_{l}^{(2,2,1)})$, $(\D_{l+1},\B_{l}^{(2,1)})$, $(\E_{6},\F_{4}^{(2,1)})$ is
similar to the case of $(\A_{2l-1},\C_{l}^{(2,1)})$.
\end{proof}

\subsection{Extremal weights in an irreducible restriction character}\label{SS:leading}

Put \[d=o(\theta).\] Then, the map \[\lambda\mapsto\frac{1}{d}\sum_{0\leq k\leq d-1}\theta^{k}\lambda,
\ \forall\lambda\in X^{\ast}(T)\] gives an isomorphism \[X^{\ast}((T^{\theta})^{0})\otimes_{\mathbb{Z}}\mathbb{Q}
\cong X^{\ast}(T)^{\theta}\otimes_{\mathbb{Z}}\mathbb{Q}.\] For this reason, we could regard $p(\Phi)$ as
a subset of $X^{\ast}((T^{\theta})^{0})\otimes_{\mathbb{Z}}\mathbb{Q}$.

\begin{lemma}\label{L:RS1}
With the above identification of $p(\Phi)$ as a subset of $X^{\ast}(T)^{\theta}\otimes_{\mathbb{Z}}
\mathbb{Q}$, we have $W_{p(\Phi)}=W_{\Phi'}^{\theta}$ when both acting on $X^{\ast}(T)^{\theta}
\otimes_{\mathbb{Z}}\mathbb{R}$.
\end{lemma}

\begin{proof}
Write \[\phi: W_{\Phi'}^{\theta}\rightarrow\GL(X^{\ast}(T)^{\theta}\otimes_{\mathbb{Z}}\mathbb{R})\] for the
action of $W_{\Phi'}^{\theta}$ on $X^{\ast}(T)^{\theta}\otimes_{\mathbb{Z}}\mathbb{R}$. Put \[\rho=\frac{1}{2}
\sum_{\alpha\in\Phi^{'+}}\alpha.\] Then, $\rho$ is regular dominant with respect to $\Phi^{'+}$. Hence, \[\Stab_{W_{\Phi'}^{\theta}}(\rho)=\Stab_{W_{\Phi'}}(\rho)=1.\] It is clear that $\rho\in X^{\ast}(T)^{\theta}
\otimes_{\mathbb{Z}}\mathbb{Q}$. Hence, $\phi$ is an injective map.

For any restricted root $\bar{\alpha}\in p(\Phi)^{\nd}$, there exists a root $\alpha\in\Phi'$ such that
$\alpha|_{(T^{\theta})^{0}}=\bar{\alpha}$ or $2\bar{\alpha}$ and any two roots in the set
$\{\theta^{j}\alpha:0\leq j\leq d-1\}$ are either equal or orthogonal. Let $k_{0}$ be the least $k\geq 1$ such
that $\theta^{k}\alpha=\alpha$. Then, the roots $\{\theta^{j}\alpha:0\leq j\leq k_{0}-1\}$ are orthogonal to
each other and we have \begin{equation}\label{Eq:theta-reflection}s_{\bar{\alpha}}=
\phi(\prod_{0\leq j\leq k_{0}-1}s_{\theta^{j}\alpha})\end{equation} while both acting on
$X^{\ast}(T)^{\theta}\otimes_{\mathbb{Z}}\mathbb{R}$. Thus, $s_{\bar{\alpha}}\in W_{\Phi'}^{\theta}$.
Hence, $W_{p(\Phi)}\subset W_{\Phi'}^{\theta}$. It is clear that the action of $W_{\Phi'}^{\theta}$ on
$X^{\ast}(T)^{\theta}\otimes_{\mathbb{Z}}\mathbb{R}$ stabilizes $p(\Phi)$. Let $W'$ be the group of
elements $w\in W_{\Phi'}^{\theta}$ which stabilize $p(\Phi)^{+}$. Then, \[W_{\Phi'}^{\theta}=W_{p(\Phi)}
\rtimes W'.\] It suffices to show that $W'=1$. To show this without loss of generality we assume that
$p(\Phi)$ is irreducible. When $p(\Phi)$ is non-simply laced, it follows from the fact that there is
no nontrivial automorphism of $p(\Phi)$ stabilizing $p(\Phi)^{+}$. When $p(\Phi)$ is simply-laced,
$\Phi'$ is isomorphic to the orthogonal union of $d$ copies of $p(\Phi)$ and $\theta$ acts on it by a
permutation on irreducible factors of it. Then, it is clear that $W'=1$.
\end{proof}

For any two vectors $\lambda,\mu\in X^{\ast}(T)\otimes_{\mathbb{Z}}\mathbb{R}$, we define $\lambda>\mu$ if
$\lambda\neq\mu$ and \[\lambda-\mu=\sum_{\alpha\in\Phi^{'+}}a_{\alpha}\alpha\] where all coefficients
$a_{\alpha}\in\mathbb{R}_{\geq 0}$. For any two vectors $\lambda,\mu\in X^{\ast}(T)^{\theta}
\otimes_{\mathbb{Z}}\mathbb{R}$, we define $\lambda>\mu$ if $\lambda\neq\mu$ and \[\lambda-\mu=
\sum_{\bar{\alpha}\in p(\Phi)^{+}}a_{\bar{\alpha}}\bar{\alpha}\] where all coefficients $a_{\alpha}\in
\mathbb{R}_{\geq 0}$. The following lemma is clear.

\begin{lemma}\label{L:RS2}
If $\lambda>\mu$, then $\lambda|_{(T^{\theta})^{0}}>\mu|_{(T^{\theta})^{0}}$.
\end{lemma}

\begin{lemma}\label{L:RS3}
Let $V_{\Phi',\lambda}$ be an irreducible character of $\Phi'$ with highest weight $\lambda$. Write
$\supp V_{\Phi',\lambda}$ for the set of weights appearing in $V_{\Phi',\lambda}$. Then, longest vectors among
$\{\mu|_{(T^{\theta})^{0}}:\mu\in\supp V_{\Phi',\lambda}\}$ are in the set $\{p(w\lambda):w\in W_{\Phi'}^{\theta}\}$.
\end{lemma}

\begin{proof}
It suffices to show that: for any weight $\mu\in\supp V_{\Phi',\lambda}$, we have $|\mu|_{(T^{\theta})^{0}}|\leq
|\lambda|_{(T^{\theta})^{0}}|$ and the equality holds only when $\mu=w\lambda$ for some $w\in W_{\Phi'}^{\theta}$.

Write $\bar{\lambda}=\lambda|_{(T^{\theta})^{0}}$ and $\bar{\mu}=\mu|_{(T^{\theta})^{0}}$. By Lemma \ref{L:RS1},
replacing $\mu$ by $w^{-1}\mu$ for some $w\in W_{\Phi'}^{\theta}$ if necessary, we may assume that $\bar{\mu}$ is
dominant with respect to $p(\Phi)^{+}$. Since $\lambda\geq\mu$, we have $\bar{\lambda}\geq\bar{\mu}$ by Lemma
\ref{L:RS2}. Then, \begin{eqnarray*}&&(\bar{\lambda},\bar{\lambda})-(\bar{\mu},\bar{\mu})\\&=&
(\bar{\lambda}-\bar{\mu},\bar{\lambda}-\bar{\mu})+2(\bar{\mu},\bar{\lambda}-\bar{\mu})\\&\geq&0.\end{eqnarray*}
Thus, $|\bar{\lambda}|\geq|\bar{\mu}|$ and the equality holds only when $\bar{\mu}=\bar{\lambda}$. By Lemma
\ref{L:RS2} again, $\bar{\mu}=\bar{\lambda}$ if and only if $\mu=\lambda$.
\end{proof}

\subsection{Degree 1 irreducible restriction characters}\label{SS:degree-one}

We introduce some notations. \begin{itemize}
\item[(1)]For an integer $a\geq 0$, write $\tau_{a}=\chi_{\B_{1}}(\frac{a}{2})$ for an irreducible character of
the root system $\B_{1}$ with dimension $a+1$. Then, \[\tau_{a}=\sum_{0\leq j\leq a}[\frac{a}{2}-j]=[\frac{a}{2}]
+[\frac{a}{2}-1]+\cdots+[-\frac{a}{2}].\]
\item[(2)]For integers $a,b\geq 0$, write \[\tau_{a,b}=\chi_{\B_{1}^{(2)}}((\frac{a}{2}),(\frac{b}{2}))=
\tau_{a}\cdot\tau_{b}=\sum_{0\leq j\leq\min\{a,b\}}\tau_{a+b-2j}.\]
\item[(3)]For $a,b\geq 0$, write \[\tau'_{a,b}=\chi_{\BC_{1}^{(2,1)}}(a,0,-b),\] where $(a,0,-b)$ means a highest
weight for the root system $\A_{2}$. Then, $\tau'_{a,b}=\tau'_{b,a}$.
\end{itemize}

\begin{lemma}\label{L:tau'}
Let $a\geq b\geq 0$. Then \begin{equation}\label{Eq:tau3}\tau'_{a,b}=\sum_{0\leq k\leq a+b}m_{k}\tau_{2a+2b-2k}
\end{equation} and the multiplicities $m_{k}$ are give as follows. \begin{itemize}
\item[(i)]When $k$ is even, \begin{equation}\label{Eq:mk-1}m_{k}=\left\{\begin{array}{cccc}\frac{k}{2}+1\textrm{ if }
0\leq\frac{k}{2}\leq\lfloor\frac{b}{2}\rfloor\\\lfloor\frac{b}{2}\rfloor+1\textrm{ if }\lfloor\frac{b}{2}\rfloor
\leq\frac{k}{2}\leq\lfloor\frac{a}{2}\rfloor\\\lfloor\frac{a}{2}\rfloor+\lfloor\frac{b}{2}\rfloor-\frac{k}{2}+1
\textrm{ if }\lfloor\frac{a}{2}\rfloor\leq\frac{k}{2}\leq\lfloor\frac{a}{2}\rfloor+\lfloor\frac{b}{2}\rfloor\\0
\quad\textrm{otherwise}.\end{array}\right.\end{equation}
\item[(ii)]When $k$ is odd, \begin{equation}\label{Eq:mk-2}m_{k}=\left\{\begin{array}{ccc}\frac{k-1}{2}+1\textrm{ if }
0\leq\frac{k-1}{2}\leq\lfloor\frac{b-1}{2}\rfloor\\\lfloor\frac{b-1}{2}\rfloor+1\textrm{ if }\lfloor\frac{b-1}{2}
\rfloor\leq\frac{k-1}{2}\leq\lfloor\frac{a-1}{2}\rfloor\\\lfloor\frac{a-1}{2}\rfloor+\lfloor\frac{b-1}{2}\rfloor
-\frac{k-1}{2}+1\textrm{ if }\lfloor\frac{a-1}{2}\rfloor\leq\frac{k-1}{2}\leq\lfloor\frac{a-1}{2}\rfloor+
\lfloor\frac{b-1}{2}\rfloor\\0\quad\textrm{otherwise}.\end{array}\right.\end{equation}
\end{itemize}
\end{lemma}

\begin{proof}
Note that \begin{equation}\label{Eq:Symk}\tau'_{a,0}=\Sym^{a}\tau_{2}=\sum_{0\leq j\leq\lfloor\frac{a}{2}\rfloor}
\tau_{2a-4j}\end{equation} for any $a\geq 0$. When $a\geq b\geq 1$, we have \[\tau'_{a,b}=\Sym^{a}\tau_{2}\Sym^{b}
\tau_{2}-\Sym^{a-1}\tau_{2}\Sym^{b-1}\tau_{2}.\] When $a$ and $b$ are both odd, using \eqref{Eq:Symk} we get \begin{eqnarray*}&&\tau'_{a,b}\\&=&\sum_{0\leq i\leq\lfloor\frac{a}{2}\rfloor}\tau_{2a-4i}\sum_{0\leq j\leq
\lfloor\frac{b}{2}\rfloor}\tau_{2b-4j}-\sum_{0\leq i\leq\lfloor\frac{a-1}{2}\rfloor}\tau_{2a-2-4i}
\sum_{0\leq j\leq\lfloor\frac{b-1}{2}\rfloor}\tau_{2b-2-4j}\\&=&\sum_{0\leq i\leq\lfloor\frac{a-1}{2}\rfloor,
0\leq j\leq\lfloor\frac{b-1}{2}\rfloor}(\tau_{2a+2b-4(i+j)}+\tau_{2a+2b-4(i+j)-2})\\&=&\sum_{0\leq k\leq a+b}
m_{k}\tau_{2a+2b-2k},\end{eqnarray*} where the multiplicities $m_{k}$ are given as in \eqref{Eq:mk-1} and
\eqref{Eq:mk-2}. When $(a,b)\equiv(0,0)$, $(1,0)$ or $(0,1)\pmod{2}$, we make similar calculation using
\eqref{Eq:Symk}.
\end{proof}

\begin{lemma}\label{L:tau'-1}
\begin{enumerate}
\item[(i)]The irreducible restriction character $\tau'_{a,b}$ is a multiplicity-free sum of $\tau_{k}$ if and only if
$\min\{a,b\}=0$ or 1.
\item[(ii)]In order that $\tau'_{a,b}=\tau_{c}\tau_{d}$ for some $c,d\in\mathbb{Z}_{\geq 0}$ it is necessary and
sufficient that $(a,b)=(0,0)$, $a=1$ or $b=1$.
\item[(iii)]Let $c\geq d\geq 0$. In order that $\tau_{c}\tau_{d}=\tau'_{a,b}$ for some $a,b\in\mathbb{Z}_{\geq 0}$ it is
necessary and sufficient that $(c,d)=(0,0)$ or $c-d=2$.
\end{enumerate}
\end{lemma}

\begin{proof}
(i)Without loss of generality we assume that $a\geq b$. When $b=0$, by \eqref{Eq:Symk} we see that $\tau'_{a,b}$ is a
multiplicity-free sum of $\tau_{k}$. When $b=1$, by \eqref{Eq:mk-1} and \eqref{Eq:mk-2} we have \begin{equation}
\label{Eq:Symk-2}\tau_{a,1}=\sum_{0\leq j\leq a}\tau_{2a+2-2j}.\end{equation} Thus, $\tau'_{a,1}$ is a multiplicity-free
sum of $\tau_{k}$. When $b\geq 2$, by \eqref{Eq:mk-1} and \eqref{Eq:mk-2} we have \[\max\{m_{k}\}=\lfloor\frac{b}{2}
\rfloor+1\geq 2.\] Thus, $\tau'_{a,b}$ is not a multiplicity-free sum of $\tau_{k}$.

(ii)Suppose that $a\geq b$ and $\tau'_{a,b}=\tau_{c}\tau_{d}$. Note that \[\tau_{c}\tau_{d}=\sum_{0\leq j\leq\min\{c,d\}}
\tau_{c+d-2j}\] is a multiplicity-free sum of $\tau_{k}$. By (1) we get $b=0$ or $1$. When $b=0$, by \eqref{Eq:Symk}
we must have $a=0$ or 1. When $b=1$, by \eqref{Eq:Symk-2} we have \[\tau'_{a,1}=\tau_{a+2}\tau_{a}.\] Then, we get the
conclusion.

(iii)It follows from (ii).
\end{proof}

\subsection{Fully decomposable irreducible restriction characters}\label{SS:fully}

Fix a reduced quasi root system $\Psi$ on a quasi torus $\tilde{S}$ and let $\chi\in\mathbb{Z}[X^{\ast}(\tilde{S}^{0})]$
be a fixed $W_{p(\Psi)}$ invariant character on $\tilde{S}^{0}$ which is not a scalar. Choose a positive system
$p(\Psi)^{+}$ of $p(\Psi)$. We intend to study pairs $(\Phi,\lambda)$ with the following properties:
\begin{enumerate}
\item[(1)]$\Phi$ is a quasi sub-root system of $\Psi$ which is obtained from the action of a diagram automorphism
$\theta$ of a root system $\Phi'$;
\item[(2)]$p(\Phi)\supset p(\Psi)^{\circ}$;
\item[(3)]$\lambda$ is a dominant integral weight of the root system $\Phi'$ such that
\begin{equation}\label{Eq:TA-chi}\chi_{\Phi}(\lambda)=\chi.\end{equation}
\end{enumerate}
Moreover, we make the following assumptions. \begin{assumption}\label{A:TA-rho}
\begin{enumerate}
\item[(4)]$\Psi$ is isomorphic to $\B_{n}^{(2)}$ or $\BC_{n}^{(2,2,1)}$ where $n\geq 2$.
\item[(5)]$\theta^{\ast}\chi_{\Phi'}(\lambda)=\chi_{\Phi'}(\lambda)^{\ast}$, where $\chi_{\Phi'}(\lambda)^{\ast}$
is the dual (contragredient) character of $\chi_{\Phi'}(\lambda)$.
\end{enumerate}
\end{assumption}

Clearly the condition in (5) is equivalent to \begin{equation}\label{Eq:TA-lambda}\theta^{\ast}\lambda
=-w_{l}\lambda,\end{equation} where $w_{l}$ is the longest element in the Weyl group of $\Phi'$. We identify
$p(\Psi)$ with \[\B_{n}=\{\pm{e_{i}}\pm{e_{j}},\pm{e_{k}}:1\leq i<j\leq n,1\leq k\leq n\}\] or \[\BC_{n}=
\{\pm{e_{i}}\pm{e_{j}},\pm{e_{k}},\pm{2e_{k}}:1\leq i<j\leq n, 1\leq k\leq n\},\] where $e_{1},\dots,e_{n}$ is an
orthonomal basis. Set \[\B_{n}^{+}=\{e_{i}\pm{e_{j}},e_{k}:1\leq i<j\leq n,1\leq k\leq n\}\] and \[\BC_{n}^{+}=
\{e_{i}\pm{e_{j}},e_{k},2e_{k}:1\leq i<j\leq n, 1\leq k\leq n\}.\] Write $I_{0}=\{1,2,\dots,n\}$. For a subset $I$
of $I_{0}$, let $p(\Psi)_{I}$ be the sub-root system of $p(\Psi)$ consisting of roots
$\pm{e_{i}}\pm{e_{j}},\pm{e_{k}},\pm{2e_{k}}$ ($i<j$) such that $i,j,k\in I$. For a partition
$I_{0}=I_{1}\sqcup I_{2}\sqcup\cdots\sqcup I_{s}$ and characters $\chi_{i}$ of $p(\Psi)_{I_{i}}$ ($1\leq i\leq s$),
we can form a {\it tensor product} $\chi_{1}\otimes\chi_{2}\otimes\cdots\otimes\chi_{s}$. Given a nonzero character
$\chi$ of $p(\Psi)$, we call $\chi=\chi_{1}\otimes\chi_{2}\otimes\cdots\otimes\chi_{s}$ a maximal decomposition
if there exists no proper refinement $I_{0}=I'_{1}\sqcup I'_{2}\sqcup\cdots\sqcup I'_{t}$ and corresponding
decomposition $\chi=\chi'_{1}\otimes\chi'_{2}\otimes\cdots\otimes\chi'_{t}$ such that each $\chi_{i}$ is a tensor
product of those $\chi'_{j}$ where $I'_{j}\subset I_{i}$.

\begin{definition}\label{D:chi-support}
For a nonzero character $\chi$ of $p(\Psi)$, we call the minimal index set $I\subset I_{0}$ such that each weight
appearing in $\chi$ is an integral linear combination of $e_{i}$ ($i\in I$) the {\it support of $\chi$}.
\end{definition}

\begin{lemma}\label{L:decom1}
Let $\chi$ be a non-scalar character on $\tilde{S}^{0}$. Then there exists a unique maximal decomposition
$\chi=\chi_{1}\otimes\cdots\otimes\chi_{s}$ up to scalars.
\end{lemma}

\begin{proof}
Without loss of generality we assume that $\supp\chi=I_{0}$. Let \[\chi=\chi_{1}\otimes\cdots\otimes\chi_{s}\] be
a maximal decomposition. Then, each $\chi_{i}$ is not a scalar. It suffices to show that: for any decomposition
$\chi=\chi'_{1}\otimes\chi'_{2}$ with $\chi'_{1}$ and $\chi'_{2}$ non-scalars, each $\chi'_{j}$ ($j=1,2$) is a
scalar multiple of the tensor product of some $\chi_{i}$. Write $I_{i}=\supp\chi_{i}$ ($1\leq i\leq s$) and
$I'_{j}=\supp\chi'_{j}$ ($j=1,2$).

When $I'_{1}$ is a union of some $I_{i}$. Without loss of generality assume that $I'_{1}=\bigsqcup_{1\leq i\leq s'}
I_{i}$. Choose a nonzero weight $\mu$ appearing in $\chi'_{2}$ and let $\chi'$ be the sum of terms in $\chi$
of the form $[\nu+\mu]$ where $\supp\nu\subset I'_{1}$. Then, there exists $c,c'\in\mathbb{Q}^{\times}$ such
that \[\chi'=\chi_{1}\otimes\cdots\otimes\chi_{s'}\otimes(c[\mu])\] and \[\chi'=\chi'_{1}\otimes(c'[\mu]).\]
Then, \[\chi'_{1}=\frac{c}{c'}\chi_{1}\otimes\chi_{2}\otimes\cdots\otimes\chi_{s'}.\] Similarly, one shows that \[\chi'_{2}=c''\chi_{s'+1}\otimes\chi_{2}\otimes\cdots\otimes\chi_{s}\] for some $c''\in\mathbb{Q}^{\times}$.
Then, $c''=\frac{c'}{c}$ and the conclusion follows.

When $I'_{1}$ is not a union of any tuple $\{I_{i}\}$. Then, there exits $i$ such that $I_{i}$ is contained in
neither $I'_{1}$ nor $I'_{2}$. Without loss of generality we assume that $I_{s}$ is contained in neither
$I'_{1}$ nor $I'_{2}$. For each $i$ ($1\leq i\leq s-1$), choose a nonzero weight $[\mu_{i}]$ appearing in
$I_{i}$. Write $\mu=\mu_{1}+\cdots+\mu_{s-1}$. Let $\chi'$ be the sum of terms in $\chi$ of the form
$[\mu+\nu]$ where $\supp\nu\subset I_{s}$. Then, there exists $c\in\mathbb{Q}^{\times}$ such that \[\chi'=
c[\mu_{1}]\otimes\cdots\otimes[\mu_{s-1}]\otimes\chi_{s}.\] For $j=1,2$, let $\chi''_{j}$ be the sum of terms
$[\lambda]$ in $\chi'_{j}$ such that the projection of $\lambda$ to each $I_{i}\cap I'_{j}$ is equal to the
projection of $\mu_{i}$ to $I_{i}\cap I'_{j}$ ($1\leq i\leq s-1$). Then, $\chi'=\chi''_{1}\otimes\chi''_{2}$.
Write $\mu'_{j}$ for the projection of $\mu$ to $I'_{j}$ ($j=1,2$). Then, $\chi''_{j}$ is the form $\chi''_{j}
=\chi'''_{j}\otimes[\mu'_{j}]$ for some character $\chi''_{j}$ supported on $I'_{j}\cap I_{s}$. Then,
$\chi_{s}=\frac{1}{c}\chi'''_{1}\otimes\chi'''_{2}$. This is in contradiction with the maximality of the
decomposition $\chi=\chi_{1}\otimes\chi_{2}\otimes\cdots\otimes\chi_{s}$.
\end{proof}

\begin{lemma}\label{L:decom2}
Let $\chi=\chi_{1}\otimes\cdots\otimes\chi_{s}$ be a maximal decomposition with $I_{i}=\supp\chi_{i}$
($1\leq i\leq s$). Assume that $\chi$ is $S_{n}$ invariant and is not a scalar. Then, either $s=1$; or $s=n$
and each $I_{i}$ consists of a single index, and $\chi_{i}$ ($1\leq i\leq s$) are $S_{n}$ conjugates of each
other up to scalars.
\end{lemma}

\begin{proof}
For any $w\in S_{n}$, we also have $\chi=(w\cdot\chi_{1})\otimes\cdots\otimes(w\cdot\chi_{s})$ by the $S_{n}$
invariance of $\chi$. By the uniqueness of maximal decomposition, we see that $\{w\cdot\chi_{1},\cdots,w\cdot
\chi_{s}\}$ only differs from $\{\chi_{1},\cdots,\chi_{s}\}$ by a permutation. Note that $\supp(w\cdot I_{i})=
w\cdot\supp I_{i}$. Then, $\{w\cdot I_{1},\dots,w\cdot I_{s}\}=\{I_{1},\dots,I_{s}\}$. This implies that: either
$s=1$; or $s=n$ and each $I_{i}$ consists of a single index. In the latter case, $\chi_{i}$ ($1\leq i\leq s$)
are $S_{n}$ conjugates of each other up to scalars.
\end{proof}

\begin{definition}\label{D:fully}
We call an $S_{n}$ invariant character $\chi$ {\it fully decomposable} if $s=n$ in Lemma \ref{L:decom2}; and call
$\chi$ in-decomposable if $s=1$ in Lemma \ref{L:decom2}.
\end{definition}

\begin{lemma}\label{L:decom4}
Let $\chi_{\Phi}(\lambda)$ be a non-trivial irreducible restriction character with $\Phi$ a non-irreducible
quasi sub-root system of $\Psi$ such that $p(\Phi)\supset p(\Psi)^{\circ}$. If $\chi$ is $S_{n}$
invariant, then there exists a degree 1 character $\chi_{0}$ with leading weight of coefficient 1 such
that \[\chi_{\Phi}(\lambda)=\chi_{0}^{\otimes n}.\]
\end{lemma}

\begin{proof}
Since $\Phi$ is non-irreducible, $\chi_{\Phi}(\lambda)$ is not in-decomposable. Then, it is fully decomposable
by Lemma \ref{L:decom2}. Note that the leading weight of $\chi=\chi_{\Phi}(\lambda)$ has coefficient 1 by Lemma
\ref{L:RS3}. By Lemma \ref{L:decom2} again we get $\chi=\chi_{0}^{\otimes n}$ where $\chi_{0}$ is a degree
one non-scalar character with leading weight of coefficient 1.
\end{proof}

\begin{lemma}\label{L:fully4}
Let $\Phi$ be a reduced irreducible quasi sub-root system of $\Psi=\B_{m}^{(2)}$ or $\BC_{m}^{(2,2,1)}$ such that
$p(\Phi)\supset p(\Psi)^{\circ}$, and be the quasi root system obtained from the action of a diagram automorphism
$\theta$ on a root system $\Phi'$. Let $\lambda$ be an integral dominant weight of $\Phi'$ such that
$\chi_{\Phi}(\lambda)=\chi_{0}^{\otimes m}$ where $\chi_{0}$ is a degree 1 non-scalar character with a leading
weight of coefficient 1. Then the triple $(\Phi,\lambda,\chi_{0})$ falls into the following list:
\begin{enumerate}
\item[(i)]$\Phi=\B_{m}$, $\lambda=(\underbrace{\frac{1}{2},\dots,\frac{1}{2}}_{m})$ and $\chi_{0}=\tau_{1}$.
\item[(ii)]$\Phi=\B_{m}^{(2)}$, $\lambda=(\lambda_{1},\lambda_{2})$ with \[\lambda_{1}=\lambda_{2}=
(\underbrace{\frac{1}{2},\dots,\frac{1}{2}}_{m})\textrm{ and }\chi_{0}=\tau_{2}+\tau_{0}\] or
\[\{\lambda_{1},\lambda_{2}\}=\{(\underbrace{\frac{1}{2},\dots,\frac{1}{2}}_{m}),0\}\textrm{ and }\chi_{0}=
\tau_{1}.\]
\item[(iii)] $\Phi=\B_{m}^{(2,1)}$, $\lambda=(\underbrace{\frac{1}{2},\dots,\frac{1}{2}}_{m+1})$ and $\chi_{0}=
\tau_{1}$.
\item[(iv)] $\Phi=\B_{2}^{(2,1)}$, $\lambda=(\frac{3}{2},\frac{3}{2},\frac{1}{2})$ and $\chi_{0}=\tau_{3}+
\tau_{1}$.
\item[(v)] $\Phi=\BC_{1}^{(2,1)}$, $\lambda=(a,0,-b)$ and $\chi_0=\tau'_{a,b}$, where $a,b\in\mathbb{Z}$ and
$a\geq b\geq 0$.
\item[(vi)] $\Phi=\B_{1}^{(2)}$, $\lambda=((\frac{a}{2}),(\frac{b}{2}))$ and $\chi_0=\tau_{a}\tau_{b}$, where
$a,b\in\mathbb{Z}_{\geq 0}$.
\item[(vii)] $\Phi=\B_{1}$, $\lambda=a$ and $\chi_0=\tau_{a}$, where $a\in\mathbb{Z}_{>0}$.
\end{enumerate}
\end{lemma}

\begin{proof}
Since it is assumed that $\Phi$ is irreducible and $p(\Phi)\supset p(\Psi)^{\circ}$, we have $\Phi\cong\B_{m}$,
$\B_{m}^{(2,1)}$, $\B_{m}^{(2)}$ or $\BC_{m}^{(2,2,1)}$.

(1)when $\Phi=\B_{m}$ and $m>1$, let $\frac{k}{2}$ ($k\in\mathbb{Z}_{>0}$) be the highest weight of $\chi_0$. Then
we have $\lambda=(\underbrace{\frac{k}{2},\dots,\frac{k}{2}}_{m})$. Reducing to a Levi factor of type $\B_{1}$ using
Lemma \ref{L:Levi2}, we get $\chi_{0}=\tau_{k}=\sum_{0\leq j\leq k}[\frac{k}{2}-j]$. Reducing to a Levi factor of
type $\B_{2}$ using Lemma \ref{L:Levi2}, we get \[\chi_{\B_{2}}(\frac{k}{2},\frac{k}{2})=\tau_{k}^{\otimes 2}.\]
Comparing dimensions using Lemma \ref{L:Levi3} we get \[\frac{1}{6}(k+1)(k+2)(k+3)=(k+1)^{2}.\] Then, $k=1$.
In this case, $\chi_{\Phi}(\lambda)$ is the degree $m$ spinor character and it is equal to $([\frac{1}{2}]+
[-\frac{1}{2}])^{\otimes m}$. This leads to (i).

(2)when $\Phi=\B_{m}^{(2)}$ and $m>1$, we have $\Phi'=\B_{m}\sqcup\B_{m}$. Write $\lambda=(\lambda_{1},\lambda_{2})$
where $\lambda_{1},\lambda_{2}$ are integral dominant weights of $\B_{m}$. Let $\frac{k}{2}$ ($k\in\mathbb{Z}_{>0}$)
be the highest weight of $\chi_0$. Then, there exists $k_{1},k_{2}\in\mathbb{Z}_{\geq 0}$ with $k_{1}+k_{2}=k$
such that $\lambda_{1}=(\underbrace{\frac{k_{1}}{2},\dots,\frac{k_{1}}{2}}_{m})$ and $\lambda_{2}=
(\underbrace{\frac{k_{2}}{2},\dots,\frac{k_{2}}{2}}_{m})$. Without loss of generality we assume that
$k_{1}\geq k_{2}$. Reducing to a Levi factor of type $\B_{1}$ using Lemma \ref{L:Levi2}, we get $\chi_{0}=
\tau_{k_{1}}\tau_{k_{2}}$. Reducing to a Levi factor of type $\B_{2}$ using Lemma \ref{L:Levi2}, we get \[\chi_{\B_{2}}(\frac{k_{1}}{2},\frac{k_{1}}{2})\chi_{\B_{2}}(\frac{k_{2}}{2},\frac{k_{2}}{2})=
(\tau_{k_{1}}\tau_{k_{2}})^{\otimes 2}.\] Comparing dimensions using Lemma \ref{L:Levi3}, we get  \[\frac{1}{6^{2}}(k_{1}+1)(k_{1}+2)(k_{1}+3)(k_{2}+1)(k_{2}+2)(k_{2}+3)=(k_{1}+1)^{2}(k_{2}+1)^{2}.\] Then,
$(k_{1},k_{2})=(1,0)$ or $(1,1)$. When $(k_{1},k_{2})=(1,0)$, $\chi_{\Phi}(\lambda)$ is the degree $m$ spinor
character and it is equal to $([\frac{1}{2}]+[-\frac{1}{2}])^{\otimes m}$; when $(k_{1},k_{2})=(1,1)$,
$\chi_{\Phi}(\lambda)$ is the square of the degree $m$ spinor character and it is equal to
$(([\frac{1}{2}]+[-\frac{1}{2}])^{2})^{\otimes m}$. This leads to (ii).

(3)when $\Phi=\B_{m}^{(2,1)}$ and $m>1$, we have $\Phi'=\D_{m+1}$. Let $\frac{k}{2}$ ($k\in\mathbb{Z}_{>0}$)
be the highest weight of $\chi_0$. Then, there exist $k'\in\mathbb{Z}$ with $|k'|\leq k$ and $2|k-k'$ such
that $\lambda=(\underbrace{\frac{k}{2},\dots,\frac{k}{2},\frac{k'}{2}}_{m+1})$. Without loss of generality
we assume that $k'\geq 0$. Write $a=\frac{k+k'}{2}$ and $b=\frac{k-k'}{2}$. Reducing to a Levi factor of type
$\B_{1}^{(2)}$ using Lemma \ref{L:Levi2}, we get $\chi_{0}=\tau_{a}\tau_{b}$. Reducing to a Levi factor of
type $\B_{2}^{(2,1)}$ using Lemma \ref{L:Levi2}, we get \[\chi_{\B_{2}^{(2,1)}}(\frac{k}{2},\frac{k}{2},
\frac{k'}{2})=(\tau_{a}\tau_{b})^{\otimes 2}.\] Comparing dimensions using Lemma \ref{L:Levi3} we get \[\frac{1}{12}(a+1)(b+1)(a+2)(b+2)(a+b+3)=(a+1)^{2}(b+1)^{2}.\] Then, $(a,b)=(1,0)$ or $(2,1)$. When $(a,b)
=(1,0)$, we have $k=k'=1$. In this case $\chi_{\Phi}(\lambda)$ is the degree $m$ spinor character. This leads
to (iii). When $(a,b)=(2,1)$, we have $(k,k')=(3,1)$ and $\chi_{0}=\tau_{2}\tau_{1}=\tau_{3}+\tau_{1}$. When
$m=2$, by analyzing weight multiplicities using $W_{\B_{m}}=\{\pm{1}\}^{m}\rtimes S_{m}$ invariance, we show
that $\chi_{\Phi}(\lambda)=\chi_{0}^{\otimes 2}$. This leads to (iv). When $m\geq 3$, reducing to a Levi
factor of type $\B_{3}^{(2)}$ using Lemma \ref{L:Levi2} and comparing dimensions using Lemma \ref{L:Levi3},
we get $224=6^{3}$, which is a contradiction.

(4)when $\Phi=\BC_{m}^{(2,2,1)}$ and $m>1$, we may assume that \[\lambda=(\underbrace{a,\dots,a}_{m},0,
\underbrace{-b,\dots,-b}_{m})\] where $a,b\in\mathbb{Z}$ and $a\geq b\geq 0$. Reducing to a Levi factor of type
$\B_{1}^{(2)}$ using Lemma \ref{L:Levi2}, we get $\chi_{0}=\tau'_{a,b}$ and $\dim\chi_{0}=\frac{1}{2}(a+1)(b+1)
(a+b+2)$. Reducing to a Levi factor of type $\B_{2}^{(2,1)}$ using Lemma \ref{L:Levi2} and comparing dimensions
using Lemma \ref{L:Levi3}, we get an equation \[(a+2)(b+2)(a+b+3)^{2}(a+b+4)=72(a+1)(b+1)(a+b+2),\] which has no
solution except $a=b=0$.

(5)when $\rank\Phi=1$, it leads to (v)-(vii).
\end{proof}

By imposing an extra condition $\theta\lambda=-w_{l}\lambda$ in Lemma \ref{L:fully4}, we reach the following lemma.

\begin{lemma}\label{L:fully2}
Let $\Phi$ be an irreducible quasi sub-root system of $\Psi=\B_{m}^{(2)}$ or $\BC_{m}^{(2,2,1)}$ such that
$p(\Phi)\supset p(\Psi)^{\circ}$, and be the quasi root system obtained from the action of a diagram
automorphism $\theta$ on a root system $\Phi'$. Let $\lambda$ be a dominant integral weight of $\Phi'$ such that
$\theta\lambda=-w_{l}\lambda$, where $w_{l}$ is the longest element in the Weyl group of $\Phi'$. Suppose that
the irreducible restriction character $\chi_{\Phi}(\lambda)=\chi_{0}^{\otimes m}$ for a degree 1 non-scalar
character $\chi_{0}$ with a leading weight of coefficient 1. Then the triple $(\Phi,\lambda,\chi_{0})$ falls
into the following list:
\begin{enumerate}
\item[(i)] $\Phi=\B_{m}$, $\lambda=(\underbrace{\frac{1}{2},\dots,\frac{1}{2}}_{m})$ and $\chi_{0}=\tau_{1}$.
\item[(ii)] $\Phi=\B_{m}^{(2)}$, $\lambda=(\lambda_{1},\lambda_{2})$ with $\lambda_{1}=\lambda_{2}=
(\underbrace{\frac{1}{2},\dots,\frac{1}{2}}_{m})$, and $\chi_{0}=\tau_{2}+\tau_{0}$.
\item[(iii)] $\Phi=\B_{m}^{(2,1)}$ with $m$ being even, $\lambda=(\underbrace{\frac{1}{2},\dots,\frac{1}{2}}_{m+1})$
and $\chi_{0}=\tau_{1}$.
\item[(iv)] $\Phi=\B_{2}^{(2,1)}$, $\lambda=(\frac{3}{2},\frac{3}{2},\frac{1}{2})$ and $\chi_{0}=\tau_{3}+\tau_{1}$.
\item[(v)] $\Phi=\BC_{1}^{(2,1)}$, $\lambda=(a,0,-b)$ and $\chi_0=\tau'_{a,b}$, where $a,b\in\mathbb{Z}$ and
$a\geq b\geq 0$.
\item[(vi)] $\Phi=\B_{1}^{(2)}$, $\lambda=((\frac{a}{2}),(\frac{a}{2}))$ and $\chi_0=\tau_{a}^{2}$, where
$a\in\mathbb{Z}_{>0}$.
\item[(vii)] $\Phi=\B_{1}$, $\lambda=a$ and $\chi_0=\tau_{a}$, where $a\in\mathbb{Z}_{>0}$.
\end{enumerate}
\end{lemma}

\begin{lemma}\label{L:fully3}
Let $\Psi=\B_{n}^{(2)}$ or $\BC_{n}^{(2,2,1)}$, $\Phi_{1}$ and $\Phi_{2}$ be two quasi sub-root systems of $\Psi$,
and $\chi_{0}$ be a degree 1 non-scalar character with a leading weight of coefficient 1. For $j=1,2$, let
$\Phi_{j}$ be the quasi root system obtained from the action of a diagram automorphism $\theta_{j}$ on a root
system $\Phi'_{j}$. Let $\lambda_{j}$ be a dominant integral weight of $\Phi'_{j}$ such that $\theta_{j}\lambda_{j}
=-w_{l,j}\lambda_{j}$, where $w_{l,j}$ is the longest element in the Weyl group of $\Phi'_{j}$. Suppose that
$\Phi_{1}$ and $\Phi_{2}$ are non-conjugate, $p(\Phi_{j})\supset p(\Psi)^{\circ}$ and   $\chi_{\Phi_{1}}(\lambda_{1})
=\chi_{\Phi_{2}}(\lambda_{2})=\chi_{0}^{\otimes n}$. Then, $\chi_{0}=\tau_{1}$, $\tau_{1}^{2}$ or $\tau_{2}$.
Moreover, all possible such quasi sub-root systems $\Phi_{j}$ ($j=1,2$) are as follows:
\begin{enumerate}
\item[(i)]when $\chi_{0}=\tau_{1}$, we have \[\Phi_{j}=(\bigsqcup_{1\leq i\leq s}
\B_{n_{i}}^{(2,1)})\bigsqcup(\bigsqcup_{1\leq i'\leq t}\B_{m_{i'}})\] ($j=1,2$) where $\sum_{1\leq i\leq s}n_{i}+
\sum_{1\leq i'\leq t}m_{i'}=n$, $n_{i}$ are all even.
\item[(ii)]when $\chi_{0}=\tau_{1}^{2}=\tau_{2}+\tau_{0}$, we have \[\Phi_{j}=
(\bigsqcup_{1\leq i\leq s}\B_{n_{i}}^{(2)})\] ($j=1,2$) where $\sum_{1\leq i\leq s}n_{i}=n$.
\item[(iii)]when $\chi_{0}=\tau_{2}$, we have \[\Phi_{1}=(\bigsqcup_{1\leq i\leq s}\BC_{1}^{(2,1)})
\bigsqcup(\bigsqcup_{1\leq i'\leq t}\B_{1})\] ($j=1,2$) where $s+t=n$.
\end{enumerate}
\end{lemma}

\begin{proof}
Suppose that $\chi_{0}\neq\tau_{1}$ or $\tau_{1}^{2}$. Then by Lemma \ref{L:fully2}, each irreducible factor of
$\Phi_{1}$ or $\Phi_{2}$ has rank 1 and $\chi_{0}$ coincides with the degree 1 character happening in at least
two cases among (v)-(vii). Then, (1)$\chi_{0}=\tau'_{a,b}=\tau_{c}^{2}$ where $a\geq b\geq 0$ and $a+b=c>0$;
(2)or $\chi_{0}=\tau'_{a,b}=\tau_{c}$ where $a\geq b\geq 0$ and $2a+2b=c>0$; (3)or $\chi_{0}=\tau_{a}^{2}=
\tau_{b}$ where $2a=b>0$. Case (3) is impossible. In cases (1)-(2), $\tau'_{a,b}$ must be multiplicity-free.
By Lemma \ref{L:tau'-1} we get $b=0$ or 1. Comparing the expansions of $\tau'_{a,0}$ and $\tau'_{a,1}$ with
$\tau_{c}^{2}$ and $\tau_{c}$, we get $(a,b)=(1,0)$ and $\chi_0=\tau'_{1,0}=\tau_{2}$.

When $\chi_{0}$ is equal to each of $\tau_{1}$, $\tau_{1}^{2}$ or $\tau_{2}$, we get all possible such quasi
sub-root systems $\Phi_{j}$ from Lemma \ref{L:fully2}.
\end{proof}

\section{S-subgroups and the functorial source problem}\label{S:S-group}

Let $G$ be a given quasi-connected compact Lie group. Choose an element $g_{0}$ generating $G/G^{0}$ and let $e$ be the
order of \[\Aut(g_{0})|_{\mathfrak{g}}\in\Out(\mathfrak{g})=\Aut(\mathfrak{g})/\Int(\mathfrak{g}).\] Then, $e$ divides
$|G/G^{0}|$ and $g_{0}^{e}\in Z_{G}(G^{0})G^{0}$. We propose the following definition of S-subgroups of $G$.

\begin{definition}\label{D:Sgroup}
We call a closed subgroup $H$ of $G$ an {\it $S$-subgroup} if $H$ is generated by $H^{0}$ and an element
$h_{0}\in g_{0}G^{0}$ such that $h_{0}^{e}\in Z_{G}(G^{0})H^{0}$ and $Z_{G}(H)=Z(G)$.
\end{definition}

Let's explain the origination of the definition of S-subgroups. When $G$ is connected, a closed subgroup $H$
is an S-subgroup if and only the complexified Lie algebra $\mathfrak{h}_{\mathbb{C}}$ is an S-subalgebra of
$\mathfrak{g}_{\mathbb{C}}$ in the sense of \cite{Dynkin} (see also \cite{Minchenko}). Motivated by the
primitization process in \cite[Section 2]{Arthur-Howe} and \cite[Section 4]{Arthur-functoriality}, we
generalize the definition of S-subgroup to the case when $G$ has a cyclic component group. We will see that
the condition $h_{0}^{e}\in Z_{G}(G^{0})H^{0}$ is convenient while comparing S-subgroups of isogenous groups.

Let's explain the motivation and the application of studying S-subgroups and their dimension data. Let
$\mathbb{G}/F$ be a connected reductive algebraic group over a number field $F$. Take a minimal Galois
extension $K$ of $F$ such that $\Gal(\overline{\mathbb{Q}}/K)$ acts trivially on the Langlands dual group
$\hat{G}$. Put $^{L}G=\hat{G}\rtimes\Gal(K/F)$ for the Langlands L-group of $\mathbb{G}/F$. The goal of this
paper is to help determine a conjectural subgroup $H_{\pi}$ of $^{L}G$ associated to an automorphic
representation $\pi$ of $G(\mathbb{A}_{F})$ of Ramanujan type. The group $H_{\pi}$ is called the functorial
source of $\pi$ (\cite{Langlands}) and is the key to show general functoriality principle. With the theory
of Einstein series and endoscopy, it reduces to the case that $\pi$ neither comes from Eisenstein series nor
comes from endoscopy. If $\pi$ has minimal possible ramification (the meaning of this requires a clarification
in future) and $\Gal(K/F)$ is cyclic, then a maximal compact subgroup $H$ of $H_{\pi}$ as a subgroup
of a maximal compact subgroup $\tilde{G}$ of $^{L}G$ should be an S-subgroup. We know that the dimension datum
of $H$ is predicted by the pole orders at $s=1$ of Langlands L-functions of $\pi$ (\cite{Langlands}). Then, if
dimension data of non-conjugate (resp. non-element conjugate) S-subgroups are non-equal, then the group
$H_{\pi}$ is determined up to conjugacy (resp. element conjugacy) by pole orders at $s=1$ of Langlands
L-functions of $\pi$. Therefore, by Theorems \ref{T2}, \ref{T4} and \ref{T6}, we have the following
assertion: let $\pi$ be an automorphic representation of $G(\mathbb{A}_{F})$ of Ramanujan type that neither
comes from Eisenstein series nor comes from endoscopy and has minimal possible ramification. When $G$ is
a split group (resp. a unitary group or an orthogonal group) the functorial source group $H_{\pi}$ of $\pi$ is
determined up to element-conjugacy (resp. conjugacy) uniquely by pole orders at $s=1$ of Langlands L-functions
of $\pi$. This shows an affirmative answer to the ``refined question" in \cite[p. 17]{Arthur-functoriality}.
We also remark that: when $\Gal(K/F)$ is not cyclic, we have no yet a good definition of
S-subgroups and a well study of dimension data of them. If $\pi$ has wilder ramification, we expect that one
must use ramification information besides dimension data in order to determine $H_{\pi}$. We wish
to reach a thorough study of $H_{\pi}$ in future when $\Gal(K/F)$ is not cyclic or when $\pi$ has
wilder ramification.

We have two remarks concerning Definition \ref{D:Sgroup}. First, let $f$ be the order of \[\Aut(g_{0})|_{\mathfrak{g}_{\der}}\in\Out(\mathfrak{g}_{\der})=
\Aut(\mathfrak{g}_{\der})/\Int(\mathfrak{g}_{\der}).\] Then, $f$ is a factor of $e$. For all cases working in
this paper, we have a further condition: $g_{0}^{f}$ acts trivially on $z(\mathfrak{g})$. Under this condition,
we have $e=f$. Second, when $G$ is connected, an S-subgroup $H$ of it is not necessarily connected.

\begin{lemma}\label{L:S1}
Let $H$ be an S-subgroup as in Definition \ref{D:Sgroup} and $S$ be a maximal commutative connected subset in
$h_{0}H^{0}$. Then there exists an element $x_{0}\in S$ such that $x_{0}^{e}\in\tilde{S}\cap Z(G)$, where $\tilde{S}$
is the quasi torus generated by $S$.
\end{lemma}

\begin{proof}
Write $(H^{0})_{\der}=[H^0,H^0]$ for the derived subgroup of $H^{0}$, which has Lie algebra $\mathfrak{h}_{\der}=
[\mathfrak{h},\mathfrak{h}]$. Let \[\phi: H\rightarrow H/Z_{H}((H^{0})_{\der})\subset\Aut(\mathfrak{h}_{\der})\]
be the adjoint homomorphism. Write $H'=\Ima\phi\subset\Aut(\mathfrak{h}_{\der})$. Then, $H'$ is generated by
$\phi(h_0)$ and $H^{'0}$, $\phi(h_0)^{e}\in H^{'0}$ and $\phi(S)$ is a maximal commutative connected subset of
$\phi(h_0)H^{'0}$. Since $H'$ is of adjoint type, there exists $x_{0}\in S$ such that $\phi(x_{0})$ is a pinned
automorphism of $(\mathfrak{h}_{\mathbb{C}})_{\der}$. Then, $\phi(x_{0})^{e}=1$. Thus, $x_{0}^{e}\in Z_{H}((H^{0})_{\der})$.
On the other hand, $x_{0}^{e}\in Z_{G}(G^0)H^{0}$ by the definition of S-subgroup. Hence, \[x_{0}^{e}\in
Z_{G}(G^{0})H^{0}\cap Z_{H}((H^{0})_{\der})=Z_{G}(G^{0})Z(H^0)\cap H.\] Since $H$ is generated by $x_{0}$ and
$H^{0}$ by the definition of S-subgroup, then $x_{0}^{e}\in Z_{G}(G^{0})Z(H^0)$ implies that $x_{0}^{e}$ commutes
with all elements of $H$. Thus, $x_{0}^{e}\in\tilde{S}\cap Z(H)$. By the definition of S-subgroup, $Z(H)\subset Z(G)$.
Then, \[x_{0}^{e}\in\tilde{S}\cap Z(H)=\tilde{S}\cap Z(G).\qedhere\]
\end{proof}

\begin{lemma}\label{L:S4}
The group $\tilde{S}/\tilde{S}\cap Z_{G}(G^{0})$ has an order $e$ cyclic component group.
\end{lemma}

\begin{proof}
First, the component group of $\tilde{S}/\tilde{S}\cap Z_{G}(G^{0})$ is generated by the coset $[x_0]$ and
Lemma \ref{L:S1} indicates that $[x_0]^{e}=1$. Second, since $\Ad(x_{0})|_{\mathfrak{g}}\in\Out(\mathfrak{g})$
has order $e$, $x_{0}^{k}\in Z_{G}(G^{0})G^{0}$ only if $e|k$. Thus, the component group of
$\tilde{S}/\tilde{S}\cap Z_{G}(G^0)$ is a cyclic group of order $e$.
\end{proof}

Lemma \ref{L:S4} also indicates that the conclusion of Lemma \ref{L:S1} is optimal. Note that we have
$\Psi'_{\tilde{S}}\subset X^{\ast}(\tilde{S})$, where $\Psi'_{\tilde{S}}$ is the enveloping quasi root
system associated to the quasi torus $\tilde{S}$ in $G$ as in \S \ref{SS:dim-ARS}. Write $\Psi_{\tilde{S}}$
for the set of roots $\alpha\in\Psi'_{\tilde{S}}$ such that $\alpha\neq 2\beta$ for any other root
$\beta\in\Psi'_{\tilde{S}}$. By definition, $\Psi_{\tilde{S}}$ is a reduced quasi root system on $\tilde{S}$.
Suppose that $e$ is not a multiple of 4. By the following lemma we have \[R(H,\tilde{S})\subset\Psi_{\tilde{S}}\]
whenever $H$ is an S-subgroup of $G$ with $\tilde{S}$ a quasi Cartan subgroup.

\begin{lemma}\label{L:S3}
If $H$ is an S-subgroup of $G$ with $\tilde{S}$ a quasi Cartan subgroup, then \begin{equation}\label{Eq:S3-3}
\bigcap_{\bar{\alpha}\in p(R(H,\tilde{S}))}\ker\bar{\alpha}=\bigcap_{\bar{\beta}\in p(\Psi'_{\tilde{S}})}
\ker\bar{\beta}=\tilde{S}^{0}\cap Z(G)\end{equation} and \begin{equation}\label{Eq:S3-1}p(R(H,\tilde{S}))
\supset p(\Psi_{\tilde{S}})^{\circ}.\end{equation} Suppose that $4\not|e$. Then,   \begin{equation}\label{Eq:S3-2}R(H,\tilde{S})\subset\Psi_{\tilde{S}}.\end{equation}
\end{lemma}

\begin{proof}
Let $x\in\bigcap_{\bar{\alpha}\in p(R(H,\tilde{S}))}\ker\bar{\alpha}$. Then $x$ commutes with $H^{0}$. Since $H$ is
generated by $H^{0}$ and $\tilde{S}$. Then, $x\in Z(H)$. As $H$ is an S-subgroup. By Definition \ref{D:Sgroup}
we get \[x\in Z(H)\subset Z_{G}(H)=Z(G).\] Thus, $x\in\tilde{S}^{0}\cap Z(G)$. Therefore,
\begin{equation*}\bigcap_{\bar{\alpha}\in p(R(H,\tilde{S}))}\ker\bar{\alpha}=
\bigcap_{\bar{\beta}\in p(\Psi'_{\tilde{S}})}\ker\bar{\beta}=\tilde{S}^{0}\cap Z(G).\end{equation*}
By \eqref{Eq:S3-3}, the root system $p(\Psi'_{\tilde{S}})$ and its sub-root system $p(R(H,\tilde{S}))$ generate
the same lattice. By Lemma \ref{L:root2}, we get \[ p(R(H,\tilde{S}))\supset p(\Psi'_{\tilde{S}})^{\circ}=
p(\Psi_{\tilde{S}})^{\circ}.\] Note that $R(H,\tilde{S})$ is a reduced quasi root system. By Lemma \ref{L:ARS0},
we get \[R(H,\tilde{S})\subset\Psi_{\tilde{S}}\] if $4\not|e$.
\end{proof}

Let $\phi: G'\rightarrow G$ be a finite covering (=surjective homomorphism with finite kernel) of compact Lie
groups, $g'_{0}$ be an element generating $G'/G^{'0}$. Write $g_0=\phi(g'_0)$. Then, both $G'$ and $G$ are
quasi-connected.

\begin{lemma}\label{L:S2}
The positive integer $e$ for $G'$ and $G$ are the the same.

Assume that $H$ is a closed subgroup of $G$, generated by $H^{0}$ and an element $h_{0}\in g_{0}G^{0}$ such that
$h_{0}^{e}\in Z_{G}(G^{0})H^{0}$. Choose an element $h'_{0}\in g'_{0}G^{'0}\cap\phi^{-1}(h_0)$ and let $H'$ be
a closed subgroup of $G'$ generated by $\phi^{-1}(H)^{0}$ and $h'_{0}$. If $H$ is an S-subgroup of $G$, then
$H'$ is an S-subgroup of $G'$.
\end{lemma}

\begin{proof}
For the first statement, by assumption the Lie algebras of $G'$ and $G$ are the same. Let $\mathfrak{g}$ denote
the complexification of them. For any element $g'\in G'$, we have $\Ad(g')|_{\mathfrak{g}}=
\Ad(\phi(g'))|_{\mathfrak{g}}$. Thus, the numbers $e$ for $G'$ and $G$ are equal.

For the second statement, $h_{0}^{e}\in Z_{G}(G^{0})H^{0}$ implies that \[h_{0}^{'e}\in\phi^{-1}(Z_{G}(G^{0})H^{0})
=Z_{G'}(G^{'0})H^{'0}.\] Let $y\in Z_{G'}(H')$. Then, $\phi(y)\in Z_{G}(H)=Z(G)\subset Z_{G}(G^{0})$. Thus,
$y\in Z_{G'}(G^{'0})$. Since $H'$ contains an element in $g'_{0}G^{'0}$. Then, $y$ commutes will all elements
in $G'$. Thus, $y\in Z(G')$ and hence $Z_{G'}(H')=Z(G')$. Therefore, $H'$ is an S-subgroup of $G'$.
\end{proof}

\subsection{S-subgroups of an isogenous form of $\U(N)\rtimes\langle\tau\rangle$}\label{SS:TA1}

Let $G$ be a compact Lie group with Lie algebra $\mathfrak{g}\cong\mathfrak{u}(N)$ ($N\geq 3$) and be generated by
$G^{0}$ which is a finite quotient of $\U(N)$ and an element $g_{0}$ such that $\Ad(g_{0})|_{\mathfrak{g}}=
\textrm{complex conjugation}$ and $g_{0}^{2}\in G^0$. Then, $g_{0}^{2}\in Z(G)$. Assume that $g_{0}^{2}=1$.
Put \[\bar{G}=\U(N)\rtimes\langle\tau\rangle\] where \[\tau X\tau^{-1}=\overline{X},\ \forall X\in\U(N)
\textrm{ and }\tau^{2}=1.\] Then, there exists a surjection $\pi:\bar{G}\rightarrow G$ with $\pi(\tau)=
g_{0}$. Let $H$ be a closed subgroup of $G$, generated by $H^0$ and an element $h_0\in G-G^{0}$ such that
$h_{0}^{2}\in Z_{G}(G^0)H^{0}$. Choose an element $\bar{h}_{0}\in\pi^{-1}(h_0)$ and let $\bar{H}\subset
\bar{G}$ be generated by $\pi^{-1}(H)^{0}$ and $\bar{h}_{0}$. By Lemma \ref{L:S2}, $\bar{H}$ is an
S-subgroup of $\bar{G}$. Hence, it suffices to study S-subgroups of $\bar{G}$. For this reason we assume
that $G=\U(N)\rtimes\langle\tau\rangle$ below. Write $\rho=\mathbb{C}^{N}$ for the natural representation
of $G^{0}=\U(N)$. The following lemma characterizes S-subgroups of $\U(N)\rtimes\langle\tau\rangle$.

\begin{lemma}\label{L:TA4}
Let $H$ be a closed subgroup of $G=\U(N)\rtimes\langle\tau\rangle$ generated by $H^0$ and an element $h_0\in G-
G^{0}$ such that $h_{0}^{2}\in Z_{G}(G^0)H^0$. In order that $H$ is an S-subgroup of $G$ it is necessary and
sufficient that $\rho|_{H^{0}}$ is irreducible and non-self dual.
\end{lemma}

\begin{proof}
The sufficiency is clear. We show the necessarity. Assume that $H$ is an S-subgroup of $G$. Let $\rho|_{H^{0}}=
\bigoplus_{1\leq i\leq s}\rho_{i}^{\oplus m_{i}}$ be the decomposition of $\rho|_{H^{0}}$ into irreducible
representations of $H^0$, where $\rho_{i}\not\cong\rho_{j}$ whenever $i\neq j$ and $m_{i}$ is the multiplicity
of $\rho_{i}$. Put $n_{i}=\dim\rho_{i}$ ($1\leq i\leq s$). We may assume that $H^0$ is block diagonal. Then, \[Z_{G^{0}}(H^{0})=\{(\Delta_{n_{1}}(X_{1}),\dots,\Delta_{n_{s}}(X_{s})): X_{i}\in\U(m_{i})\}.\] By the definition
of S-subgroup, $(Z_{G^{0}}(H^{0}))^{h_0}=Z(G)=\{\pm{I}\}$. Since $h_{0}^{2}\in Z_{G}(G^0)H^0$, then
$(\Ad(h_0)|_{Z_{G^{0}}(H^{0})})^{2}=1$. Thus, we have all $m_{i}=1$ and $s=1$. Hence, $\rho|_{H^{0}}$ is irreducible.

If $\rho|_{H^{0}}$ is irreducible and self-dual, we may assume that $H^{0}\subset\SO(N)$ or $H^{0}\subset\Sp(N/2)$.
When $H^{0}\subset\SO(N)$, write $h_{0}=y_{0}\tau$ where $y_0\in G^{0}$. Then, there exists $\sigma\in\Aut(H^{0})$
such that \[\sigma(x)=y_{0}\tau x(y_{0}\tau)^{-1}=y_{0}xy_{0}^{-1},\quad\forall x\in H^{0}.\] Taking conjugation,
we get $\sigma(x)=\bar{y}_{0}x\bar{y}_{0}^{-1}$ ($\forall x\in H^{0}$). Then, $y_{0}^{-1}\bar{y_{0}}\in
Z_{G^{0}}(H^{0})=Z(G^{0})$. Thus, there exists $y\in\O(N)$ and $\eta\in\mathbb{C}^{\times}$ such that $y_0=\eta y$.
In this case $(\eta I)\tau$ commutes with $H$, which is in contradiction with $Z_{G}(H)=Z(G)$. When $H^{0}\subset
\Sp(N/2)$, the proof proceeds similarly by writing $h_{0}=y_{0}\Ad(J_{N/2})\tau$ where $y_0\in G^{0}$.
\end{proof}

\subsection{S-subgroups of an isogenous form of $\O(2N)$}\label{SS:Dn}

Let $G$ be a compact Lie group with Lie algebra $\mathfrak{g}\cong\mathfrak{so}(2N)$ ($N\geq 4$) and be generated by
$G^{0}$ and an element $g_{0}$ such that $\Ad(g_{0})|_{\mathfrak{g}}=\Ad(I_{1,2n-1})$ and $g_{0}^{2}\in G^{0}$.
Then, $g_{0}^{2}\in Z(G)$. Assume that $g_{0}^{2}=1$. Put \[\bar{G}=\O(2N).\] Let $H$ be a closed subgroup of
$G$, generated by $H^0$ and an element $h_0\in G-G^{0}$ such that $h_{0}^{2}\in Z_{G}(G^0)H^{0}$. Then, $G^0$
is isomorphic to one of the following: $\Spin(2N)$, $\SO(2N)$, $\PSO(2N)$.   \begin{enumerate}
\item[(i)]When $G^{0}=\PSO(2N)$, there is a double covering $\pi:\bar{G}\rightarrow G$ with $\ker\pi=Z(\bar{G})$.
Choose an element $\bar{h}_{0}\in\bar{G}\in\pi^{-1}(h_{0})$ and let $\bar{H}$ be generated by $\pi^{-1}(H)^{0}$
and $\bar{h}_{0}$.
\item[(ii)]When $G^{0}\cong\SO(2N)$, we identity $G$ with $\bar{G}$. Let $\bar{H}=H$.
\item[(iii)]When $G^{0}\cong\Spin(2N)$, there is a double covering $\pi':G\rightarrow\bar{G}$ with $\ker\pi'=Z(G)=
\{\pm{1}\}$. Let $\bar{H}=\pi'(H)$.
\end{enumerate}

\begin{lemma}\label{L:TD-reduction}
If $H$ is an S-subgroup of $G$, then $\bar{H}$ is an S-subgroup of $\bar{G}$.
\end{lemma}

\begin{proof}
When $G^{0}=\PSO(2N)$, $\bar{H}$ is an S-subgroup of $\bar{G}$ by Lemma \ref{L:S2}. When $G^{0}=\SO(2N)$, the conclusion
is trivial.

When $G^0=\Spin(2N)$, put $S=Z_{\bar{G}}(\bar{H})$. First we have $\{\pm{I}\}=Z(\bar{G})\subset S$. For any $x\in H$
and any $s\in S$, choose an element $s'\in\pi'^{-1}(s)$ and define \[f(x,s)=xs'x^{-1}s'^{-1}.\] It is clear that
$f(x,s)$ does not depend on the choice of $s'$. As $s$ commutes with $\pi'(x)\in\bar{H}$, we have $f(x,s)\in
\ker\pi'=\{\pm{1}\}\subset Z(G)$. Then, $f$ is multiplicative on both variables and $f(x,s)=1$ whenever $x\in H^{0}$
or $x\in\{y^2: y\in H\}$. Let $H'$ be the subgroup of $H$ generated by $H^0$ and $\{y^2: y\in H\}$. Since $H$
is generated by $H^0$ and $h_{0}\in G-G^0$, we have $H/H'\cong C_{2}$. Then, $f$ descends to a map \[f': H/H'\times
S\rightarrow\{\pm{1}\}.\] On the other hand, for any $1\neq s\in S$, there must exist $x\in H$ such that
$xs'x^{-1}s'^{-1}\neq 1$ (since $H$ is an S-subgroup). Thus, $f'$ is a perfect pairing and $S\cong H/H'\cong C_{2}$.
Hence, $S=Z(\bar{G})$ and $\bar{H}$ is an $S$-subgroup of $\bar{G}$.
\end{proof}

By Lemma \ref{L:TD-reduction}, it suffices to study S-subgroups of $\bar{G}$. Assume that $G=\O(2N)$ below. Write
$\rho=\mathbb{C}^{2N}$ for the natural representation of $G=\O(2N)$. The following lemma characterizes S-subgroups
of $\O(2N)$.

\begin{lemma}\label{L:TD1}
Let $H$ be a closed subgroup of $G$, generated by $H^0$ and an element $h_0\in G-G^{0}$ such that $h_{0}^{2}\in
Z_{G}(G^0)H^{0}$. Then $Z_{G}(H)=Z(G)$ if and only if one of the following conditions holds true:
\begin{enumerate}
\item[(i)]$\rho|_{H^{0}}$ is irreducible.
\item[(ii)]$\rho|_{H^{0}}=\rho_{1}\oplus\rho_{2}$ where $\rho_{1},\rho_{2}$ are odd-dimensional self-dual
irreducible representations of orthogonal type and $h_0$ maps $\rho_{1}$ (resp. $\rho_{2}$) to $\rho_{2}$
(resp. $\rho_{1}$).
\item[(iii)]$\rho|_{H^{0}}=\rho_{1}\oplus\rho_{2}$ where $\rho_{1},\rho_{2}$ are odd-dimensional non-self dual
irreducible representations of $H^{0}$, $\rho_2\cong\rho_{1}^{\ast}$ and $h_0$ maps $\rho_{1}$ (resp. $\rho_{2}$)
to $\rho_{2}$ (resp. $\rho_{1}$).
\end{enumerate}
\end{lemma}

\begin{proof}
The sufficiency is clear. We show the necessarity. The condition $Z_{G}(H)=Z(G)$ just means that $H$ is an
irreducible subgroup. Suppose that $\rho|_{H^{0}}$ is not irreducible. Since $h_{0}^{2}\in Z_{G}(G^{0})H^{0}=
\{\pm{I}\}H^{0}$, then $\rho|_{H^{0}}\cong\rho_{1}\oplus\rho_{2}$ is the direct sum of two non-isomorphic
irreducible representations $\rho_1,\rho_{2}$ and the action of $h_0$ permutes them. Then, $\dim\rho_{1}=
\dim\rho_{2}=N$. Since $\rho|_{H}$ is irreducible and of orthogonal type, either $\rho_1,\rho_2$ are non-self
dual and are conjugate, or $\rho_{1},\rho_{2}$ are self-dual and of orthogonal type. It remains to show that
$N$ is odd.

In the first case, we may assume that \[H^{0}\subset\U(N)=\{\left(\begin{array}{cc}A&B\\-B&A\\\end{array}\right)
\in\SO(2N)\}\] and $h_0=\left(\begin{array}{cc}A_0&B_0\\-B_0&A_0\\\end{array}\right)\diag\{I,-I\}.$ Since $-1=
\det h_0=(-1)^{N}$, then $N$ is odd. In the second case, we may assume that \[H^{0}\subset\SO(N)\times\SO(N)\]
and $h_0=\diag\{A_{0},B_{0}\}\left(\begin{array}{cc}&I_{n}\\I_{n}&\\\end{array}\right)$ with $A_{0},B_{0}\in
\O(N)$. Since \[\diag\{A_{0}B_{0},B_{0}A_{0}\}=h_{0}^{2}\in\{\pm{I}\}H^{0},\] we have $\det A_{0}\det B_{0}=1$
or $(-1)^{N}$. On the other hand, $-1=\det h_{0}=(-1)^{N}\det A_{0}\det B_{0}$. Then, $\det A_{0}\det B_{0}=
1$ and $N$ is odd.
\end{proof}

\subsection{S-subgroups of an isogenous form of $\E_{6}\rtimes\langle\tau\rangle$}\label{SS:E6}

Let $G$ be a compact Lie group with Lie algebra $\mathfrak{g}\cong\mathfrak{e}_{6}$ and be generated by $G^{0}$
and an element $g_{0}$ such that $\mathfrak{g}^{\Ad(g_{0})}=\mathfrak{f}_{4}$ and $g_{0}^{2}\in G^{0}$.
Then, $Z(G)=1$ and $g_{0}^{2}\in Z(G)$. Thus, $g_{0}^{2}=1$. Put \[\bar{G}=\E_{6}\rtimes\langle\tau\rangle,\]
where $\tau^2=1$ and $\E_{6}^{\tau}=\F_{4}$. Let $H$ be a closed subgroup of $G$, generated by $H^0$ and an
element $h_0\in G-G^{0}$ such that $h_{0}^{2}\in Z_{G}(G^0)H^{0}=Z(G^{0})H^{0}$.

We have $G^{0}\cong\E_{6}$ or $\E_{6}^{\ad}$. When $G^{0}\cong\E_{6}$, we may assume that $G=\bar{G}$ and let
$\bar{H}=H$. When $G^{0}\cong\E_{6}^{\ad}$, there is a 3-covering $\pi:\bar{G}\rightarrow G$ with $\ker\pi=
Z(\bar{G}^{0})$. Choose $\bar{h}_{0}\in\pi^{-1}(h_{0})$ and let $\bar{H}$ be generated by $\pi^{-1}(H)^{0}$
and $\bar{h}_{0}$. By Lemma \ref{L:S2}, $\bar{H}$ is an S-subgroup of $\bar{G}$. Hence, it suffices to study
S-subgroups of $\bar{G}$. For this reason we let $G=\bar{G}=\E_{6}\rtimes\langle\tau\rangle$ below.

\begin{lemma}\label{L:FixedPoint2}
Let $\phi$ be an order 2 automorphism of a compact Lie group $L$. If $L^{\phi}=1$, then $L$ is a finite abelian
group of odd order and \[\phi(x)=x^{-1},\quad\forall x\in L.\qedhere\]
\end{lemma}

\begin{proof}
Since $L^{\phi}=1$, we have $\mathfrak{l}^{\phi}=0$. By $\phi^{2}=1$, $\phi$ acts by $-1$ on $\mathfrak{l}$ and
$\mathfrak{l}$ is abelian. Then, $\phi$ fixes order 2 elements in $L^{0}$. Thus, $L^{\phi}=1$ indicates that
$L^{0}=1$. Hence, $L$ is a finite group. By $L^{\phi}=1$ again, $L-\{1\}$ is a disjoint union of some subsets
$\{x,\phi(x)\}$. Hence, $L$ has odd order. By this $L$ is solvable. We prove by induction on $|L|$ that $L$ is
abelian and \[\phi(x)=x^{-1},\quad\forall x\in L.\]

When $L$ is abelian, we have $\phi(x\phi(x))=x\phi(x)$ ($\forall x\in L$). Thus, $x\phi(x)=1$ by $L^{\phi}=1$.
Then, $\phi(x)=x^{-1}$ ($\forall x\in L$). In general, by the solvability of $L$, $[L,L]$ is a proper characteristic
subgroup of $L$. By induction $[L,L]$ is abelian and $\phi$ acts on it by $\phi(x)=x^{-1}$ ($\forall x\in [L,L]$).
Consider the induced action of $\phi$ on $L/[L,L]$. Suppose that it has a nontrivial
fixed element, i.e., there exists $x\in L-[L,L]$ such that $\phi(x)=xy$ for some $y\in [L,L]$. Since $[L,L]$ is
a finite abelian group of odd order, there exists $z\in [L,L]$ such that $y=z^{2}$. Then \[\phi(xz)=\phi(x)\phi(z)
=xyz^{-1}=xz.\] This is in contradiction with $L^{\phi}=1$. Hence, the induced action of $\phi$ on $L/[L,L]$ has
no nontrivial fixed element. By induction the induced action of $\phi$ on $L/[L,L]$ is $x\mapsto x^{-1}$. Take
a maximal subgroup $L'$ of $L$ containing $[L,L]$. Then, $L/L'$ is a cyclic group of odd prime order. Now the
action of $\phi$ stabilizes $L'$. By induction $L'$ is abelian and $\phi(y)=y^{-1}$ ($\forall y\in L'$). Choose
an element $x$ generating $L/L'$. Write $\phi(x)=x^{-1}z$ for some $z\in L'$. For any $y\in L'$, we have
\[xy^{-1}x^{-1}=\phi(xyx^{-1})=\phi(x)\phi(y)\phi(x)^{-1}=x^{-1}y^{-1}x.\] Then, $x^2$ commutes with $y$.
Since $x^2$ and $L'$ generate $L$, then $L$ is abelian.
\end{proof}

\begin{lemma}\label{L:TE6-2}
Let $H$ be an S-subgroup of $G$. Then $H^{0}$ is semisimple and $|Z(H^0)|$ is odd. Moreover, $Z_{G^0}(H^0)$ is
a finite abelian group of odd order and \[\Ad(h_0)|_{Z_{G^0}(H^0)}=-1.\]
\end{lemma}

\begin{proof}
Since $H$ is an S-subgroup of $G$, we have $Z_{G}(H)=Z(G)=1$. Thus, $Z_{G^0}(H^0)^{h_0}=1$. Since $h_{0}^{2}\in
Z_{G}(G^{0})H^{0}$, we have $(\Ad(h_0)|_{Z_{G^0}(H^0)})^{2}=1$. By Lemma \ref{L:FixedPoint2}, $Z_{G^0}(H^0)$ is a
finite abelian group of odd order and $\Ad(h_0)|_{Z_{G^0}(H^0)}=-1$. Since $Z(H^0)\subset Z_{G^0}(H^0)$, then
$Z(H^0)$ is a finite abelian group of odd order. Thus, $H^0$ is semisimple.
\end{proof}

\begin{lemma}\label{L:TE6-3}
If $H$ is an S-subgroup of $G$ and $H^{0}$ is not an S-subgroup of $G^0$, then the Lie algebra $\mathfrak{h}$
has type $3\A_{2}$.
\end{lemma}

\begin{proof}
Since $H^{0}$ is not an S-subgroup, there exists $y\in Z_{G^0}(H^0)-Z(G^0)$. Then, $H^0\subset (G^0)^{y}$. Thus,
$Z((G^0)^{y})\subset Z_{G^0}(H^0)$. By Lemma \ref{L:TE6-2}, $Z_{G^0}(H^0)$ is a finite abelian group of odd order,
and so is $Z((G^0)^{y})$. This implies that $(G^0)^{y}$ is semisimple and hence $\Ad(y)|_{\mathfrak{g}}$ is an isolated
automorphism (cf. \cite[p. 213]{Yu-Geom16}). By \cite[Example 2.8]{Yu-Geom16}, one shows that $\mathfrak{g}^{y}$
is of type $3\A_{2}$. A little more argument shows that $Z(G^0)$ and $y$ generates $Z_{G^0}(H^0)$ and $H^0=(G^0)^{y}$.
Thus, $\mathfrak{h}$ has type $3\A_{2}$.
\end{proof}

Continue the argument in the above proof. By \cite[Example 2.8]{Yu-Geom16}, we have \[H^0=(G^0)^{y}=\SU(3)^{3}/\langle
(\omega I,\omega I,\omega I)\rangle\] where $Z(G^0)$ is generated by $(I,\omega I,\omega^{-1}I)$ and $y=(\omega I,I,I)$.
Then, \[Z((G^0)^{y})=\langle Z(G^0),y\rangle,\] and $Z((G^0)^{y})$ and $(G^0)^{y}$ are centralizers of each other in
$G^{0}$. Moreover, we have \[N_{G}((G^0)^{y})=N_{G}(\langle Z(G^0),y\rangle)=(\SU(3)^{3}/\langle(\omega I,\omega I,\omega I)
\rangle)\rtimes\langle\eta,\sigma,\tau\rangle,\] where $\eta,\sigma\in G^{0}$, $\tau\in G-G^0$,
\[\eta(X_{1},X_{2},X_{3})\eta^{-1}=(X_{2},X_{3},X_{1}),\] \[\sigma(X_{1},X_{2},X_{3})\sigma^{-1}=
(\overline{X_{1}},\overline{X_{3}},\overline{X_{2}}),\] \[\tau(X_{1},X_{2},X_{3})\tau^{-1}=(X_{1},X_{3},X_{2}),\] $\eta^{3}=\sigma^2=\tau^{2}=\sigma\tau\sigma^{-1}\tau^{-1}=1$ and $\sigma\eta\sigma^{-1}=\tau\eta\tau^{-1}=\eta^{-1}$.

\begin{lemma}\label{L:TE6-4}
With the above identification, we have \[H\sim(\SU(3)^{3}/\langle(\omega I,\omega I,\omega I)\rangle)\rtimes\langle
\tau\sigma\rangle\] where $\tau\sigma(X_{1},X_{2},X_{3})(\tau\sigma)^{-1}=
(\overline{X_{1}},\overline{X_{2}},\overline{X_{3}})$.
\end{lemma}

\begin{proof}
We have $H\subset N_{G}(H^0)$. Since $Z_{G}(G^0)=Z(G^0)\subset H^{0}$, we get $H/H^{0}\cong C_{2}$ by the condition
$h_{0}^{2}\in Z_{G}(G^0)H^{0}$. Conjugating by an element in $\langle\eta\rangle$ if necessary, we see that $H/H^0$
is generated by $\tau$ or $\tau\sigma$. Since $y$ commutes with $H^{0}$ and $\tau$, the group generated by $H^0$ and
$\tau$ is not an S-subgroup. Therefore, \[H\sim\langle H^{0},\tau\sigma\rangle=(\SU(3)^{3}/\langle(\omega I,\omega I,
\omega I)\rangle)\rtimes\langle\tau\sigma\rangle\] where $\tau\sigma(X_{1},X_{2},X_{3})(\tau\sigma)^{-1}=
(\overline{X_{1}},\overline{X_{2}},\overline{X_{3}})$.
\end{proof}

\begin{lemma}\label{L:TE6-5}
If $H$ is an S-subgroup of $G$ and $H^{0}$ is an S-subgroup of $G^0$, then the complexified Lie algebra
$\mathfrak{h}_{\mathbb{C}}$ is one of the S-subalgebras \[\A_{2}^{(9)},\quad\G_{2}^{(1)}+\A_{2}^{(2)},\quad\A_{1}^{(28)}
+\A_{2}^{(2)}\] of $\mathfrak{g}_{\mathbb{C}}=\mathfrak{e}_{6}(\mathbb{C})$.
\end{lemma}

\begin{proof}
Since $Z(H)\subset Z_{G}(H)=1$, $h_0$ must act on $\mathfrak{h}$ by an outer automorphism. Thus, $\Out(\mathfrak{h})
:=\Aut(\mathfrak{h})/\Int(\mathfrak{h})$ has even order. Since $H^{0}$ is an S-subgroup of $G^0$, the complexified
Lie algebra $\mathfrak{h}_{\mathbb{C}}$ is an S-subalgebra of $\mathfrak{g}$. By the classification of S-subalgebras
of $\mathfrak{e}_{6}(\mathbb{C})$ by Dynkin (cf. \cite{Dynkin},\cite{Minchenko}), $\mathfrak{h}$ is conjugate to one
of the subalgebras $\A_{2}^{(9)}, \G_{2}^{(1)}+\A_{2}^{(2)}, \A_{1}^{(28)}+\A_{2}^{(2)}$ of
$\mathfrak{g}_{\mathbb{C}}=\mathfrak{e}_{6}(\mathbb{C})$.
\end{proof}

\begin{lemma}\label{L:TE6-6}
When $\mathfrak{h}$ is conjugate to any of the subalgebras $\A_{2}^{(9)}$, $\G_{2}^{(1)}+\A_{2}^{(2)}$,
$\A_{1}^{(28)}+\A_{2}^{(2)}$ of $\mathfrak{g}_{\mathbb{C}}=\mathfrak{e}_{6}(\mathbb{C})$, there is at most
one conjugacy class of S-subgroups $H$.
\end{lemma}

\begin{proof}
Assume that $H$ is an S-subgroup of $G$ with $\mathfrak{h}$ equal to a given S-subalgebra of type $\A_{2}^{(9)}$,
$\G_{2}^{(1)}+\A_{2}^{(2)}$ or $\A_{1}^{(28)}+\A_{2}^{(2)}$. Then, the conjugacy class of $H^0$ is given. Since $H$
is an S-subgroup, any element $h$ generating $H/H^0$ must act on $\mathfrak{h}$ as an outer automorphism. Since $\Aut(\mathfrak{h})/\Int(\mathfrak{h})\cong C_{2}$ in any of these three cases, we may choose an element $h_0$
generating $H/H^0$ such that $\Ad(h_0)|_{\mathfrak{h}}$ is a fixed outer automorphism. Suppose that $H'$ is another
S-subgroup of $G$ with Lie algebra equal to $\mathfrak{h}$. By the same argument, there exists an element $h'_0$
generating $H'/H^{'0}$ such that $\Ad(h'_0)|_{\mathfrak{h}}=\Ad(h_0)|_{\mathfrak{h}}$. Then, $h_{0}^{-1}h'_{0}\in
Z_{G^{0}}(\mathfrak{h})=Z(G^{0})$ (as $\mathfrak{h}$ is an S-subalgebra). Then, $H'$ is equal to the conjugate
of $H$ by an element in $Z(G^{0})$.
\end{proof}

\begin{corollary}\label{C:TE6-7}
If $H_{1},H_{2}$ are two non-conjugate S-subgroups of $G$ with $h_{j}$ generating $H_{j}/H_{j}^{0}$, then any two
maximal commutative connected subsets of $h_{1}H_{1}^{0},h_{2}H_{2}^{0}$ respectively (or any two maximal tori
of $H_{j}^{0}$ respectively) are not conjugate.
\end{corollary}

\begin{proof}
Write $L_{j}$ ($j=1,2,3,4$) for an S-subgroup of $G$ (it it exists) with $\mathfrak{l}$ of type $3\A_2$, $\A_{2}^{(9)}$,
$\G_{2}^{(1)}+\A_{2}^{(2)}$, $\A_{1}^{(28)}+\A_{2}^{(2)}$, respectively. For each $j$ ($j=1,2,3,4$), choose an element
$l_{j}$ generating $L_{j}/L_{j}^{0}$, choose a maximal commutative connected subset $S_{j}$ of $l_{j}L_{j}^{0}$ and
choose a maximal torus $T_{j}$ of $L_{j}^{0}$. By Lemmas \ref{L:TE6-4}-\ref{L:TE6-6}, we have \[\dim S_{1}=3,\quad
\dim S_{2}=1,\quad \dim S_{3}=3,\quad \dim S_{4}=2\] and \[\dim T_{1}=6,\quad\dim T_{2}=2,\quad \dim T_{3}=4,\quad
\dim T_{4}=3.\] Then, only $S_{1}$ and $S_{3}$ might be conjugate. Note that the restricted root system
$R(L_{1},\tilde{S}_{1}^{0})$ has type $3\A_1$ (roots having equal length) and the restricted root system
$R(L_{3},\tilde{S}_{3}^{0})$ has type $\G_2+\A_1$. Since there is no rank 3 root system containing
both $3\A_1$ (roots having equal length) and $\G_2+\A_1$, $\tilde{S}_{1}^{0}$ and $\tilde{S}_{3}^{0}$ are not
conjugate by Lemma \ref{L:Psi'}. Then, $S_{1}$ and $S_{3}$ are not conjugate.
\end{proof}

\begin{remark}\label{R:S1}
While considering isogenous forms of $\U(N)\rtimes\langle\tau\rangle$, $\O(2N)$, $\E_{6}\rtimes\langle\tau\rangle$,
we always assume that $g_{0}^{2}=1$. This excludes some interesting groups like $\Pin(2N)$ while $N>4$ is odd.
We may loose this condition to $g_{0}^{2}\in Z(G)$. The arguments in this section might still be able to
treat S-subgroups of these more general isogenous forms.
\end{remark}

\subsection{S-subgroups of an isogenous form of $\Spin(8)\rtimes\langle\tau\rangle$}\label{SS:D4}

Let $G$ be a compact Lie group with Lie algebra $\mathfrak{g}\cong\mathfrak{so}(8)$ and be generated by $G^{0}$ and
an element $g_{0}$ such that $\mathfrak{g}^{\Ad(g_{0})}=\mathfrak{g}_{2}$ and $g_{0}^{3}\in G^{0}$. Then $Z(G)=1$
and $g_{0}^{3}\in Z(G)$. Thus, $g_{0}^{3}=1$.

\begin{lemma}\label{L:D4}
There is no proper S-subgroup of $G$.
\end{lemma}

\begin{proof}
Suppose that $H$ is a proper S-subgroup of $G$. That is to say, $H$ is a proper closed subgroup of $G$, generated by
$H^0$ and an element $h_0\in g_{0}G^{0}$ such that $h_{0}^{3}\in Z_{G}(G^0)H^{0}=Z(G^{0})H^{0}$, and $Z_{G}(H)=
Z(G)=1$.

First, we show that $H^{0}$ is semisimple. Write $T=Z(H^0)^{0}$ and $\phi=\Ad(h_{0})|_{T}$. Since $h_{0}^{3}\in
Z(G^{0})H^{0}$, we have $\phi^{3}=\Ad(h_{0}^{3})|_{T}=1$. Since $Z_{G}(H)=1$, then $T^{\phi}=1$. By Lemma
\ref{L:FixedPoint}, 1 is not a root of $f(t)$ and $f(1)=1$. Then, $f(t)=1$ due to $\phi^{3}=1$. Hence, $T=1$ and
$H^{0}$ is semisimple.

Second, we deduce a contradiction. Since $Z_{G}(H)=1$ again, the conjugation action of some $x_0\in h_{0}H^{0}$
on $H^0$ must be an outer automorphism of order $3$. This forces $\mathfrak{h}_{0}$ to be of type $4\A_{1}$ or
$3\A_{1}$, and $H^0$ is either a full rank subgroup of type $4\A_{1}$ or is a normal subgroup of it of type
$3\A_{1}$. We observe that $Z(H^0)\cong C_{2}^{3}$ or $C_{2}$. Then, $Z(H^{0})^{x_{0}}\neq 1$, which is
in contradiction with $Z_{G}(H)=1$.
\end{proof}

\section{Dimension data of S-subgroups of a connected compact Lie group}\label{S:connected}

Let $G$ be a given compact Lie group. Two homomorphisms $\phi_{1},\phi_{2}: H\rightarrow G$ are said to be element
conjugate if \[\phi_{1}(x)\sim\phi_{2}(x),\ \forall x\in H.\] Element-conjugacy of homomorphisms is defined and
studied by Michael Larsen (\cite{Larsen2}). We call two closed subgroups $H_{1},H_{2}$ of $G$ element conjugate if
there exists an isomorphism $\phi: H_{1}\rightarrow H_{2}$ such that \[\phi(x)\sim x,\ \forall x\in H_{1}.\]
Clearly element conjugate subgroups must have the same dimension datum.

\begin{lemma}\label{L:root3}
Let $G$ be a connected compact Lie group and $T\subset G$ be a torus, and $H$ be a connected closed subgroup with
$T$ a maximal torus. If $Z(H)\subset Z(G)$, then $\Psi_{T}^{'\circ}\subset R(H,T)$. Particularly, if $H$ is
an S-subgroup of $G$, then $\Psi_{T}^{'\circ}\subset R(H,T)$.
\end{lemma}

\begin{proof}
We have $R(H,T)\subset\Psi_{T}^{'}\subset X^{\ast}(T)$. If all roots $\alpha\in R(H,T)$ taking value 1 at an element
$t\in T$, then $t\in Z(H)\subset Z(G)$. Thus, all roots $\alpha\in\Psi_{T}^{'}$ taking value 1 at $t\in T$. Hence,
the subgroup of $X^{\ast}(T)$ generated by roots in $R(H,T)$ and that generated by roots in $\Psi_{T}^{'}$ are equal.
By Lemma \ref{L:root2}, we get $\Psi_{T}^{'\circ}\subset R(H,T)$. If $H$ is an S-subgroup of $G$, then $Z(H)\subset
Z_{G}(H)=Z(G)$. Thus, $\Psi_{T}^{'\circ}\subset R(H,T)$.
\end{proof}

\begin{proof}[Proof of Theorem \ref{T1}.]
Let $\{H_{1},\dots,H_{s}\}$ be a tuple of pairwise non-element conjugate connected S-subgroups of $G$. Considering the
supports of Sato-Tate measures, we may further assume that $H_{1},\dots,H_{s}$ share a maximal torus, denoted by $T$.
Put \[\Gamma^{\circ}=N_{G}(T)/Z_{G}(T).\] By Lemma \ref{L:root3}, we have $\Psi_{T}^{'\circ}\subset R(H_{i},T)$
($1\leq i\leq s$). For simplicity we write $\Psi=\Psi'_{T}$ and $\Phi_{i}=R(H_{i},T)$ ($1\leq i\leq s$). Since
$H_{1},\dots,H_{s}$ are pairwise non-element conjugate, the root systems $\Phi_{1},\dots,\Phi_{s}$ are pairwise
non-conjugate with respect to the action of $\Gamma^{\circ}$ on $\Psi$ (\cite[Proposition 4.1]{Yu-dimension}).
Then, the characters $F_{\Phi_{i},\Gamma^{\circ}}$ ($1\leq i\leq s$) have distinct leading terms. Thus, the
characters $F_{\Phi_{i},\Gamma^{\circ}}$ ($1\leq i\leq s$) are linearly independent. Hence, the dimension data
$\mathscr{D}_{H_{i}}$ are linearly independent by Proposition \ref{P:dim-char}.
\end{proof}

Theorem \ref{T1} is a common generalization of \cite[Theorem 2]{Larsen-Pink} and \cite[Theorem 1.3]{Yu-dimension}. The
above proof is much simpler than that in \cite{Larsen-Pink} and \cite{Yu-dimension}.

\begin{proof}[Proof of Theorem \ref{T2}.]
Suppose that $\mathscr{D}_{H_{1}}=\mathscr{D}_{H_{2}}$. Choose a maximal torus $T_{i}$ of $H_{i}^{0}$ ($i=1,2$).
By Lemma \ref{L:S1}, we have $H_{i}=(Z(G)\cap H_{i})H_{i}^{0}$. Any maximal commutative connected subset $S_{i}$
of $H_{i}$ is of the form $S_{i}=zT_{i}$ ($z\in Z(G)\cap H_{i}$). Considering supports of the Sato-Tate measures
$\st_{H_{1}}$ and $\st_{H_{2}}$, we see that $T_{1}$ and $T_{2}$ are conjugate. Replacing $H_{2}$ by a conjugate
one we may assume that $T_{1}=T_{2}$. For any $z\in Z(G)\cap H_{1}$, considering the supports of Sato-Tate measures
$\st_{H_{1}}$ and $\st_{H_{2}}$ again we see that: there exists $z'\in Z(G)\cap H_{2}$ such that $zT_{1}$ and
$z'T_{2}$ are conjugate. Choosing $g\in G$ such that $g(zT_{1})g^{-1}=z'T_{2}$. Then, \[z\in zT_{1}\cap Z(G)=
g(zT_{1}\cap Z(G))g^{-1}=z'T_{2}\cap Z(G)\subset Z(G)\cap H_{2}.\] Thus, $Z(G)\cap H_{1}\subset Z(G)\cap H_{2}$.
Similarly we have $Z(G)\cap H_{2}\subset Z(G)\cap H_{1}$. Then, $Z(G)\cap H_{1}=Z(G)\cap H_{2}$. By Theorem
\ref{T1}, $H_{1}^{0}$ and $H_{2}^{0}$ are element-conjugate, i.e., there exists an isomorphism
$\phi: H_{1}^{0}\rightarrow H_{2}^{0}$ such that $\phi(x)\sim x$ ($\forall x\in H_{1}^{0}$). Define
$\phi': H_{1}\rightarrow H_{2}$ by \[\phi'(zx)=z\phi(x),\quad\forall(z,x)\in(Z(G)\cap H_{1})\times H_{1}^{0}.\]
Then, $\phi'$ is an isomorphisms and $\phi'(x)\sim x$ ($\forall x\in H_{1}$). Hence, $H_{1}$ and $H_{2}$ are
element-conjugate.
\end{proof}

The following theorem generalizes Theorems \ref{T1} and \ref{T2}. Its proof shows the power of quasi root system.

\begin{theorem}\label{T:S3}
Let $G$ be a connected compact Lie group. Then the dimension data of any tuple of pairwise non-element conjugate
S-subgroups of $G$ are linearly independent.
\end{theorem}

\begin{proof}
Let $\{H_{1},\dots,H_{s}\}$ be a tuple of pairwise non-element conjugate S-subgroups of $G$. Suppose that the
dimension data \[\mathscr{D}_{H_{1}},\dots,\mathscr{D}_{H_{s}}\] are linearly dependent. For each $i$
($1\leq i\leq s$), choose a maximal commutative connected subset $S_{i}$ of $H_{i}$ generating $H_{i}/H_{i}^{0}$.
Write $d_{i}$ for the minimal order of elements in $S_{i}$. Without loss of generality we assume that
$\dim S_{i}\geq\dim S_{1}$ for each $i$ and $d_{i}\leq d_{1}$ whenever $\dim S_{i}=\dim S_{1}$. Moreover, we
assume that $(\dim S_{i},d_{i})=(\dim S_{1},d_{1})$ holds when $1\leq i\leq s'$ for some integer $s'\leq s$;
for some integer $t\leq s'$, $S_{i}=S_{1}$ when $1\leq i\leq t$ and any maximal commutative connected subset
$S'_{i}$ of $H_{i}$ generating $H_{i}/H_{i}^{0}$ is not conjugate to $S_{1}$ whenever $t+1\leq i\leq s'$. Then,
by Lemma \ref{L:MCC1}, any maximal commutative connected subset $S'_{i}$ of $H_{i}$ is not conjugate to $S_{1}$
whenever $t+1\leq i\leq s$. Put $S=S_{1}$. Let $\tilde{S}$ be the quasi torus generated by $S$. Put
\[\Gamma^{\circ}=N_{G}(\tilde{S})/Z_{G}(\tilde{S})\] and $\Psi=\Psi'_{\tilde{S}}$. Write $\Phi_{i}=
R(H_{i},\tilde{S})$ ($1\leq i\leq t$). Then, the characters \[F_{\Phi_{1},\Gamma^{\circ}},\dots,
F_{\Phi_{t},\Gamma^{\circ}}\] are linearly dependent by Proposition \ref{P:dim-char}. Note that $S$ is of the
form $S=zT$ with $T$ a torus and $z\in Z(G)$. Then, each root $\alpha\in\Psi$ takes value 1 at $z$ and is
determined by the restricted root $\alpha|_{T}$. The remaining proof is the same as that for Theorem \ref{T1}.
\end{proof}

Using the condition $Z(G)\cap H_{1}=Z(G)\cap H_{2}$ shown in the proof of Theorem \ref{T2}, one reduces the
element-conjugacy against conjugacy question for S-subgroups in a connected compact Lie group $G$ to connected
S-subgroups case. Note that it is equivalent to study the element-conjugacy against conjugacy question for
connected S-subgroups in any isogeny form of $G$. Furthermore, this reduces to the case that $G$ is a torus of or a
connected compact simple Lie group. When $G$ is a torus, trivially element-conjugacy and conjugacy are equivalent.
When $G\cong\SU(N)$, $\Sp(N)$, $\SO(2N+1)$ or $\G_{2}$, it is strongly acceptable (\cite{Griess}, \cite{Larsen2})
and hence element-conjugate subgroups are conjugate. When $G$ is a connected compact exceptional simple Lie group,
by Dynkin's classification (\cite{Dynkin},\cite{Minchenko}), it follows that element-conjugate connected
S-subgroups are conjugate. When $G=\SO(2N)$, an element-conjugacy conjugacy class of connected S-subgroups may
split into one or two conjugacy classes.

\section{Dimension data of S-subgroups of $\U(N)\rtimes\langle\tau\rangle$}\label{S:TA1}

Put \[G=\U(N)\rtimes\langle\tau\rangle.\] Write $\rho=\mathbb{C}^{N}$ for the natural representation of $\U(N)$.

\subsection{Quasi root systems of S-subgroups: properties and constraints}\label{SS:TA1-ARS}

Fix a closed commutative connected subset $S$ in $\tau\U(N)$, and let $\tilde{S}$ be the quasi torus generated by 
$S$. Write $\Psi'_{\tilde{S}}\subset X^{\ast}(\tilde{S})$ for the enveloping quasi root system associated to 
$\tilde{S}$ as in Subsection \ref{SS:dim-ARS} and let $\Psi_{\tilde{S}}$ the set of roots $\alpha\in\Psi'_{\tilde{S}}$ 
such that $\alpha\neq 2\beta$ for any root $\beta\in\Psi'_{\tilde{S}}$. For simplicity we write $\Psi=\Psi_{\tilde{S}}$. 
Let \[\Psi=\bigsqcup_{1\leq i\leq r}\Psi_{i}\] be the decomposition of $\Psi$ into an orthogonal union of irreducible 
quasi root systems. Put \[\Gamma^{\circ}=N_{G^{0}}(\tilde{S})/Z_{G^{0}}(\tilde{S}).\] 

\begin{lemma}\label{L:TA1-Psi0}
Let $H$ be an S-subgroup of $G$ with $\tilde{S}$ a quasi Cartan subgroup. Then we have \begin{equation*}
\bigcap_{\bar{\alpha}\in p(R(H,\tilde{S}))}\ker\bar{\alpha}=\bigcap_{\bar{\alpha}\in p(\Psi)}\ker\bar{\alpha}
=\tilde{S}^{0}\cap Z(G),\end{equation*} \[R(H,\tilde{S})\subset\Psi\] and
\[p(R(H,\tilde{S}))\supset p(\Psi)^{\circ}.\]
\end{lemma}

\begin{proof}
We have $e=2$. By Lemma \ref{L:S3} we get the conclusion of the lemma.
\end{proof}

Fix an element $x_{0}\in S$ such that $x_{0}^{2}=\pm{I}$. Define $\beta_{0}\in X^{\ast}(\tilde{S})$ by
\[\beta_{0}|_{\tilde{S}^{0}}=1\textrm{ and }\beta_{0}(x_{0})=-1.\] For an element $y\in\tilde{S}\cap G^{0}$
with $y^2=I$, define $t_{y}\in\Aut(\tilde{S})$ by \[t_{y}|_{\tilde{S}^{0}}=\id\textrm{ and }t_{y}(x_{0})
=yx_{0}.\]

\begin{lemma}\label{L:TA1-Gamma1}
For each element $y\in\tilde{S}\cap G^{0}$ with $y^{2}=I$, we have $t_{y}\in\Gamma^{\circ}$.
\end{lemma}

\begin{proof}
Since $(yx_{0})^{2}=x_{0}^{2}=\pm{I}$, there exists $g\in\U(N)$ such that $gx_{0}g^{-1}=yx_{0}$. Then,
$g^{-1}\tilde{S}^{0}g\subset\U(N)^{x_{0}}$ and $\tilde{S}^{0}\subset\U(N)^{x_{0}}$. Since $\U(N)^{x_{0}}$ ($\sim\O(N)$
or $\Sp(N/2)$) strongly controls its fusion in $\U(N)$, there exists $h\in\U(N)^{x_{0}}$ such that \[g^{-1}y'g=
hy'h^{-1},\ \forall y\in\tilde{S}^{0}.\] Put $g'=gh$. Then, $g'x_{0}g'^{-1}=yx_{0}$ and $\Ad(g')|_{\tilde{S}^{0}}=\id$.
Thus, $g'\in N_{G^{0}}(\tilde{S})$ and $t_{y}=\Ad(g')|_{\tilde{S}}\in\Gamma^{\circ}$.
\end{proof}

\begin{lemma}\label{L:TA1-Psi1}
Let $\alpha\in\Psi'_{\tilde{S}}$ be a root. Except when $\alpha|_{\tilde{S}^{0}}=2\lambda$ for some
$\lambda\in X^{\ast}(\tilde{S}^{0})$, we have $\alpha+\beta_{0}\in\Psi'_{\tilde{S}}$.
\end{lemma}

\begin{proof}
Suppose that $\alpha|_{\tilde{S}^{0}}\neq 2\lambda$ for any $\lambda\in X^{\ast}(\tilde{S}^{0})$. Then, there exists
$I\neq y\in\tilde{S}^{0}$ with $y^2=I$ such that $\alpha(y)=-1$. By Lemma \ref{L:TA1-Gamma1} we have
$t_{y}\in\Gamma^{\circ}$. Then, \[\alpha+\beta_{0}=t_{y}(\alpha)\in\Psi'_{\tilde{S}}.\qedhere\]
\end{proof}

\begin{lemma}\label{L:TA1-Psi2}
Let $\Psi_{i}$ be an irreducible factor of $\Psi$. \begin{enumerate}
\item[(i)]When $p(\Psi_{i})$ is not of type $\mathbf{C}$ or $\mathbf{BC}$, we have $\Psi_{i}\cong p(\Psi_{i})^{(2)}$.
\item[(ii)]When $p(\Psi_{i})\cong\C_{m}$ ($m\geq 2$), we have $\Psi_{i}\cong\C_{m}^{(2)}$ or $\C_{m}^{(2,1)}$.
\item[(iii)]When $p(\Psi_{i})\cong\C_{1}\cong\A_{1}$, we have $\Psi_{i}\cong\A_{1}^{(2)}$ or $\A_{1}$.
\item[(iv)]When $p(\Psi_{i})\cong\BC_{m}$ ($m\geq 2$), we have $\Psi_{i}\cong\BC_{m}^{(2,2,1)}$.
\item[(v)]When $p(\Psi_{i})\cong\BC_{1}$, we have $\Psi_{i}\cong\BC_{1}^{(2,1)}$.
\end{enumerate}
\end{lemma}

\begin{proof}
Suppose that $\Psi_{i}$ is not of type $\mathbf{C}$ or $\mathbf{BC}$. Then, for any root $\alpha\in\Psi_{i}$, there
exists $\beta\in\Psi_{i}$ such that $\frac{2(\alpha,\beta)}{(\beta,\beta)}=1$. Thus, $\alpha|_{\tilde{S}^{0}}
\neq 2\lambda$ for any $\lambda\in X^{\ast}(\tilde{S}^{0})$. By Lemma \ref{L:TA1-Psi1}, we get
$\alpha+\beta_{0}\in\Psi_{\tilde{S}}$. Then, the conclusion in (i) follows.

Suppose that $p(\Psi_{i})\cong\C_{m}$ ($m\geq 2$). Then, for any root $\alpha\in\Psi_{i}$ such that $p(\alpha)$ is a
short root in $p(\Psi_{i})$, there exists $\beta\in\Psi_{i}$ such that $\frac{2(\alpha,\beta)}{(\beta,\beta)}=1$.
Hence, $p(\alpha)\neq 2\lambda$ for any $\lambda\in X^{\ast}(\tilde{S}^{0})$. By Lemma \ref{L:TA1-Psi1}, we get
$\alpha+\beta_{0}\in\Psi_{\tilde{S}}$. Then, the conclusion in (ii) follows.

The conclusion in each item (iii)-(v) is clear.
\end{proof}

Fix a positive system $p(\Psi)^{+}$ of $p(\Psi)$. Let $H$ be an S-subgroup with $\tilde{S}$ a quasi
Cartan subgroup. Write $\Phi=R(H,\tilde{S})$ for the quasi root system of $H$ with respect to $\tilde{S}$.
Choose a pinned element $x\in S$ of $\Phi$. Then, $\alpha(x)=\pm{1}$ for each root $\alpha\in\Phi$. Thus,
$x^{2}\in Z(H)\subset Z(G)=\{\pm{I}\}$. Write $x^{2}=\epsilon_{\Phi}I$ where $\epsilon_{\Phi}=\pm{1}$.

\begin{lemma}\label{L:TA1-fracFactor3}
Let $\Phi\subset\Psi_{\tilde{S}}$ be the quasi root system of an S-subgroup $H$ with $\tilde{S}$ a generalized
Cartan subgroup. Suppose that the isomorphism type of each quasi root system $\Phi\cap\Psi_{i}$ ($1\leq i\leq r$)
is given. Then the sign $\epsilon_{\Phi}$ and the $\Gamma^{\circ}$ orbit containing $\Phi$ are uniquely determined.
\end{lemma}

\begin{proof}
Replacing $\Phi$ by a $W_{\Psi}$ conjugate of it if necessary, we may assume that $p(\Phi)\cap p(\Psi_{i})$
is given for each $i$ ($1\leq i\leq r$). Choose a pinned element $x$($\in S$) of $\Phi$. Write $H_{x}=(H^{x})^{0}$.
Then, $\tilde{S}^{0}$ is a maximal torus of $H_{x}$ and the root system $R(H_{x},\tilde{S}^{0})$ is a given subset
of $X^{\ast}(\tilde{S}^{0})$. Moreover, we have \[H_{x}\subset\U(N)^{x}\sim\left\{\begin{array}{cc}
\O(N)\textrm{ if }\epsilon_{\Phi}=1\\\Sp(N/2)\textrm{ if }\epsilon_{\Phi}=-1.\end{array}\right.\] Hence,
$\rho|_{H_{x}}$ is of orthogonal type if $\epsilon_{\Phi}=1$, and is of symplectic type if $\epsilon_{\Phi}=-1$.
On the other hand, the character $\chi_{\rho}|_{\tilde{S}^{0}}$ is given. By Lemmas \ref{L:RS1} and \ref{L:RS3},
longest weights appearing in $\chi_{\rho}|_{\tilde{S}^{0}}$ consist in a $W_{p(\Phi)}$ orbit.
Hence, $\rho|_{H_{x}}$ contains a unique irreducible representation with highest weight of maximal length. Then,
the symplectic/orthogonal type of $\rho|_{H_{x}}$ is determined by the Schur index of this irreducible factor of
$\rho|_{H_{x}}$. Therefore, the sign $\epsilon_{\Phi}$ is uniquely determined.

Given the sign $\epsilon_{\Phi}$, the pinned element $x$ is determined up to a translation $t_{y}$ where $y\in
\tilde{S}\cap G^{0}$ having $y^{2}=I$. Since all such $t_{y}$ are contained in $\Gamma^{\circ}$ by Lemma
\ref{L:TA1-Gamma1}, we may assume that the pinned element $x$ is given. Then, a simple system of $\Phi$ is
determined by the isomorphism type of $p(\Phi)$. Furthermore, the quasi root system $\Phi$ is determined by its
folding indices.
\end{proof}

The following lemma says that: once the semisimplicity or non-semisimplicity of an S-subgroup $H$ is given, the
quasi root system of $H$ determines it up to conjugacy. This is an instance showing the power of quasi root
system.

\begin{lemma}\label{L:ARS-conjugacy-TA1}
Let $H_{1}$ and $H_{2}$ be two S-subgroups of $G=\U(N)\rtimes\langle\tau\rangle$ with $\tilde{S}$ a common quasi
Cartan subgroup and with quasi root systems $\Phi_{j}=R(H_{i},\tilde{S})$ ($j=1,2$). Assume that $H_{1}^{0}$
and $H_{2}^{0}$ are both semisimple or both non-semisimple. Then, the following three conditions are equivalent:
\begin{enumerate}
\item[(1)]$H_{1}$ and $H_{2}$ are conjugate.
\item[(2)]$H_{1}$ and $H_{2}$ are element-conjugate.
\item[(3)]$\Phi_{1}$ and $\Phi_{2}$ are conjugate with respect to $\Gamma^{\circ}$.
\end{enumerate}
\end{lemma}

\begin{proof}
The implication $(1)\Rightarrow(2)$ is trivial.

$(2)\Rightarrow(3)$. Suppose that $\phi: H_{1}\rightarrow H_{2}$ is an isomorphism such that $\phi(x)\sim x$
($\forall x\in H_{1}$). Then, $\phi(\tilde{S})$ is also a quasi Cartan subgroup of $H_{2}$. Replacing $\phi$
by a conjugate of it if necessary we may assume that $\phi(\tilde{S})=\tilde{S}$. Since $\phi(x_{0})\sim x_{0}$,
there exists $g_{1}\in G^{0}$ such that $x_{0}=g_{1}\phi(x_{0})g_{1}^{-1}$. Now, both $\tilde{S}^{0}$ and
$\Ad(g_{1})(\phi(\tilde{S}^{0}))$ are contained in $\U(N)^{x_{0}}$ and $\Ad(g_{1})\circ\phi: \tilde{S}^{0}
\rightarrow\Ad(g_{1})(\phi(\tilde{S}^{0}))$ is an element-conjugate homomorphism. Since $\U(N)^{x_{0}}$
($\cong\O(N)$ or $\Sp(N/2)$) strongly controls its fusion in $\U(N)$ (\cite{Griess}), there exists
$g_{2}\in\U(N)^{x_{0}}$ such that $g_{1}\phi(y)g_{1}^{-1}=g_{2}^{-1}yg_{2}$ ($\forall y\in\tilde{S}^{0}$).
Put $g=(g_{2}g_{1})^{-1}$. Then, $g\in N_{G^{0}}(\tilde{S})$ and $\phi|_{\tilde{S}}=\Ad(g)|_{\tilde{S}}$.
Put $\gamma=[g]\in\Gamma^{\circ}$. Then, $\Phi_{1}=\phi^{\ast}\Phi_{2}=\gamma^{\ast}\Phi_{2}$.

$(3)\Rightarrow(1)$. Suppose that $\Phi_{1}$ and $\Phi_{2}$ are $\Gamma^{\circ}$ conjugate. Replacing $H_{2}$
by a conjugate one if necessary we may assume that $\Phi_{1}=\Phi_{2}$, and let it be denoted by $\Phi$. Choose
a pinned element $x\in S$ of $\Phi$. Let $T_{j}=Z_{H_{j}^{0}}(\tilde{S}^{0})$ ($j=1,2$). Then, $T_{j}$ is a
maximal torus of $(H_{j})^{0}$ by Lemma \ref{L:T}. Define a root $\beta\in R(H_{j},T_{j})$ to be positive if
$\beta|_{\tilde{S}^{0}}>0$. Let $\lambda_{j}$ be the highest weight of $\rho|_{T_{j}}$ with respect to this
positive system. By Lemma \ref{L:ARS-isogeny2} there is an isomorphism $\sigma:\mathfrak{h}_{1}\rightarrow
\mathfrak{h}_{2}$ such that $\sigma|_{\mathfrak{s}}=\id$ and $\sigma\circ\Ad(x)=\Ad(x)\circ\sigma$. Then,
the restriction characters on $\tilde{S}^{0}$ of two irreducible representations of $\mathfrak{h}_{1}$ with
highest weights $\sigma^{\ast}\lambda_{2}$ and $\lambda_{1}$ are both equal to $\chi_{\rho}|_{\tilde{S}^{0}}$.
From the decomposition of $\Phi$ into irreducible factors, we get direct product decompositions $\mathfrak{h}_{j}
=\bigoplus_{1\leq i\leq s}\mathfrak{h}_{j,i}$ of $\mathfrak{h}_{j}$ ($j=1,2$). Then, $\sigma(\mathfrak{h}_{1,i})
=\mathfrak{h}_{2,i}$ ($1\leq i\leq s$). Write $\lambda_{j,i}$ for the component of $\lambda_{j}$ on the $i$-th
factor of $\mathfrak{h}_{j}$ ($j=1,2$, $1\leq i\leq s$). By Lemma \ref{L:Levi4}, for each $i$ ($1\leq i\leq s$)
we have $\sigma^{\ast}\lambda_{2,i}=\lambda_{1,i}$ or $\Ad(x)^{\ast}\lambda_{1,i}$. Replacing the action of
$\sigma$ on the $i$-th factor of $\mathfrak{h}_{1}$ by that of $\Ad(x)\circ\sigma$ for some $i$ ($1\leq i\leq s$)
if necessary, we may assume that $\sigma^{\ast}\lambda_{2,i}=\lambda_{1,i}$ for all $i$ ($1\leq i\leq s$). Then,
$\sigma^{\ast}\lambda_{2}=\lambda_{1}$ and $V_{H_{1}^{0},\lambda_{1}}\cong\sigma^{\ast}V_{H_{2}^{0},\lambda_{2}}$.
Thus, there exists $g\in\U(N)$ such that $\sigma=\Ad(g)|_{\mathfrak{h}_{1}}$. Hence, $gH_{1}^{0}g^{-1}=H_{2}^{0}$
and $[g,x]$ commutes with $H_{1}^{0}$. By the irreducibility of $H_{1}^{0}$ (cf. Lemma \ref{L:TA4}), there
exists $z\in\mathbb{C}$ with $|z|=1$ such that $[g,x]=z^{2}I$. Put $g'=(z^{-1}I)g$. Then, $\Ad(g')|_{H_{1}^{0}}=
\Ad(g)|_{H_{1}^{0}}$ and $g'xg'^{-1}=x$. Thus, $g'H_{1}g'^{-1}=H_{2}$ and $\Ad(g')|_{\tilde{S}}=\id$. Hence,
$H_{1}$ and $H_{2}$ are conjugate.
\end{proof}

\subsection{Distinction and linear independence of dimension data}\label{SS:TA1-linear}

For each $i$ ($1\leq i\leq r$), let $\tilde{S}_{i}^{0}$ be the neutral subgroup of
\[\{y\in\tilde{S}^{0}:\alpha(y)=1,\forall\alpha\in\Psi_{j},j\neq i\}.\] Then, the natural map \[\tilde{S}_{1}^{0}
\times\cdots\times\tilde{S}_{r}^{0}\rightarrow\tilde{S}^{0},\ (y_{1},\dots,y_{r})\mapsto y_{1}\cdots y_{r}\]
is a finite covering. From this decomposition, it is associated a canonical character $\chi_{i}$ on
$\tilde{S}_{i}^{0}$ such that: \[\chi=\chi_{1}\otimes\cdots\otimes\chi_{r}.\]

\begin{lemma}\label{L:TA1-char1}
Suppose that each $(\Psi_{i},\chi_{i})$ ($1\leq i\leq r$) is not isomorphic to $(\B_{n}^{(2)},([\frac{1}{2}]+
[-\frac{1}{2}])^{\otimes n})$ ($n\geq 2$). Let $H_{1},H_{2},\dots,H_{s}$ be a set of pairwise non-conjugate
S-subgroups of $G$ with $\tilde{S}$ a common quasi Cartan subgroup, which are all semisimple or all non-semisimple.
Then the characters $F_{\Phi_{1},\Gamma^{\circ}},\dots,F_{\Phi_{s},\Gamma^{\circ}}$ are linearly independent, where
$\Phi_{j}=R(H_{j},\tilde{S})\subset\Psi$ is the quasi root system of $H_{j}$ with respect to $\tilde{S}$
($1\leq j\leq s$).
\end{lemma}

\begin{proof}
Write \[L_{\Psi}=\span_{\mathbb{Z}}p(\Psi),\] which is a sub-lattice of $X^{\ast}(\tilde{S}^{0})$. Since the action of
$\Gamma^{\circ}$ on $X^{\ast}(\tilde{S})$ preserves $\Psi$, its induced action on $X^{\ast}(\tilde{S}^{0})$ preserves
$p(\Psi)$. Then, $L_{\Psi}$ is preserved by $\Gamma^{\circ}$. For each $j$ ($1\leq j\leq s$), write $F'_{\Phi_{j},
\Gamma^{\circ}}$ for the sum of terms $[\lambda]$ in $F_{\Phi_{j},\Gamma^{\circ}}$ such that \[\lambda|_{\tilde{S}^{0}}
\in 2L_{\Psi}.\] Note that, each $t_{y}$ ($y\in\tilde{S}\cap G^{0}$, $y^{2}=I$) fixes all such linear characters.
Then, $F'_{\Phi_{j},\Gamma^{\circ}}$ is preserved by $\Gamma^{\circ}$ and is fixed by all $t_{y}$ ($y\in\tilde{S}\cap G^{0}$,
$y^{2}=I$). It suffices to show that: the characters \[F'_{\Phi_{1},\Gamma^{\circ}},\dots,F'_{\Phi_{s},\Gamma^{\circ}}\]
are linearly independent. By Lemma \ref{L:ARS-conjugacy-TA1}, the quasi root systems $\Phi_{1},\dots,\Phi_{s}$ are
pairwise $\Gamma^{\circ}$ non-conjugate as the S-subgroups $H_{1},H_{2},\dots,H_{s}$ are pairwise non-conjugate and are
all semisimple or all non-semisimple. We show that: the leading weights of $\{F'_{\Phi_{j},\Gamma^{\circ}}:1\leq j\leq s\}$
are distinct, which implies the above assertion.

Without loss of generality we suppose that $F'_{\Phi_{1},\Gamma^{\circ}}$ and $F'_{\Phi_{2},\Gamma^{\circ}}$ have equal
leading weights. Replacing $\Phi_{2}$ by a $\Gamma^{\circ}$ conjugate of it if necessary, we may assume that $A_{\Phi_{1}}$
and $A_{\Phi_{2}}$ have equal actual leading weights (cf. Lemma \ref{L:muPhi}), i.e., $\bar{\mu}_{\Phi_{1}}=
\bar{\mu}_{\Phi_{2}}$. For $j=1,2$, write $\Phi_{j,i}=\Phi_{j}\cap\Psi_{i}$ for each $i$ ($1\leq i\leq r$). We have
\[\bar{\mu}_{\Phi_{j}}=(\bar{\mu}_{\Phi_{j,1}},\dots,\bar{\mu}_{\Phi_{j,r}}).\] Then, $\bar{\mu}_{\Phi_{1,i}}=
\bar{\mu}_{\Phi_{2,i}}$ for each $i$ ($1\leq i\leq r$). By assumption, each $(\Psi_{i},\chi_{i})$ ($1\leq i\leq r$) is
not isomorphic to $(\B_{n}^{(2)},([\frac{1}{2}]+[-\frac{1}{2}])^{\otimes n})$ ($n\geq 2$). By Lemmas \ref{L:TA1-Psi0},
\ref{L:TA1-Psi2} and \ref{L:fully2}, the triple $(\Psi_{i},\chi_{i},\Phi_{j,i})$ ($j=1,2$) fall into the following list:
\begin{enumerate}
\item[(i)]When $p(\Psi_{i})$ is simply-laced, we have $\Psi_{i}\cong p(\Psi_{i})^{(2)}$ and $\Phi_{j,i}\cong
 p(\Psi_{i})^{(2)}$ or $p(\Psi_{i})$.
\item[(ii)]When $p(\Psi_{i})\cong\B_{n}$ or $\BC_{n}$ and $\chi_{i}$ is in-decomposable ($n\geq 2$), we have
$\Psi_{i}\cong\B_{n}^{(2)}$ or $\BC_{n}^{(2,2,1)}$ and $\Phi_{j,i}\cong\B_{n}$, $\B_{n}^{(2,1)}$, $\B_{n}^{(2)}$ or
$\BC_{n}^{(2,2,1)}$.
\item[(iii)]When $p(\Psi_{i})\cong\B_{n}$ or $\BC_{n}$ and $\chi_{i}=(\tau_{1}^{2})^{\otimes n}$, we have
$\Psi_{i}\cong\B_{n}^{(2)}$ or $\BC_{n}^{(2,2,1)}$ and $\Phi_{j,i}\cong\bigsqcup_{1\leq k\leq u}\B_{n_{k}}^{(2)}$
where $\sum_{1\leq k\leq u}n_{k}=n$.
\item[(iv)]When $p(\Psi_{i})\cong\B_{n}$ or $\BC_{n}$ and $\chi_{i}=\tau_{2}^{\otimes n}$, we have $\Psi_{i}\cong
\B_{n}^{(2)}$ or $\BC_{n}^{(2,2,1)}$ and \[\Phi_{j,i}\cong(\sqcup_{1\leq k\leq u}\B_{1})\bigsqcup
(\sqcup_{1\leq k'\leq u'}\BC_{1}^{(2,1)})\] where $u+u'=n$.
\item[(v)]When $p(\Psi_{i})\cong\C_{n}$ ($n\geq 3$), we have $\Psi_{i}\cong\C_{n}^{(2,1)}$ or $\C_{n}^{(2)}$ and
$\Phi_{j,i}\cong\D_{n}$, $\C_{n}$, $\D_{n}^{(2)}$, $\C_{n}^{(2,1)}$ or $\C_{n}^{(2)}$.
\item[(vi)]When $p(\Psi_{i})\cong\F_{4}$, we have $\Psi_{i}\cong\F_{4}^{(2)}$ and $\Phi_{j,i}\cong\D_{4}$,
$\C_{4}$, $\F_{4}$, $\D_{4}^{(2)}$, $\C_{4}^{(2,1)}$, $\F_{4}^{(2,1)}$, $\C_{4}^{(2)}$ or $\F_{4}^{(2)}$.
\item[(vii)]When $p(\Psi_{i})\cong\G_{2}$, we have $\Psi_{i}\cong\G_{2}^{(2)}$ and $\Phi_{j,i}\cong\A_{2}^{S}$,
$\G_{2}$, $\A_{2}^{S(2)}$ or $\G_{2}^{(2)}$.
\end{enumerate}
By the formulas of $\bar{\mu}_{\Phi}$ given in the proof of Lemma \ref{L:muPhi}, we see that in each case (i)-(vii)
the condition $\bar{\mu}_{\Phi_{1,i}}=\bar{\mu}_{\Phi_{2,i}}$ forces that: $\Phi_{1,i}\cong\Phi_{2,i}$. By Lemma
\ref{L:TA1-fracFactor3}, it follows that $\Phi_{1}$ and $\Phi_{2}$ are $\Gamma^{\circ}$ conjugate, which gives a
contradiction.
\end{proof}

\begin{proof}[Proof of Theorem \ref{T3}.]
Let $\{H_{1},\dots,H_{s}\}$ be a tuple of pairwise non-conjugate S-subgroups of $G$ satisfying conditions in the
assumption. Suppose that the dimension data \[\mathscr{D}_{H_{1}},\dots,\mathscr{D}_{H_{s}}\] are linearly dependent.
For each $j$ ($1\leq j\leq s$), choose a maximal commutative connected subset $S_{j}$ of $H_{j}$ generating
$H_{j}/H_{j}^{0}$. Write $d_{j}$ for the minimal order of elements in $S_{j}$. Without loss of generality we assume
that $\dim S_{j}\geq\dim S_{1}$ for each $j$; $d_{j}\leq d_{1}$ whenever $\dim S_{j}=\dim S_{1}$; and
$(\dim S_{j},d_{j})=(\dim S_{1},d_{1})$ holds exactly when $1\leq j\leq s'$ for some integer $s'\leq s$. Moreover,
for some integer $t\leq s'$, $S_{j}$ is conjugate to $S_{1}$ when $1\leq j\leq t$ and any maximal commutative
connected subset $S'_{j}$ of $H_{j}$ generating $H_{j}/H_{j}^{0}$ is not conjugate to $S_{1}$ whenever
$t+1\leq j\leq s'$. Then, by Lemma \ref{L:MCC1}, any maximal commutative connected subset $S'_{j}$ of $H_{j}$ is
not conjugate to $S_{1}$ whenever $t+1\leq j\leq s$. Write $S=S_{1}$ and let $\tilde{S}$ be the quasi torus
generated by $S$. Put $\Gamma^{\circ}=N_{G^{0}}(\tilde{S})/Z_{G^{0}}(\tilde{S})$ and write
$\Phi_j=R(H_j,\tilde{S})$ ($1\leq j\leq t$). Then, the characters
\[F_{\Phi_{1},\Gamma^{\circ}},\dots,F_{\Phi_{t},\Gamma^{\circ}}\] are linearly dependent by Proposition
\ref{P:dim-char} and Corollary \ref{C:dim-char2}, which is in contradiction with Lemma \ref{L:TA1-char1}.
\end{proof}

\begin{lemma}\label{L:TA1-char2}
Let $H_{1},H_{2}$ be two non-conjugate S-subgroups of $G$ both with $\tilde{S}$ a quasi Cartan subgroup. Suppose
that $H_{1}$ and $H_{2}$ are both semisimple or both non-semisimple. Let $\Phi_{1}\subset\Psi_{\tilde{S}}$ (resp.
$\Phi_{2}\subset\Psi_{\tilde{S}}$) be the quasi root system of $H_{1}$ (resp. $H_{2}$) with respect to $\tilde{S}$.
Then \[F_{\Phi_{1},\Gamma^{\circ}}\neq F_{\Phi_{2},\Gamma^{\circ}}.\]
\end{lemma}

\begin{proof}
Suppose that $F_{\Phi_{1},\Gamma^{\circ}}=F_{\Phi_{2},\Gamma^{\circ}}$. By Lemma \ref{L:fully3}, there is at
most one index $i$ such that \[(\Psi_{i},\chi_{i})\cong(\B_{n}^{(2)},([\frac{1}{2}]+[-\frac{1}{2}])^{\otimes n}).\]
When there is no such index $i$, $F_{\Phi_{1},\Gamma^{\circ}}$ and $F_{\Phi_{2},\Gamma^{\circ}}$ are linearly
independent by Lemma \ref{L:TA1-char1}, which gives a contradiction. When there is such an index
$i$, without loss of generality we assume that $(\Psi_{r},\chi_{r})\cong(\B_{n}^{(2)},([\frac{1}{2}]+
[-\frac{1}{2}])^{\otimes n})$. Write \[\Psi'=\bigsqcup_{1\leq i\leq r-1}\Psi_{i}.\] Then, \begin{equation}
\label{Eq:split}\Psi=\Psi'\bigsqcup\Psi_{r}\end{equation} and both $\Psi'$ and $\Psi_{r}$ are stable under the
action of $\Gamma^{\circ}$. Thus, we have natural homomorphisms $\Gamma^{\circ}\rightarrow\Aut(\Psi')$ and
$\Gamma^{\circ}\rightarrow\Aut(\Psi_{r})$. Let $\Gamma'$ and $\Gamma_{r}$ be the images of them, respectively.
Then, $F_{\Phi_{1},\Gamma^{\circ}}=F_{\Phi_{2},\Gamma^{\circ}}$ indicates that \begin{equation}\label{Eq:TA1-char1}
F_{\Phi_{1},\Gamma'\times\Gamma_{r}}=F_{\Phi_{2},\Gamma'\times\Gamma_{r}}.\end{equation} Write $\Phi_{j,i}=
\Phi_{j}\cap\Psi_{i}$ ($j=1,2$, $1\leq i\leq r$) and \[\Phi'_{j}=\Phi_{j}\cap\Psi'=\bigsqcup_{1\leq i\leq r-1}
\Psi_{j,i},\ j=1,2.\] Then, \begin{equation}\label{Eq:TA1-char2}F_{\Phi_{j},\Gamma'\times\Gamma_{r}}=
F_{\Phi'_{j},\Gamma'}\otimes F_{\Phi_{j,r}\Gamma_{r}}.\end{equation} By \eqref{Eq:TA1-char1} and
\eqref{Eq:TA1-char2}, we get \begin{equation}\label{Eq:TA1-char3}F_{\Phi'_{1},\Gamma'}=F_{\Phi'_{2},\Gamma'}
\end{equation} and \begin{equation}\label{Eq:TA1-char4}F_{\Phi_{1,r},\Gamma_{s}}=F_{\Phi'_{2,r},\Gamma_{r}}.
\end{equation}

By \eqref{Eq:TA1-char3}, replacing $H_{2}$ by a conjugate one if necessary we may assume that
$\bar{\mu}_{\Phi'_{1}}=\bar{\mu}_{\Phi'_{2}}$. By the same argument as in the proof of Lemma \ref{L:TA1-char1}
one shows that \begin{equation}\label{Eq:TA1-Phi1}\Phi_{1,i}\cong\Phi_{2,i},\ \forall i, 1\leq i\leq s-1.
\end{equation} By Lemma \ref{L:fully3}, we have \[\Phi_{1,r}\cong(\bigsqcup_{1\leq i\leq t_{1}}\B_{n_{i}}^{(2,1)})
\bigsqcup(\bigsqcup_{1\leq i'\leq t'_{1}}\B_{n'_{i'}})\] and \[\Phi_{2,r}\cong(\bigsqcup_{1\leq j\leq t_{2}}
\B_{m_{j}}^{(2,1)})\bigsqcup(\bigsqcup_{1\leq j'\leq t'_{2}}\B_{m'_{j'}}),\] where $n_{i},m_{j}$ are even and
\[\sum_{1\leq i\leq t_{1}}n_{i}+\sum_{1\leq i'\leq t'_{1}}n'_{i'}=\sum_{1\leq j\leq t_{2}}m_{j}+
\sum_{1\leq j'\leq t'_{2}}m'_{j'}=n.\] By Lemma \ref{L:Ave2}, we get \begin{equation}\label{Eq:TA1-multi}
\prod_{1\leq i\leq t_{1}}(b_{n_{i}/2}c_{n_{i}/2})\prod_{1\leq i'\leq t'_{1}}a_{n'_{i'}}=\prod_{1\leq j\leq t_{2}}
(b_{m_{j}/2}c_{m_{j}/2})\prod_{1\leq j'\leq t'_{2}}a_{m'_{j'}}\end{equation} from \eqref{Eq:TA1-char4}, where
$a_{k}$, $b_{k}$, $c_{k}$ (and $b'_{k}=\sigma(b_{k})$ and $d_{k+1}$ below) are polynomials as defined in
\cite[\S 7]{Yu-dimension}. By \cite[Proposition 7.2]{Yu-dimension}, we have \[a_{2k}=b_{k}b'_{k}\] and \[a_{2k+1}=
c_{k}d_{k+1}\] for any $k\in\mathbb{Z}_{\geq 1}$. By \cite[Proposition 7.3]{Yu-dimension}, there is no nontrivial
multiplicative relation among the polynomials \[\{b_{k},b'_{k},c_{k},d_{k+1}:k\geq 1\}.\] Then, \eqref{Eq:TA1-multi}
indicates that $t_{1}=t_{2}$ and $\{n_{i}:1\leq i\leq t_{1}\}=\{m_{j}:1\leq j\leq t_{2}\}$, $t'_{1}=t'_{2}$ and
$\{n'_{i}:1\leq i'\leq t'_{1}\}=\{m'_{j}:1\leq j'\leq t'_{2}\}$. Thus, \begin{equation}\label{Eq:TA1-Phi2}
\Phi_{1,r}\cong\Phi_{2,r}.\end{equation} By Lemma \ref{L:TA1-fracFactor3}, \eqref{Eq:TA1-Phi1} and
\eqref{Eq:TA1-Phi2} imply that $\Phi_{1}$ and $\Phi_{2}$ are $\Gamma^{\circ}$ conjugate. By Lemma
\ref{L:ARS-conjugacy-TA1} and the assumption that $H_{1}$ and $H_{2}$ are both semisimple or both non-semisimple,
it follows that $H_{1}$ and $H_{2}$ are conjugate, which contradicts to the assumption of the lemma.
\end{proof}

\begin{proof}[Proof of Theorem \ref{T4}.]
Suppose that $\mathscr{D}_{H_{1}}=\mathscr{D}_{H_{2}}$. Choose a maximal torus $T_{j}$ of $H_{j}^{0}$ and a maximal
commutative connected subset $S_{j}$ of $H_{j}$ generating $H_{j}/H_{j}^{0}$ ($j=1,2$). Then, $T_{1}$ is $G$ conjugate
to $T_{2}$ and we may assume that $S_{1}$ is $G^{0}$ conjugate to $S_{2}$. Without loss of generality we assume that
$S_{1}=S_{2}$ and let it be denoted by $S$. Let $\tilde{S}$ be the quasi torus generated by $S$. Put \[\Gamma^{\circ}
=N_{G^{0}}(S)/Z_{G^{0}}(S).\] For each $j$ ($j=1,2$), write $\Phi_{j}=R(H_{j},\tilde{S})$ for the quasi root system
of $H_{j}$ with respect to $\tilde{S}$.

Since $T_{1}\sim T_{2}$, they either both contain the center of $\U(N)$, or neither contains that. By Lemma \ref{L:TA4},
both $H_{1}^{0}$ and $H_{2}^{0}$ are irreducible subgroups of $\U(N)$. When $T_{1}$ and $T_{2}$ do not contain the
center of $\U(N)$, both $H_{1}$ and $H_{2}$ are semisimple; when they neither contain that, $H_{1}$ and $H_{2}$ are
both non-semisimple. By Lemma \ref{L:TA1-char2} we get \[F_{\Phi_{1},\Gamma^{\circ}}\neq F_{\Phi_{2},\Gamma^{\circ}}.\]
Then, $\mathscr{D}_{H_{1}}\neq\mathscr{D}_{H_{2}}$ by Proposition \ref{P:dim-char} and Corollary \ref{C:dim-char2}.
\end{proof}

In Theorem \ref{T3}, we impose two conditions: (1)for any quasi Cartan subgroup $\tilde{S}$ of any of
$H_{1},\dots,H_{s}$, $(\Psi_{\tilde{S}},\chi_{\rho}|_{\tilde{S}^{0}})$ does not contain an irreducible factor
$(\Psi_{i},\chi_{i})$ isomorphic to \[(\B_{n}^{(2)},([\frac{1}{2}]+[-\frac{1}{2}])^{\otimes n})\] ($n\geq 2$).
(2)$H_{1},\dots,H_{s}$ are all semisimple or all non-semisimple. In the next two subsections, we construct two
kinds of examples of S-subgroups of $G=\U(N)\rtimes\langle\tau\rangle$ with linearly dependent dimension data.
This indicates that both conditions are necessary for showing linear independence of dimension data.

\subsection{S-subgroups with linearly dependent dimension data, I}\label{SS:counter1}

When $n\geq 3$ is odd, the spin group $\Spin(n)$ has a spinor module $M_{n}$ of dimension $2^{\frac{n-1}{2}}$,
which is a fundamental irreducible representation. When $n\geq 4$ is even, the spin group $\Spin(n)$ has two
spinor modules $M_{+,n}$ and $M_{-,n}$ both of dimension $2^{\frac{n-2}{2}}$, which are both fundamental
irreducible representations. The spinor modules $M_{+,n}$ and $M_{-,n}$ are distinguished by their central
characters and they are permuted by an outer involutive automorphism of $\Spin(n)$. The self-duality and
orthogonal/symplectic type of spinor modules as described in the following lemma are well-known.

\begin{lemma}\label{L:spinor1}
\begin{enumerate}
\item[(i)]When $n\equiv -1$ or $1$ $\pmod{8}$, the spin module $M_{n}$ is self-dual and of orthogonal type.
\item[(ii)]When $n\equiv 3$ or $5$ $\pmod{8}$, the spin module $M_{n}$ is self-dual and of symplectic type.
\item[(iii)]When $n\equiv 0$ $\pmod{8}$, the spin modules $M_{+,n}$ and $M_{-,n}$ are self-dual and of orthogonal type.
\item[(iv)]When $n\equiv 4$ $\pmod{8}$, the spin modules $M_{+,n}$ and $M_{-,n}$ are self-dual and of symplectic type.
\item[(v)]When $n\equiv 2$ or $6$ $\pmod{8}$, the spin modules $M_{+,n}$ and $M_{-,n}$ are non-self dual and they
are dual to each other.
\end{enumerate}
\end{lemma}

Recall that the {\it Schur index} of a self-dual irreducible finite-dimensional complex linear representation $\rho$
of a compact Lie group $K$ is defined as follows: it is 1 if $\rho$ is of orthogonal type; it is $-1$ if $\rho$ is
of symplectic type.

\begin{lemma}\label{L:spinor2}
Let $n_1,\dots,n_{s}$ be positive integers such that $\sum_{1\leq i\leq s}n_{i}=n$ and let $k$ be the number
of indices $i$ ($1\leq i\leq s$) such that $n_{i}\equiv 1$ or $2$ ($\pmod{4}$). Then, the Schur index of
$M_{2n_{1}+1}\otimes\cdots\otimes M_{2n_{s}+1}$ is equal to $(-1)^{k}$.
\end{lemma}

\begin{proof}
By Lemma \ref{L:spinor1}, the Schur index of $M_{2n_{i}+1}$ is equal to \[\left\{\begin{array}{cc}-1\textrm{ if }
n_{i}\equiv 1\textrm{ or }2\pmod{4}\\1\textrm{ if }n_{i}\equiv 0\textrm{ or }3\pmod{4}.\end{array}\right.\] Then,
the Schur index of $M_{2n_{1}+1}\otimes\cdots\otimes M_{2n_{s}+1}$ is equal to $(-1)^{k}$ by multiplicativity of
Schur index.
\end{proof}

Let $n$ be a positive integer. Put $N=2^{n}$ and \[G_{n}=\U(N)\rtimes\langle\tau\rangle.\] We construct a quasi
torus $\tilde{S}_{n}$ of $G_{n}$ as follows. Take a torus $T_{n}=\U(1)^{n}$. Write $(e^{\mathbf{i}\theta_{1}},\dots,
e^{\mathbf{i}\theta_{n}})$ ($\theta_{i}\in\mathbb{R}$) for a general element of $T_{n}$. For each $i$ ($1\leq i\leq n$),
define a 2-dimensional orthogonal representation $\sigma_{i}$ of $T_{n}$ by \[\sigma_{i}(e^{\mathbf{i}\theta_{1}},
\dots,e^{\mathbf{i}\theta_{n}})=\left(\begin{array}{cc}\cos\theta_{i}&\sin\theta_{i}\\-\sin\theta_{i}&
\cos\theta_{i}\\\end{array}\right).\] Write $\sigma=\sigma_{1}\otimes\cdots\otimes\sigma_{n}$. Then, $\sigma(T_{n})
\subset\SO(N)$. Put \[S_{n}=\sigma(T_{n})\tau\] and let $\tilde{S}_{n}$ be the quasi torus generated by $S_{n}$.
Write \[\Psi_{n}=\Psi_{\tilde{S}_{n}}\] and \[\Gamma_{n}^{\circ}=N_{G_{n}^{0}}(S_{n})/Z_{G_{n}^{0}}(S_{n}).\]

\begin{lemma}\label{L:Bn-Phi}
Let $n_1,\dots,n_{s}$ be even positive integers and $n'_{1},\dots,n'_{t}$ be positive integers such that
$\sum_{1\leq i\leq s}n_{i}+\sum_{1\leq j\leq t}n'_{j}=n$. Then there exists an S-subgroup $H$ of $G_{n}$ with
$\tilde{S}_{n}$ a quasi Cartan subgroup and with quasi root system \[R(H,\tilde{S}_{n})\cong
(\bigsqcup_{1\leq i\leq s}\B_{n_{i}}^{(2,1)})\bigsqcup(\bigsqcup_{1\leq j\leq t}\B_{n'_{j}}).\]
\end{lemma}

\begin{proof}
Put \[K\!=\!(\!\Spin(2n_{1}\!+\!2)\!\times\!\cdots\!\times\!\Spin(2n_{s}\!+\!2)\!)\!\times\!(\!\Spin(2n'_{1}\!+\!1)\!
\times\!\cdots\!\times\!\Spin(2n'_{t}\!+\!1)\!)\] and \[\rho=(M_{2n_{1}+2,+}\otimes\cdots\otimes M_{2n_{s}+2,+})\times
(M_{2n'_{1}+1}\otimes\cdots\otimes M_{2n'_{t}+1}).\] Take a subgroup of $K$: \[K'=\!(\!\Spin(2n_{1}\!+\!1)\!\times\!
\cdots\!\times\!\Spin(2n_{s}\!+\!1)\!)\!\times\!(\!\Spin(2n'_{1}\!+\!1)\!\times\!\cdots\!\times\!\Spin(2n'_{t}\!+\!1)\!.\]
Then, there exists an involutive automorphism $\theta$ of $K$ such that $K^{\theta}=K'$. Choose an element $g_{1}\in
\tau\U(N)$. Define \[\rho'(x)=g_{1}^{-1}\rho(\theta(x))g_{1}.\] Then, $\rho'\cong(\theta\rho)^{\ast}\cong\rho$ by Lemma
\ref{L:spinor1}. Thus, there exists $g_{2}\in\U(N)$ such that \[\rho'(x)=g_{2}\rho(x)g_{2}^{-1},\ \forall x\in K.\] Put
$g=g_{1}g_{2}$. Then, \begin{equation}\label{Eq:Bn-Phi1}\rho(\theta(x))=g\rho(x)g^{-1},\ \forall x\in K.\end{equation}
Inputting $x\in K'$, we get $g\in Z_{\U(N)\tau}(\rho(K'))$. Since $\rho|_{K'}$ is irreducible, this determines $g$ up
to multiplying a scalar matrix. Since $g^{2}$ commutes with $\rho(K')$ and $g\in\U(N)\tau$, we must have $g^{2}=\pm{I}$.
Choose a maximal torus $T$ of $K'$ and let \[\tilde{S}'=\langle\rho(T),g\rangle.\] Put $H'=\langle\rho(K),g\rangle$. Then,
$H'$ is an S-subgroup of $G$ with $\tilde{S}'$ a quasi Cartan subgroup by Lemma \ref{L:TA4}. By \eqref{Eq:Bn-Phi1},
we have \[R(H',\tilde{S}')\cong(\bigsqcup_{1\leq i\leq s}\B_{n_{i}}^{(2,1)})\bigsqcup(\bigsqcup_{1\leq j\leq t}\B_{n'_{j}}).\]
By the argument in the proof of Lemma \ref{L:TA1-Gamma1} one shows that $\tilde{S}'$ and $\tilde{S}_{n}$ are conjugate
in $G_{n}$. Then, the conclusion of this lemma follows.
\end{proof}

\begin{lemma}\label{L:Bn-Psi}
Let $n\geq 2$. Then we have \[\Psi_{n}\cong\B_{n}^{(2)}\textrm{ and }\Gamma_{n}^{\circ}=W_{\Psi_{n}}.\]
\end{lemma}

\begin{proof}
(1)For each $j$ ($1\leq j\leq n$), define a linear character $e_{j}$ of $T_{n}$ by \[e_{j}(e^{\mathbf{i}\theta_{1}},
\dots,e^{\mathbf{i}\theta_{n}})=e^{2\mathbf{i}\theta_{j}}.\] Then, each $e_{j}$ descends to a character of
$\tilde{S}_{n}^{0}$, still denoted by $e_{j}$. Define an inner product $(\cdot,\cdot)$ on $X^{\ast}(T_{n})$ (and on
$X^{\ast}(\tilde{S}_{n}^{0})$) by taking $\{e_{i}:1\leq i\leq n\}$ an orthonormal basis. By our construction, the
character $\chi_{\sigma}$ is equal to \[([\frac{e_{1}}{2}]+[-\frac{e_{1}}{2}])\otimes\cdots\otimes([\frac{e_{n}}{2}]
+[-\frac{e_{n}}{2}]).\] From this, we see that $(\cdot,\cdot)$ coincides with the inner product induced by the
conjugation invariant inner product $(\cdot,\cdot)'$ on $\mathfrak{u}(N)$ up to a scalar, where \[(X,Y)'=-\tr(XY),
\ \forall X,Y\in\mathfrak{u}(N).\] Let $\bar{\alpha}=\sum_{1\leq i\leq n}a_{i}e_{i}$ be a character in $p(\Psi_{n})$.
Then $\bar{\alpha}(-I)=1$. Hence, $\bar{\alpha}$ descends to a linear character of $\tilde{S}_{n}^{0}/\langle-I\rangle
\cong T_{n}/\{\pm{1}\}^{n}$. Hence, all $a_{i}$ are integers. For any weight $\lambda$ appearing in $\sigma$, we have
\[\frac{2(\lambda,\bar{\alpha})}{(\bar{\alpha},\bar{\alpha})}\in\mathbb{Z}.\] Since all $\sum_{1\leq i\leq n}\pm{}
\frac{e_{i}}{2}$ appear as weights of $\sigma$, we have \[\frac{\sum_{1\leq i\leq n}\pm{a_{i}}}{\sum_{1\leq i\leq n}
a_{i}^{2}}\in\mathbb{Z}.\] Then, each $a_{i}\in\{-1,0,-1\}$ and $\sum_{1\leq i\leq n}a_{i}^{2}\leq 2$. Hence,
\[ p(\Psi_{n})\subset\B_{n}=\{\pm{e}_{i}\pm{e}_{j},\pm{e}_{k}:1\leq i<j\leq n,1\leq k\leq n\}.\] By Lemma
\ref{L:Bn-Phi}, there exists a quasi sub-root system of $\Phi$ of $\Psi_{n}$ such that \[\Phi\cong\left\{
\begin{array}{cc}\B_{n}^{(2,1)}\textrm{ if }n\textrm{ is even}\\\B_{n-1}^{(2,1)}\sqcup\B_{1}\textrm{ if }n
\textrm{ is odd}.\\\end{array}\right.\] Then, $p(\Psi_{n})=\B_{n}$. By Lemma \ref{L:TA1-Psi1}, we get
$\Psi_{n}\cong\B_{n}^{(2)}$.

(2)Apparently, $W_{\Psi_{n}}\subset\Gamma_{n}^{\circ}$. Write $\phi: \Gamma_{n}^{\circ}\rightarrow
\Aut(\tilde{S}_{n}^{0})$ for the natural homomorphism induced from the action of $\Gamma_{n}^{\circ}$ on
$\tilde{S}_{n}^{0}$. It is clear that \[\phi(W_{\Psi_{n}})=\phi(\Gamma_{n}^{\circ})=W_{\B_{n}}.\] We know that
\[\ker\phi=\{t_{y}:y\in\tilde{S}_{n}^{0},y^{2}=I\}.\] For any root $\alpha\in\Psi_{n}$, we have
$s_{\alpha+\beta_{0}}s_{\alpha}=t_{\check{\alpha}(-1)}\in\ker\phi$. The coroot of $e_{1}-e_{2}$ is equal to
\[\sigma(\mathbf{i},-\mathbf{i},1,\dots,1)\] and the coroot of $e_{1}+e_{2}$ is equal to
\[\sigma(\mathbf{i},\mathbf{i},1,\dots,1),\] similarly for coroots of $e_{i}-e_{j}$ and $e_{i}+e_{j}$. Thus,
$\{\check{\alpha}(-1):\alpha\in\Psi_{n}\}$ generate $\{y\in\tilde{S}^{0},y^{2}=I\}$. Hence,
$\ker\phi\subset W_{\Psi_{n}}$. Therefore, $\Gamma_{n}^{\circ}=W_{\Psi_{n}}$.
\end{proof}

We lift each $e_{i}$ ($\leq i\leq n$) to a linear character of $\tilde{S}_{n}$ by defining $e_{i}(\tau)=1$,
and let it be still denoted by $e_{i}$. Then, \[\Psi_{n}=\{\pm{e}_{i}\pm{e}_{j},\pm{e}_{i}\pm{e}_{j}+\beta_{0},
\pm{e}_{k},\pm{e}_{k}+\beta_{0}:1\leq i<j\leq n,1\leq k\leq n\}\] by Lemma \ref{L:Bn-Psi}. Put \[\Psi_{n,0}=
\{\pm{e}_{i}\pm{e}_{j},\pm{e}_{k}:1\leq i<j\leq n,1\leq k\leq n\}\subset\Psi_{n}.\]

\begin{definition}\label{D:Schur2}
Let $\Phi$ be a quasi sub-root system of $\Psi_{n}$ such that $p(\Phi)\supset p(\Psi_{n})^{\circ}$. If $\Phi$
is $\Gamma_{n}^{\circ}$ conjugate to a quasi sub-root system with a simple system contained in $\Psi_{n,0}$,
then we define the Schur index of $\Phi$ to to 1; otherwise we define the Schur index of $\Phi$ to be $-1$.
\end{definition}

\begin{lemma}\label{L:Bn-2}
A quasi sub-root system $\Phi$ of $\Psi_{n}$ occurs as the quasi root system of an S-subgroup of $G_{n}$ with
$\tilde{S}_{n}$ a quasi Cartan subgroup if and only if \[\Phi\cong(\bigsqcup_{1\leq i\leq s}\B_{n_{i}}^{(2,1)})
\bigsqcup(\bigsqcup_{1\leq j\leq t}\B_{n'_{j}})\] for some even positive integers $n_1,\dots,n_{s}$ and positive
integers $n'_{1},\dots,n'_{t}$ such that $\sum_{1\leq i\leq s}n_{i}+\sum_{1\leq j\leq t}n'_{j}=n$ and the Schur
index of $\Phi$ is equal to $(-1)^{k}$, where $k$ is the sum of the number of indices $i$ ($1\leq i\leq s$)
such that $n_{i}\equiv 2\pmod{4}$ and the number of indices $j$ ($1\leq j\leq t$) such that
$n'_{j}\equiv 1,2\pmod{4}$.
\end{lemma}

\begin{proof}
Suppose that $\Phi$ is the quasi root system of an S-subgroup of $G_{n}$ with $\tilde{S}_{n}$ a generalized
Cartan subgroup. By Lemma \ref{L:fully3}, we get \[\Phi\cong(\bigsqcup_{1\leq i\leq s}\B_{n_{i}}^{(2,1)})\bigsqcup
(\bigsqcup_{1\leq j\leq t}\B_{n'_{j}})\] for some even positive integers $n_1,\dots,n_{s}$ and positive integers
$n'_{1},\dots,n'_{t}$ such that $\sum_{1\leq i\leq s}n_{i}+\sum_{1\leq j\leq t}n'_{j}=n$. From the proof of
Lemma \ref{L:TA1-fracFactor3}, we see that the Schur index of $\Phi$ is equal to $(-1)^{k}$. By Lemma \ref{L:Bn-Phi},
each such quasi sub-root system of $\Psi_{n}$ occurs as the quasi root system of an S-subgroup of $G_{n}$ with
$\tilde{S}$ a quasi Cartan subgroup.
\end{proof}

Put \[W_{n}=W_{\Psi_{n,0}}.\] Write $Y_{n,2}$ for the algebra of the restriction to $\tau\tilde{S}_{n}^{0}$ of
$W_{n}$ invariant rational coefficient characters on $\tilde{S}_{n}/\langle-I\rangle$. It is clear that $Y_{n,2}$
possesses a set of basis $\chi^{\ast}_{a_{1}e_{1}+\cdots+a_{n}e_{n}+k\beta_{0},W_{n}}$, where $a_{1},\dots,a_{n}$
are integers with $a_{1}\geq\cdots a_{n}\geq 0$ and $k\in\{0,1\}$. The following lemma is a special case of Lemma
\ref{L:Sn-Wn}.

\begin{lemma}\label{L:Sn-Wn2}
By mapping $\chi^{\ast}_{a_{1}e_{1}+\cdots+a_{n}e_{n}+k\beta_{0},W_{n}}$ to $(-1)^{k}x_{a_{1}}\cdots x_{a_{n}}$, we
get an algebra injection \[E:\lim_{\longrightarrow_{n}}Y_{n,2}\hookrightarrow\mathbb{C}[x_{0},x_{1},\dots].\]
\end{lemma}

Note that $W_{\B_{n}^{(2)}}$ is not equal to $W_{\BC_{n}^{(2,2,1)}}$, but an index 2 subgroup of it. For this
reason we can not use results in Subsection \ref{SS:BCn} directly. For each $p\geq 1$, write $b_{p}^{'+}$
(resp. $c_{p}^{'+}$) for the sum of terms in $b_{p}$ (resp. $c_{p}$) of the form $x_{2}^{t_{2}}x_{4}^{t_{4}}
\cdot x_{2q}^{t_{2q}}\cdots$; write $b_{p}^{'-}$ (resp. $c_{p}^{'-}$) for the sum of terms in $b_{p}$
(resp. $c_{p}$) of the form $x_{1}^{t_{1}}x_{3}^{t_{3}}\cdot x_{2q-1}^{t_{2q-1}}\cdots$. Put \[b_{p}^{+}=
b_{p}^{'+},\quad b_{p}^{-}=(-1)^{\frac{p(p+1)}{2}}b_{p}^{'-}\] and \[c_{p}^{+}=c_{p}^{'+},\quad c_{p}^{-}=
(-1)^{\frac{p(p+1)}{2}}c_{p}^{'-}.\]

\begin{lemma}\label{L:poly-Phi}
Let $\Phi$ be the quasi root system of an S-subgroup of $G_{n}$ with $\tilde{S}_{n}$ a quasi Cartan subgroup.
Suppose that \[\Phi\cong(\bigsqcup_{1\leq i\leq s}\B_{n_{i}}^{(2,1)})\bigsqcup(\bigsqcup_{1\leq j\leq t}\B_{n'_{j}})\]
for some even positive integers $n_1,\dots,n_{s}$ and positive integers $n'_{1},\dots,n'_{t}$ with $\sum_{1\leq i\leq s}
n_{i}+\sum_{1\leq j\leq t}n'_{j}=n$. Then, \[E(F_{\Phi,W_{\Psi_{n}}})=\prod_{1\leq i\leq s}c_{n_{i}}^{+}
\prod_{1\leq j\leq t}b_{n'_{j}}^{+}+\prod_{1\leq i\leq s}c_{n_{i}}^{-}\prod_{1\leq j\leq t}b_{n'_{j}}^{-}.\]
\end{lemma}

\begin{proof}
When $0\leq i\leq s$, write $m_{i}=\sum_{1\leq i'\leq i}n_{i}$. When $0\leq j\leq t$, write $m'_{j}=m_{s}+
\sum_{1\leq j'\leq j}n'_{j}$. Let $k$ be the sum of the number of indices $i$ ($1\leq i\leq s$) such that $n_{i}
\equiv 2\pmod{4}$ and the number of indices $j$ ($1\leq j\leq t$) such that $n'_{j}\equiv 1,2\pmod{4}$. By Lemma
\ref{L:Bn-2}, we can take the following simple systems of irreducible factors of $\Phi$: \begin{enumerate}
\item[(i)]When $k$ is even, $\{e_{m_{i-1}+1}-e_{m_{i-1}+2},\dots,e_{m_{i}-1}-e_{m_{i}},e_{m_{i}}\}$ (type
$\B_{n_{i}}^{(2,1)}$, $1\leq i\leq s$), $\{e_{m'_{j-1}+1}-e_{m'_{j-1}+2},\dots,e_{m'_{j}-1}-e_{m'_{j}},
e_{m'_{j}}\}$ (type $\B_{n'_{j}}$, $1\leq j\leq t$).
\item[(ii)]When $k$ is odd and $s\geq 1$, $\{e_{1}-e_{2}+\beta_{0},\dots,e_{m_{1}-1}-e_{m_{1}},e_{m_{1}}\}$ (type
$\B_{n_{1}}^{(2,1)}$), $\{e_{m_{i-1}+1}-e_{m_{i-1}+2},\dots,e_{m_{i}-1}-e_{m_{i}},e_{m_{i}}\}$ (type
$\B_{n_{i}}^{(2,1)}$, $2\leq i\leq s$), $\{e_{m'_{j-1}+1}-e_{m'_{j-1}+2},\dots,e_{m'_{j}-1}-e_{m'_{j}},
e_{m'_{j}}\}$ (type $\B_{n'_{j}}$, $1\leq j\leq t$).
\item[(iii)]When $k$ is odd and $s=0$, $\{e_{1}-e_{2}+\beta_{0},\dots,e_{m'_{1}-1}-e_{m'_{1}},e_{m'_{1}}\}$ (type
$\B_{n'_{1}}$), $\{e_{m'_{j-1}+1}-e_{m'_{j-1}+2},\dots,e_{m'_{j}-1}-e_{m'_{j}},e_{m'_{j}}\}$ (type $\B_{n'_{j}}$,
$2\leq j\leq t$).
\end{enumerate}

Note that $F_{\Phi,W_{n}}$ is a finite linear combination of terms of the form \[\chi^{\ast}_{a_{1}e_{1}+\dots+
a_{n}e_{n}+a_{0}\beta_{0},W_{n}}=\frac{1}{|W_{n}|}\sum_{\gamma\in W_{n}}\gamma([a_{1}e_{1}+\dots+a_{n}e_{n}+
a_{0}\beta_{0}])\] where $a_{1}\geq a_{2}\geq\cdots\geq a_{n}\geq 0$ and $a_{0}\in\{0,1\}$ and \[F_{\Phi,
W_{\Psi_{n}}}=\frac{1}{|W_{\Psi}|}\sum_{\gamma\in W_{\Psi}}\gamma F_{\Phi,W_{n}}.\] When $a_{1},a_{2},\dots,
a_{n}$ are all even or all odd, we have \[\gamma\chi^{\ast}_{a_{1}e_{1}+\dots+a_{n}e_{n}+a_{0}\beta_{0},W_{n}}
=\chi^{\ast}_{a_{1}e_{1}+\dots+a_{n}e_{n}+a_{0}\beta_{0},W_{n}},\ \forall\gamma\in W_{\Psi}.\] Thus,
\[\frac{1}{|W_{\Psi_{n}}|}\sum_{\gamma\in W_{\Psi_{n}}}\gamma\chi^{\ast}_{a_{1}e_{1}+\dots+a_{n}e_{n}+
a_{0}\beta_{0},W_{n}}=\chi^{\ast}_{a_{1}e_{1}+\dots+a_{n}e_{n}+a_{0}\beta_{0},W_{n}}.\] When $a_{1},a_{2},\dots,
a_{n}$ are not all even or all odd, there exists $t_{y}\in\Gamma_{n}^{\circ}=W_{\Psi_{n}}$ ($y\in\tilde{S}^{0}$,
$y^2=I$) such that \[t_{y}\chi^{\ast}_{a_{1}e_{1}+\dots+a_{n}e_{n}+a_{0}\beta_{0},W_{n}}=-\chi^{\ast}_{a_{1}e_{1}
+\dots+a_{n}e_{n}+a_{0}\beta_{0},W_{n}}.\] Then, \[\frac{1}{|W_{\Psi_{n}}|}\sum_{\gamma\in W_{\Psi_{n}}}\gamma
\chi^{\ast}_{a_{1}e_{1}+\dots+a_{n}e_{n}+a_{0}\beta_{0},W_{n}}=0.\] Therefore, $F_{\Phi,W_{\Psi_{n}}}$ is the
sum of terms $\chi^{\ast}_{a_{1}e_{1}+\dots+a_{n}e_{n}+a_{0}\beta_{0},W_{n}}$ in $F_{\Phi,W_{n}}$ such that
$a_{1},a_{2},\dots,a_{n}$ are all even or all odd.

When $k$ is even, we have $E(F_{\Phi,W_{n}})=\prod_{1\leq i\leq s}c_{n_{i}}\prod_{1\leq j\leq t}b_{n'_{j}}$.
Then, \begin{eqnarray*}&&E(F_{\Phi,W_{\Psi_{n}}})\\&=&\prod_{1\leq i\leq s}c_{n_{i}}^{'+}\prod_{1\leq j\leq t}
b_{n'_{j}}^{'+}+\prod_{1\leq i\leq s}c_{n_{i}}^{'-}\prod_{1\leq j\leq t}b_{n'_{j}}^{'-}\\&=&\prod_{1\leq i\leq s}
c_{n_{i}}^{+}\prod_{1\leq j\leq t}b_{n'_{j}}^{+}+\prod_{1\leq i\leq s}c_{n_{i}}^{-}\prod_{1\leq j\leq t}
b_{n'_{j}}^{-}\end{eqnarray*} due to \[\sum_{1\leq i\leq s}\frac{n_{i}(n_{i}+1)}{2}+
\sum_{1\leq j\leq t}\frac{n'_{j}(n'_{j}+1)}{2}\equiv k\equiv 0\pmod{2}.\]

When $k$ is odd, define $\Phi'$ to be a quasi sub-root system by substituting the simple root $e_{1}-e_{2}+\beta_{0}$
appearing in (ii) and (iii) by $e_{1}-e_{2}$. Then, $F_{\Phi,W_{\Psi_{n}}}$ can be obtained from
$F_{\Phi',W_{\Psi_{n}}}$ by substituting each term $\chi^{\ast}_{a_{1}e_{1}+\dots+a_{n}e_{n}+a_{0}\beta_{0},W_{n}}$
by $\chi^{\ast}_{a_{1}e_{1}+\dots+a_{n}e_{n}+(a_{0}+a_{1})\beta_{0},W_{n}}$. When $a_{1}$ is even, this has no
affection; when $a_{1}$ is even, this effects by multiplying -1 for the image under the homomorphism $E$. We have
\[E(F_{\Phi',W_{\Psi_{n}}})=\prod_{1\leq i\leq s}c_{n_{i}}^{'+}\prod_{1\leq j\leq t}b_{n'_{j}}^{'+}+
\prod_{1\leq i\leq s}c_{n_{i}}^{'-}\prod_{1\leq j\leq t}b_{n'_{j}}^{'-}.\] Then, \begin{eqnarray*}&&
E(F_{\Phi,W_{\Psi_{n}}})\\&=&\prod_{1\leq i\leq s}c_{n_{i}}^{'+}\prod_{1\leq j\leq t}b_{n'_{j}}^{'+}-
\prod_{1\leq i\leq s}c_{n_{i}}^{'-}\prod_{1\leq j\leq t}b_{n'_{j}}^{'-}\\&=&\prod_{1\leq i\leq s}
c_{n_{i}}^{+}\prod_{1\leq j\leq t}b_{n'_{j}}^{+}+\prod_{1\leq i\leq s}c_{n_{i}}^{-}\prod_{1\leq j\leq t}
b_{n'_{j}}^{-}\end{eqnarray*} due to \[\sum_{1\leq i\leq s}\frac{n_{i}(n_{i}+1)}{2}+
\sum_{1\leq j\leq t}\frac{n'_{j}(n'_{j}+1)}{2}\equiv k\equiv 1\pmod{2}.\qedhere\]
\end{proof}

Consider weighted homogeneous polynomials of the variables $\{b_{k},c_{2k}:k\geq 1\}$ where the weight of $b_{k}$
(resp. $c_{2k}$) is assigned to be $k$ (resp. $2k$). Note that each of $b_{k}^{+}$, $b_{k}^{-}$, $c_{k}^{+}$,
$c_{k}^{-}$ is a homogeneous polynomial of degree $k$. The highest index variables appearing in $b_{k}^{+}$,
$b_{k}^{-}$, $c_{2k}^{+}$, $c_{2k}^{-}$ are $x_{2k-2}$, $x_{2k-1}$, $x_{4k}$, $x_{4k-1}$, respectively. Then,
there exist a nonzero integral coefficients weighted homogeneous polynomial $f_{1}$ of $b_1,b_2,b_3,c_2$ and
a nonzero integral coefficients weighted homogeneous polynomial $f_{2}$ of $b_1,b_2,c_2$ such that
\[f_{1}(b_{1}^{+},b_{2}^{+},b_{3}^{+},c_{2}^{+})=0\] and \[f_{2}(b_{1}^{-},b_{2}^{-},c_{2}^{-})=0.\] Put
$f=c_{4k_{1}}f_{1}f_{2}$ for some $k_{1}\geq 1$. Let $n'$ be the degree of $f$. By Lemma \ref{L:poly-Phi},
monomial terms of $f$ give $\Gamma_{n'}^{\circ}$ pairwise non-conjugate quasi sub-root systems
$\Phi_{1},\dots,\Phi_{s}$ of $\Psi_{n'}$ all occurring as quasi root systems of S-subgroups of $G_{n'}$ with
$\tilde{S}_{n'}$ a quasi Cartan subgroup such that \begin{equation}\label{Eq:counter1}\sum_{1\leq i\leq s}
c_{i}F_{\Phi_{i},\Gamma_{n'}^{\circ}}=0\end{equation} for some non-zero integers $c_{1},\dots,c_{s}$. Let
$H_{i}$ ($1\leq i\leq s$) be an S-subgroup of $G_{n'}$ with $\tilde{S}_{n'}$ a quasi Cartan subgroup and with
quasi root system \[R(H_{i},\tilde{S}_{n'})=\Phi_{i}.\] Choose a pinned element $g'_{i}\in\tilde{S}_{n'}$ of
$\Phi_{i}$. Write $g_{i}^{'2}=\epsilon_{i}I$ where $\epsilon_{i}=\pm{1}$.

Let $k_{2},k_{3}$ be two positive integers. Put $n=n'+4k_{2}+4k_{3}$. For each $i$ ($1\leq i\leq s$) and $j=1,2$,
form subgroup $H_{i,j}$ of $G_{n}$ generated by \[H_{i,j}^{0}=H_{i}^{0}\times K_{k_{2},j}\times K_{k_{3},j}\]
and an element $g_{i,j}\in\U(2^{n})\tau$, where \[K_{k_{p},1}=K_{k_{p},2}=\Spin(8k_{p}+2)\] are embedded into
$\U(2^{4k_{p}})$ via the spinor representation $M_{+,8k_{p}+2}$ ($p=2,3$), \[g_{i,j}^2=\epsilon_{i}I
\textrm{ and }\Ad(g_{i,j})|_{H_{i}^{0}}=\Ad(g'_{i})|_{H_{i}^{0}},\] $\Ad(g_{i,1})|_{K_{k_{p},1}}$ is a pinned
automorphism with fixed point subgroup equal to $\Spin(8k_{p}+1)$ ($p=2,3$), $\Ad(g_{i,2})|_{K_{k_{p},2}}=\id$
($p=2,3$).

\begin{theorem}\label{T:counter1}
All $H_{i,j}$ ($1\leq i\leq s$, $j=1,2$) are S-subgroups of $G_{n}$ and we have \[\sum_{1\leq i\leq s}
c_{i}(\mathscr{D}_{H_{i,1}}-\mathscr{D}_{H_{i,2}})=0.\]
\end{theorem}

\begin{proof}
All $H_{i,1}$ ($1\leq i\leq s$) share a common maximal commutative connected subset $S$ containing $g_{i,1}$
and all $H_{i,2}$ ($1\leq i\leq s$) share a common maximal commutative connected subset $S'$ containing
$g_{i,2}$; $H_{i,1}^{0}$ and $H_{i,2}^{0}$ ($1\leq i\leq s$) share a common maximal torus $T_{i}$. Let
$\tilde{S}$ (resp. $\tilde{S}'$) be the quasi torus generated by $S$ (resp. $S'$). Put \[\Gamma=
N_{G_{n}^{0}}(\tilde{S})/Z_{G_{n}^{0}}(\tilde{S});\] \[\Gamma'=N_{G_{n}^{0}}(\tilde{S}')/Z_{G_{n}^{0}}(\tilde{S}');\]
\[\Gamma^{(i)}=N_{G_{n}^{0}}(T_{i})/Z_{G_{n}^{0}}(T_{i}).\] Write \[\Phi_{i,1}=R(H_{i,1},\tilde{S}),
\ 1\leq i\leq s;\] \[\Phi_{i,2}=R(H_{i,2},\tilde{S}'),\ 1\leq i\leq s;\] \[\Phi'_{i,j}=R(H_{i,j},T_{i}),
\ 1\leq i\leq s,j=1,2.\] By \eqref{Eq:counter1}, we get \begin{equation}\label{Eq:counter1-2}\sum_{1\leq i\leq s}
c_{i}F_{\Phi_{i,1},\Gamma}=0\end{equation} and \begin{equation}\label{Eq:counter1-3}\sum_{1\leq i\leq s}
c_{i}F_{\Phi_{i,2},\Gamma'}=0\end{equation} By our construction, $\Phi'_{i,1}=\Phi'_{i,2}$. Then,
\begin{equation}\label{Eq:counter1-4}F_{\Phi'_{i,1},\Gamma^{(i)}}-F_{\Phi'_{i,2},\Gamma^{(i)}}=0.\end{equation}
By Proposition \ref{P:dim-char} and Corollary \ref{C:dim-char2}, the equalities \eqref{Eq:counter1-2},
\eqref{Eq:counter1-3}, \eqref{Eq:counter1-4} together indicate that \[\sum_{1\leq i\leq s}
c_{i}(\mathscr{D}_{H_{i,1}}-\mathscr{D}_{H_{i,2}})=0.\qedhere\]
\end{proof}

\begin{remark}\label{R:counter1}
We take $f=c_{4k_{1}}f_{1}f_{2}$ instead of $f_{1}f_{2}$ to make sure that each $H_{i}^{0}$ is non-self dual so
that $H_{i}$ is an S-subgroup. We multiply two groups $K_{p,j}\cong\Spin(8k_{p}+2)$ ($p=2,3$, $j=1,2$) while
defining $H_{i,j}$ instead of only one so that the Schur index of some subgroup of $H_{i,j}$ ($j=1,2$)
is equal to that of the Schur index of the respective subgroup of $H_{i}^{0}$ .
\end{remark}

\subsection{S-subgroups with linearly dependent dimension data, II}\label{SS:counter2}

Let $N>1$ and $K$ be a connected semisimple subgroup $K$ of $\U(N)$ containing $-I$ and with two outer involutive
automorphisms $\theta_{1}$ and $\theta_{2}$. Write $\rho$ for the restriction to $K$ of the natural representation
of $\U(N)$.

\begin{assumption}\label{A:counter2}
Suppose the following conditions hold: \begin{enumerate}
\item[(i)]$\rho$ is irreducible and non-self dual.
\item[(ii)]$\theta_{1}^{\ast}\rho\cong\theta_{2}^{\ast}\rho\cong\rho^{\ast}$.
\item[(iii)]The images of $\theta_{1}$ and $\theta_{2}$ in $\Out(K):=\Aut(K)/\Int(K)$ are non-conjugate.
\end{enumerate}
\end{assumption}

\begin{example}\label{E:counter2}
Let $k\geq 4$ be an even integer. Write $c$ for a generator of the center $Z(\SU(k))$ and write $\rho'$ for the
natural representation of $\SU(k)$. Take \[K=(\SU(k)\times\SU(k))/\langle(c,c^{-1})\rangle\] and
\[\rho=\rho'\otimes\rho'.\] Write $\sigma$ for conjugation on $\SU(k)$. Define \[\theta_{1}([x,y])=
([\sigma(x),\sigma(y)]),\ \forall x,y\in\SU(k)\] and \[\theta_{2}([x,y])=([\sigma(y),\sigma(x)]),
\ \forall x,y\in\SU(k).\] Put $N=k^{2}$. Then, the tuple $(K,\rho,\theta_{1},\theta_{2})$ satisfies the Assumption
\ref{A:counter2}.
\end{example}

Put $G=\U(N)\rtimes\langle\tau\rangle$, where $\tau^{2}=I$ and \[\tau X\tau^{-1}=\overline{X},\ \forall X\in\U(N).\]

\begin{lemma}\label{L:counter2}
Suppose the tuple $(K,\rho,\theta_{1},\theta_{2})$ satisfies the Assumption \ref{A:counter2}. Then for each $i=1,2$,
there exists a subgroup $H_{i}$ of $G$ generated $K$ and an element $g_{i}\in\U(N)\tau$ such that $g_{i}^{2}=\pm{}I$
and $\Ad(g_{i})|_{K}=\theta_{i}$.
\end{lemma}

\begin{proof}
Define \[\rho_{i}'(x)=\tau^{-1}\rho(\theta_{i}(x))\tau,\ \forall x\in K.\] Then, \[\rho'_{i}\cong
(\theta_{i}^{\ast}\rho)^{\ast}\cong(\rho^{\ast})^{\ast}=\rho\] by Assumption \ref{A:counter2}. Thus, there exist
$y_{i}\in\U(N)$ such that \[\rho_{i}'(x)=y_{i}\rho(x)y_{i}^{-1},\ \forall x\in K.\] Put $g_{i}=\tau y_{i}$. Then,
\begin{equation}\label{Eq:counter2-1}\rho(\theta_{i}(x))=g_{i}\rho(x)g_{i}^{-1},\ \forall x\in K.\end{equation} This
just means $\Ad(g_{i})|_{K}=\theta_{i}$. By iteration we get \[g_{i}^{2}\rho(x)g_{i}^{-2}=g_{i}\rho(\theta_{i}(x))
g_{i}^{-1}=\rho(\theta_{i}^{2}(x))=\rho(x),\ \forall x\in K\] from \eqref{Eq:counter2-1}. Since $K$ is an irreducible
subgroup of $\U(N)$, we have $g_{i}^{2}=\lambda I$ for some $\lambda\in\mathbb{C}$ with $|\lambda|=1$. Since $g_{i}^{2}$
commutes with $g_{i}\in\U(N)\tau$, we have \[\lambda I=g_{i}(\lambda I)g_{i}^{-1}=\bar{\lambda} I.\] Thus, $\lambda=
\pm{1}$. Hence, $g_{i}^{2}=\pm{}I$.
\end{proof}

Put \[H_{3}=\langle H_{1},Z(\U(N))\rangle\] and \[H_{4}=\langle H_{2},Z(\U(N))\rangle.\]

\begin{theorem}\label{T:counter2}
For each $i=1,2,3,4$, $H_{i}$ is an S-subgroup of $G$. The groups $H_{1}$, $H_{2}$, $H_{3}$, $H_{4}$ are pairwise
non-isomorphic. We have \begin{equation}\label{Eq:counter2}(\mathscr{D}_{H_{1}}-\mathscr{D}_{H_{2}})-
(\mathscr{D}_{H_{3}}-\mathscr{D}_{H_{4}})=0.\end{equation}
\end{theorem}

\begin{proof}
By Assumption \ref{A:counter2}(i), $\rho|_{H_{i}^{0}}$ is irreducible and non-self dual for each $i=1,2,3,4$. Then,
$H_{i}$ is an S-subgroup of $G$ by Lemma \ref{L:TA4}. Note that it is assumed that $-I\in K$. By Lemma \ref{L:counter2},
each $H_{i}$ ($1\leq i\leq 4$) has two connected components.

We show that $H_{1}$, $H_{2}$, $H_{3}$, $H_{4}$ are pairwise non-isomorphic. Then, they are pairwise non-conjugate.
Since $H_{1},H_{2}$ are semisimple and $H_{3},H_{4}$ are non-semisimple, it suffices to show that $H_{1}\not\cong
H_{2}$ and $H_{3}\not\cong H_{4}$. By construction we have \[(H_{3}^{0})_{\der}=(H_{4}^{0})_{\der}=H_{1}^{0}=
H_{2}^{0}=K.\] Then, $H_{1}\not\cong H_{2}$ and $H_{3}\not\cong H_{4}$ by Assumption \ref{A:counter2}(iii) and
Lemma \ref{L:counter2}.

Choose a quasi Cartan subgroup $\tilde{S}_{i}$ of $H_{i}$ containing $g_{i}$ ($i=1,2$). Choose a maximal torus
$T$ of $H_{1}^{0}$ and put $T'=\langle T,Z(\U(N))\rangle$. Then, $\tilde{S}_{1}$ is a quasi Cartan subgroup of
$H_{3}$, $\tilde{S}_{2}$ is a quasi Cartan subgroup of $H_{4}$, $T$ is a maximal torus of $H_{2}^{0}$, $T'$ is
a maximal torus of $H_{3}^{0}$ (and $H_{4}^{0}$). Moreover by the construction we have \[R(H_{1},T)=R(H_{2},T),\]
\[R(H_{3},T')=R(H_{4},T'),\] \[R(H_{1},\tilde{S}_{1})=R(H_{3},\tilde{S}_{1}),\] \[R(H_{2},\tilde{S}_{2})=
R(H_{4},\tilde{S}_{2}).\] By Proposition \ref{P:dim-char} and Corollary \ref{C:dim-char2}, \eqref{Eq:counter2}
follows.
\end{proof}

\section{Dimension data of S-subgroups of $\O(2N)$}\label{S:TD}

Put \[G=\O(2N)\] and write $\rho=\mathbb{C}^{2N}$ for the natural representation of it. By Lemma \ref{L:TD1}, 
there are three classes of S-subgroups $H$ of $\O(2N)$. \begin{enumerate} 
\item[(i)]$\rho|_{H^{0}}$ is irreducible.
\item[(ii)]$\rho|_{H^{0}}$ is the direct sum of two odd-dimensional self-dual irreducible representations of orthogonal
type. 
\item[(iii)]$\rho|_{H^{0}}$ is the direct sum of two odd-dimensional non-self dual irreducible representations which
are dual to each other. 
\end{enumerate}
Choose a commutative connected closed subset $S$ of $\O(2N)$ consisting of elements of determinant $-1$ and let
$\tilde{S}$ be the quasi torus generated by $S$. Write $\chi=\chi_{\rho}|_{\tilde{S}}$. Put \[\chi_{o}=
\chi_{\rho}|_{\tilde{S}^{0}}\quad\textrm{ and }\quad\chi_{t}=\chi_{\rho}|_{S},\] which are called the restriction 
character and the twisted character of $\rho$, respectively. Put \[\Gamma^{\circ}=
N_{G^{0}}(\tilde{S})/N_{G^{0}}(\tilde{S}).\] By Lemma \ref{L:S1}, choose and fix an element $x_{0}\in S$ such that 
$x_{0}^{2}=\pm{I}$. Define a character $\beta_{0}\in X^{\ast}(\tilde{S})$ by \[\beta_{0}|_{\tilde{S}^{0}}=
1\textrm{ and }\beta_{0}(x_{0})=-1.\] For any element $y\in\tilde{S}\cap G^{0}$ with $y^2=I$, define $t_{y}\in
\Aut(\tilde{S})$ by \[t_{y}|_{\tilde{S}^{0}}=\id\textrm{ and }t_{y}(x_{0})=yx_{0}.\] Define $\Psi_{\tilde{S}}$ to 
be the set of roots $\alpha\in\Psi'_{\tilde{S}}$ such that there exists no root $\beta\in\Psi'_{\tilde{S}}$ such 
that $\alpha=2\beta$. Write $\Psi=\Psi'_{\tilde{S}}$ for simplicity. Then, $\Psi$ is stable under $\Gamma^{\circ}$ 
and $W_{\Psi}=W_{\Psi'_{\tilde{S}}}$. Let \[\Psi=\bigsqcup_{1\leq i\leq r}\Psi_{i}\] be the decomposition of $\Psi$ 
into an orthogonal union of irreducible quasi root systems. 

\begin{lemma}\label{L:ot}
Assume that $\tilde{S}$ is a quasi Cartan subgroup of an S-subgroup $H$ of $\O(2N)$. If $H$ is in Class (i), then
$\chi_{\rho}|_{S}\neq 0$; if $H$ is in Classes (ii)-(iii), then $\chi_{\rho}|_{S}=0$.
\end{lemma}

\begin{proof}
When $H$ is in Class (i), by the twining character formula \eqref{Eq:chi'2} we have $\chi_{\rho}|_{S}\neq 0$.
When $H$ is in Classes (ii)-(iii), we have $\rho|_{H^{0}}=\rho_{1}\oplus\rho_{2}$ where $\rho_{1},\rho_{2}$
are two irreducible representations of $H^{0}$ and the action of any $x\in S$ maps $\rho_{1}$ (resp. $\rho_{2}$)
into $\rho_{2}$ (resp. $\rho_{1}$). Thus, $\chi_{\rho}|_{S}=0$. 
\end{proof}

\begin{lemma}\label{L:TD-Psi0}
Let $H$ be an S-subgroup of $G$ with $\tilde{S}$ a quasi Cartan subgroup. Then we have \begin{equation*}
\bigcap_{\bar{\alpha}\in p(R(H,\tilde{S}))}\ker\bar{\alpha}=\bigcap_{\bar{\alpha}\in p(\Psi)}\ker\bar{\alpha}
=\tilde{S}^{0}\cap Z(G),\end{equation*} \[R(H,\tilde{S})\subset\Psi\] and \[p(R(H,\tilde{S}))
\supset p(\Psi)^{\circ}.\]
\end{lemma}

\begin{proof}
We have $e=2$. By Lemma \ref{L:S3} we get the conclusion of the lemma.
\end{proof}

The following lemma characterizes $\Gamma^{\circ}$ using the character $\chi$.

\begin{lemma}\label{L:TD-Gamma2}
We have \[\Gamma^{\circ}=\{\gamma\in\Aut(\tilde{S}):\gamma^{\ast}\chi=\chi\}.\]
\end{lemma}

\begin{proof}
It is clear that $\Gamma^{\circ}\subset\{\gamma\in\Aut(\tilde{S}):\gamma^{\ast}\chi=\chi\}$. For any $\gamma\in
\Aut(\tilde{S})$ such that $\gamma^{\ast}\chi=\chi$, by the character theory of representations of $\tilde{S}$
there exists $g\in\U(2N)$ such that $\Ad(g)|_{\tilde{S}}=\gamma$. Since $\O(2N)$ strongly controls its fusion in
$\U(2N)$ (\cite{Griess}), there exists $g'\in\O(2N)$ such that $\Ad(g)|_{\tilde{S}}=\Ad(g')|_{\tilde{S}}$.
Then, $\gamma=\Ad(g')|_{\tilde{S}}\in\Gamma^{\circ}$.
\end{proof}

The following lemma characterizes the conjugacy class of an S-subgroup of $\O(2N)$ using quasi root system while
the semisimplicity (or non-semisimplicity) is given.

\begin{lemma}\label{L:ARS-conjugacy-TD}
Let $H_{1}$ and $H_{2}$ be two S-subgroups of $\O(2N)$ with $\tilde{S}$ a common quasi Cartan subgroup and with
quasi root system $\Phi_{j}=R(H_{i},\tilde{S})$ ($j=1,2$). Assume that $H_{1}^{0}$ and $H_{2}^{0}$ are both
semisimple or both non-semisimple. Then, the following three conditions are equivalent:
\begin{enumerate}
\item[(1)]$H_{1}$ and $H_{2}$ are conjugate.
\item[(2)]$H_{1}$ and $H_{2}$ are element-conjugate.
\item[(3)]$\Phi_{1}$ and $\Phi_{2}$ are conjugate with respect to $\Gamma^{\circ}$.
\end{enumerate}
\end{lemma}

\begin{proof}
$(1)\Leftrightarrow(2)$: by the fact that $\O(2N)$ is strongly acceptable (\cite{Griess}, \cite{Larsen2}).

$(1)\Rightarrow(3)$. Suppose that $H_{1}$ and $H_{2}$ are conjugate. Then, there exists $g_{1}\in G$ such that
$g_{1}H_{1}g_{1}^{-1}=H_{2}$. Thus, $\Ad(g_{1})(\tilde{S})$ are $\tilde{S}$ are both quasi Cartan subgroup of
$H_{2}$. Hence, there exists $h_{2}\in H_{2}$ such that $\Ad(g_{1})(\tilde{S})=\Ad(h_{1})(\tilde{S})$.
Put $g=h_{1}^{-1}g_{1}$. Then, $\Ad(g)(H_{1})=H_{2}$ and $\Ad(g)(\tilde{S})=\tilde{S}$. Put $\gamma=
\Ad(g)|_{\tilde{S}}\in\Gamma^{\circ}$. Then, $\Phi_{1}=\gamma^{\ast}\Phi_{2}$.

$(3)\Rightarrow(1)$. Suppose that $\Phi_{1}$ and $\Phi_{2}$ are $\Gamma^{\circ}$ conjugate. Replacing $H_{2}$ by 
a conjugate one if necessary we may assume that $\Phi_{1}=\Phi_{2}$, and let it be denoted by $\Phi$. Choose a 
pinned element $x\in S$ of $\Phi$. Let $T_{j}=Z_{H_{j}^{0}}(\tilde{S}^{0})$ ($j=1,2$). Then, $T_{j}$ is a maximal 
torus of $(H_{j})^{0}$ by Lemma \ref{L:T}. Define a root $\beta\in R(H_{j},T_{j})$ to be positive if 
$\beta|_{\tilde{S}^{0}}>0$. Let $\lambda_{j}$ be the highest weight of $\rho|_{T_{j}}$ with respect to this positive 
system. By Lemma \ref{L:ARS-isogeny2} there is an isomorphism $\sigma:\mathfrak{h}_{1}\rightarrow\mathfrak{h}_{2}$ 
such that $\sigma|_{\mathfrak{s}}=\id$ and $\sigma\circ\Ad(x)=\Ad(x)\circ\sigma$. Then, the restriction characters 
on $\tilde{S}^{0}$ of two irreducible representations of $\mathfrak{h}_{1}$ with highest weights 
$\sigma^{\ast}\lambda_{2}$ and $\lambda_{1}$ are both equal to $\chi_{\rho}|_{\tilde{S}^{0}}$. From the decomposition 
of $\Phi$ into irreducible factors, we get direct product decompositions $\mathfrak{h}_{j}=\bigoplus_{1\leq i\leq s}
\mathfrak{h}_{j,i}$ of $\mathfrak{h}_{j}$ ($j=1,2$). Then, $\sigma(\mathfrak{h}_{1,i})=\mathfrak{h}_{2,i}$ 
($1\leq i\leq s$). Choose a highest weight $\lambda_{j}$ of $\rho|_{H_{j}^{0}}$ and write $\lambda_{j,i}$ for the 
$i$-th component of $\lambda_{j}$ ($j=1,2$, $1\leq i\leq s$). By Lemma \ref{L:TD1} we have \[\chi_{\Phi}(\lambda_{1})
=\chi_{\Phi}(\lambda_{2}).\] By Lemma \ref{L:Levi4}, for each $i$ ($1\leq i\leq s$) we have $\sigma^{\ast}\lambda_{2,i}
=\lambda_{1,i}$ or $\Ad(x)^{\ast}\lambda_{1,i}$. Replacing the action of $\sigma$ to the $i$-th factor of 
$\mathfrak{h}_{1}$ for some indices $i$ ($1\leq i\leq s$) by $\Ad(x)\circ\sigma$ if necessary, we are led to 
$\sigma^{\ast}\lambda_{2,i}=\lambda_{1,i}$ for all $i$ ($1\leq i\leq s$). Then, $\sigma^{\ast}\lambda_{2}=
\lambda_{1}$. Note that, when $\chi_{\rho}|_{S}\neq 0$, $H_{1}$ and $H_{2}$ are in Class (i). Thus, 
$\mathfrak{h}_{j}$ ($j=1,2$) acts irreducibly on $\rho$. When $\chi_{\rho}|_{S}=0$, $H_{1}$ and $H_{2}$ are in 
Classes (ii)-(iii). Thus, $\rho|_{\mathfrak{h}_{j}}$ ($j=1,2$) is the direct sum of two irreducible representations 
and the other highest weight of $\rho|_{H_{j}^{0}}$ is $\Ad(x)\cdot\lambda_{j}$. Then, 
$\sigma^{\ast}(\rho|_{\mathfrak{h}_{2}})\cong\rho|_{\mathfrak{h}_{1}}$ in any case. Thus, there exists $g\in\U(2N)$
such that $\Ad(g)|_{\mathfrak{h}_{1}}=\sigma$. Hence, $gH_{1}^{0}g^{-1}=H_{2}^{0}$. Moreover,
$\sigma\circ\Ad(x)=\Ad(x)\circ\sigma$ indicates that $[g,x]$ commutes with $H_{1}^{0}$.

When $H_{1}$ is in Class (i), by the irreducibility of $H_{1}^{0}$, we have $gxg^{-1}=(\epsilon I)x$
where $\epsilon=\pm{1}$. Then, $\epsilon\tr(x)=\tr((\epsilon I)x)=\tr(x)$. By the twining
character formula, we have $\tr(x)\neq 0$. Then, $\epsilon=1$. Thus, $gH_{1}g^{-1}=H_{2}$ and
$\Ad(g)|_{\tilde{S}}=\id$.

When $H_{1}$ is in Class (ii), we may assume that $H_{1}$ is contained in $(\SO(N)\times\SO(N))\rtimes
\langle J'_{N}\rangle$. Then, $gxg^{-1}=yx$ where \[y\in Z_{\SO(2N)}(H_{1}^{0})=\langle I_{N,N},
-I_{N}\rangle.\] Since $\det(y)=\det(yx)\det(x)^{-1}=\det(gxg^{-1})\det(x)^{-1}=1$, we must have $y=
\pm{I}$. Put \[g'=\left\{\begin{array}{cc}gI_{N,N}\textrm{ if }y=-I;\\g\quad\textrm{ if }y=I.\end{array}
\right.\] Then, $g'xg'^{-1}=x$. Thus, $g'H_{1}g'^{-1}=H_{2}$ and $\Ad(g')|_{\tilde{S}}=\id$.

When $H_{1}$ is in Class (iii), we may assume that $H_{1}$ is contained in $\U(N)\rtimes\langle\tau\rangle$
where \[\U(N)=\{\left(\begin{array}{cc}A&B\\-B&A\end{array}\right)\in\SO(2N):A,B\in M_{N}(\mathbb{R})\}\]
and $\tau=I_{N,N}$. The proof reduces to the case for $\U(N)\rtimes\langle\tau\rangle$.
\end{proof}

\subsection{Quasi root systems of S-subgroups of $\O(2N)$ in Class (i)}\label{SS:TD1}

\begin{lemma}\label{L:TDi-1}
Let $\Phi=R(H,\tilde{S})$ be the quasi root system of an S-subgroup $H$ in Class (i) with $\tilde{S}$ a quasi
Cartan subgroup. \begin{enumerate} 
\item[(1)]If $\alpha\in\Phi$ is a root such that $p(\alpha)\in p(\Phi)^{\nd}$, then
there exist two weights $\mu$ and $\mu'$ appearing in the twisted character $\chi_{t}$ such that
\[\mu'-\mu=m_{\alpha}\alpha.\]
\item[(2)]For any two weights $\mu$ and $\mu'$ appearing in the twisted character $\chi_{t}$, we have
\[\mu'-\mu\in\span_{\mathbb{Z}}\{m_{\alpha}\alpha: \alpha\in\Phi, p(\alpha)\in p(\Phi)^{\nd}\}.\]  
\end{enumerate} 
\end{lemma}

\begin{proof}
We first show (1). By assumption, $\alpha$ appears in a simple system $\{\alpha_{j}:1\leq j\leq l\}$ of $\Phi$.
Let $[\lambda]$ ($\lambda\in X^{\ast}(\tilde{S})$) be a leading term of $\chi_{t}$. For each weight
$\mu\in X^{\ast}(\tilde{S})$ appearing in $\chi_{t}$, by Lemma \ref{L:twining3}, we have 
\[\lambda-\mu=\sum_{1\leq j\leq l}b_{j}m_{\alpha_{j}}\alpha_{j}\] where $b_{j}\in\mathbb{Z}_{\geq 0}$. Define 
\[\supp(\lambda-\mu)=\{\alpha_{j}:b_{j}\neq 0\}\] and \[\depth(\mu)=\sum_{1\leq j\leq l}b_{j}.\] Among those 
weights $\mu$ appearing in $\chi_{t}$ and such that $\alpha\in\supp(\lambda-\mu)$, choose a weight $\mu_{0}$ 
such that $\depth(\mu_{0})$ is minimal. If $(\mu_{0},\alpha_{j})<0$ for a simple root $\alpha_{j}\neq\alpha$, 
then $s_{\alpha_{j}}(\mu_0)=\mu_{0}-\frac{2(\mu_{0},\alpha_{j})}{(\alpha_{j},\alpha_{j})}\alpha_{j}$ appears 
in $\chi_{t}$ and such that $\alpha\in\supp(\lambda-\mu)$. We have \[\depth(s_{\alpha_{j}}(\mu_0))-\depth(\mu_{0})=
\frac{2(\mu_{0},\alpha_{j})}{(\alpha_{j},\alpha_{j})}<0,\] which is in contradiction with the minimality of
$\depth(\mu_{0})$. Then, $\alpha$ is the only simple root such that $(\mu_{0},\alpha)<0$. Thus, $\mu':=\mu+
m_{\alpha}\alpha$ appears in $\chi_{t}$ as well. It is clear that $\mu'-\mu=m_{\alpha}\alpha$.  

(2)follows from the twining character formula directly.
\end{proof}

\begin{lemma}\label{L:TDi-2}
Let $\Phi=R(H,\tilde{S})$ be the quasi root system of an S-subgroup $H$ in Class (i) with $\tilde{S}$ a quasi
Cartan subgroup. Suppose that the restricted root system $p(\Phi)$ is given. Then, the quasi root system $\Phi$
is uniquely determined. 
\end{lemma}

\begin{proof}
Let $\alpha\in\Phi$ be a root in $\Phi$ such that $p(\alpha)=\alpha|_{\tilde{S}^{0}}\in p(\Phi)^{\nd}$. By Lemma
\ref{L:TDi-1}, the linear character $m_{\alpha}\alpha$ is determined by the twisted character $\chi_{t}$. When 
$2p(\alpha)\not\in p(\Phi)$, we have $m'_{\alpha}=m_{\alpha}$; when $2p(\alpha)\not\in p(\Phi)$, we have 
$m'_{\alpha}=2$. Then, both the folding index $m'_{\alpha}$ and the fractional factor of $\alpha$ are determined. 
Thus, the quasi root system $\Phi$ is determined. 
\end{proof}

\begin{lemma}\label{L:TDi-3}
Let $\Phi=R(H,\tilde{S})$ be the quasi root system of an S-subgroup $H$ in Class (i). Then
$\span_{\mathbb{Z}}\{m_{\alpha}\alpha|_{\tilde{S}^{0}}:\alpha\in\Phi\}$ is a fixed sub-lattice of
$X^{\ast}(\tilde{S}^{0})$ determined by the twisted character $\chi_{t}$, and it is $\Gamma^{\circ}$ invariant.
\end{lemma}

\begin{proof}
By Lemma \ref{L:TDi-1}, we see that  \[\span_{\mathbb{Z}}\{m_{\alpha}\alpha_{\tilde{S}^{0}}:\alpha\in\Phi\}=
\span_{\mathbb{Z}}\{\mu'_{\tilde{S}^{0}}-\mu_{\tilde{S}^{0}}:\mu,\mu'\textrm{ appears in }\chi_{t}\}.\] Since
$\chi_{t}$ is a fixed $\Gamma^{\circ}$ invariant character, $\span_{\mathbb{Z}}\{\mu'_{\tilde{S}^{0}}-
\mu_{\tilde{S}^{0}}:\mu,\mu'\textrm{ appears in }\chi_{t}\}$ is determined by $\chi_{t}$ and is a
$\Gamma^{\circ}$ invariant sub-lattice of $X^{\ast}(\tilde{S}^{0})$. Then, so is
$\span_{\mathbb{Z}}\{m_{\alpha}\alpha_{\tilde{S}^{0}}\}$.
\end{proof}

\begin{lemma}\label{L:TDi-4}
Let $\Phi=R(H,\tilde{S})$ be the quasi root system of an S-subgroup $H$ in Class (i). Let $\alpha\in\Phi$ and
$\beta\in\Psi$ be two roots such that $p(\alpha)=\alpha_{\tilde{S}^{0}}\in p(\Phi)^{\nd}$ and $p(\beta)=
\beta_{\tilde{S}^{0}}\in p(\Phi)^{\nd}$. If $m_{\beta}=2k$ ($k=1$ or 2) in $\Psi$ and $m_{\alpha}=1$ in $\Phi$,
then \[\frac{(\alpha,\beta)}{k(\beta,\beta)}\in \mathbb{Z}.\]
\end{lemma}

\begin{proof}
If $m_{\beta}=2$ in $\Psi$, then $\beta\in\Psi$ and $\beta+\beta_{0}\in\Psi$. Thus, $s_{\beta}\chi_{t}=\chi_{t}$
and $s_{\beta+\beta_{0}}\chi_{t}=\chi_{t}$. Hence, \[s_{\beta+\beta_{0}}s_{\beta}\chi_{t}=\chi_{t}.\] By Lemma
\ref{L:TDi-1} and the assumption of $m_{\alpha}=1$, we get \[s_{\beta+\beta_{0}}s_{\beta}\alpha=\alpha,\]
which is equivalent to \[\frac{(\alpha,\beta)}{(\beta,\beta)}\in\mathbb{Z}.\]

If $m_{\beta}=4$ in $\Psi$, then $\beta\in\Psi$ and $2\beta+\beta_{0}\in\Psi$. Thus, $s_{\beta}\chi_{t}=\chi_{t}$
and $s_{2\beta+\beta_{0}}\chi_{t}=\chi_{t}$. Hence, \[s_{2\beta+\beta_{0}}s_{\beta}\chi_{t}=\chi_{t}.\] By Lemma
\ref{L:TDi-1} and the assumption of $m_{\alpha}=1$, we get \[s_{2\beta+\beta_{0}}s_{\beta}\alpha=\alpha,\]
which is equivalent to \[\frac{(\alpha,\beta)}{2(\beta,\beta)}\in\mathbb{Z}.\]
\end{proof}

\begin{lemma}\label{L:TD-char1}
Let $H_{1}$, $H_{2}$, $\dots$, $H_{s}$ be a set of pairwise non-conjugate S-subgroups of $G=\O(2N)$ in Class (i)
with $\tilde{S}$ a common quasi Cartan subgroup. Then the characters \[F_{\Phi_{1},\Gamma^{\circ}},\dots,
F_{\Phi_{s},\Gamma^{\circ}}\] are linearly independent, where $\Phi_{j}=R(H_{j},\tilde{S})\subset\Psi$ is the
quasi root system of $H_{j}$ with respect to $\tilde{S}$ ($1\leq j\leq s$). 
\end{lemma}

\begin{proof}
Write \[L_{\Psi}=\span_{\mathbb{Z}}p(\Psi),\] which is a sub-lattice of $X^{\ast}(\tilde{S}^{0})$. Since the action of
$\Gamma^{\circ}$ on $X^{\ast}(\tilde{S})$ preserves $\Psi$, its induced action on $X^{\ast}(\tilde{S}^{0})$ preserves
$p(\Psi)$. Then, $L_{\Psi}$ is preserved by $\Gamma^{\circ}$. For each $j$ ($1\leq j\leq s$), write $F'_{\Phi_{j},
\Gamma^{\circ}}$ for the sum of term $[\lambda]$ in $F_{\Phi_{j},\Gamma^{\circ}}$ such that \[\lambda|_{\tilde{S}^{0}}
\in 2L_{\Psi}.\] Note that, each $t_{y}$ ($y\in\tilde{S}\cap G^{0}$, $y^{2}=I$) fixes all linear characters in 
$2L_{\Psi}$. Then, $F'_{\Phi_{j},\Gamma^{\circ}}$ is preserved by $\Gamma^{\circ}$ and is fixed by all $t_{y}$ 
($y\in\tilde{S}\cap G^{0}$, $y^{2}=I$). It suffices to show that: the characters \[F'_{\Phi_{1},\Gamma^{\circ}},\dots,
F'_{\Phi_{s},\Gamma^{\circ}}\] are linearly independent. By Lemma \ref{L:ARS-conjugacy-TD}, the quasi root systems
$\Phi_{1},\dots,\Phi_{s}$ are pairwise $\Gamma^{\circ}$ non-conjugate as the S-subgroups $H_{1},H_{2},\dots,H_{s}$
are pairwise non-conjugate and are all semisimple. We show that: the leading weights of
$\{F'_{\Phi_{j},\Gamma^{\circ}}:1\leq j\leq s\}$ are distinct, which implies the above assertion.

Without loss of generality we suppose that $F'_{\Phi_{1},\Gamma^{\circ}}$ and $F'_{\Phi_{2},\Gamma^{\circ}}$ have equal 
leading weights. Replacing $\Phi_{2}$ by a $\Gamma^{\circ}$ conjugate one of it if necessary, we may assume that
$A_{\Phi_{1}}$ and $A_{\Phi_{2}}$ have equal actual leading weights (cf. Lemma \ref{L:muPhi}), i.e.,
$\bar{\mu}_{\Phi_{1}}=\bar{\mu}_{\Phi_{2}}$. For $j=1,2$ and each $i$ ($1\leq i\leq r$), write $\Phi_{j,i}=
\Phi_{j}\cap\Psi_{i}$. We have\[\bar{\mu}_{\Phi_{j}}=(\bar{\mu}_{\Phi_{j,1}},\dots,\bar{\mu}_{\Phi_{j,r}}).\]
Then, $\bar{\mu}_{\Phi_{1,i}}=\bar{\mu}_{\Phi_{2,i}}$ for each $i$ ($1\leq i\leq r$). By Lemmas \ref{L:TD-Psi0}
and \ref{L:TDi-1}-\ref{L:TDi-4}, while $\Psi_{i}$ is fixed, the possible $\Phi_{j,i}$ ($j=1,2$) fall into the
following list: \begin{enumerate}
\item[(i-1)]When $p(\Psi_{i})$ is simply-laced and $\Psi_{i}\not\cong\A_{1}^{(2)}$, we have $\Phi_{j,i}=\Psi$.
\item[(i-2)]When $\Psi_{i}\cong\A_{1}^{(2)}$, we have $\Phi_{j,i}\cong\A_{1}$ or $\A_{1}^{(2)}$.
\item[(ii-1)]When $\Psi_{i}=\B_{n}$, $\Phi_{j,i}\cong\bigsqcup_{1\leq i\leq t}\B_{n_{i}}$ where
$\sum_{1\leq i\leq t}n_{i}=n$.
\item[(ii-2)]When $\Psi_{i}=\B_{n}^{(2,1)}$ ($n\geq 2$), $\Phi_{j,i}\cong\bigsqcup_{1\leq i\leq t}\B_{n_{i}}$
(where $\sum_{1\leq i\leq t}n_{i}=n$) or $\B_{n}^{(2,1)}$. If $\B_{n}^{(2,1)}$ occurs, then any
$\bigsqcup_{1\leq i\leq t}\B_{n_{i}}$ can not occur by Lemma \ref{L:TDi-3}.
\item[(ii-3)]When $\Psi_{i}=\B_{n}^{(2)}$ ($n\geq 3$), $\Phi_{j,i}\cong\bigsqcup_{1\leq i\leq t}\B_{n_{i}}^{(2)}$
where $\sum_{1\leq i\leq t}n_{i}=n$.
\item[(ii-4)]When $\Psi_{i}=\B_{2}^{(2)}$, $\Phi_{j,i}\cong\B_{2}^{(2)},\B_{1}^{(2)}\bigsqcup\B_{1}^{(2)}$ or
$\B_{2}^{(2,1)}$. If $\B_{2}^{(2,1)}$ occurs, then $\B_{2}^{(2)}$ and $\B_{1}^{(2)}\bigsqcup\B_{1}^{2}$ can not
occur by Lemma \ref{L:TDi-3}.
\item[(iii-1)]When $\Psi_{i}=\BC_{n}^{(2,2,1)}$ ($n\geq 2$), $\Phi_{j,i}\cong\bigsqcup_{1\leq i\leq t}
\B_{n_{i}}^{(2)}$ (where $\sum_{1\leq i\leq t}n_{i}=n$) or $\BC_{n}^{(2,2,1)}$. If $\BC_{n}^{(2,2,1)}$ occurs,
then any $\bigsqcup_{1\leq i\leq t}\B_{n_{i}}^{(2)}$ can not occur by Lemma \ref{L:TDi-3}.
\item[(iii-2)]When $\Psi_{i}=\BC_{1}^{(2,1)}$, $\Phi_{j,i}\cong\B_{1}^{(2)}$ or $\BC_{1}^{(2,1)}$, but only one
can occur by Lemma \ref{L:TDi-3}.
\item[(iv-1)]When $\Psi_{i}=\C_{n}$ ($n\geq 3$), $\Phi_{j,i}\cong\D_{n}$ or $\C_{n}$.
\item[(iv-2)]When $\Psi_{i}=\C_{n}^{(2,1)}$ ($n\geq 3$), $\Phi_{j,i}\cong\D_{n}^{(2)}$ or $\C_{n}^{(2,1)}$, but
only one can occur by Lemma \ref{L:TDi-3}.
\item[(iv-3)]When $\Psi_{i}=\C_{n}^{(2)}$ ($n\geq 3$), $\Phi_{j,i}\cong\D_{n}^{(2)}$, $\C_{n}^{(2,1)}$ or
$\C_{n}^{(2)}$. If $\C_{n}^{(2,1)}$ occurs, then $\D_{n}^{(2)}$ and $\C_{n}^{(2)}$ can not occur by Lemma
\ref{L:TDi-3}.
\item[(v-1)]When $\Psi_{i}=\F_{4}$, $\Phi_{j,i}\cong\D_{4}^{S}$, $\C_{4}$ or $\F_{4}$.
\item[(v-2)]when $\Psi_{i}=\F_{4}^{(2,1)}$, $\Phi_{j,i}\cong\D_{4}^{S,(2)}$, $\C_{4}^{(2,1)}$ or
$\F_{4}^{(2,1)}$, but only one can occur by Lemma \ref{L:TDi-3}.
\item[(v-3)]When $\Psi_{i}=\F_{4}^{(2)}$, $\Phi_{j,i}\cong\D_{4}^{S,(2)}$, $\C_{4}^{(2)}$ or $\F_{4}^{(2)}$.
\item[(vi-1)]when $\Psi_{i}=\G_{2}$, $\Phi_{j,i}\cong\A_{2}^{S}$ or $\G_{2}$.
\item[(vi-2)]when $\Psi_{i}=\G_{2}^{(2)}$, $\Phi_{j,i}\cong\A_{2}^{S,(2)}$ or $\G_{2}^{(2)}$.
\end{enumerate}

By the formulas of $\bar{\mu}_{\Phi}$ given in the proof of Lemma \ref{L:muPhi}, we see that in each case
(i-1)-(vi-2) the condition $\bar{\mu}_{\Phi_{1,i}}=\bar{\mu}_{\Phi_{2,i}}$ forces that: $\Phi_{1,i}\cong
\Phi_{2,i}$. By Lemma \ref{L:TDi-2}, it follows that $\Phi_{1}$ and $\Phi_{2}$ are $\Gamma^{\circ}$ conjugate,
which gives a contradiction.
\end{proof}

\subsection{Quasi root systems of S-subgroups of $\O(2N)$ in Classes (ii)-(iii)}\label{SS:TD2}

\begin{lemma}\label{L:TD-Gamma1}
Suppose that $\tilde{S}$ is a quasi Cartan subgroup of an S-subgroup $H$ of $G$ in Classes (ii)-(iii). Then for each
element $y\in\tilde{S}\cap G^{0}$ with $y^{2}=I$, we have $t_{y}\in\Gamma^{\circ}$.
\end{lemma}

\begin{proof}
Since it is assumed that $H$ is in Classes (ii)-(iii), we have $\chi_{t}=0$. Then, $t_{y}^{\ast}\chi=\chi$. Thus,
$t_{y}\in\Gamma^{\circ}$ by Lemma \ref{L:TD-Gamma2}.
\end{proof}

\begin{lemma}\label{L:TD-Psi1}
Suppose that $\tilde{S}$ is a quasi Cartan subgroup of an S-subgroup $H$ of $G$ in Classes (ii)-(iii). Let
$\alpha\in\Psi$ be a root. Except when $\alpha|_{\tilde{S}^{0}}=2\beta|_{\tilde{S}^{0}}$ for another root
$\beta\in\Psi$, we have $\alpha+\beta_{0}\in\Psi$.
\end{lemma}

\begin{proof}
Since it is assumed that $H$ is in Classes (ii)-(iii), then $N$ is odd and $-I\not\in\tilde{S}^{0}$. By Lemma
\ref{L:TD-Psi0}, $\{\alpha|_{\tilde{S}^{0}}:\alpha\in\Psi\}$ generate the character group $X^{\ast}(\tilde{S}^{0})$.
When $\alpha|_{\tilde{S}^{0}}\neq 2\beta|_{\tilde{S}^{0}}$ for any other root $\beta\in\Psi$, $\alpha|_{\tilde{S}^{0}}$
lies in a simple system of $p(\Psi)$. Then, there exists $y\in\tilde{S}^{0}$ with $y^2=I$ such that $\alpha(y)=-1$.
Thus, \[\alpha+\beta_{0}=t_{y}\cdot\alpha\in\Psi\] by Lemma \ref{L:TD-Gamma1}.
\end{proof}

The following Lemma \ref{L:TD-Psi2} is a direct consequence of Lemma \ref{L:TD-Psi1}.

\begin{lemma}\label{L:TD-Psi2}
Suppose that $\tilde{S}$ is a quasi Cartan subgroup of an S-subgroup $H$ of $G$ in Classes (ii)-(iii).
Let $\Psi_{i}$ be an irreducible factor of $\Psi_{\tilde{S}}$. \begin{enumerate}
\item[(i)]when $p(\Psi_{i})$ is reduced, we have $\Psi_{i}\cong p(\Psi_{i})^{(2)}$.
\item[(ii)]when $p(\Psi_{i})\cong\BC_{m}$ ($m\geq 2$), we have $\Psi_{i}\cong\BC_{m}^{(2,2,1)}$.
\item[(iii)]when $p(\Psi_{i})\cong\BC_{1}$, we have $\Psi_{i}\cong\BC_{1}^{(2,1)}$.
\end{enumerate}
\end{lemma}

\begin{lemma}\label{L:TD-fracFactor}
Let $\Phi=R(H,\tilde{S})\subset\Psi$ be the quasi root system of an S-subgroup $H$ in Classes (ii)-(iii) with
$\tilde{S}$ a quasi Cartan subgroup. Suppose that the isomorphism type of each quasi root system
$\Phi\cap\Psi_{i}$ ($1\leq i\leq r$) is given. Then the $\Gamma^{\circ}$ orbit containing $\Phi$ is uniquely
determined.
\end{lemma}

\begin{proof}
Since $H$ is in Classes (ii)-(iii), then $N$ is odd and $-I\not\in\tilde{S}^{0}$. Choose a simple system
$\{\alpha_{i}:1\leq i\leq l\}$ of $\Phi$. Then, $\{\alpha_{i}|_{\tilde{S}^{0}}:1\leq i\leq l\}$ is a basis of
$X^{\ast}(\tilde{S}^{0})$ by Lemma \ref{L:TD-Psi0}. For any tuple $(t_{1},\dots,t_{l})\in\{0,1\}^{l}$, there
exists $y\in\tilde{S}^{0}$ with $y^2=I$ and $\alpha_{i}(y)=(-1)^{t_{i}}$ ($1\leq i\leq l$). Then,
\[t_{y}\cdot\alpha_{i}=\alpha_{i}+t_{i}\beta_{0},\ \forall i, 1\leq i\leq l.\] By Lemma \ref{L:TD-Gamma1},
we have $t_{y}\in\Gamma^{\circ}$. Then, the $\Gamma^{\circ}$ orbit containing $\Phi$ is uniquely determined.
\end{proof}

The following lemma follows from Lemma \ref{L:fully4} directly.

\begin{lemma}\label{L:TD-fully5}
Let $\Phi$ be an irreducible quasi sub-root system of $\Psi=\B_{m}^{(2)}$ or $\BC_{m}^{(2,2,1)}$ such that
$p(\Phi)^{\circ}\supset p(\Psi)^{\circ}$ and $\Phi$ be the quasi root system obtained from the action
of an involutive (or trivial) diagram automorphism $\theta$ on a root system $\Phi'$. Let $\lambda$ be an
integral dominant weight of $\Phi'$ such that the irreducible restriction character $\chi_{\Phi}(\lambda)=
\chi_{0}^{\otimes m}$ for a degree 1 character $\chi_{0}$ with $\dim\chi_{0}$ odd. Then the triple
$(\Phi,\lambda,\chi_{0})$ falls into the following list:
\begin{enumerate}
\item[(i)] $\Phi=\BC_{1}^{(2,1)}$, $\lambda=(a,0,-b)$ and $\chi_0=\tau'_{a,b}$, where $a,b\in\mathbb{Z}_{\geq 0}$
and $a+b>0$, at least one of $a$ and $b$ is even and neither of $a+1$, $b+1$, $a+b+2$ is a multiple of 4.
\item[(ii)] $\Phi=\B_{1}^{(2)}$, $\lambda=(\lambda_{1},\lambda_{2})$ with $\lambda_{1}=a$, $\lambda_{2}=b$, and
$\chi_0=\tau_{a,b}=\tau_{a}\tau_{b}$, where $a,b\in 2\mathbb{Z}_{\geq 0}$ and $a+b>0$.
\item[(iii)] $\Phi=\B_{1}$, $\lambda=a$ and $\chi_0=\tau_{a}$, where $a\in 2\mathbb{Z}_{>0}$.
\end{enumerate}
\end{lemma}

\begin{lemma}\label{L:TD-char2}
Let $H_{1}$, $H_{2}$, $\dots$, $H_{s}$ be a set of pairwise non-conjugate S-subgroups of $G=\O(2N)$ in Classes
(ii)-(iii) with $\tilde{S}$ a common quasi Cartan subgroup, and they are all semisimple or all non-semisimple.
Then the characters \[F_{\Phi_{1},\Gamma^{\circ}},\dots,F_{\Phi_{s},\Gamma^{\circ}}\] are linearly independent,
where $\Phi_{j}=R(H_{j},\tilde{S})\subset\Psi$ is the quasi root system of $H_{j}$ with respect to $\tilde{S}$
($1\leq j\leq s$).
\end{lemma}

\begin{proof}
Write \[L_{\Psi}=\span_{\mathbb{Z}}p(\Psi),\] which is a sub-lattice of $X^{\ast}(\tilde{S}^{0})$. Since the action of
$\Gamma^{\circ}$ on $X^{\ast}(\tilde{S})$ preserves $\Psi$, its induced action on $X^{\ast}(\tilde{S}^{0})$ preserves
$p(\Psi)$. Then, $L_{\Psi}$ is preserved by $\Gamma^{\circ}$. For each $j$ ($1\leq j\leq s$), write $F'_{\Phi_{j},
\Gamma^{\circ}}$ for the sum of term $[\lambda]$ in $F_{\Phi_{j},\Gamma^{\circ}}$ such that \[\lambda|_{\tilde{S}^{0}}
\in 2L_{\Psi}.\] Note that, each $t_{y}$ ($y\in\tilde{S}\cap G^{0}$, $y^{2}=I$) fixes all linear characters in
$2L_{\Psi}$. Then, $F'_{\Phi_{j},\Gamma^{\circ}}$ is preserved by $\Gamma^{\circ}$ and is fixed by all $t_{y}$
($y\in\tilde{S}\cap G^{0}$, $y^{2}=I$). It suffices to show that: the characters
\[F'_{\Phi_{1},\Gamma^{\circ}},\dots,F'_{\Phi_{s},\Gamma^{\circ}}\] are linearly independent. By Lemma
\ref{L:ARS-conjugacy-TD}, the quasi root systems $\Phi_{1},\dots,\Phi_{s}$ are pairwise $\Gamma^{\circ}$ non-conjugate
as the S-subgroups $H_{1},H_{2},\dots,H_{s}$ are pairwise non-conjugate and are all semisimple or all non-semisimple.
We show that: the leading weights of $\{F'_{\Phi_{j},\Gamma^{\circ}}:1\leq j\leq s\}$ are distinct, which implies the
above assertion.

Without loss of generality we suppose that $F'_{\Phi_{1},\Gamma^{\circ}}$ and $F'_{\Phi_{2},\Gamma^{\circ}}$ have equal
leading leading weights. Replacing $\Phi_{2}$ by a $\Gamma^{\circ}$ conjugate one of it if necessary, we may assume 
that $A_{\Phi_{1}}$ and $A_{\Phi_{2}}$ have equal actual leading weights, i.e., $\bar{\mu}_{\Phi_{1}}=
\bar{\mu}_{\Phi_{2}}$. For $j=1,2$, write $\Phi_{j,i}=\Phi_{j}\cap\Psi_{i}$ for each $i$ ($1\leq i\leq r$). We have
\[\bar{\mu}_{\Phi_{j}}=(\bar{\mu}_{\Phi_{j,1}},\dots,\bar{\mu}_{\Phi_{j,r}}).\] Then, $\bar{\mu}_{\Phi_{1,i}}=
\bar{\mu}_{\Phi_{2,i}}$ for each $i$ ($1\leq i\leq r$). By Lemmas \ref{L:TD-Psi0}, \ref{L:TD-Psi2} and
\ref{L:TD-fully5}, the triple $(\Psi_{i},\chi_{i},\Phi_{j,i})$ ($j=1,2$) fall into the following list:
\begin{enumerate}
\item[(i)]When $p(\Psi_{i})$ is simply-laced, we have $\Phi_{j,i}\cong p(\Psi_{i})^{(2)}$ or $p(\Psi_{i})$.
\item[(ii)]When $p(\Psi_{i})\cong\B_{n}$ or $\BC_{n}$ and $\chi_{i}$ is in-decomposable ($n\geq 2$), we have
$\Phi_{j,i}\cong\B_{n}$, $\B_{n}^{(2,1)}$, $\B_{n}^{(2)}$, $\BC_{n}^{(2,2,1)}$.
\item[(iii)]When $p(\Psi_{i})\cong\B_{n}$ or $\BC_{n}$ and $\chi_{i}$ is fully decomposable, $\Phi_{j,i}$ is a
product of $\BC_{1}^{(2,1)}$, $\B_{1}^{(2)}$, $\B_{1}$.
\item[(iv)]When $p(\Psi_{i})\cong\C_{n}$ ($n\geq 3$), we have $\Phi_{j,i}\cong\D_{n}$, $\C_{n}$, $\D_{n}^{(2)}$,
$\C_{n}^{(2,1)}$ or $\C_{n}^{(2)}$.
\item[(v)]When $p(\Psi_{i})\cong\F_{4}$, we have $\Phi_{j,i}\cong\D_{4}$, $\C_{4}$, $\F_{4}$, $\D_{4}^{(2)}$,
$\C_{4}^{(2,1)}$, $\F_{4}^{(2,1)}$, $\C_{4}^{(2)}$ or $\F_{4}^{(2)}$.
\item[(vi)]When $p(\Psi_{i})\cong\G_{2}$, we have $\Phi_{j,i}\cong\A_{2}^{S}$, $\G_{2}$, $\A_{2}^{S(2)}$ or
$\G_{2}^{(2)}$.
\end{enumerate}
By the formulas of $\bar{\mu}_{\Phi}$ given in the proof of Lemma \ref{L:muPhi}, we see that in each case (i)-(vi)
the condition $\bar{\mu}_{\Phi_{1,i}}=\bar{\mu}_{\Phi_{2,i}}$ forces that: $\Phi_{1,i}\cong\Phi_{2,i}$. By Lemma
\ref{L:TD-fracFactor}, it follows that $\Phi_{1}$ and $\Phi_{2}$ are $\Gamma^{\circ}$ conjugate, which gives a
contradiction.
\end{proof}

\subsection{Distinction and linear independence of dimension data}\label{SS:TD-linear}

\begin{proof}[Proof of Theorem \ref{T5}.]
Let $\{H_{1},\dots,H_{s}\}$ be a tuple of pairwise non-element conjugate S-subgroups of $G$. Suppose that the
dimension data \[\mathscr{D}_{H_{1}},\dots,\mathscr{D}_{H_{s}}\] are linearly dependent. For each $i$
($1\leq i\leq s$), choose a maximal commutative connected subset $S_{i}$ of $H_{i}$ generating $H_{i}/H_{i}^{0}$.
Write $d_{i}$ for the minimal order of elements in $S_{i}$. Without loss of generality we assume that
$\dim S_{i}\geq\dim S_{1}$ for each $i$ and $d_{i}\leq d_{1}$ whenever $\dim S_{i}=\dim S_{1}$;
$(\dim S_{i},d_{i})=(\dim S_{1},d_{1})$ holds exactly when $1\leq i\leq s'$ for some integer $s'\leq s$; for
some integer $t\leq s'$, $S_{i}$ is conjugate to $S_{1}$ when $1\leq i\leq t$ and any maximal commutative
connected subset $S'_{i}$ of $H_{i}$ generating $H_{i}/H_{i}^{0}$ is not conjugate to $S_{1}$ whenever
$t+1\leq i\leq s'$. Then, by Lemma \ref{L:MCC1}, any maximal commutative connected subset $S'_{i}$ of $H_{i}$
is not conjugate to $S_{1}$ whenever $t+1\leq i\leq s$. Write $S=S_{1}$ and let $\tilde{S}$ be the quasi torus
generated by $S$. Put \[\Gamma^{\circ}=N_{G^{0}}(\tilde{S})/Z_{G^{0}}(\tilde{S}).\] Then, the
characters \[F_{\Phi_{1},\Gamma^{\circ}},\dots,F_{\Phi_{t},\Gamma^{\circ}}\] are linearly dependent by Proposition
\ref{P:dim-char} and Corollary \ref{C:dim-char2}. When $\chi_{\rho}|_{S}\neq 0$, $H_{1},\dots,H_{t}$ are all
in Class (i). Then, the characters $F_{\Phi_{1},\Gamma^{\circ}},\dots,F_{\Phi_{t},\Gamma^{\circ}}$ are linearly
independent by Lemma \ref{L:TD-char1}. When $\chi_{\rho}|_{S}=0$, $H_{1},\dots,H_{s}$ are all in Classes
(ii)-(iii). Then, the characters $F_{\Phi_{1},\Gamma^{\circ}},\dots,F_{\Phi_{t},\Gamma^{\circ}}$ are linearly
independent by Lemma \ref{L:TD-char2}. Thus, there is a contradiction in any case. 
\end{proof}

\begin{proof}[Proof of Theorem \ref{T6}.]
Suppose that $\mathscr{D}_{H_{1}}=\mathscr{D}_{H_{2}}$. Choose a maximal commutative subset $\tilde{S}_{j}$ of 
$H_{j}$ generating $H_{j}/H_{j}^{0}$ ($j=1,2$). Without loss of generality we assume that $S_{1}=S_{2}$ and let 
it be denoted $S$. Let $\tilde{S}$ be the quasi torus generated by $S$. Put \[\Gamma^{\circ}=
N_{G^{0}}(\tilde{S})/Z_{G^{0}}(\tilde{S})\] and write $\Phi_{j}=R(H_{j},\tilde{S})$ ($j=1,2$) for the quasi root 
system of $H_{j}$ with respect to $\tilde{S}$.

If $H_{1}$ and $H_{2}$ are both semisimple or non-semisimple, we get a contradiction by Theorem \ref{T5}. Otherwise, 
without loss of generality we assume that $H_{1}$ is non-semisimple and $H_{2}$ is semisimple. Then, $H_{1}$ is in 
Class (iii) and $H_{2}$ is in Classes (ii)-(iii). By Proposition \ref{P:dim-char} and Corollary \ref{C:dim-char2}, 
we get $F_{\Phi_{1},\Gamma^{\circ}}=F_{\Phi_{1},\Gamma^{\circ}}$ from $\mathscr{D}_{H_{1}}=\mathscr{D}_{H_{2}}$. 
By the proof of Lemma \ref{L:TD-char2}, it follows that $\Phi_{1}$ and $\Phi_{2}$ are $\Gamma^{\circ}$ conjugate. 
Then, $\dim H_{1}=\dim H_{2}+1$, which contradicts to the assumption that $\mathscr{D}_{H_{1}}=\mathscr{D}_{H_{2}}$ 
since dimension datum determines the dimension of a closed subgroup.
\end{proof}

\subsection{S-subgroups with linearly dependent dimension data}\label{SS:counter3}

Let $N>1$ be an odd integer. Put $G=\O(2N)$. Define an imbedding \[\phi: \U(N)\rtimes\langle\tau\rangle\rightarrow
\O(2N)\] by \[\phi(A+B\mathbf{i})=\left(\begin{array}{cc}A&B\\-B&A\\\end{array}\right),\ \forall (A+B\mathbf{i})\in
\U(N)\ (A,B\in M_{n}(\mathbb{R}))\] and \[\phi(\tau)=\left(\begin{array}{cc}I_{N}&0_{N}\\0_{N}&-I_{N}\\\end{array}
\right).\] As in Theorem \ref{T:counter2}, we have S-subgroups $H_{1},H_{2},H_{3},H_{4}$ of $\U(N)\rtimes\langle
\tau\rangle$.

\begin{theorem}\label{T:counter3}
For each $i=1,2,3,4$, $\phi(H_{i})$ is an S-subgroup of $G=\O(2N)$. The subgroups $\phi(H_{1})$, $\phi(H_{2})$,
$\phi(H_{3})$, $\phi(H_{4})$ are pairwise non-conjugate and we have \begin{equation}\label{Eq:counter3}
(\mathscr{D}_{\phi(H_{1})}-\mathscr{D}_{\phi(H_{2})})-(\mathscr{D}_{\phi(H_{3})}-\mathscr{D}_{\phi(H_{4})})=0.
\end{equation}
\end{theorem}

\begin{proof}
Write $\sigma$ for the natural representation of $\U(N)$ and write $\rho$ for the natural representation of
$\O(2N)$. Then, \[\rho|_{\U(N)}\cong\sigma\oplus\sigma^{\ast}.\] By the construction in Theorem \ref{T:counter2},
we know that $\sigma|_{H_{i}^{\circ}}$ is irreducible and $\sigma|_{H_{i}^{\circ}}\not\cong
\sigma^{\ast}|_{H_{i}^{\circ}}$. Then, \[Z_{\O(2N)}(H_{i}^{0})=Z(\U(N)).\] By construction, $H_{i}$ contains an
element $g_{i}\in\tau\U(N)$. Then, \[Z_{\O(2N)}(\phi(H_{i}))=(Z_{\O(2N)}(H_{i}^{0}))^{g_{i}}=Z(\U(N))^{g_{i}}=
\{\pm{I}\}.\] On the other hand, we have $\det\phi(g_{i})=-1$ due to $N$ is odd. Thus, each $\phi(H_{i})$ is an
S-subgroup of $G=\O(2N)$.

In the proof of Theorem \ref{T:counter3}, it is shown that the groups $H_{1}$, $H_{2}$, $H_{3}$, $H_{4}$ are
pairwise non-isomorphic. Then, the subgroups $\phi(H_{1})$, $\phi(H_{2})$, $\phi(H_{3})$, $\phi(H_{4})$ of
$\O(2N)$ are pairwise non-conjugate.

The equality \eqref{Eq:counter3} follows from the equality \eqref{Eq:counter2} directly.
\end{proof}

\appendix

\section{Affine root systems and quasi root systems}

\medbreak
\begin{center}
{\scshape Jiu-Kang Yu\footnote{Chinese University of Hong Kong}} 
\end{center}

\subsection{Introduction}

\ \medbreak\indent
After quasi root systems were introduced and classified in Section 1 of
the paper, we observed that the classification is remarkably similar
to the classification of affine root systems by Macdonald \cite{Macdonald}.

In this appendix, we define a notion called the {\it
  \(A\)-refinement} of a root datum, where \(A\) is an abelian group.
When \(A=\bbR\) or \(\bbZ\), \(A\)-refinements are essentially the
affine root systems.  When \(A\) is a finite cyclic groups,
\(A\)-refinements are exactly the quasi root systems.

We will give a classification of \(A\)-refinements for any cyclic
group \(A\) (finite or not).  When \(A=\bbZ\), we will recover
Macdonald's classification.  This will allows us 
to make a precise connection between affine root systems and quasi
root systems.  The connection essentially says that the quasi root systems are
``affine root systems modulo \(d\)''.  For the precise statement, see Theorem~\ref{T:affine-vs-quasi}.

\subsection{Definitions}

\begin{adefinition}
By a {\it normed root datum}, we mean a pair \((\Phi,(\cdot,\cdot))\), where
\(\Phi=(X,R,\check X,\check R)\) is a root datum and \((\cdot,\cdot)\)
is an inner product on \(X\otimes\bbR\).  It is well-known that the
datum \((X,R,\check X,\check R,(\cdot,\cdot))\) can be reduced to a
triple \((X,R,(\cdot,\cdot))\).
\end{adefinition}

\begin{adefinition}
Let \(\Phi\) be either a root datum or a normed root datum, and let
\(A\) be an abelian group.  By an \(A\)-refinement of \(\Phi\), we
mean a pair \((\tilde X,\tilde R)\) where \(\tilde X\) is an extension of \(X\) by
\(A\), i.e.~there is a short exact sequence of abelian groups
\[
0\to A\to \tilde X\to X\to 0,
\]
and \(\tilde R \subset \tilde X\) is a subset satisfying
\begin{itemize}
  \item the projection \(p:\tilde X\to X\) takes \(\tilde R\) onto
    \(R\);
   \item \(s_a(\tilde R)=\tilde R\) for all \(a \in \tilde R\), where
     \(s_a(x)=x-p(a)^\vee(p(x))a\). Notice that \(s_a\) is an order 2 element in
\(\Aut(\tilde X)\).
\end{itemize}
\end{adefinition}

It is clear that when \(A\) is a finite cyclic group, an
\(A\)-refinement of a normed root datum is exactly a quasi root
system.

When \(A=\bbR\), a discrete non-degenerate \(\bbR\)-refinement of a normed semisimple root
datum is exactly an affine root system as defined by
Macdonald.  Here we say that an \(\bbR\)-refinement is {\it discrete}
(resp.~{\it non-degenerate}) if  \(\{a \in \tilde R : p(a) = \bar
a\}\) is discrete in \(p^{-1}(\bar a)\simeq\bbR\) (resp.~is infinite) for all \(\bar a \in R\).

It is clear that all non-degenerate \(\bbZ\)-refinements of \(\Phi\) are
affine root systems.  In practice it is enough to use affine root
systems that are \(\bbZ\)-refinements.

Fix an extension \(\tilde X\) of \(X\) by \(A\).
Let \(\scrR_{\tilde X}(\Phi)\) be the set of all  \(\tilde R\) such
that \((\tilde X,\tilde R)\) is an \(A\)-refinements of
\(\Phi\).  The automorphism group of \(\tilde X\) as an extension
 is isomorphic to \(\Hom(X,A)\), and it acts on
\(\scrR_{\tilde X}(\Phi)\) naturally.  This action factors through the
action given by the next lemma, whose proof is easy and omitted.

\begin{alemma} Let \(\langle R\rangle\) be the \(\bbZ\)-submodule of \(X\)
generated by \(R\).  If \(\tilde R \in \scrR_{\tilde X}(\Phi)\),  then \(\tilde R+c:=\{a+c(a) : a \in \tilde R\}\)
is also an element of \(\scrR_{\tilde X}(\Phi)\) for any \(c \in
\Hom(p^{-1}\langle R\rangle,A)\).  This gives an action of the subgroup \(\Hom(\langle R\rangle
,A) \subset \Hom(p^{-1}\langle R\rangle ,A)\) on \(\scrR_{\tilde A}(\Phi)\).\qed
\end{alemma}

The main result of this appendix will be the classification of
\(\scrR_{\tilde X}(\Phi)\) modulo \(\Hom(\langle R\rangle ,A)\).

\subsection{The classification}

\ \medbreak\indent
Fix an extension \(0\to A\to \tilde X\to X\to 0\).  It splits since \(X\) is a
free abelian group.  From now on, we also fix a splitting and identify
\(\tilde X\) with \(X \oplus A\).  This amounts to fixing a section as
in Lemma~\ref{L:ARS5}.

Let \(\tilde R\in \scrR_{\tilde X}(\Phi)\).  Put \(A_a=\{u\in A: a+u \in \tilde
R\}\) for each \(a \in R\).  Obviously \(\tilde R\) and
\(\{A_a\}_{a\in R}\) determine each other.  The proofs of the next
two lemmas are easy and omitted.

\begin{alemma}
  For each \(a\), the set
\(A_a\) has to be closed under reflection in the sense that \(u,v\in
A_a\) implies \(2u-v\in A_a\).  If \(A\) is cyclic, this implies that
\(A_a\) is of the form \(u+A'\), where \(u \in A\) and \(A'\) is a
subgroup of \(A\).\qed
\end{alemma}

\begin{alemma}\label{Aa-condition}
Let \(\{A_a\}_{a \in R}\) be non-empty subsets of \(A\) and put
\(\tilde R=\{a+u: a\in R, u\in A_a\}\).  Then \((\tilde X,\tilde R)\) is an
\(A\)-refinement of \(\Phi\) if and only if \(A_{s_a(b)}=A_b-\check
a(b)u\) for all \(a,b\in R\), \(u \in A_a\).\qed
\end{alemma}

\begin{alemma} Assume that \(A\) is cyclic.  After replacing \(\tilde R\) by \(\tilde R+c\) if
necessary, we may assume that \(A_a\) is a subgroup of \(A\) for all
\(a\) non-divisible.
\end{alemma}

\begin{proof}
We can choose a system \(\Delta\) of simple roots for
\(R^{\rm nondiv}\), the root subsystem of non-divisible roots in \(R\).  Then \(\Delta \) is a \(\bbZ\)-basis of \(\langle
R\rangle \).  Write \(A_a = u_a + A'_a\) for each \(a \in R\),
where \(u_a\in A\) and \(A'_a\) is a subgroup of \(A\).  Let
\(c:\langle R\rangle \to A\) be given by \(c(a)=-u_a\) for all \(a \in
\Delta\).  Repalce \(\tilde R\) by \(\tilde R+c\).  Clearly, now we
may assume that \(A_a\) is a subgroup for all \(a \in \Delta\).

By the preceding lemma, we can show easily by induction on length that
\(A_{wa}\) is a subgroup of \(A\) for all \(a \in \Delta\) and all \(w\) in the Weyl group of
\(\Phi\).  Since the union of the Weyl group orbits of simple roots
are exactly the non-divisible roots, the lemma is proved.
\end{proof}

Suppose that \(R=R_1\coprod \cdots\coprod R_k\) is a decomposition
into irreducible subsystems.  Put \(\Phi_i=(X,R_i,\check X,\check
R_i)\) where \(\check R_i=\{\check a: a\in R_i\}\).  Then we have:

\begin{alemma}  Let \(\tilde R \subset \tilde X\).  Then \((\tilde X,\tilde R)\)
  is an \(A\)-refinement of \(\Phi\) if and only if
\[
\bigl(\tilde X,\{a \in\tilde R: p(a) \in R_i\}\bigr)
\]
is an \(A\)-refinement of \(\Phi_i\) for \(i=1,\ldots,k\).\qed
\end{alemma}

The two preceding lemmas allow us to reduce our classification problem
by working under the following hypothesis: {\it from now on till the
  end of this section, we assume: \(A\) is cyclic, \(\Phi\) is
  irreducible, and \(A_a\) is a subgroup for all non-divisible \(a\).}

\begin{atheorem} \label{cls-reduced}
 Assume that \(\Phi\) is reduced.  If \(\Phi\) is simply-laced,
then there is a subgroup \(A'\) of \(A\) such that \(A_a=A'\) for \(a
\in R\).

If \(\Phi\) has two root lengths such that \(|\mbox{long
  root}|^2/|\mbox{short root}|^2=r\) (so \(r=2\) or \(3\)).  Then
there are subgroups \(A_\ell\) and \(A_\rms\) of \(A\) such that \(A_a=A_\ell\)
for all long root \(a\) and \(A_a=A_\rms\) for all short root \(a\).
Furthermore, either \(A_\ell=A_\rms\) or \(A_\ell=rA_\rms\).
\end{atheorem}

\begin{proof}  Lemma~\ref{Aa-condition} implies that if \(a,a'\in R\)
  are in the same Weyl group orbit, then the subgroups \(A_a\) and
  \(A_{a'}\) are the same.  This implies the existence of \(A'\) in
  the simply laced case, and the existence of \(A_\ell\) and
  \(A_\rms\) in the other case.

Lemma~\ref{Aa-condition} also implies: \(A_b\supset \check
  a(b)A_a\) for all \(a, b \in R\).  Now assume that there are two
  root length.  We can find a long root \(a\) and a short root \(b\)
  such that \(\check a(b)=-1\) and \(\check b(a)=-r\).  This implies
  \(A_\rms \supset A_\ell\supset rA_\rms\).
\end{proof}

\bigbreak
Now assume that \(\Phi\) is non-reduced (hence of type \(BC_n\) for some
\(n\geq 1\)).  The preceding
proof implies that there are
two subgroups \(A_\rmm,A_{\rm ndnm}\) (if \(n\geq 2\)) of \(A\) such
that \(A_a=A_\rmm\) for all \(a\) multipliable, \(A_a=A_{\rm ndnm}\) for all
\(a\) non-divisible and non-multipliable.  If \(b,b'\) are divisible
and \(b\neq b'\), we can
  find a non-divisible \(a\in R\) such that \(b'=s_a(b)\).
  Lemma~\ref{Aa-condition} implies \(A_b=A_{b'}\) (because we can take
  \(u=0\)).  Therefore, there is a subset \(A_\rmd\) of \(A\) (which
  is a coset of certain subgroup) such that
  \(A_b=A_\rmd\) for all divisible \(b\).

\begin{atheorem}\label{cls-nonreduced}
If \(n=1\), then \(A_\rmd=\) \(A_\rmm\), \(2A_\rmm\),
  \(A_\rmm\setminus 2A_\rmm\), \(4A_\rmm\), or \(2A_\rmm\setminus
  4A_\rmm\).

  If \(n\geq 2\), there are at most six (non-exclusive) cases:
\begin{itemize}
  \item[--] if \(A_{\rm ndnm}=A_\rmm\), then \(A_\rmd=\) \(A_\rmm\), \(2A_\rmm\), or
    \(A_\rmm\setminus 2A_\rmm\).
   \item[--] if \(A_{\rm ndnm}=2A_\rmm\), then \(A_\rmd=\)
     \(2A_\rmm\), \(4A_\rmm\) or \(2A_\rmm\setminus 4A_\rmm\).
\end{itemize}
\end{atheorem}

\begin{proof} Assume \(n\geq 2\).  Lemma~\ref{Aa-condition} implies
  that all the restrictions on \(A_\rmm\), \(A_{\rm nmnd}\), \(A_\rmd\)
  are the following: \(A_\rmm\supset A_{\rm nmnd}\supset 2A_\rmm\cup A_\rmd\),
  and \(A_\rmd\) is invariant under
  translation by elements of  \(4A_\rmm \cup 2A_{\rm nmnd}\cup
  2A_\rmd\).  The result then follows easily.  The case \(n=1\) is similar.
\end{proof}

\begin{aremark}  If \(A_\rmm=2A_\rmm\), then the case \(A_\rmd=A_\rmm\setminus
  2A_\rmm\) can not occur and our list of cases may have repetitions.
  A similar remark applies when \(2A_\rmm=4A_\rmm\).
\end{aremark}

\begin{aremark} Assume \(2A_\rmm\neq 4A_\rmm\).  The case
  \(A_\rmd=4A_\rmm\) and the case \(A_\rmd=2A_\rmm\setminus 4A_\rmm\)
  are in the same \(\Hom(\langle R\rangle ,A)\)-orbit.  Therefore, we
  may omit one of them for the purpose of classification.
\end{aremark}

This completes our classification of \(A\)-refinements when \(A\) is
cyclic.  When \(A=\bbZ\), our theory recovers
Macdonald's classification easily.  Let us make a quick comparison and
establish some notation.  Assume that \(\Phi\) is semisimple and
irreducible of adjoint type.

If \(\Phi\) is reduced, the
\(\bbZ\)-refinement with \(A_a=k\bbZ\) for all \(a\) is denoted by \(\Phi[k]\).

If \(\Phi\) is reduced, \(R\) has two root lengths and \(r\) is as in
Theorem~\ref{cls-reduced}, the \(\bbZ\)-refinement with \(A_a=k\bbZ\)
for all short \(a\) and \(A_a=rk\bbZ\) for all long \(a\) is denoted
by \(\Phi^\vee[k]\).

Assume \(\Phi\) is non-reduced of type \(BC_n\).  We denote by the
following symbols the \(\bbZ\)-refinements with indicated
\((A_\rmm,A_{\rm nmnd},A_\rmd)\):
{\setlength{\leftmargini}{70pt}
\begin{itemize}
  \item[{\(BC_n[k]\)}:] \((A_\rmm,A_{\rm
      nmnd},A_\rmd)=(k\bbZ,k\bbZ,k\bbZ\setminus 2k\bbZ)\), \(n\geq 1\).
  \item[{\(BCC_n[k]\)}:] \((A_\rmm,A_{\rm
      nmnd},A_\rmd)=(k\bbZ,k\bbZ,k\bbZ)\), \(n\geq 1\).
  \item[{\(C^\vee BC_n[k]\)}:] \((A_\rmm,A_{\rm
      nmnd},A_\rmd)=(k\bbZ,2k\bbZ,4k\bbZ)\), \(n \geq 1\).
  \item[{\(BB^\vee _n[k]\)}:] \((A_\rmm,A_{\rm
      nmnd},A_\rmd)=(k\bbZ,k\bbZ,2k\bbZ)\), \(n\geq 2\).
  \item[{\(C^\vee C_n[k]\)}:] \((A_\rmm,A_{\rm
      nmnd},A_\rmd)=(k\bbZ,2k\bbZ,2k\bbZ)\), \(n \geq 1\).
\end{itemize}}

For all the symbols of the form \(Y[k]\) introduced above, \(Y[1]\)
represents an affine root system of type \(Y\) in Macdonald's
notation, and they exhaust all irreducible affine root systems.  We will write \(Y\) for \(Y[1]\) if there is no
danger of confusion.

\subsection{The reduction theorem}

\ \medbreak\indent
Let \(\pi:A \to A'\) be a homomorphism of abelian groups.  If
\((\tilde X,\tilde R)\) is an \(A\)-refinement of \(\Phi\), then
\((\pi_* \tilde X,\pi_*\tilde R)\) is an \(A'\)-refinement of \(\Phi\), where
\(\pi_*\tilde X\) is the pullout of \(\tilde X\) by \(\pi\), and
\(\pi_*\tilde R\) is the image of \(\tilde R\) under the natural map
\(\tilde X\to \pi_*\tilde X\).

\begin{adefinition}
When \(A=\bbZ\) and
\(A'=\bbZ/d\bbZ\), we say that \((\pi_*\tilde X,\pi_*\tilde R)\) is obtained from
\((\tilde X,\tilde R)\) by {\it reduction modulo \(d\)}.
\end{adefinition}

\begin{atheorem}\label{T:affine-vs-quasi}
  Let \(\Phi\) be a root datum and let \(d\geq 0\).  Every
  \(\bbZ/d\bbZ\)-refinement of \(\Phi\) is obtained from certain
  \(\bbZ\)-refinement of \(\Phi\) by reduction modulo \(d\).
\end{atheorem}

\begin{proof}  Let \((\tilde X,\tilde R)\) be a \(\bbZ/d\bbZ\)-refinement of
\(\Phi\).  We need to show that it is obtained by reduction modulo
\(d\).  We may and do assume that \(\tilde X=X\oplus \bbZ/d\bbZ\) so that we
can use notation and ideas from the preceding section.

It is clear that \(\Hom(\langle R\rangle ,\bbZ)\to\Hom(\langle
R\rangle ,\bbZ/d\bbZ)\) is surjective, and it is compatible with the
action of these groups on \(\scrR_{X\oplus\bbZ}(\Phi)\) and
\(\scrR_{X\oplus\bbZ/d\bbZ}(\Phi)\) respectively.  This allows us to assume
that \(A_a\) is a subgroup for all non-divisible \(a\).  It is also
clear that we can reduce to the case of an irreducible \(\Phi\).

Under these hypothesis, all possible \(\tilde R\)'s are given by
either Theorem~\ref{cls-reduced} or Theorem~\ref{cls-nonreduced}.
It is clear from these that \((\tilde X,\tilde R)\) can be obtained from
reduction modulo \(d\).
\end{proof}

The following is a dictionary of the terminologies used in the
main paper in terms of those in the appendix.  The root datum \(\Phi\) is
assumed to be irreducible.
\begin{itemize}
\item[--]folding index of \(a\): \(|A_a|\).
\item[--] total folding index: \(|A_\rms|\) when \(\Phi\) reduced,
  \(|A_\rmm|\) when \(\Phi\) is non-reduced.
\item[--] ordinary folding index: \(|A_\ell|\) when \(\Phi\) reduced,
  \(|A_{\rmd}|\) when \(\Phi\) is non-reduced.
\item[--] ndm folding index: \(|A_{\rm ndnm}|\).
\item[--] twisted number: \(|A_\rms|/|A_\ell|\) when \(\Phi\) reduced,
  \(|A_\rmm|/|A_{\rmd}|\) when \(\Phi\) non-reduced.
\item[--] ndm twisted number: \(|A_{\rm ndnm}|/|A_\rmd|\).
\end{itemize}

\begin{aexample} Let \(\Phi\) be a symbol denoting an irreduible reduced
  root system.  Then the notations introduced after Lemma~\ref{L:ARS5} are as
  follows:

{\setlength{\leftmargini}{45pt}
  \begin{itemize}
    \item[(i)] \(\Phi\) is \(\Phi [1]\bmod 1\), and \(\Phi^{(2)}\) is
      \(\Phi[1] \bmod 2\).
    \item[(ii)-(iv)] \(\Phi^{(2,1)}\) is \(\Phi^\vee[1] \bmod 2\).
    \item[(v)] \(BC_1^{'(2,1)}\) is \(C^\vee C_1[1]\bmod 2\), and
      \(BC_1^{(2,1)}\) is \(BC_1[1] \bmod 2\).
     \item[(vi)] \(BC_n^{(2,1,1)}\) is \(C^\vee C_n[1] \bmod 2\),
       \(BC_n^{'(2,2,1)}\) is \(BB^\vee[1]\bmod 2\), and
       \(BC_n^{(2,2,1)}\) is \(BC_n[1]\bmod 2\).
  \end{itemize}}
\end{aexample}

\end{document}